\documentclass[oneside,a4paper, 10pt]{report}
\usepackage[T1]{fontenc}
\usepackage[latin1]{inputenc}
\usepackage{amsmath}
\usepackage{amssymb}
\usepackage{amsthm}
\usepackage[dvips]{graphics}
\usepackage{textcomp}
\usepackage{graphicx}

\newtheorem{ex}{Example}[chapter]
\newtheorem{Def}{Definition}[chapter]
\newtheorem{Not}{Notational convention}[chapter]
\newtheorem{rem}{Remark}[chapter]
\newtheorem{Th}{Theorem}[chapter]
\newtheorem{lem}{Lemma}[chapter]
\newtheorem{aux}{Auxiliary Lemma}[chapter]
\newtheorem{cor}{Corollary}[chapter]

\newcommand{\PP}{\mathbb{P}}
\newcommand{\EE}{\mathbb{E}}
\newcommand{\II}{\mathbb{I}}
\newcommand{\CC}{\mathbb{C}}
\newcommand{\RR}{\mathbb{R}}

\newcommand{\NN}{\mathbb{N}}
\newcommand{\ZZ}{\mathbb{Z}}

\newcommand{\polm}{{\RR_m [X]}}
\newcommand{\cB}{{\cal B}}

\newcommand{\cK}{{\cal K}}
\newcommand{\cU}{{\cal U}}

\newcommand{\cR}{{\cal R}}

\newcommand{\cE}{{\cal E}}
\newcommand{\cI}{{\cal I}}
\newcommand{\cS}{{\cal S}}
\newcommand{\cF}{{\cal F}}

\newcommand{\cG}{{\cal G}}
\newcommand{\cH}{{\cal H}}
\newcommand{\cD}{{\cal D}}

\newcommand{{{\cadlag}}}{c\`adl\`ag}

\newcommand{\diag}{\mathrm{diag} \ }
\newcommand{\esssup}{\mathrm{ess \ sup} \ }
\newcommand{\id}{\mathrm{id}}

\renewcommand{\baselinestretch}{1} 

\begin{document}

\pagestyle{plain}

\title{Recent approaches to high-dimensional American and Bermudan option pricing}
\author{by Frederik S Herzberg\\Mathematical Institute and Merton College\\University of Oxford}
\date{Trinity Term 2005}

\renewcommand{\baselinestretch}{1} 

\maketitle

\begin{abstract} 

\begin{centering}
of ``Recent approaches to high-dimensional American and Bermudan option pricing'',\\by Frederik S Herzberg, Merton College, University of Oxford\\
\end{centering}

A number of Bermudan option pricing methods that are applicable to options on multiple assets are studied in this thesis, one of the dominating questions being the natural scaling needed to extrapolate from Bermudan to American (both approximate and ``exact'') option prices. Among the Bermudan option pricing techniques discussed in more detail will be (1) the use of cubature formulae for symmetric measures to price non-perpetual Bermudan options, and (2) r\'eduite-based approximation of subharmonic functions (corresponding to piecewise harmonic interpolation in the one-dimensional setting). These (iterative) algorithms shall be proven to be {\em sound} in a sense that is yet to be made precise. Also, for each iteration sequence resulting from these algorithms, monotone convergence to the least fixed point of the {iteration procedure} will be shown. (As part of an appendix, the fixed points of a Bermudan option pricing algorithm based on polynomial interpolation shall be characterised as well.)

If the iteration procedure is based on cubature, a linear convergence rate of the iteration sequence can be derived. Moreover, for a couple of practically relevant settings one can find $L^1$ estimates for a non-perpetual American option pricing algorithm based on cubature.

At the outset of this dissertation, the existence of exercise regions for multi-dimensional Bermudan options is established; afterwards one can proceed to prove bounds on the American-Bermudan barrier put option price difference (``continuity correction'') when the argument of this function -- as a function of the logarithmic start price -- approaches the exercise boundary. In particular, results of Feller's shall be generalised to show that an extrapolation from the exact Bermudan prices to the American price cannot be polynomial in the exercise mesh size in the setting of many common market models, and more specific bounds on the natural scaling exponent of the non-polynomial extrapolation for a number of (both one- and multi-dimensional) market models will be deduced.

Finally, three approximate $\Delta$-hedging algorithms for high-dimensional derivative securities are proposed and implemented, alongside with a measure of comparing their effectiveness.

\end{abstract}

{\bf Acknowledgements.} The author is highly indebted to Professor Terry Lyons for numerous helpful discussions. Furthermore, he gratefully acknowledges a post-graduate scholarship of the German National Academic Foundation ({\em Studienstiftung des deutschen Volkes}) which funded the first year of work on this Dissertation, as well as a pre-doctoral research grant of the German Academic Exchange Service ({\em Doktorandenstipendium des Deutschen Akademischen Austauschdienstes}).

\tableofcontents

\renewcommand{\baselinestretch}{1} 

\noindent

\part{Introduction}

\chapter{Definitions and basic facts on Bermudan and American options}

\label{formalintro}

In order to clarify terminology, we start by introducing the mathematical notions corresponding to the financial concepts that we shall allude to.

Our first definition is a notational convention.

\begin{Def}\label{lieexp} Let $d\in\NN$. By $\exp:\RR^d\rightarrow {\RR_{>0}}^d$ and $\ln:{\RR_{>0}}^d\rightarrow\RR^d$ we denote componentwise exponentiation and taking natural logarithms componentwise, respectively.
\end{Def}
\begin{rem} For any $d\in\NN$, $\RR^d$ is a Lie group with respect to componentwise multiplication $\cdot:(x,y)\mapsto (x_iy_i)_{i\in\{1,\dots,d\}}$. Its Lie algebra is the vector space $\RR^d$ with its usual (componentwise) addition. The exponential map from the Lie algebra $(\RR^d,+)$ into the Lie group $(\RR^d,\cdot)$ is componentwise exponentiation $\exp:x\mapsto\left(e^{x_i}\right)_{i\in\{1,\dots,d\}}$. Therefore the abbreviation introduced in Definition \ref{lieexp} is consistent with standard notation.
\end{rem}

\begin{Def} Let $T$ be a positive real number. Consider a real-valued stochastic process $X:=(X_t)_{t\in[0,T]}$, adapted to a filtered probability space $\left(\Omega,(\cF_t)_{t\in[0,T]},P\right)$. We will call $X$ a {\em logarithmic price process for a non-dividend paying asset} (for short, a {\em logarithmic price process} or simply {\em log-price process}), if and only if there exists a probability measure $Q$ equivalent to $P$ on $\cF_T$ and a constant $r>0$ such that the stochastic process $\exp\left(X_t-rt\right)_{t\in[0,T]}$ is a martingale with respect to the filtration $\cF:=(\cF_t)_{t\in[0,T]}$ and the probability measure $Q$. In this case, such a $Q$ is called a {\em martingale measure} and $r$ a {\em market price of risk} or a {\em dicsount rate} for the stochastic process $X$ and the probability measure $P$.
\end{Def}

\begin{Def} Let $d\in\NN$. A {\em $d$-dimensional basket} is a $d$-tuple of logarithmic price processes such that there exists a probability measure $Q$ and a market price of risk $r>0$ such that $Q$ is a martingale measure and $r$ a market price of risk for all components of the $d$-tuple.
\end{Def}

For the rest of this Chapter, we will adopt the terminology and the notation for Markov processes of Revuz and Yor \cite{RY}.

In particular, for all probability measures $\nu$ on $\cB\left(\RR^d\right)$, $Q_\nu$ is the probability measure induced by the transition function $\left(Q_s\right)_{s\geq 0}$ via the Ionescu-Tulcea-Kolmogorov projective limit construction, cf Revuz and Yor \cite[Theorem 1.5]{RY}).

For any $d\in\NN$, we will denote the $\sigma$-algebra of Borel subsets of $\RR^d$ by $\cB\left(\RR^d\right)$.

\begin{Def} \label{Markovbasket}Let again $d\in\NN$. A family $Y:=\left(Y^x\right)_{x\in\RR^d}$ of $\RR^d$-valued homogeneous Markov processes $Y^x$ adapted to a filtered probability space $\left(\Omega,\cF,\tilde Q\right)$ with respect to $\cF$, with transition function $\left(P_s\right)_{s\geq 0}$ and initial measure $\delta_x$, is called a {\em $d$-dimensional Markov basket} if and only if there is a homogeneous transition function $\left(Q_s\right)_{s\geq 0}$ on the measurable space $\left(\RR^d,\cB\left(\RR^d\right)\right)$ and a constant $r>0$ such that the following three assertions hold:
\begin{enumerate}
\item The process $Y^x$ is a Markov process with transition function $\left(Q_s\right)_{s\geq 0}$ with respect to $\cF$ for all $x\in\RR^d$.
\item The process $\exp\left({Y^x}_t-rt\right)_{t\in[0,T]}$ is a martingale with respect to $\cF$ and $Q_{\delta_x}$.
\item The measures $P_{\delta_x}$ and $\PP^x:=Q_{\delta_x}$ are equivalent for all $x\in\RR^d$.
\end{enumerate}

In this case, $\PP$ is called a {\em family of martingale ({\rm or:} risk-neutral) measures associated with $Y$}, and $r$ is called the {\em discount rate for $Y$}. 

The expectation operator for the probability measure $\PP^x$ will be denoted by $\EE^x$ for all $x\in\RR^d$.

If the transition function $P$ is a Feller semigroup, then we shall refer to $Y$ as a {\em Feller basket}.

If $P$ is a translation-invariant Feller semigroup, we shall call $Y$ a {\em L\'evy basket}.

\end{Def}

\begin{rem} A priori, it is not clear if there are logical connections between the three assertions in the previous Definition \ref{Markovbasket}, in particular the author does not know whether the third assertion implies the first one.
\end{rem}

\begin{Not} If no ambiguity can arise, we will drop the superscript of a Markov basket. Thus, in the notation of Definition \ref{Markovbasket}, we set $$\EE^x\left[\left.f\left(Y_{\tau_1},\dots,Y_{\tau_n}\right)\right|\cF_s\right]:=\EE^x\left[\left.f\left(Y_{\tau_1}^x,\dots,Y_{\tau_n}^x\right)\right|\cF_s\right]$$ for all $s\geq 0$, $n\in\NN$ and $n$-tuples of stopping times $\vec \tau=(\tau_1,\dots,\tau_n)$ whenever $f:\RR^n\rightarrow\RR$ is nonnegative or $f\left(Y_{\tau_1}^x,\dots,Y_{\tau_n}^x\right)\in L^1\left(\PP^x\right)$. Here we are using the term  ``stopping time'' as a synonym for $\RR_+$-valued stopping time, that is a stopping time with values in $[0,+\infty]$. 

Also, since we are explicitly allowing stopping times (with respect to the filtration generated by a process $X$) to attain the value $+\infty$, we stipulate that the random variable $f\left(X_\tau\right)$ (for any Lebesgue-Borel measurable function $f$) should be understood to be multiplied by the characteristic function of the event $\left\{\tau<+\infty\right\}$. Formally, this can be done by introducing a constant $\Delta\not\in \RR^d$, called {\em cemetery}, and stipulating that $X_\tau=\Delta$ on $\left\{\tau=+\infty\right\}$ and $f(\Delta)=0$ for all measurable functions $f$ (cf eg Revuz and Yor \cite[pp 84,102]{RY}).
\end{Not}

We will not formally define what we mean by an option itself, but we will rather define what {\em expected payoffs} and {\em prices} of some classes of financial derivatives are.

\begin{Def}

Consider a $d$-dimensional Markov basket $Y$ with an associated family $\PP^\cdot$ of martingale measures and discount rate $r>0$. 

The expected payoff of a {\em Bermudan option with (log-price) payoff function $g:\RR^d\rightarrow \RR_{\geq 0}$ on the underlying Markov basket $Y$ with exercise times in $J\subset[0,+\infty)$, log start-price $x$ and maturity $T\in[0,+\infty]$} is defined to be $$U^J(T)(x)=\sup_{\tau\text{ stopping time}, \ \tau(\Omega)\subseteq J\cup\{+\infty\} }\EE^x\left[e^{-r(\tau\wedge T)}g\left(X_{\tau\wedge T}\right)\right].$$ 

The expected payoff of a {\em perpetual Bermudan option} is the expected payoff of a Bermudan option of maturity $+\infty$.

The expected payoff of a {\em Bermudan option with exercise mesh size $h>0$} is the expected payoff of a Bermudan option with exercise times in $h\cdot\NN_0$ .

The expected payoff of an {\em American option} is the expected payoff of a Bermudan option with exercise times in $[0,+\infty)$.

We shall call the expected payoff of a Bermudan option (or an American option) a {\em Bermudan option price} (or an {\em American option price}) if and only if the martingale measures associated with the underlying basket are unique (that is, if the market model described by $P$, $\cF$ and $X$ is {\em complete}).

\end{Def}

In recent years, there has been increasing interest in incomplete market models that are governed by general L\'evy processes as log-price processes, as is not only witnessed by a tendency in research papers to focus on L\'evy process settings (for instance Boyarchenko and Levendorskii \cite{BL}; Asmussen, Avram and Pistorius \cite{AAP}; \O ksendal and Proske \cite{OP}, to take a random sample). Even textbooks, such as Karatzas' \cite{K} and Mel'nikov's \cite{M} introductory works, are putting considerable emphasis on incomplete markets. Finally, ``L\'evy finance'' has already been treated in survey articles intended for a general mathematical audience, e g Applebaum's article \cite{A04}. We will try not to deviate too much from this consensus that tries to accomplish as much mathematical generality as possible, while stopping short of studying Markov process models in their full generality. Instead we note that a substantial proportion of our results is concerned with perpetual Bermudan and American options, and it is precisely the medium and long-term risk theory where L\'evy finance seems to be applied most frequently. As a last remark on this issue, we consider it as beyond the scope of this thesis to question whether it is reasonable from an economist's point of view to study incomplete markets. 

Whilst there are some points to be made about market failures on stock markets that might entail arbitrage opportunities (for example, when assets are traded simultaneously on several stock exchanges, or in the event of insider trading), the transaction costs to exploit these arbitrage opportunities usually tend to be close to the actual gain that can be achieved through taking advantage of the arbitrage. Therefore we shall, for the sake of mathematical simplicity, merely refer to the works of Corcos et al \cite{CEMMS} as well as Imkeller et al \cite{I,IPW}, and impose a strict no-arbitrage assumption (which under certain regularity conditions on the basket is equivalent to the existence of a martingale measure, cf Karatzas \cite[Theorem 0.2.4]{K}).

\begin{ex}[A few common examples] 
\begin{enumerate} 
\item The price of a European call option on a single asset with maturity $T$ and strike price $K$ is the price of a Bermudan option with the set of exercise times being the singleton $\{T\}$ and the (log-price) payoff function $\left(\exp(\cdot)-K\right)\vee 0$.
\item The price of a perpetual American put of exercise mesh size $h>0$ on the arithmetic average of two assets in an underlying basket with strike price $K$ is the price of a Bermudan put option with the set of exercise times being the whole of the half-line $[0,+\infty)$, the maturity being $T$ and the payoff function $\left(K-\frac{\exp\left((\cdot)_1\right)+\exp\left((\cdot)_2\right)}{2}\right)\vee 0$.
\item Consider a perpetual Bermudan call option on a single asset that continuously pays dividends at a rate $\delta$ and whose logarithm follows a Markov process $Z$ adapted to some probability space $\left(\Omega,(\cF_t)_{t\geq 0},P\right)$. Then, in order to exclude arbitrage, we will have to require the existence of a family of measures $\PP^\cdot$ such that each $\PP^x$ is equivalent to $P^x$ (in particular, $\PP^{x}_{Z_0}=\delta_x$) and such that $\left(e^{-rt+\delta t+ Z_t}\right)_{t\geq 0}$ is a $\PP^x$-martingale for all $x\in \RR^d$. The expected payoff of the option will then be $$\tilde U^{h\cdot \NN_0}(\cdot)=\sup_{\tau\text{ stopping time}, \ \tau(\Omega)\subseteq h\NN_0\cup\{+\infty\} }\EE^\cdot\left[e^{-(r-\delta)\tau}\left(e^{Z_{\tau}}-K\right)\vee 0\right]$$
\end{enumerate}
\end{ex}

As an auxiliary result, let us remark  

\begin{lem}[Lower semi-continuity of $\sup$] If $I$ is a set and $\left(a_{k,\ell}\right)_{\ell\in I,k\in\NN_0}$ is a family of real numbers, then $$\sup_{\ell\in I} \liminf_{k\rightarrow\infty} a_{k,\ell}\leq \liminf_{k\rightarrow\infty} \sup_{\ell\in I} a_{k,\ell}.$$
\end{lem}
\begin{proof} We have $$\sup_{\ell}a_{k,\ell}\geq a_{k,\ell_0}$$ for all $k\in\NN_0$ and $\ell_0\in I$, therefore for all $n\in\NN$ and $\ell_0\in I$, $$\inf_{k\geq n }\sup_{\ell}a_{k,\ell}\geq \inf_{k\geq n} a_{k,\ell_0},$$ thus $$\inf_{k\geq n }\sup_{\ell}a_{k,\ell}\geq \sup_{\ell_0}\inf_{k\geq n} a_{k,\ell_0},$$ hence $$\sup_n\inf_{k\geq n }\sup_{\ell}a_{k,\ell}\geq \sup_n\sup_{\ell_0}\inf_{k\geq n} a_{k,\ell_0}=\sup_{\ell_0}\sup_n\inf_{k\geq n} a_{k,\ell_0}.$$ This is the assertion.
\end{proof}

This estimate enables us to prove the following Lemma that is asserting the approximability of expected payoffs or prices of American options by sequences of expected payoffs or prices of Bermudan options, respectively.

\begin{lem} \label{dyadic limits generic}Let $d\in\NN$, $T> 0$, $x\in\RR^d$, consider a bounded continuous function $g\geq 0$ (the {\em payoff function}), and a $d$-dimensional basket $X$ having a modification with continuous paths. If the expected payoff of an American option of maturity $T$, log start-price $x$ and payoff function $g$ on this basket $X$ is less than infinity, then the limit $$\lim_{h\downarrow 0}U^{h\NN_0}(T)(x)=\sup_{k\in\NN} U^{2^{-k}\NN_0}(T)(x)$$ exists and equals the American expected payoff.
\end{lem}
\begin{proof} Consider a sequence $(h_k)_{k\in\NN_0}\in\left(\RR_{>0}\right)^{\NN_0}$ such that $h_k\downarrow 0$ as $k\rightarrow\infty$. Choose a sequence of stopping times $(\tau_\ell)_{\ell\in\NN_0}$ such that for all $x\in\RR^d$, $$\sup_{\tau\text{ stopping time}}\EE^x\left[e^{-r(\tau\wedge T)}g\left(X_{\tau\wedge T}\right)\right]=\sup_{\ell}\EE^x\left[e^{-r(\tau_\ell\wedge T)}g\left(X_{\tau_\ell\wedge T}\right)\right]$$ and define $$\tau_{\ell,k}:= \inf\left\{t\in h_k\NN_0 \ : \ t\geq \tau_\ell\right\}.$$ Then, due to the continuity conditions we have imposed on $g$ and on the paths of (a modification of) the basket $X$, we get $$\sup_{\ell} e^{-r\left(\tau_\ell\wedge T\right)}g\left(X_{\tau_\ell\wedge T}\right)= \sup_{\ell} \lim_{k\rightarrow\infty} e^{-r\left(\tau_{\ell,k}\wedge T\right)}g\left(X_{\tau_{\ell,k}\wedge T}\right)$$ and hence by the lower semi-continuity of $\sup$, one obtains $$\sup_{\ell} e^{-r\left(\tau_\ell\wedge T\right)}g\left(X_{\tau_\ell\wedge T}\right)\leq \liminf_{k\rightarrow\infty} \sup_{\ell} e^{-r\left(\tau_{\ell,k}\wedge T\right)}g\left(X_{\tau_{\ell,k}\wedge T}\right).$$ Now we can use the Montone Convergence Theorem and Lebesgue's Dominated Convergence Theorem (this is applicable because of the boundedness of $g$) to swap limits/suprema with the expectation operator. Combining this with the specific choice of the sequence $\left(\tau_\ell\right)_{\ell\in\NN_0}$, this yields for all $x\in\RR^d$, \begin{eqnarray*} &&\sup_{\tau\text{ stopping time}}\EE^x\left[e^{-r(\tau\wedge T)}g\left(X_{\tau\wedge T}\right)\right]=\sup_{\ell}\EE^x\left[e^{-r(\tau_\ell\wedge T)}g\left(X_{\tau_\ell\wedge T}\right)\right] \\ &=& \EE^x\left[\sup_{\ell}e^{-r(\tau_\ell\wedge T)}g\left(X_{\tau_\ell\wedge T}\right)\right] \leq \EE^x\left[\liminf_{k\rightarrow\infty} \sup_{\ell} e^{-r\left(\tau_{\ell,k}\wedge T\right)}g\left(X_{\tau_{\ell,k}\wedge T}\right)\right]\\ &\leq& \liminf_{k\rightarrow\infty}\EE^x\left[ \sup_{\ell} e^{-r\left(\tau_{\ell,k}\wedge T\right)}g\left(X_{\tau_{\ell,k}\wedge T}\right)\right] \\ &=& \liminf_{k\rightarrow\infty}\sup_{\ell} \EE^x\left[ e^{-r\left(\tau_{\ell,k}\wedge T\right)}g\left(X_{\tau_{\ell,k}\wedge T}\right)\right] \\&\leq& \liminf_{k\rightarrow\infty}\sup_{\tau(\Omega)\subseteq h_k\NN_0\cup\{+\infty\}} \EE^x\left[ e^{-r\left(\tau\wedge T\right)}g\left(X_{\tau\wedge T}\right)\right] \\&\leq& \limsup_{k\rightarrow\infty}\sup_{\tau(\Omega)\subseteq h_k\NN_0\cup\{+\infty\}} \EE^x\left[ e^{-r\left(\tau\wedge T\right)}g\left(X_{\tau\wedge T}\right)\right] \\&\leq& \sup_{\tau\text{ stopping time}}\EE^x\left[e^{-r(\tau\wedge T)}g\left(X_{\tau\wedge T}\right)\right]. \end{eqnarray*} This finally gives \begin{eqnarray*}&&\sup_{\tau\text{ stopping time}}\EE^x\left[e^{-r(\tau\wedge T)}g\left(X_{\tau\wedge T}\right)\right]\\&=&\lim_{k\rightarrow\infty}\sup_{\tau(\Omega)\subseteq h_k\NN_0\cup\{+\infty\}} \EE^x\left[ e^{-r\left(\tau\wedge T\right)}g\left(X_{\tau\wedge T}\right)\right]=U^{h_k\NN_0}(T)(x). \end{eqnarray*} Since the left hand side does not depend on $(h_k)_k$, we conclude that $\lim_{h\downarrow 0}U^{h\NN_0}(T)(x)$ exists and is equal to $\sup_{\tau\text{ stopping time}}\EE^x\left[e^{-r(\tau\wedge T)}g\left(X_{\tau\wedge T}\right)\right]$.
\end{proof}

\section{Some classes of Bermudan option pricing algorithms}

Let $C^0\left(\RR^d,[0,+\infty)\right)$ and $L^0\left(\RR^d,[0,+\infty)\right)$, as usual, denote the spaces of nonnegative continuous functions defined on $\RR^d$, and of nonnegative measurable functions defined on $\RR^d$, respectively.

The purpose of the following definitions is merely to introduce a {\it fa\c{c}on de parler} which will allow us to quickly describe desirable properties of approximative Bermudan pricing algorithms in the later parts of this thesis.

\begin{Def}\label{algorithmclasses} A map $D:L^0\left(\RR^d,[0,+\infty)\right)\rightarrow L^0\left(\RR^d,[0,+\infty)\right)$ is said to be a {\em sound iterative Bermudan option pricing algorithm} (for short, a {\em sound algorithm}) for a payoff function $g\in C^0\left(\RR^d,[0,+\infty)\right)$ if and only if $Df\geq g$ for all $f\in L^0\left(\RR^d,[0,+\infty)\right)$ and the map $D$ is pointwise monotone, that is \begin{eqnarray*}&& \forall f_0,f_1\in L^0\left(\RR^d,[0,+\infty)\right)\\ && \left(\forall x\in\RR^d \quad f_0(x)\leq f_1(x)\Rightarrow \forall x\in\RR^d \quad Df_0(x)\leq Df_1(x) \right).\end{eqnarray*} A sound iterative Bermudan option pricing algorithm $D$ is said to {\em have a perpetual limit} if and only if $$\sup_{n\in\NN_0} D^{\circ n}g\in L^0\left(\RR^d,[0,+\infty)\right)$$ (rather than this supremum being allowed to equal $+\infty$ on a subset of positive measure of its range). In that very case, the function in the last line is simply referred to as the {\em perpetual limit} of the algorithm. Finally, $D$ is said to {\em converge linearly in $L^\infty$ to the perpetual limit} if and only if there exists a $c\in(0,1)$ such that $$\forall n\in\NN\quad \left\|\left(D^{n+1}-D^n\right)g\right\|_{L^\infty(\RR^d,\RR)}\leq c\cdot \left\|\left(D^{n}-D^{n-1}\right)g\right\|_{L^\infty(\RR^d,\RR)},$$ $D^n$ being shorthand for $D^{\circ n}$ for all $n\in\NN_0$.
\end{Def}

\begin{rem} The elements of $L^0\left(\RR^d,[0,+\infty)\right)$ should be conceived of assigning the value -- that is, the expected payoff -- of an option to the vector of logarithmic start prices of the components of the basket (at least on the complement of a Lebesgue null set).
\end{rem}
\begin{rem} \label{soundmonotone} The monotonicity condition imposed on sound iterative Bermudan pricing algorithms entail that the sequence of functions $\left(D^ng\right)_{n\in\NN_0}=\left(D^{\circ n}g\right)_{n\in\NN_0}$ is always pointwise increasing. Thus, this sequence has a limit: $$g\leq \sup_{n\in\NN} D^ng= \lim_{n\rightarrow\infty} D^ng.$$
\end{rem}

The infimum of all $D$-fixed points is always an upper bound for the perpetual limit: 

\begin{lem} \label{inffixedpoint}Let $D$ be a sound iterative Bermudan pricing algorithm for $g$ with a perpetual limit $u$. Then the function $u$ is smaller than any fixed point of $D$; moreover, $Du\geq u$.
\end{lem}
\begin{proof} Any fixed point $h$ of $D$ is in the image of $D$ and therefore, due to our assumptions on sound algorithms, pointwise greater or equal $g$. Now, as $D$ (and thus $D^n$) is pointwise monotone, $$\forall n\in\NN_0 \quad h=D^nh\geq D^ng,$$ therefore $$h=\sup_m D^mh\geq \sup_mD^mg,$$ where the right hand side is just the perpetual limit. Hence, any fixed point of $D$ is greater or equal the perpetual limit. Furthermore, observe that due to the pointwise monotonicity of $D$, $$D\left(\sup_nD^ng\right)\geq D^{m+1}g,$$ therefore for all $m\in\NN_0$, $$D\left(\sup_nD^ng\right)\geq \sup_m D^{m+1}g=\sup_mD^mg.$$ 
\end{proof}

Later on, it will turn out that if $D$ is based on either cubature or piecewise harmonic interpolation or the r\'eduite, the perpetual limit is, in fact the minimal fixed point (cf Theorem \ref{fixexist} and Lemma \ref{Lemma7.6} for piecewise harmonic interpolation, Theorem \ref{Theorem7.2} for r\'eduite-based approximation, and Theorem \ref{monoconv} for a result on a map $D$ which is based on cubature).

Moreover, we shall show that the algorithm based on cubature converges linearly in the sense of the definition above.

\section{Outline of this thesis}

We will postpone giving a more informal account of our motivation to use cubature formulae for Bermudan option pricing until we study the numerical implementation of some  Bermudan pricing algorithms. In this thesis, we will first of all prepare the derivation of bounds on the natural scaling of the difference between an American and a Bermudan perpetual barrier option price (conceived of as a function of the Bermudan's exercise mesh size). Later on, we will show that this is sufficient to obtain bounds on the natural scaling for the difference of certain non-perpetual American and Bermudan barrier options. For the one-dimensional setting, analogous results have been obtained by Broadie, Glasserman and Kou \cite{BGK}.

Furthermore, we will prove soundness and existence of a perpetual limit for a number of Bermudan pricing algorithms, including pricing based on cubature. In particular, as was previously remarked, we will obtain a linear convergence rate for the latter class of algorithms.

Later on, we shall prove convergence bounds for a non-perpetual American pricing algorithm in which one is computing non-perpetual Bermudan prices of a certain exercise mesh size via cubature and successively halves their exercise mesh size while leaving the number of paths at which the option is evaluated constant.

The natural scaling for which bounds are derived in the first chapters of this report, can be used to consistently extrapolate from a finite number of Bermudan barrier prices to an approximation for the American price. En passant, we will sketchily explain how some of the features of object-oriented (C++) programming can be exploited to enhance the efficiency of a Bermudan pricing algorithm based on cubature (that is to say, how to circumvent exponential complexity by achieving recombination through the map class template).

\section{Notation}

We are following largely standard probabilistic notation, as can be found for instance in the works by It\^o and McKean jr. \cite{IK} or Revuz and Yor \cite{RY}.

Both $A\subseteq B$ and $A\subset B$ for sets $A$ and $B$ will mean that $A$ is a subset of $B$ (possibly $A=B$).
\chapter{Exercise regions}

In this Chapter, we will give a rigorous proof for the fact that an American/Bermudan option price coincides with the payoff that is expected if one exercises at the first possible entry of the log-price process into the immediate exercise region ($G\subset \RR^d$). Characterisations of such regions for special cases have been proven in recent years \cite{BD,P}. 

\begin{Def} Given a countable subset $\cI\subset [0,+\infty)$ and a Lebesgue-Borel measurable set $G\subset\RR^d$, often referred to as {\em exercise region}, we define the stopping time $$\tau_G^\cI:=\min\left\{t\in \cI \ : \ X_t\in G \right\},$$ (the superscript will be dropped when no ambiguity can arise) which is just the first (nonnegative) entry time in $\cI$ into $G$. If $\cG$ is a subset of space-time, that is $\cG\subset\RR^d\times[0,+\infty)$ rather than space (ie $\RR^d$) itself, we use the space-time process rather than just the process itself to give an analogous definition: $$\tau_\cG^\cI:=\min\left\{t\in t_0+\cI\ : \ (X_t,t)\in \cG \right\},$$ where $t_0$ is the time-coordinate at which the space-time process was started. Also, for $h>0$ we set $$\tau^h_G:=\tau_G^{h\NN},\quad \tau^h_\cG:=\tau_\cG^{h\NN}$$ to denote the first positive entry time in $h\NN_0$ into $G$ or $\cG$, respectively, whilst finally $\bar\tau_\cG^h:=\tau_\cG^{h\NN_0}$ and $\bar\tau_G^h:=\tau_G^{h\NN_0}$ will denotes the first nonnegative entry time into $\cG$ and $G$, respectively.

\end{Def}

For convenience, we will also adopt the following convention for this Chapter:

\begin{Def} Let $\cI\subset[0,+\infty)$. A stopping time $\tau$ is called {\em $\cI$-valued} if the range of $\tau$, denoted by $\mathrm{ran \ } \tau$, is a subset of $\cI\cup\{+\infty\}$.
\end{Def}

\begin{lem}\label{exerciseboundaryLemma} Consider a countable subset $\cI\subset [0,+\infty)$. Let $X$ be a $d$-dimensional basket with an associated risk-neutral measure $\PP$ and discount rate $r>0$. Suppose $g=(K-f)\vee 0$, $s>0$ and $e^{-r\cdot}f(X_{\cdot})$ is a $\PP$-submartingale. For all $\cI$-valued stopping times $\tau$ there is a space-time region $B=\bigcup_{u\in \cI} \{u\}\times B_u$ such that $$\EE\left[e^{-r\left(\tau\wedge T\right)}g(X_{\tau\wedge T})\right]\leq \EE\left[e^{-r\left(\tau_B^\cI\wedge T\right)}g\left(X_{\tau_B^\cI\wedge T}\right)\right]$$
for all $T\in[0,+\infty)$ where $$\tilde\tau:=\tau_B^\cI=\inf\left\{u\in \cI\ : \ X_u\in B_u\right\}.$$ If $\sup\left\{e^{-r\tau}g\left(X_\tau\right)\ : \ \tau\text{ stopping time} \right\}$ is $\PP$-integrable, then the latter inequality will also hold for $T=+\infty$. 

The Lemma holds in particular for $\cI=s\NN_0$ for arbitrary $s>0$.
\end{lem}
\begin{proof} Firstly, we will treat the case of $T<+\infty$. Define $$\forall t\in\cI\quad B_t:=X_{\tau}\left(\left\{\tau=t\right\}\right)\subset \RR^d.$$ Let us first of all assume that \begin{equation}\label{gXtau>0onCt}\forall t\in \cI\quad g(X_t)>0 \text{ a.s. on }\left\{X_t\in B_t\right\},\end{equation} and let us also for the moment suppose \begin{equation}\label{gXtau>0onAt}\forall t\in \cI\quad\{\tilde\tau=t\}\cap\left\{\tau>T\right\}=\emptyset.\end{equation} Both of these assumptions will be dropped at the end of the proof for the case $T<+\infty$ in order to show the Lemma in its full strength. Now, from equations (\ref{gXtau>0onAt}) and (\ref{gXtau>0onCt}) one may derive \begin{eqnarray} &&\EE\left[e^{-r\left(\tau\wedge T\right)}g(X_{\tau\wedge T}),\left\{\tilde \tau=t\right\}\right] \nonumber \\  &=& \EE\left[e^{-r\tau}g\left(X_{\tau}\right),\left\{\tilde \tau=t\right\}\right] \nonumber \\ &=& \EE\left[e^{-r\tau}\left(K-f(X_{\tau})\right),\left\{\tilde \tau=t\right\}\right] \nonumber\\ &=& \EE\left[e^{-r\left(\tau\wedge T\right)}\left(K-f(X_{\tau\wedge T})\right),\left\{\tilde \tau=t\right\}\right]  \label{gbyK-f}\end{eqnarray} for all $t\in[0,T]\cap \cI$.

Furthermore, observe that $\tau\geq \tilde\tau $ a.s. Using Doob's Optional Stopping Theorem (see eg Varadhan \cite[Theorem 5.11]{Va}), we infer from our assumption of $e^{-r\cdot}f(X_{\cdot})$ being a $\PP$-submartingale with respect to the canonical filtration $\cF$ the assertion that $\left(e^{-r\upsilon}f\left(X_\upsilon\right)\right)_{\upsilon\in \left\{\tilde\tau\wedge T, \tau\wedge T\right\}}$ is a $\PP$-submartingale with respect to the filtration $\left\{\cF_{\tilde\tau\wedge T},\cF_{\tau\wedge T}\right\}$. Hence, if we combine this with equation (\ref{gbyK-f}) and note that $\left\{\tilde \tau=t\right\}=\left\{\tilde \tau\wedge T=t\right\}\in\cF_{\tau\wedge T}$ for all $t\in[0,T)\cap \cI$, we obtain for every $t\in[0,T)\cap \cI$, \begin{eqnarray*}&&\EE\left[e^{-r\left(\tau\wedge T\right)}g(X_{\tau\wedge T}),\left\{\tilde \tau=t\right\}\right]\\&= &  \EE\left[e^{-r\left(\tau\wedge T\right)}\left(K-f(X_{\tau\wedge T})\right),\left\{\tilde \tau=t\right\}\right]\\ &\leq & K\cdot \EE\left[e^{-r\left(\tau\wedge T\right)},\left\{\tilde \tau=t\right\}\right]-\EE\left[e^{-r\left(\tilde\tau\wedge T\right)}f\left(X_{\tilde\tau\wedge T})\right),\left\{\tilde \tau=t\right\}\right]\\ &\leq & K\cdot \EE\left[e^{-r\left(\tilde\tau\wedge T\right)},\left\{\tilde \tau=t\right\}\right]-\EE\left[e^{-r\left(\tilde\tau\wedge T\right)}f\left(X_{\tilde\tau\wedge T})\right),\left\{\tilde \tau=t\right\}\right] \\ &= & \EE\left[e^{-r\left(\tilde\tau\wedge T\right)}g\left(X_{\tilde\tau\wedge T})\right),\left\{\tilde \tau=t\right\}\right] .\end{eqnarray*} On the other hand, since $\tilde \tau\leq \tau$, if $\tilde\tau\geq T$, then also $\tau\geq T$, entailing $$\tilde \tau \wedge T=T=\tau\wedge T \text{ on } \left\{\tilde\tau \geq T\right\}.$$ Summarising these last two remarks, one concludes \begin{eqnarray*}&&\EE\left[e^{-r\left(\tau\wedge T\right)}g\left(X_{\tau\wedge T}\right)\right]\\ &=&\sum_{t\in \cI}\EE\left[e^{-r\left(\tau\wedge T\right)}g\left(X_{\tau\wedge T}\right),\left\{\tilde \tau=t\right\}\right] \\ &=& \sum_{t\in \cI\cap [0,T)}\EE\left[e^{-r\left(\tau\wedge T\right)}g\left(X_{\tau\wedge T}\right),\left\{\tilde \tau=t\right\}\right] \\ &&+\sum_{t\in \cI\cap [T,+\infty)}\EE\left[e^{-r\left(\tau\wedge T\right)}g\left(X_{\tau\wedge T}\right),\left\{\tilde \tau=t\right\}\right]  \\&\leq & \sum_{t\in \cI\cap [0,T)}\EE\left[e^{-r\left(\tilde\tau\wedge T\right)}g\left(X_{\tilde\tau\wedge T}\right),\left\{\tilde \tau=t\right\}\right] \\ && + \sum_{t\in \cI\cap [T,+\infty)}\EE\left[e^{-r\left(\tilde\tau\wedge T\right)}g\left(X_{\tilde\tau\wedge T}\right),\left\{\tilde \tau=t\right\}\right]\\ &=& \sum_{t\in \cI}\EE\left[e^{-r\left(\tilde\tau\wedge T\right)}g\left(X_{\tilde\tau\wedge T}\right),\left\{\tilde \tau=t\right\}\right] \\ &=&  \EE\left[e^{-r\left(\tilde\tau\wedge T\right)}g\left(X_{\tilde\tau\wedge T}\right)\right].\end{eqnarray*} 

In order to complete the proof for the case of $T<+\infty$, let us show that the assumptions (\ref{gXtau>0onCt}) and (\ref{gXtau>0onAt}) are dispensable. 

If the assertion (\ref{gXtau>0onCt}) failed to hold, we would simply define the stopping time $$\tau':= \chi_{\complement\left(\bigcup_{t\in \cI}\left\{\tau =t\right\}\cap\left\{g(X_t)>0\right\}\right)}\cdot\infty+\sum_{t\in \cI}\chi_{\left\{\tau=t\right\}\cap\left\{g(X_t)>0\right\}}\cdot t$$ and based on this definition, we would set $$\forall t\in\cI\quad B_t':=X_{\tau'}\left(\left\{\tau'=t\right\}\right)=X_{\tau}\left(\left\{\tau'=t\right\}\right)\subset \RR^d.$$ Then we would have on the one hand (\ref{gXtau>0onCt}) for $B'$ instead of $B$ which, according to what we have been able to show under the assumption of (\ref{gXtau>0onCt}), yields $$\EE\left[e^{-r\left(\tau'\wedge T\right)}g(X_{\tau'\wedge T})\right]\leq \EE \left[e^{-r\left(\tau_{B'}^\cI\wedge T\right)}g(X_{\tau_{B'}^\cI\wedge T})\right],$$ provided the condition (\ref{gXtau>0onAt}) is satisfied. However, in any case $$\EE\left[ e^{-r\left(\tau'\wedge T\right)}g(X_{\tau'\wedge T})\right]=\EE\left[ e^{-r\left(\tau\wedge T\right)}g(X_{\tau\wedge T})\right].$$ Thus, with $B'$ we have found a set that can play the r\^ole of $B$ in the Lemma's statement, under the assumption (\ref{gXtau>0onAt}).

Now suppose the condition (\ref{gXtau>0onAt}) was not satisfied, and (\ref{gXtau>0onCt}) may or may not hold (in the former case, one may even replace $\tau'$ by $\tau$ in what follows). In this situation we consider the stopping time $$\tau''=\chi_{\left\{\tau'\leq T\right\}\cup\left(\left\{\tau'>T\right\}\cap \left\{\tau_{B'}^\cI\geq T\right\}\right)}\cdot\tau'+\chi_{\complement\left(\left\{\tau'\leq T\right\}\cup\left(\left\{\tau'>T\right\}\cap \left\{\tau_{B'}^\cI\geq T\right\}\right)\right)}\cdot\infty.$$ If one now defines $$\forall t\in\cI\quad B_t'':=X_{\tau''}\left(\left\{\tau''=t\right\}\right)$$ then 
$$\forall t\in\cI\cap [0,T]\quad B_t''=X_{\tau''}\left(\left\{\tau'=\tau''=t\right\}\right)=X_{\tau'}\left(\left\{\tau'=t\right\}\right)= B_t'$$ and \begin{eqnarray*}\forall t\in\cI\cap (T,+\infty)\quad B_t''&=&X_{\tau''}\left(\left\{\tau''=t\right\}\right)=X_{\tau''}\left(\left\{\tau''=\tau'=t\right\}\cap \left\{\tau_{B'}^\cI\geq T\right\}\right)\\ &\subset& X_{\tau'}\left(\left\{\tau'=t\right\}\right)= B_t',\end{eqnarray*} hence (\ref{gXtau>0onCt}) holds for $B''$ instead of $B$ (since (\ref{gXtau>0onCt}) holds for $B'$ instead of $B$ and as we have just seen $B''_t\subset B_t'$ for all $t\in\cI$). Furthermore,
\begin{eqnarray*}\tau_{B''}^\cI &=&\chi_{\left\{\tau'\leq T\right\}\cup\left(\left\{\tau'>T\right\}\cap \left\{\tau_{B'}^\cI\geq T\right\}\right)}\cdot\tau_{B'}^\cI+\chi_{\complement\left(\left\{\tau'\leq T\right\}\cup\left(\left\{\tau'>T\right\}\cap \left\{\tau_{B'}^\cI\geq T\right\}\right)\right)}\cdot\infty\\ &=& \chi_{\left\{\tau''\leq T\right\}}\cdot \tau_{B'}^\cI + \chi_{\left\{+\infty>\tau''>T\right\}\cap \left\{\tau_{B'}^\cI\geq T\right\}}\cdot\tau_{B'}^\cI\\ && +\chi_{\complement\left(\left\{\tau'\leq T\right\}\cup\left(\left\{\tau'>T\right\}\cap \left\{\tau_{B'}^\cI\geq T\right\}\right)\right)}\cdot\infty\end{eqnarray*} (the first line because of $X_{\tau''}\left(\complement\left(\left\{\tau'\leq T\right\}\cup\left(\left\{\tau'>T\right\}\cap \left\{\tau_{B'}^\cI\geq T\right\}\right)\right)\right)=\emptyset$). Therefore $$\tau_{B''}^\cI\geq T\quad \text{ on }\left\{\tau''>T\right\},$$ thus (\ref{gXtau>0onAt}) holds for $\tau_{B''}^\cI$ instead of $\tilde\tau$ and $\tau''$ instead of $\tau$. But we have already seen that (\ref{gXtau>0onCt}) holds for $B''$ instead of $B$. Therefore, using what we have proven under the assumption of both (\ref{gXtau>0onCt}) and (\ref{gXtau>0onAt}), we get $$\EE\left[ e^{-r\left(\tau''\wedge T\right)}g(X_{\tau''\wedge T})\right]\leq \EE \left[e^{-r\left(\tau_{B''}^\cI\wedge T\right)}g(X_{\tau_{B''}^\cI\wedge T})\right].$$ On the other hand, however, $$\EE\left[ e^{-r\left(\tau''\wedge T\right)}g(X_{\tau''\wedge T})\right]=\EE\left[ e^{-r\left(\tau'\wedge T\right)}g(X_{\tau'\wedge T})\right]$$ (as $\tau''=\tau'$ on $\{\tau< T\}$, as well as $\tau''\geq T$ on $\{\tau\geq T\}$, thus $\tau''\wedge T=T=\tau\wedge T $ on $\{\tau\geq T\}$) and we have already seen that $$\EE\left[ e^{-r\left(\tau'\wedge T\right)}g(X_{\tau'\wedge T})\right]=\EE\left[ e^{-r\left(\tau\wedge T\right)}g(X_{\tau\wedge T})\right].$$ Finally, $$\EE\left[ e^{-r\left(\tau\wedge T\right)}g(X_{\tau\wedge T})\right]\leq \EE \left[e^{-r\left(\tau_{B''}^\cI\wedge T\right)}g(X_{\tau_{B''}^\cI\wedge T})\right]$$ whence with $B''$ we have found a set that can play the r\^ole of $B$ in the Lemma's statement.

Finally, we need to consider the case where $T=+\infty$. The random variable $\sup\left\{e^{-r\tau}g\left(X_\tau\right)\ : \ \tau\text{ stopping time} \right\}$ is an upper bound on $e^{-r\tau'\wedge T}g\left(X_{\tau'\wedge T}\right)$ for all $T>0$. Therefore, as soon as $\sup_{\tau\text{ stopping time} }e^{-r\tau}g\left(X_\tau\right)$ is $\PP$-integrable, we are allowed to apply Lebesgue's Dominated Convergence Theorem and the assertion for $T=+\infty$ follows by letting $T$ tend to infinity in $$\EE\left[e^{-r\left(\tau\wedge T\right)}g(X_{\tau\wedge T})\right]\leq \EE\left[e^{-r\left(\tau_B^\cI\wedge T\right)}g\left(X_{\tau_B^\cI\wedge T}\right)\right]$$ (which we have already proven for all $T\in (0,+\infty)$).

\end{proof}

\begin{cor}[Formula for an option price using hitting times] \label{exerciseboundaryCorollary}Let $X$ be a $d$-dimensional basket with an associated risk-neutral measure $\PP$ and discount rate $r>0$. Consider a countable subset $\cI\subset[0,+\infty)$. Suppose $g=(K-f)\vee 0$, $s>0$ and assume that the process $e^{-r\cdot}f(X_{\cdot})$ is a $\PP$-submartingale. Then one has \begin{eqnarray*}&&\sup_{\cG\subset \RR^d\times [0,T] \text{ measurable}}\EE \left[e^{-r\tau_\cG^\cI}g\left(X_{\tau_\cG^\cI}\right)\right] \\&=&\sup_{\tau \ \cI \cap [0,T]\text{-valued stopping time}}\EE \left[e^{-r\tau}g\left(X_\tau\right)\right]\end{eqnarray*} for all $T<+\infty$. If the random variable $\sup_{\tau\text{ stopping time} }e^{-r\tau}g\left(X_\tau\right)$ is $\PP$-integrable, then the equation \begin{eqnarray*}&& \sup_{\cG\subset \RR^d\times [0,+\infty) \text{ measurable}}\EE \left[e^{-r\tau_\cG^\cI}g\left(X_{\tau_\cG^\cI}\right)\right] \\&=& \sup_{\tau \ \cI\text{-valued stopping time}}\EE \left[e^{-r\tau}g\left(X_\tau\right)\right]\end{eqnarray*} holds.

\end{cor}

\begin{Def} Let $\cI\subset[0,+\infty)$ be countable, $\cG\subset \RR^d\times [0,+\infty)$ and $G\subset \RR^d$ measurable, and $X$ a $d$-dimensional Markov basket with an associated family of risk-neutral measures $\PP^\cdot$ and discount rate $r>0$. We define 
$$V_{\cG,X}^\cI:(x,t)\mapsto\left\{\begin{array}{*{2}{c}} e^{rt}\EE^{(x,t)}\left[e^{-r\tau_\cG^\cI}g\left(X_{\tau_\cG^\cI}\right)\right],&(x,t)\notin \cG, \\ g(x), & (x,t)\in \cG.\end{array}\right.$$ as well as $$V_{G,X}^\cI:x\mapsto\left\{\begin{array}{*{2}{c}}\EE^x\left[e^{-r\tau_G^\cI}g\left(X_{\tau_G^\cI}\right)\right], &x\notin G, \\ g(x), & x\in G.\end{array}\right.$$ Instead of $V_{\cG,X}^\cI(x,0)$, we shall often simply write $V_{\cG,X}^\cI(x)$. Also, the subscript $X$ will be dropped when no ambiguity can arise. Also, $V_\cG^h$ and $V_G^h$ will be shorthand for $V_\cG^{h\NN_0}$ and $V_G^{h\NN_0}$, respectively.
\end{Def}

As another notational convention, let us from now on use $\sup_{\cG\subset \RR^d\times [0,T] }$ and $\sup_{G\subset \RR^d}$ to denote $\sup_{\cG\subset \RR^d\times [0,T] \text{ measurable}}$ and $\sup_{G\subset \RR^d\text{ measurable}}$, respectively.

\begin{Th}[Optimality of the immediate exercise region]\label{exerciseboundary} Let $X$ be a $d$-dimensional Feller basket with $\PP^\cdot$ being an associated family of risk-neutral measures and $r>0$ being the discount rate belonging to $\PP^\cdot$. Suppose $g=(K-f)\vee 0$, $\cI\subset[0,+\infty)$ is countable, and $T\in[0,+\infty]$. Assume, moreover, that $e^{-r\cdot}f(X_{\cdot})$ is a $\PP^x$-submartingale for all $x\in\RR^d$. Define $$F^{\cI,T} =\left\{(x,t)\in\RR^d\times[0,T] \ : \sup_{\cG\subset \RR^d\times[0,T]}V_\cG^\cI(x,t) \leq g(x)\right\} $$ if $T<+\infty$ (we may drop the superscript $T$ wherever this is unambiguous) and else $$F^{\cI,+\infty} =\left\{(x,t)\in\RR^d\times[0,\infty) \ : \sup_{\cG\subset \RR^d\times[0,+\infty)}V_\cG^\cI(x,t) \leq g(x)\right\} $$  Then $$\forall x\in\RR^d \quad V_{F^{\cI,T}}^\cI(x,0)=U^{s\NN_0}(T)(x)$$ if $T<+\infty$, and $V_{F^{\cI,+\infty}}^\cI(x,0)=U^{s\NN_0}(+\infty)(x)$ for all $x\in\RR^d$ such that the random variable $\sup_{\tau\text{ stopping time} }e^{-r\tau}g\left(X_\tau\right)$ is $\PP^x$-integrable.
\end{Th}
\begin{proof} Let $T<+\infty$. Using Corollary \ref{exerciseboundaryCorollary} and recalling the definition of $U^{s\NN_0}$, all we have to show is $$\forall x\in\RR^d \quad V_{F^{\cI,T}}^\cI(x,0)= \sup_{\cG\subset \RR^d\times [0,T] }\EE^{(x,0)} \left[e^{-r\bar\tau_\cG^\cI}g\left(\bar\tau_\cG^\cI\right)\right]$$ (where we recall that $\bar\tau_\cG^\cI\leq \tau_\cG^\cI$ denotes the first nonnegative entry time into $\cG$). However, after exploiting the special particular shape of $F^\cI$, we can -- due to the boundedness of $g\leq K$ which yields $V_\cG^\cI\leq K$ for all $\cG$ which allows us to apply Lebesgue's Dominated Convergence Theorem -- swap $\sup$ and $\EE$ to get for all $x\in\RR^d$, \begin{eqnarray*}V_{F^\cI}^\cI(x,0) &=& \EE^{(x,0)}\left[e^{-r\bar\tau_{F^\cI}^\cI} g\left(X_{\bar\tau_{F^\cI}^\cI}\right) \right] \\ &=& \EE^{(x,0)}\left[e^{-r\bar\tau_{F^\cI}^\cI} \sup_{\cG\subset \RR^d\times[0,T]}V_\cG^\cI\left(X_{\bar\tau_{F^\cI}^\cI},{\bar\tau_{F^\cI}^\cI}\right) \right] \\ &=& \EE^{(x,0)}\left[e^{-r\bar\tau_{F^\cI}^\cI} \sup_{\cG\subset \RR^d\times[0,T]} e^{r\bar\tau_{F^\cI}^\cI}\EE^{\left(X_{\bar\tau_{F^\cI}^\cI},{\bar\tau_{F^\cI}^\cI}\right)} e^{-r\bar\tau_{\cG}^\cI} g\left(X_{\bar\tau_{\cG}^\cI}\right)\right] \\ &=& \sup_{\cG\subset \RR^d\times [0,T]} \EE^{(x,0)}\left[\EE^{\left(X_{\bar\tau_{F^\cI}^\cI},{\bar\tau_{F^\cI}^\cI}\right)}e^{-r\bar\tau_{\cG}^\cI}g\left(X_{\bar\tau_\cG^\cI}\right) \right]\end{eqnarray*} (where for notational convenience $\bar\tau_\cG^\cI$ should denote the first nonnegative entry time into $\cG$). Now, let us use the strong Markov property of the Feller process $X$, and for this purpose, let $\bar\theta$ denote the shift operator on the space-time path space $D\left([0,+\infty),\RR^d\times [0,+\infty)\right)$ (which is the set of all c\`adl\`ag functions from $[0,+\infty)$ into $\RR^d\times [0,+\infty)$ -- recall that all Feller processes have a c\`adl\`ag modification). We obtain \begin{eqnarray*}V_{F^\cI}^\cI(x,0) &=& \sup_{\cG\subset \RR^d\times [0,T]} \EE^{(x,0)}\left[\EE^{\left(X_{\bar\tau_{F^\cI}^\cI},{\bar\tau_{F^\cI}^\cI}\right)}e^{-r\bar\tau_{\cG}^\cI}g\left(X_{\bar\tau_\cG^\cI}\right) \right]\\ &=& \sup_{\cG\subset \RR^d\times [0,T]} \EE^{(x,0)}\left[\EE^{\left(X_{0},{0}\right)}\left[\left.e^{-r\bar\tau_{\cG}^\cI\circ \bar\theta_{\bar\tau_{F^\cI}^\cI}}g\left(X_{\bar\tau_\cG^\cI}\circ \bar\theta_{\bar\tau_{F^\cI}^\cI}\right)\right|\cF_{\bar\tau_{F^\cI}^\cI}\right] \right] \\ &=& \sup_{\cG\subset \RR^d\times [0,T]} \EE^{(x,0)}\left[e^{-r\bar\tau_\cG^\cI\circ \bar\theta_{\bar\tau_{F^\cI}^\cI}}g\left(X_{\bar\tau_\cG^\cI\circ \bar\theta_{\bar\tau_{F^\cI}^\cI}}\right) \right]\end{eqnarray*} 

But $e^{-r\cdot}g(X_{\cdot})$ is a $\PP^x$-supermartingale for all $x\in\RR^d$, therefore by Doob's Optional Stopping Theorem, $\left(e^{-r\upsilon}g\left(X_{\upsilon}\right)\right) _{\upsilon\in\left\{\bar\tau_\cG^\cI\wedge T',\bar\tau_\cG^\cI\circ \bar\tau_{F^\cI}^\cI\wedge T'\right\}}$ must also be a $\PP^x$-submartingale for all $x\in\RR^d$ and $T'\in(0,\infty)$ (note that $\bar\tau_\cG^\cI\leq \bar\tau_\cG^\cI\circ \bar\tau_{F^\cI}^\cI$ a.s. because of the fact that $\bar\theta$ is the shift operator for the space-time process $(t,X_t)_{t\geq 0}$, rather than simply for $X$). Letting $T'$ tend to infinity, we can employ Lebesgue's Dominated Convergence Theorem (as $g\leq K$ yields $e^{-r\upsilon}g\left(X_{\upsilon}\right)\leq K\in L^1(\PP^x)$ for $\upsilon\in\left\{\bar\tau_\cG^\cI\wedge T',\bar\tau_\cG^\cI\circ \bar\theta_{\bar\tau_{F^\cI}^\cI}\wedge T'\right\}$ for all $T'\in(0,+\infty)$ and $x\in\RR^d$) in order to get that the expected value of $e^{-r\bar\tau_\cG^\cI\circ \bar\theta_{\bar\tau_{F^\cI}^\cI}}g\left(X_{\bar\tau_\cG^\cI\circ \bar\theta_{\bar\tau_{F^\cI}^\cI}}\right)$ is always greater or equal than the expectation of $e^{-r\bar\tau_\cG^\cI}g\left(X_{\bar\tau_\cG^\cI}\right)$.
Hence \begin{eqnarray*}V_{F^\cI}^\cI(x,0) &=& \sup_{\cG\subset \RR^d\times [0,T]} \EE^{(x,0)}\left[e^{-r\bar\tau_\cG^\cI\circ \bar\theta_{\bar\tau_{F^\cI}^\cI}}g\left(X_{\bar\tau_\cG^\cI\circ \bar\theta_{\bar\tau_{F^\cI}^\cI}}\right) \right] \\ &\geq & \sup_{\cG\subset \RR^d\times [0,T]} \EE^{(x,0)}\left[e^{-r\bar\tau_\cG^\cI}g\left(X_{\bar\tau_\cG^\cI}\right) \right]\\ &=& \sup_{\cG\subset \RR^d\times [0,T]} \EE^{(x,0)}V_\cG^\cI(x,0).\end{eqnarray*} 

The case $T=+\infty$ can be dealt with analogously.

\end{proof}

\begin{lem}[Time-stationarity of immediate exercise regions for perpetual Bermudans]\label{time-stationarity} Let $X$ be a L\'evy basket with $\PP^\cdot$ being an associated family of probability measures and discount rate $r>0$. Then for all $s>0$ we have $$U^{s\NN_0}(+\infty)(x)=V_{\left\{x\in\RR^d\ : \ U^{s\NN_0}(+\infty)(x)\leq g(x)\right\}}^{s\NN_0}(x)$$ for all $x\in\RR^d$ satisfying the condition that the random variable \linebreak $\sup\left\{e^{-r\tau}g\left(X_\tau\right) \ : \ \tau\text{ stopping time}\right\}$ be $\PP^x$-integrable.
\end{lem}
\begin{proof} Consider an integer $n\in \NN_0$, and an $x\in\RR^d$ such that $\sup\left\{e^{-r\tau}g\left(X_\tau\right) \ : \ \tau\text{ stopping time}\right\}$ is $\PP^x$-integrable. Then we shift the time scale by $ns$ to get\begin{eqnarray*}&& e^{rns}\EE^{(x,ns)} \left[\sup_{\cG\subset \RR^d\times [0,+\infty)}e^{-r\tau_\cG^{s}} g\left(X_{\tau_\cG^{s}}\right)\right] \\&=& e^{rns}\EE^{(x,ns)} \left[\sup_{\cG\subset \RR^d\times [ns,+\infty)}e^{-r\tau_\cG^{s}} g\left(X_{\tau_\cG^{s}}\right)\right] \\&=& e^{rns}\EE^{(x,ns)} \left[\sup_{\cG'\subset \RR^d\times [0,+\infty)}e^{-r\left(\tau_{\cG'}^{s}\circ\theta_{ns}+ns\right)} g\left(X_{\tau_{\cG'}^{s}}{\circ\theta_{ns}}\right)\right] \\&=& \EE^{(x,0)} \left[\sup_{\cG'\subset \RR^d\times [0,+\infty)}e^{-r\tau_{\cG'}^{s}} g\left(X_{\tau_{\cG'}^{s}}\right)\right] 
\end{eqnarray*} where $\theta$ denotes the shift operator on the space (as opposed to space-time) path space $D\left([0,+\infty),\RR^d\right)$. Because of the boundedness of $g\leq K$ which entitles us to apply Lebegue's Dominated Convergence Theorem, we may swap $\sup$ and $\EE$ to obtain \begin{eqnarray*} && e^{rns}\EE^{(x,ns)} \left[\sup_{\cG\subset \RR^d\times [0,+\infty)}e^{-r\tau_\cG^{s}} g\left(X_{\tau_\cG^{s}}\right)\right]\\ &=& \sup_{\cG\subset \RR^d\times [0,+\infty)}e^{rns}\EE^{(x,ns)} \left[e^{-r\tau_\cG^{s}} g\left(X_{\tau_\cG^{s}}\right)\right]\\ &=& \sup_{\cG\subset \RR^d\times [0,+\infty)}V_\cG^{s\NN_0}(x,ns)\end{eqnarray*} for all $n\in\NN_0$. Thus we conclude $$\sup_{\cG\subset \RR^d\times [0,+\infty)}V_\cG^{s\NN_0}(x,t)=\sup_{\cG\subset \RR^d\times [0,+\infty)}V_\cG^{s\NN_0}(x,0)$$ for all $t\in s\NN_0$. If we insert this equality fact into the definition of $F^{s\NN_0}$, we see that the condition determining whether a pair $(x,t)$ belongs to $F^{s\NN_0}$ does not depend on $t$. On the other hand, by Corollary \ref{exerciseboundaryCorollary}, $$\sup_{\cG\subset \RR^d\times [0,+\infty)}V_\cG^{s\NN_0}(x,0)=U^{s\NN_0}(+\infty)(x),$$ and the left hand side equals -- by our previous observations in this proof -- the term featuring in the definition of $F^{s\NN_0}$.

\end{proof}

Summarising the two previous Lemmas and applying them to a more concrete setting, we deduce that the expected payoff of a perpetual Bermudan option of mesh size $h>0$ equals $$U^{h\NN_0}(+\infty)(x)=\left\{\begin{array}{*{2}{c}} \EE^x e^{-r\tau_G}g\left(X_{\tau_G}\right), & x\not\in G \\ g(x), & x\in G \end{array}\right. ,$$ where $G=\left\{U^{h\NN_0}(+\infty)(\cdot)\leq g(\cdot)\right\}$.

\begin{lem} \label{convolveproject}Let us fix a L\'evy basket with an associated family of risk-neutral probability measures $\PP^\cdot$ and discount rate $r>0$, as well as a region $G\subset \RR^d$ and a real number $h>0$. Then we have $$\forall x\not\in G \quad V_G^h(x)=e^{-rh}\PP_{X_0-X_h}\ast V_G^h(x).$$ In particular, using Lemma \ref{time-stationarity}, one has the following equation for the expected perpetual Bermudan option payoff: $$\forall x\not\in G \quad U^{h\NN_0}(+\infty)(x)=e^{-rh}\PP_{X_0-X_h}\ast U^{h\NN_0}(+\infty)(x).$$
\end{lem}
\begin{proof} Using the Markov property of $X$, denoting by $\theta$ the shift operator on the path space $D\left([0,+\infty),\RR^d\right)$ of a L\'evy process $X$, and taking into account the fact that $\tau_G^h>0$ (i.e. $\tau_G^h\geq h$) in case $x\not\in G$, we obtain:
\begin{eqnarray*}\forall x\not\in G \quad V_G^h(x) &=& e^{-rh}\EE^x \EE^{x}\left[ e^{-r\left(\tau_G^h-h\right)}g\left(X_{\tau_G^h}\right) \cdot \chi_{\{\tau_G^h\geq h\}}\right]  \\ &=& e^{-rh}\EE^x \EE^{x}\left[ e^{-r\tau_G^h\circ \theta_h}g\left(X_{\tau_G^h}\circ\theta_h \right) \cdot \chi_{\{\tau_G^h\geq h\}}\right] \\ &=& e^{-rh} \EE ^x\EE^{x}\left[ e^{-r\tau_G^h\circ \theta_h}g\left(X_{\tau_G^h}\circ\theta_h \right) \left( \chi_{\{\tau_G^h\geq h\}} + \chi_{\{\tau_G^h<h\}}\right)\right] \\ &=& e^{-rh} \EE^x \EE^x \left[ e^{-r\tau_G^h \circ \theta_h}g\left(X_{\tau_G^h}\circ\theta_h\right) | \cF_{h} \right]\\ &=& e^{-rh} \EE^x \EE^{X_h} e^{-r\tau_G^h}g\left(X_{\tau_G^h}\right) = e^{-rh}\EE ^xV_G^h\left(X_h\right) \\ &=& e^{-rh} \PP_{x-X_h}\ast V_G^h(x) \\ &=& e^{-rh} \PP_{X_0-X_h}\ast V_G^h(x). \end{eqnarray*} 
\end{proof}

\part{Bounds on the American-Bermudan barrier option price difference}

\chapter{Scaling the difference between perpetual American and Bermudan barrier options}

Embracing the terminology of Broadie, Glasserman and Kou \cite{BGK}, we shall refer to the difference between an American and the corresponding Bermudan options (on the same basket and with the same payoff function) as ``continuity correction''.

\section{The exercise boundary and its relevance for continuity corrections}

\begin{lem} \label{phaselemma1} Consider a $d$-dimensional Feller basket $X$ with an associated family of risk-neutral probability measures $\PP^\cdot$ and discount rate $r>0$. Furthermore, let $f:\RR^d\rightarrow \RR_+$ be a nonnegative continuous function such that $e^{-r\cdot}f(X_\cdot)$ is a martingale, $K$ a  nonnegative real number, and define $g:=K-f$. Finally, let $G$ be a measurable set such that $g\geq 0$ on $G$. Then we have for all $x\in {\RR}^d$ and $s>0$, $$V_{G}^s(x)= K\cdot \sum_{n=1}^\infty  e^{-rns} \PP^{x} \left[ \bigcap_{i=1}^{n-1}\left\{ X_{is}\not\in G \right\}\cap\left\{X_{ns}\in G \right\}\right] -f(x) .$$ Similarly, if $g:=f-K$ instead and one assumes that this $g$ is nonnegative on $G$, the identity $$V_{G}^s(x)= f(x)-K\cdot \sum_{n=1}^\infty  e^{-rns} \PP^{x} \left[ \bigcap_{i=1}^{n-1}\left\{ X_{is}\not\in G \right\}\cap\left\{X_{ns}\in G \right\}\right] $$ holds for all $x\in\RR^d$ and $s>0$.
\end{lem}
\begin{proof} Let $x\in\RR^d$ and $s>0$, and let us first set consider the case of $g=K-f$. Since $g\geq 0$ on $G$ by assumption, one has the identity $$V_G^s(x)=\EE^x\left[e^{-r\tau_G^s}g\left(X_{\tau_G^s}\right)\vee 0\right]=\EE^x\left[e^{-r\tau_G^s}(K-f)\left(X_{\tau_G^s}\right)\right].$$ Moreover, $$\lim_{n\rightarrow\infty }\underbrace{e^{-r\left(\tau_G^s\wedge ns\right)}(K-f)\left(X_{\tau_G^s\wedge ns}\right)}_{\leq K}= e^{-r\tau_G^s}(K-f)\left(X_{\tau_G^s}\right)$$ $\PP^x$-almost surely, therefore by Lebesgue's Dominated Convergence Theorem, \begin{equation} \label{phaselemma1 equation}V_G^s(x)=\lim_{n\rightarrow\infty }\EE^x\left[e^{-r\left(\tau_G^s\wedge ns\right)}(K-f)\left(X_{\tau_G^s\wedge ns}\right)\right].\end{equation} But since $e^{-r\cdot}f\left(X_\cdot\right)=\left(e^{-rt}f\left(X_t\right)\right)_{t\geq 0}$ is a $\PP^x$-martingale, we may apply Doob's Optional Stopping Theorem to get that $\left(e^{-r \left({\tau_G^s\wedge ns}\right)}f\left(X_{\tau_G^s\wedge ns}\right)\right)_{n\in \NN_0}$ is a $\PP^x$-martingale, too, whence $$\forall n\in\NN_0 \quad \EE^x\left[e^{-r\left(\tau_G^s\wedge ns\right)}(K-f)\left(X_{\tau_G^s\wedge ns}\right)\right]=K\cdot \EE^x\left[e^{-r\left(\tau_G^s\wedge ns\right)} \right]-f(x).$$ This finally yields, because of equation (\ref{phaselemma1 equation}) and the monotonicity of the sequence $\left(e^{-r\left(\tau_G^s\wedge ns\right)}\right)_{n\in\NN_0}$, $$V_G^s(x)=K\cdot \lim_{n\rightarrow\infty}\EE^x\left[e^{-r\left(\tau_G^s\wedge ns\right)}\right]-f(x)=K \EE^x\left[e^{-r\tau_G^s}\right]-f(x)$$ which completes the proof as $$\EE^x\left[e^{-r\tau_G^s}\right]=\sum_{n=1}^\infty  e^{-rns} \PP^{x} \left[ \bigcap_{i=1}^{n-1}\left\{ X_{is}\not\in G \right\}\cap\left\{X_{ns}\in G \right\}\right]$$ for all $n\in\NN_0$.

The case of $g=f-K$ can be treated analogously.
\end{proof}

An important feature of the immediate exercise region for a perpetual Bermudan option with payoff function of the form $(K-f)\vee 0$ where $f$ is monotonely increasing in each component, and exercise mesh size $h$, is that 
-- owing to the fact that $U^h(+\infty)(\cdot)$, the option price as function of the logarithmic start price, is monotonely decreasing in each component -- it is south-west connected in the following sense:

\begin{Def} A set $E\subseteq \RR^d$ is called {\em north-east connected} if and only if for all $x\in\RR^d$ such that $x\geq 0$ componentwise, $E+x\subseteq E$. Likewise, any set $F\subseteq \RR^d$ is called {\em south-west connected} if and only if for all $x\in\RR^d$ such that $x\geq 0$ componentwise, $F+x\supseteq F$. 
\end{Def}

\begin{rem} If $F$ is a south-west connected subset of $\RR^d$, then $\sup F:=\left(\sup_{x\in F} x_i\right)_{i=1}^d$ is an element of the boundary of $F$. Analogously, if $E$ is a north-east connected subset of $\RR^d$, then $\inf E:=\left(\inf_{x\in E} x_i\right)_{i=1}^d$ is an element of the boundary of $E$.
\end{rem}

\begin{lem}[Characterisation of the American-Bermudan barrier difference for perpetual puts] \label{characteriseputdifference} 
Let $X$ be the logarithmic price process of the multidimensional Black-Scholes model with constant volatility and interest rate, that is $$\forall t\geq 0 \quad X_t=\left(\left(X_0\right)_i+\sigma_i\cdot (B_t)_i + \left(r-\frac{1}{2}{\sigma_i}^2\right) t\right)_{i=1}^d$$ (where $B$ is the $d$-dimensional Wiener process) for some $r>0$ and $\sigma\in{\RR_{>0}}^d$.
Let $g=K-f$, wherein $K\geq 0$ be a real number and $f\geq 0$ be a continuous function that is monotonely increasing in each component and such that $\left(e^{-rt}f(X_t)\right)_{t\geq 0}$ is a martingale.
Finally, consider a measurable set of the shape $G=\gamma-H$ for some $\gamma\in\RR^d$ and some convex north-east connected set $H\subseteq {\RR_+}^d$ (making $G$ convex and south-west connected) such that $g=K-f$ is nonnegative on $G$. 
Then we have for all $s>0$, $$ V_{G}^{\frac{s}{2}}-V_{G}^s = K\cdot\EE^\cdot\left[e^{-r\tau_{G}^\frac{s}{2}} - e^{-r\tau_{G}^s}\right] \text{ on }\complement G.$$ 
\end{lem}

\begin{proof}[Proof] Fix an $s>0$. Introduce stopping times $\tau_{G,n}^t$ for $n\in\NN$, $t>0$ through $$\forall n\in\NN \quad \forall t>0\quad \tau_{G,n}^t:=\tau_G^t\wedge ns$$ and define $$\forall n\in\NN \quad \forall t>0\quad V_{G,n}^t:=\EE^\cdot \left[e^{-r\tau_{G,n}^t}g\left(X_{\tau_{G,n}^t}\right)\right]\quad \text{ on }\complement G.$$ 
Now consider an arbitrary $n\in\NN$ and $t>0$. Due to Doob's Optional Stopping Theorem, applied to the two-component sequences of stopping times $\left(0,\tau_{G,n}^\frac{s}{2}\right)$ and $\left(0,\tau_{G,n}^\frac{s}{2}\right)$ (which both are bounded by $ns$), combined with the fact that $e^{-r\cdot}f(X_\cdot)$ is a martingale, 
\begin{eqnarray} \EE^\cdot \left[e^{-r\tau_{G,n}^{s/2}}f\left(X_{\tau_{G,n}^{s/2}}\right)\right] \nonumber &=& \EE^\cdot \left[e^{-r\cdot 0}f\left(X_{0}\right)\right]=\EE^\cdot \left[e^{-r\tau_{G,n}^s}f\left(X_{\tau_{G,n}^s}\right)\right] \label{isotone lemma stopping},\end{eqnarray} 
thus as $g=K-f\geq 0$ on $G$ and $\tau_{G}^{s/2}\leq \tau_{G}^{s}$, 
\begin{eqnarray*} && V_{G,n}^{\frac{s}{2}}-V_{G,n}^s \\ 
& =& \EE^\cdot\left[\left. e^{-r\tau_{G,n}^\frac{s}{2}}\left(K-f\right)\left(X_{\tau_{G,n}^{s/2}}\right)-e^{-r\tau_{G,n}^s}\left(K-f\right)\left(X_{\tau_{G,n}^{s}}\right)\right. ,\left\{\tau_{G}^{s}< ns\right\}\right]\\ 
&& + \EE^\cdot\left[\left. e^{-r\tau_{G,n}^\frac{s}{2}}\left(K-f\right)\left(X_{\tau_{G,n}^{s/2}}\right)\vee 0-e^{-r\tau_{G,n}^s}\left(K-f\right)\left(X_{\tau_{G,n}^{s}}\right)\vee 0\right. ,\left\{\tau_{G}^{s/2},\tau_{G}^{s}\geq ns\right\}\right]\\ 
&& + \EE^\cdot\left[\left. e^{-r\tau_{G,n}^\frac{s}{2}}\left(K-f\right)\left(X_{\tau_{G,n}^{s/2}}\right)\vee 0-e^{-r\tau_{G,n}^s}\left(K-f\right)\left(X_{\tau_{G,n}^{s}}\right)\vee 0\right. ,\left\{\tau_{G}^{s/2}\leq ns\leq \tau_{G}^{s} \right\} \right]\\ 
& =& \EE^\cdot\left[\left. e^{-rt}\left(K-f\right)\left(X_{t}\right)\right|_{t=\tau_{G,n}^s=\tau_{G}^s}^{t=\tau_{G,n}^\frac{s}{2}=\tau_{G}^\frac{s}{2}} ,\left\{\tau_{G,n}^\frac{s}{2}\leq \tau_{G}^{s}< ns\right\}\right]\\ 
&& + \underbrace{\EE^\cdot\left[\left. e^{-rt}\left(K-f\right)\left(X_{t}\right)\vee 0\right|_{t=\tau_{G,n}^s}^{t=\tau_{G,n}^\frac{s}{2}} ,\left\{\tau_{G,n}^\frac{s}{2},\tau_{G}^{s}\geq ns\right\}\right] }_{=0}\\ 
&& + \EE^\cdot\left[\left. e^{-rt}\left(K-f\right)\left(X_{t}\right)\vee 0\right|_{t=\tau_{G,n}^s}^{t=\tau_{G,n}^\frac{s}{2}} ,\left\{\tau_{G,n}^\frac{s}{2}<ns \leq \tau_{G}^{s}\right\}\right]\\ 
& =& \EE^\cdot\left[\left. e^{-rt}\left(K-f\right)\left(X_{t}\right)\right|_{t=\tau_{G}^s}^{t=\tau_{G}^\frac{s}{2}} ,\left\{\tau_{G,n}^\frac{s}{2}\leq \tau_{G}^{s}< ns\right\}\right]\\ 
&& + \underbrace{\EE^\cdot\left[\left. e^{-rt}\left(K-f\right)\left(X_{t}\right)\right|_{t=\tau_{G,n}^s=ns}^{t=\tau_{G,n}^\frac{s}{2}=ns} ,\left\{\tau_{G,n}^\frac{s}{2},\tau_{G}^{s}\geq ns\right\}\right] }_{=0}\\ 
&& + \EE^\cdot\left[\left. e^{-rt}\left(K-f\right)\left(X_{t}\right)\right|_{t=\tau_{G,n}^s}^{t=\tau_{G,n}^\frac{s}{2}} ,\left\{\tau_{G,n}^\frac{s}{2}<ns \leq \tau_{G}^{s}\right\}\right]\\ 
&& - \EE^\cdot\left[\left. e^{-rns}\left(K-f\right)\left(X_{ns}\right)\wedge 0\right.,\left\{\tau_{G,n}^\frac{s}{2}<ns \leq \tau_{G}^{s}\right\}\right]\\ 
& =& \EE^\cdot\left[\left. e^{-rt}\left(K-f\right)\left(X_{t}\right)\right|_{t=\tau_{G,n}^s}^{t=\tau_{G,n}^\frac{s}{2}} \right]\\ 
&& - \EE^\cdot\left[\left. e^{-rns}\left(K-f\right)\left(X_{ns}\right)\wedge 0\right.,\left\{\tau_{G,n}^\frac{s}{2}<ns \leq \tau_{G}^{s}\right\}\right]
\end{eqnarray*}

Note, however that the subtractor converges to zero exponentially. For, exploiting not only the Markov property, but also $f\geq 0$ and that $e^{-r\cdot f(X_\cdot)}$ is a martingale, we have 

\begin{eqnarray*}0&\geq& \EE^\cdot\left[\left. e^{-rns}\left(K-f\right)\left(X_{ns}\right)\wedge 0\right.,\left\{\tau_{G,n}^\frac{s}{2}<ns \leq \tau_{G}^{s}\right\}\right]\\ 
&=& \EE^\cdot\left[e^{-rns}\left(K-f\right)\left(X_{ns}\right)\wedge 0,\left\{\tau_{G,n}^\frac{s}{2}\leq (n-1)s\leq ns\leq \tau_{G}^{s}\right\}\right] \\ 
&&+ \EE^\cdot\left[e^{-rns}\left(K-f\right)\left(X_{ns}\right)\wedge 0,\left\{\tau_{G,n}^\frac{s}{2} =\left(n-\frac{1}{2}\right)s\leq ns \leq \tau_{G}^{s}\right\}\right]\\
&=& \EE^\cdot\left[\left. e^{-rns}\left(K-f\right)\left(X_{ns}\right)\wedge 0\right|\left\{X_{\left(n-1  \right)s}\not\in G\right\}\right]\PP^\cdot\left\{\tau_{G,n}^\frac{s}{2} \leq\left(n-1\right)s\leq ns \leq \tau_{G}^{s}\right\}\\
&&+ \EE^\cdot\left[\left. e^{-rns}\left(K-f\right)\left(X_{ns}\right)\wedge 0\right|\left\{X_{\left(n-\frac{1}{2}\right)s}\in G\right\}\right]\PP^\cdot\left\{\tau_{G,n}^\frac{s}{2} =\left(n-\frac{1}{2}\right)s\leq ns \leq \tau_{G}^{s}\right\}\\
&\geq & \EE^\cdot\left[\left. e^{-rns}\left(K-f\right)\left(X_{ns}\right) \wedge 0 \right|\left\{X_{\left(n-1  \right)s}\not\in G\right\}\right]\PP^\cdot\left[\left\{X_{\left(n-1\right)s}\not\in G\right\}\cap \bigcap_{i=1}^{n-1} \left\{X_{is}\not\in G\right\} \right]\\
&&+ \EE^\cdot\left[\left. e^{-rns}\left(K-f\right)\left(X_{ns}\right) \wedge 0 \right|\left\{X_{\left(n-\frac{1}{2}\right)s\in G}\right\}\right]\PP^\cdot\left[\left\{X_{\left(n-\frac{1}{2}\right)s} \in G\right\}\cap \bigcap_{i=1}^{n-1} \left\{X_{is}\not\in G\right\} \right]\\
&\geq & \EE^\cdot\left[\left. e^{-rns}\left(-f\right)\left(X_{ns}\right) \right|\left\{X_{\left(n-1  \right)s}\not\in G\right\}\right]\PP^\cdot\left[\left\{X_{\left(n-1  \right)s}\not\in G\right\}\cap \bigcap_{i=1}^{n-1} \left\{X_{is}\not\in G\right\} \right]\\
&&+ \EE^\cdot\left[\left. e^{-rns}\left(-f\right)\left(X_{ns}\right) \right|\left\{X_{\left(n-\frac{1}{2}\right)s\in G}\right\}\right]\PP^\cdot\left[\left\{X_{\left(n-\frac{1}{2}\right)s} \in G\right\}\cap \bigcap_{i=1}^{n-1} \left\{X_{is}\not\in G\right\} \right]\\
&\geq & \EE^\cdot\left[\left. e^{-rns}\left(-f\right)\left(X_{ns}\right) \right.,\left\{X_{\left(n-1  \right)s}\not\in G\right\}\right]\PP^\cdot\left[\left. \bigcap_{i=1}^{n-1} \left\{X_{is}\not\in G\right\} \right|\left\{X_{\left(n-1\right)s}\not\in G\right\}\right]\\
&&+ \EE^\cdot\left[\left. e^{-rns}\left(-f\right)\left(X_{ns}\right) \right.,\left\{X_{\left(n-\frac{1}{2}\right)s\in G}\right\}\right]\PP^\cdot\left[\left. \bigcap_{i=1}^{n-1} \left\{X_{is}\not\in G\right\} \right|\left\{X_{\left(n-\frac{1}{2}\right)s}\in G\right\}\right]\\
&\geq & \EE^\cdot\left[\left. e^{-rns}\left(-f\right)\left(X_{ns}\right) \right.\right]\PP^\cdot\left[ \left. \bigcap_{i=1}^{n-1} \left\{X_{is}\not\in G\right\} \right|\left\{X_{\left(n-1\right)s}\not\in G\right\}\right]\\
&&+ \EE^\cdot\left[\left. e^{-rns}\left(-f\right)\left(X_{ns}\right) \right.\right]\PP^\cdot\left[ \left. \bigcap_{i=1}^{n-1} \left\{X_{is}\not\in G\right\} \right|\left\{X_{\left(n-\frac{1}{2}\right)s}\in G\right\}\right]\\
&=& \EE^\cdot\left[\left.\left(-f\right)\left(X_0\right) \right.\right]\PP^\cdot\left[ \left. \bigcap_{i=1}^{n-1} \left\{X_{is}\not\in G\right\} \right|\left\{X_{\left(n-1\right)s}\not\in G\right\}\right]\\
&&+ \EE^\cdot\left[\left.\left(-f\right)\left(X_0\right) \right.\right]\PP^\cdot\left[ \left. \bigcap_{i=1}^{n-1} \left\{X_{is}\not\in G\right\} \right|\left\{X_{\left(n-\frac{1}{2}\right)s}\in G\right\}\right]\\
&=& -f(\cdot)\left(\PP^\cdot\left[ \left. \bigcap_{i=1}^{n-1} \left\{X_{is}\not\in G\right\} \right|\left\{X_{\left(n-1\right)s}\not\in G\right\}\right]+ \PP^\cdot\left[ \left. \bigcap_{i=1}^{n-1} \left\{X_{is}\not\in G\right\} \right|\left\{X_{\left(n-\frac{1}{2}\right)s}\in G\right\}\right]\right)
\end{eqnarray*}

on $\complement G$ -- and both of the conditioned probabilities in the last line converge to nought as $n$ tends to $+\infty$.

Therefore, due to equation (\ref{isotone lemma stopping}), 
\begin{eqnarray}&& \lim_{n\rightarrow\infty} \left(V_{G,n}^{\frac{s}{2}}-V_{G,n}^s\right)\nonumber \\ 
&=&\lim_{n\rightarrow\infty} \EE^\cdot\left[\left. e^{-rt}\left(K-f\right)\left(X_{t}\right)\right|_{t=\tau_{G,n}^s}^{t=\tau_{G,n}^\frac{s}{2}} \right]\nonumber \\ 
&=&\lim_{n\rightarrow\infty} K\cdot\EE^\cdot\left[\left. e^{-rt}\right|_{t=\tau_{G,n}^s}^{t=\tau_{G,n}^\frac{s}{2}} \right]\label{isotone lemma limiting difference}
\end{eqnarray} on $\complement G$.

However, we also know that
\begin{eqnarray*} && \left. e^{-rt}\left(K-f\right)\left(X_{t}\right)\right|_{t=\tau_{G,n}^s}^{t=\tau_{G, n}^\frac{s}{2}}\\
&\rightarrow& \left. e^{-rt}\left(K-f\right)\left(X_{t}\right)\right|_{t=\tau_{G}^s}^{t=\tau_{G}^\frac{s}{2}} \text{ as }n\rightarrow\infty\end{eqnarray*} almost surely. 
Combining this convergence assertion with Lebesgue's Dominated Convergence Theorem yields 
$$V_{G}^{\frac{s}{2}}-V_{G}^s= \lim_{n\rightarrow\infty} \left(V_{G,n}^{\frac{s}{2}}-V_{G,n}^s\right) ,$$ 
and recalling equation (\ref{isotone lemma limiting difference}), we finally arrive at 
\begin{equation}V_{G}^{\frac{s}{2}}-V_{G}^s= \lim_{n\rightarrow\infty} K\cdot\EE^\cdot\left[\left. e^{-rt}\right|_{t=\tau_{G,n}^s}^{t=\tau_{G,n}^\frac{s}{2}} \right]\text{ on }\complement G. \label{differencecharacterisation} \end{equation}

\end{proof}

Next we shall prove that the expression on the right hand side of the last equation is monotonely increasing as $\cdot\downarrow\gamma$ componentwise. 

\begin{lem} \label{putdifferenceboundary} 
Under the assumptions of Lemma \ref{characteriseputdifference}, one has for all $s>0$ and $x,y\in\RR^d$ that satisfy the relation $$x\geq y\geq\gamma \text{ componentwise}$$ the following lower bound: \begin{eqnarray*}\left(V_G^{\frac{s}{2}}-V_G^s\right)(y) &\geq & K\sum_{i=0}^\infty \PP^y \left\{\tau_G^\frac{s}{2}={\left(i+\frac{1}{2}\right)s} \right\} e^{-r\left(i+\frac{1}{2}\right)s}\\ &&\cdot\left(1- \EE^{x}\left[\left. e^{-r\tau_G^{\left(\frac{1}{2}+\NN_0\right)s}} \right| \left\{X_{\left(i+\frac{1}{2}\right)s}\in G\right\}\right]\right) 
\end{eqnarray*}
\end{lem}

In the proof of Lemma \ref{putdifferenceboundary}, we will apply the following

\begin{aux} \label{aux isotone Lemma}Consider a measurable set $A\subseteq\RR$ of positive Lebesgue measure and two continuous positive integrable functions $p,q\in L^1\left(A\right)$ on $A$, as well as a nonnegative monotonely decreasing function $f\in L^1\left(p \ d\lambda\right)\cap L^1\left(q \ d\lambda\right)$, and assume that not only is the function $\frac{p}{q}$ monotonely decreasing, but also $\int_A p \ d\lambda =\int_A q \ d\lambda $. Then $$\int_Afp \ d\lambda\geq \int_Afq \ d\lambda.$$
\end{aux}
\begin{proof} If the continuous monotonely decreasing function $p-q$ was either $>0$ or $<0$ on all of $A$, then one would get $\int_A p \ d\lambda >\int_A q \ d\lambda $ or $\int_A p \ d\lambda <\int_A q \ d\lambda $, respectively. Hence there exists a real number $b\in A$ such that $p(a)=q(a)$ as well as $p\geq q$ on $A\subset(-\infty,a)$ and $p\leq q$ on $A\subset(a,+\infty)$. But since $f$ is monotonely decreasing and nonnegative, these inequalities yield \begin{eqnarray*}&&\int_{A\cap (-\infty,a)}f\cdot(p-q)d\lambda+ \int_{A\cap(a,+\infty)} f\cdot(p-q)d\lambda \\ &\geq &\int_{A\cap (-\infty,a)} f(a)(p-q)d\lambda+ \int_{A\cap(a,+\infty)} f(a)(p-q)d\lambda \geq 0,\end{eqnarray*} thus $$\int_Afp\ d\lambda-\int_Afq \ d\lambda=\int_Af\cdot(p-q)d\lambda\geq 0.$$
\end{proof}

\begin{proof}[Proof of Lemma \ref{putdifferenceboundary}] Let $x\geq y\geq\gamma$, which in particular entails $x\geq y\geq \sup G$. 
Let us define, for all $z\in \RR^d$ and $i\in\NN_0$, a measure $Q^{(z,i)}$ on $G$ by $$Q^{(z,i)}:= \frac{\PP^z_{ X_{ \left(i+\frac{1}{2}\right)s} }\left[\cdot\cap G\right] }{\PP^z\left\{X_{ \left(i+\frac{1}{2}\right)s}\in G\right\}}.$$ 
Then, due to our choice of $X$ as being a Brownian motion with drift, $Q^{(z,i)}$ will have a positive continuous Lebesgue density, denoted by $\frac{dQ^{(z,i)}}{d\lambda^d}$.

Now whenever $z_0,z_1\in G$ such that \begin{equation}\exists\alpha>1\quad \left(\sup G- z_0\right)=\alpha\cdot \left(\sup G- z_1\right)\geq 0 \label{z0z1direction} \end{equation} (in particular $z_0\leq z_1$), note the identity $$\frac{dQ^{(z,i)}}{d\lambda^d}=\frac{\frac{d\PP^z_{ X_{ \left(i+\frac{1}{2}\right)s} }\left[\cdot\cap G\right] }{d\lambda^d}}{\PP^z\left\{X_{ \left(i+\frac{1}{2}\right)s}\in G\right\}}= \frac{ g_{0,\diag\left({\sigma_1}^2,\dots,{\sigma_d}^2\right)\left(i+\frac{1}{2}\right)s}\left(\cdot-z-\left(i+\frac{1}{2}\right)s\mu\right) }{\PP^z\left\{X_{ \left(i+\frac{1}{2}\right)s}\in G\right\}}\chi_G$$ (where $g_{\alpha,C}$ shall for every symmetric positive semidefinite $C\in\RR^{d\times d}$ and $\alpha\in\RR^d$ denote the density of $\nu_{\alpha,C}$, the Gaussian measure with covariance matrix $C$ and mean $\alpha$) for all $z\in\RR^d$. For, we shall then obtain \begin{eqnarray*}\frac{ \frac{dQ^{(z,i)}}{d\lambda^d} (z_0)}{ \frac{dQ^{(z,i)}}{d\lambda^d} (z_1) }&=&\frac{\prod_{j=1}^d e^{ -\frac{\left|z_0-z-\left(i+\frac{1}{2}\right)s\mu\right|^2}{2{\sigma_j}^2\left(i+\frac{1}{2}\right)s} } } {\prod_{j=1}^d e^{ -\frac{\left|z_1-z-\left(i+\frac{1}{2}\right)s\mu\right|^2}{2{\sigma_j}^2\left(i+\frac{1}{2}\right)s} } } = \prod_{j=1}^d e^{ -\frac{|z_0|^2-|z_1|^2-2\cdot{^t}\left( \left(i+\frac{1}{2}\right)s\mu+z\right)\left(z_0-z_1\right)}{2{\sigma_j}^2\left(i+\frac{1}{2}\right)s} }  \\ &=& \prod_{j=1}^d e^{ -\frac{|z_0|^2-|z_1|^2-2\cdot \left(i+\frac{1}{2}\right)s\cdot{^t}\mu\left(z_0-z_1\right)}{2{\sigma_j}^2\left(i+\frac{1}{2}\right)s} } \cdot\prod_{j=1}^d e^{ -\frac{2\cdot{^t}z\left(z_1-z_0\right)}{2{\sigma_j}^2\left(i+\frac{1}{2}\right)s} } \end{eqnarray*}
and therefore $$\frac{ \frac{dQ^{(y,i)}}{d\lambda^d} (z_0)}{ \frac{dQ^{(y,i)}}{d\lambda^d} (z_1) }\geq \frac{ \frac{dQ^{(x,i)}}{d\lambda^d} (z_0)}{ \frac{dQ^{(x,i)}}{d\lambda^d} (z_1) }.$$ This implies \begin{equation}\frac{ dQ^{(x,i)} }{d Q^{(y,i)} }(z_0)\leq \frac{d Q^{(x,i)} }{ dQ^{(y,i)} }(z_1)\end{equation} whence we have shown that $\frac{ dQ^{(x,i)} }{d Q^{(y,i)} }$ is monotonely increasing on the ray-segment $R(z_2):=\left(\sup G-z_2\cdot{\RR_+}^d\right)\cap G$.
But in addition, $$\EE^{z_0}\left[e^{-r\tau_G^{\left(\frac{1}{2}+\NN_0\right)s}}\right]\geq \EE^{z_1}\left[e^{-r\tau_G^{\left(\frac{1}{2}+\NN_0\right)s}}\right]$$ will hold (since $G$ was assumed to be convex and south-west connected).
Thus $\EE^{\cdot}\left[e^{-r\tau_G^{\left(\frac{1}{2}+\NN_0\right)s}}\right]$ decreases monotonely on any ray (or ray-segment) $R(z_2)=\left(\sup G-z_2\cdot\RR_+\right)\cap G$ (for $z_2\geq 0$ componentwise) if we look at this ray (-segment) as a linearly ordered set with respect to the componentwise real order relation $\leq$.
From these two sets of monotonicity assertions (for all $z_2\geq 0$ componentwise), together with the fact that $Q^{(x,i)}(G)=Q^{(y,i)}(G)=1$, one deduces by means of the above Auxiliary Lemma \ref{aux isotone Lemma}, \begin{eqnarray*} \int_{R(z_2)} \EE^{z}\left[e^{-r\tau_G^{\left(\frac{1}{2}+\NN_0\right)s}}\right] \ Q^{(x,i)}(dz) &\leq& \int_{R(z_2)} \EE^{z}\left[e^{-r\tau_G^{\left(\frac{1}{2}+\NN_0\right)s}}\right] \ Q^{(y,i)}(dz) ,\end{eqnarray*} therefore \begin{eqnarray*} && \int_{\partial B_1(0)\cap{\RR_+}^d}\int_{R(z_2)} \EE^{z}\left[e^{-r\tau_G^{\left(\frac{1}{2}+\NN_0\right)s}}\right] \ Q^{(x,i)}(dz)  \frac{dz_2 }{\lambda\left[\partial B_1(0)\cap{\RR_+}^d\right]}\\ &\leq &\int_{\partial B_1(0)\cap{\RR_+}^d} \int_{R(z_2)} \EE^{z}\left[e^{-r\tau_G^{\left(\frac{1}{2}+\NN_0\right)s}}\right] \ Q^{(y,i)}(dz) \frac{dz_2 }{\lambda\left[\partial B_1(0)\cap{\RR_+}^d\right]}\end{eqnarray*} (where $B_1(0)$ denotes the $d$-dimensional unit ball) which via the Fubini-Tonelli Theorem amounts to 
\begin{equation}\int \EE^{z}\left[e^{-r\tau_G^{\left(\frac{1}{2}+\NN_0\right)s}}\right] \ Q^{(x,i)}(dz) \leq \int \EE^{z}\left[e^{-r\tau_G^{\left(\frac{1}{2}+\NN_0\right)s}}\right] \ Q^{(y,i)}(dz) \label{$Q^{(x,i)}$ integrals}.\end{equation} 
On the other hand, the definition of $Q^{(w,i)} $ gives
\begin{eqnarray*} \int \EE^{z}\left[e^{-r\tau_G^{\left(\frac{1}{2}+\NN_0\right)s}}\right] \ Q^{(w,i)}(dz) 
&=& \EE^{w}\left[\left. e^{-r\tau_G^{\left(\frac{1}{2}+\NN_0\right)s}}\right| \left\{X_{\left(i+\frac{1}{2}\right)s}\in G\right\}\right] \end{eqnarray*}
for arbitrary $w\in\RR^d$. Hence inequality (\ref{$Q^{(x,i)}$ integrals}) becomes 
\begin{equation} \EE^{x}\left[\left. e^{-r\tau_G^{s}} \right| \left\{X_{\left(i+\frac{1}{2}\right)s}\in G\right\}\right]\leq \EE^{y}\left[\left. e^{-r\tau_G^{s}} \right| \left\{X_{\left(i+\frac{1}{2}\right)s}\in G\right\}\right].\label{isotone lemma conditioned inequality} \end{equation}

The Markov property of $X$, moreover, entitles us to state 
$$\EE^w\left[\left. e^{-r\tau_{G}^s}\right| \left\{\tau_G^\frac{s}{2}={\left(i+\frac{1}{2}\right)s} \right\}\right] =e^{-r\left(i+\frac{1}{2}\right)s}\EE^{w}\left[\left. e^{-r\tau_G^{\left(\frac{1}{2}+\NN_0\right)s}} \right| \left\{X_{\left(i+\frac{1}{2}\right)s}\in G\right\}\right]$$ 
for all $i$ and $w\in\RR^d$. Thus we can use inequality (\ref{isotone lemma conditioned inequality}) to derive the following estimate: 
\begin{eqnarray} && \EE^y\left[e^{-r\tau_{G}^\frac{s}{2}} - e^{-r\tau_{G}^s}\right] 
\\&= & \sum_{i=0}^\infty \PP^y \left\{\tau_{G}^\frac{s}{2}={\left(i+\frac{1}{2}\right)s} \right\}\cdot \EE^y\left[\left. e^{-r\tau_{G}^\frac{s}{2}} - e^{-r\tau_{G}^s}\right| \left\{\tau_G^\frac{s}{2}={\left(i+\frac{1}{2}\right)s} \right\}\right] \nonumber
\\ & =& \sum_{i=0}^\infty \PP^y \left\{\tau_G^\frac{s}{2}={\left(i+\frac{1}{2}\right)s} \right\} e^{-r\left(i+\frac{1}{2}\right)s}  \nonumber\\ &&\cdot \left(1- \EE^{y}\left[\left. e^{-r\tau_G^{\left(\frac{1}{2}+\NN_0\right)s}} \right| \left\{X_{\left(i+\frac{1}{2}\right)s}\in G\right\}\right]\right) \label{isotone Lemma conditioned identity}
\\ & \geq &  \nonumber \sum_{i=0}^\infty \PP^y \left\{\tau_G^\frac{s}{2}={\left(i+\frac{1}{2}\right)s} \right\} e^{-r\left(i+\frac{1}{2}\right)s}\\ &&\cdot \left(1- \EE^{x}\left[\left. e^{-r\tau_G^{\left(\frac{1}{2}+\NN_0\right)s}} \right| \left\{X_{\left(i+\frac{1}{2}\right)s}\in G\right\}\right]\right)  \nonumber
\end{eqnarray}

which via Lemma \ref{characteriseputdifference} gives the desired result.

\end{proof}

\begin{cor} Suppose the assumptions of the preceding Lemma \ref{putdifferenceboundary} hold. Consider $x\geq y\geq \gamma$  and assume, in addition, there is a positive lower bound for the sequence $\left(\frac{\PP^y \left\{\tau_G^\frac{s}{2}={\left(i+\frac{1}{2}\right)s} \right\}}{\PP^x \left\{\tau_G^\frac{s}{2}={\left(i+\frac{1}{2}\right)s} \right\}}\right)_{i\in\NN_0}$. Then one has $$\left(V_G^{\frac{s}{2}}-V_G^s\right)(y)\geq \left(\inf_{i\in\NN_0}\frac{\PP^y \left\{\tau_G^\frac{s}{2}={\left(i+\frac{1}{2}\right)s} \right\}}{\PP^x \left\{\tau_G^\frac{s}{2}={\left(i+\frac{1}{2}\right)s} \right\}}\right)\cdot \left(V_G^{\frac{s}{2}}-V_G^s\right)(x).$$
\end{cor}
\begin{proof} Equation (\ref{isotone Lemma conditioned identity}) and Lemma \ref{characteriseputdifference} yield \begin{eqnarray} && V_G^\frac{s}{2}(z)-V_G^s(z)\nonumber
\\ & =& K\sum_{i=0}^\infty \PP^z \left\{\tau_G^\frac{s}{2}={\left(i+\frac{1}{2}\right)s} \right\} e^{-r\left(i+\frac{1}{2}\right)s}  \nonumber\\ &&\cdot \left(1- \EE^{z}\left[\left. e^{-r\tau_G^{\left(\frac{1}{2}+\NN_0\right)s}} \right| \left\{X_{\left(i+\frac{1}{2}\right)s}\in G\right\}\right]\right) \end{eqnarray} for every $z\not\in G$. Inserting $x$ and $y$ for $z$ in this estimate, the Corollary can be deduced via the estimate of the preceding Lemma \ref{putdifferenceboundary}.

\end{proof}

We can explicitly state a partial differential equation that the said difference $V_G^{s/2}-V_G^s$ obeys:

\begin{lem} \label{VG pde}Consider a measurable set $G\subset \RR^d$. Let again $X$ be the logarithmic price process of the multidimensional Black-Scholes model, that is, $$\forall t\geq 0 \quad X_t=\left(\left(X_0\right)_i+\sigma_i\cdot (B_t)_i + \underbrace{\left(r-\frac{1}{2}{\sigma_i}^2\right)}_{=:\mu_i} t\right)_{i=1}^d,$$ and let us assume that $g$ and $f$ are as in Lemma \ref{putdifferenceboundary} Then for all $t>s>0$, the partial differential equation $$\left(\frac{1}{2}\sum_{i=1}^d{\sigma_i}^2 \partial_{i}\partial_{i} -{^t}\mu\nabla\right)\left( V_G^s-V_G^t\right) = \left.\frac{\partial}{\partial u}\right|_{u=0}\left(V_G^{s,u}-V_G^{t,u}\right)$$ holds on $\complement G$, wherein for all $u,v>0$ $$V_G^{v,u}:= e^{ru} V_G^{v\NN_0+u}= e^{ru}\EE^\cdot\left[e^{-r\tau_G^{v\NN_0+u}}g\left(X_{\tau_G^{v\NN_0+u}}\right)\right],$$ and (always following the notation introduced in Chapter \ref{formalintro}) $$\tau_G^{v\NN_0+u}=\inf\left\{w\in v\NN_0+u\ : \ X_w\in G\right\}.$$
\end{lem}
\begin{proof} Fix $s>0$. We shall prove the Lemma by studying the space-time Markov process $(X_t,s-t)_{t\in[0,s]}$ and the functions $$f^s:\complement G\rightarrow\RR,\quad (x,t)\mapsto e^{rt}\EE^x\left[e^{-r\tau_G^s(s-t)}g\left(X_{\tau_G^s(s-t)}\right)\right],$$ where $$\forall t\in[0,s)\quad\tau_G^s(t):= \tau_G^{s\NN-t}= \inf\left\{u\in s\NN-t \ : X_u\in G\right\}.$$ It is clear that $$\forall u>t\in[0,s) \quad \tau_G^s(u)\circ\theta_{u-t}=\tau_G^s(t)-(u-t), \quad X_{\tau_G^s(u)}\circ\theta_{u-t}=X_{\tau_G^s(t)},$$ hence if we employ the Markov property in two directions we can for all $u>t$ get the following $\PP^\cdot$-almost sure identities: \begin{eqnarray*} &&\EE^\cdot\left[\left.\EE^{X_u}\left[e^{-r\tau_G^s(u)}g\left(X_{\tau_G^s(u)}\right)\right]\right| \cF_t\right]\\ &=& \EE^\cdot \left[\left.\EE^{X_{u-t}\circ\theta_t}\left[e^{-r\tau_G^s(u)}g\left(X_{\tau_G^s(u)}\right)\right]\right| \cF_t\right]\\ &=& \EE^{X_t}\EE^{X_{u-t}}\left[e^{-r\tau_G^s(u)}g\left(X_{\tau_G^s(u)}\right)\right]\\ &=& \EE^{X_t}\EE^{X_{0}}\left[\left.e^{-r\tau_G^s(u)\circ\theta_{u-t}}g\left(X_{\tau_G^s(u)}\circ\theta_{u-t}\right)\right|\cF_{u-t}\right]\\ &=& \EE^{X_t}\EE^{X_{0}}\left[\left.e^{-r(u-t)}e^{-r\tau_G^s(t)}g\left(X_{\tau_G^s(t)}\right)\right|\cF_{u-t}\right]\\ &=& e^{r(u-t)}\EE^{X_{t}}\left[e^{-r\tau_G^s(t)}g\left(X_{\tau_G^s(t)}\right)\right].\end{eqnarray*}
Thus, the process $\left(f^s\left(X_t,s-t\right)\right)_{t\in[0,s]}$ is a martingale. The infinitesimal generator of the space-time Markov process $(X_t,s-t)_{t\in[0,s]}$ is $$L:=-\frac{\partial}{\partial t}+\frac{1}{2} \sum_{i=1}^d{\sigma_i}^2 \partial_{i}\partial_{i}-{^t}\mu \nabla=-{\partial _{d+1}}+\frac{1}{2} \sum_{i=1}^d{\sigma_i}^2 \partial_{i}\partial_{i}-{^t}\mu \nabla$$ due to the well-known result on the infinitesimal generator of Brownian motion with drift (see eg Revuz and Yor \cite[p. 352]{RY}). Therefore we have proven $$\forall s>0 \quad 0=Lf^s =-\frac{\partial}{\partial t}f^s + \left(\frac{1}{2} \sum_{i=1}^d{\sigma_i}^2 \partial_{i}\partial_{i}-{^t}\mu \nabla \right)f^s.$$ 

Now one observes that $$\forall s>0 \quad f^s(\cdot,0) = V_G^s \text{ on }\complement G $$ which yields 
\begin{eqnarray*} \forall s>0 \quad && \left( \frac{1}{2} \sum_{i=1}^d{\sigma_i}^2 \partial_{i}\partial_{i}-{^t}\mu \nabla\right)V_G^s \\ &= & \left( \frac{1}{2} \sum_{i=1}^d{\sigma_i}^2 \partial_{i}\partial_{i}-{^t}\mu \nabla\right)f^s(cdot,0) \\ &= & \partial_{d+1}f^s(\cdot,0) = \left.\frac{\partial}{\partial u}\right|_{u=0}V_G^{s,u}
\end{eqnarray*}
and thus brings the proof of the Lemma to a close.
\end{proof}

\begin{rem}
In the remainder of this Chapter, we will continue to largely focus on put options, thus always setting $g=K-f$ for some componentwise monotonely increasing $f$ and assuming the $G$ occurring in the defition of $V_G^s$, $s>0$, to be south-west connected. However, one can easily derive analogous results for call options on dividend-paying assets, by simply cutting the interest rate to discount the dividends and by replacing $g$ by $-g$ and $G$ by some north-east connected measurable subset of $\RR^d$ that is assumed to satisfy the condition $f-K\geq 0$ on $G$.
\end{rem}

\section{Continuity corrections in a one-dimensional setting}

\begin{Def} For measurable $G\subseteq \RR^d$ and $s>0$ define $$\bar V_G^s:=\EE^\cdot\left[e^{-r\tau_G^{s\NN}}g\left(X_{\tau_G^{s\NN}}\right)\right].$$
\end{Def}

\begin{rem}Note that in general, $\bar V_G^s\neq g = V_G^s $ on $G$, but always $\bar V_G^s=V_G^s$ on $\complement G$.
\end{rem}

\begin{Th}\label{phase1d} Suppose $d=1$, let $G=(-\infty,\gamma]$ and $g=K-\exp$, and assume $(X_t)_{t\geq 0}=\left(X_0+\sigma\cdot B_t+\left(r-\frac{\sigma^2}{2}\right)\right)_{t\geq 0}$, in words: $X$ is the logarithmic price process of the one-dimensional Black-Scholes model with constant volatility $\sigma$ and discount rate $r>0$. Set $\mu:=r-\frac{\sigma^2}{2}$. Then one has for all $s>0$ the relations \begin{eqnarray} && \lim_{t\downarrow 0}\bar V_{G}^t(\gamma) - \bar V_{G}^s(\gamma)\nonumber \\ &=& \left.K\sum_{n=1}^\infty  e^{-rnt} \PP^{\gamma} \left[ \bigcap_{i=1}^{n-1}\left\{ X_{it}> \gamma \right\}\cap\left\{X_{nt}\leq \gamma \right\}\right]\right|^{t\downarrow 0}_{t=s}\nonumber \\ &=& K\cdot \left[\exp\left(-\sum_{n=1}^\infty\frac{e^{-rnt}}{n} \PP^0\left\{X_{nt}\leq 0\right\}\right) \right]^{t=s}_{t\downarrow 0}\label{phase1d identity}.\end{eqnarray}  

Furthermore, if $\mu\geq 0$, there exist constants $c_0, C_0>0$ such that for all sufficiently small $s>0$,  \begin{eqnarray*} c_0s^\frac{1}{2}\leq \lim_{t\downarrow 0}\bar V_{G}^t(\gamma) - \bar V_{G}^s(\gamma) \leq C_0s^\frac{1}{2\sqrt{2}}\end{eqnarray*} If both $\mu\leq 0$ and $r> \frac{\mu^2}{2\sigma^2}$, there exist constants $c_1, C_1>0$ such that for all sufficiently small $s>0$,  \begin{eqnarray*}c_1s^\frac{1}{\sqrt{2}} \leq \lim_{t\downarrow 0}\bar V_{G}^t(\gamma) - \bar V_{G}^s(\gamma) \leq C_1s^\frac{1}{{2}}.\end{eqnarray*} 
\end{Th}

\begin{rem} Although computing the constants $c_0,C_0,c_1,C_1$ explicitly is possible, we refrain from it for the moment, as it is not required to find the right scaling for an extrapolation for $V_G^s$ from $s>0$ to $s=0$ and it would not provide any additional useful information for our extrapolation purposes. The same remark applies to all examples and generalisations that are studied subsequently.
\end{rem}

\begin{proof} The existence of $\lim_{t\downarrow 0}\bar V_{G}^t(\gamma)$ is a consequence of Lemma \ref{dyadic limits generic}. The first identity in the statement of the Theorem is a consequence of the previously established Lemma \ref{phaselemma1}, whereas the second equation in the statement of the Theorem follows from a result by Feller \cite[p. 606, Lemma 3]{F} on processes with stationary and independent increments. For, if we define $$\forall s>0\quad\forall q\in[0,1) \quad \xi(q,s):= \sum_{n=1}^\infty q^n \PP^{0} \left[ \bigcap_{i=1}^{n-1}\left\{ X_{is}> 0\right\}\cap\left\{X_{ns}\leq 0\right\}\right] ,$$ then Feller's identity \cite[p. 606, Lemma 3]{F} reads \begin{equation}\label{Feller} \forall s>0\quad\forall q\in[0,1)\quad -\ln \left(1-\xi(q,s)\right) = \sum_{n=1}^\infty\frac{q^n}{n} \PP^0\left\{X_{ns}\leq 0\right\} \end{equation} and holds whenever $X$ has stationary and independent increments, in particular for all L\'evy processes (note that our definition of a L\'evy process requires them to be Feller processes in addition). This entails \begin{eqnarray} \label{rhosimVG}\xi\left(e^{-rs},s\right)&=&\sum_{n=1}^\infty e^{-rns} \PP^{0} \left[ \bigcap_{i=1}^{n-1}\left\{ X_{is}> 0\right\}\cap\left\{X_{ns}\leq 0\right\}\right] \\&=& 1-\exp\left(-\sum_{n=1}^\infty\frac{e^{-rns}}{n} \PP^0\left\{X_{ns}\leq 0\right\}\right) ,\end{eqnarray} which is enough to prove the second identity (\ref{phase1d identity}) in the Theorem.
This ushers in the derivation of the estimates on $\exp\left(-\sum_{n=1}^\infty \frac{e^{-rns}}{n}\PP^{0} \left\{X_{ns}\leq 0\right\} \right)$ which are needed in order to prove the inequalities of the second half of the Theorem. We shall show that if $\mu\geq 0$, there exist constants $c_0, C_0>0$ such that for all sufficiently small $s>0$,  \begin{eqnarray*} c_0s^\frac{1}{2}\leq \exp\left(-\sum_{n=1}^\infty \frac{e^{-rns}}{n}\PP^{0} \left\{X_{ns}\leq 0\right\} \right) \leq C_0s^\frac{1}{2\sqrt{2}},\end{eqnarray*} and if both $\mu\leq 0$ and $r> \frac{\mu^2}{2\sigma^2}$, there exist constants $c_1, C_1>0$ such that for all sufficiently small $s>0$,  \begin{eqnarray*}c_1s^\frac{1}{\sqrt{2}} \leq \exp\left(-\sum_{n=1}^\infty \frac{e^{-rns}}{n}\PP^{0} \left\{X_{ns}\leq 0\right\} \right) \leq C_1s^\frac{1}{{2}}.\end{eqnarray*} Now, the scaling invariance of Brownian motion yields for all $n\in\NN$ and $s>0$: \begin{eqnarray} \label{scaling}\PP^0\left\{X_{ns}\leq 0\right\}&=& \PP^0\left\{B_{ns}\leq -\frac{\mu}{\sigma} ns\right\}=\PP^0\left\{ B_1\leq -\frac{\mu}{\sigma}(ns)^{1/2}\right\} \nonumber \\ &=& (2\pi)^{-1/2}\int_{-\infty}^{-\frac{\mu}{\sigma}(ns)^{1/2}} \exp\left(\frac{-x^2}{2}\right)dx .\end{eqnarray} We divide the remainder of the proof, which will essentially consist in finding estimates for the right hand side of the last equation, into two parts according to the sign of $\mu$.\\
{\em Case I: $\mu\geq 0$.} In this case we use the estimates $$\forall x\leq 0\forall y\leq 0\quad -y^2-x^2\leq -\frac{|x+y|^2}{2} \leq -\frac{y^2}{2}-\frac{x^2}{2}, $$ thus \begin{eqnarray*}&&\forall y\leq 0 \\ &&e^{-y^2}\int_{-\infty}^0 e^{-x^2}dx \leq\int_{-\infty}^0 \exp\left(-\frac{|x+y|^2}{2}\right) dx\leq e^{-y^2/2}\int_{-\infty}^0 e^{-x^2/2}dx ,\end{eqnarray*} hence by transformation for all $y\leq 0$ \begin{equation}\label{elestim}\frac{ \sqrt{\pi}}{2} e^{-y^2}\leq \int_{-\infty}^y \exp\left(-\frac{x^2}{2}\right) dx \leq \sqrt{\frac{\pi}{2}}e^{-\frac{y^2}{2}} .\end{equation} Due to equation (\ref{scaling}), this entails for all $n\in\NN$, $s>0$, $\mu\geq0$ (if we insert $-\frac{\mu}{\sigma}(ns)^{1/2}$ for $y$) $$\frac{e^{-\left(\frac{\mu}{\sigma}\right)^2 ns}}{2\sqrt{2}}\leq \PP^0\left\{ X_{ns}\leq 0 \right\}\leq \frac{e^{-\frac{\mu^2 ns }{2\sigma^2}}}{2}.$$ Therefore for arbitrary $r,s>0$, \begin{equation}\label{geoharmonic}\frac{1}{2\sqrt{2}}\sum_{n=1}^\infty \frac{e^{-ns\left(r+\frac{\mu^2}{\sigma^2}\right)}}{n} \leq \sum_{n=1}^\infty \frac{e^{-rns}}{n}\PP^0\left\{X_{ns}\leq 0\right\} \leq \frac{1}{2}\sum_{n=1}^\infty \frac{e^{-ns\left(r+\frac{\mu^2}{2\sigma^2}\right)}}{n} \end{equation} The sums in equation (\ref{geoharmonic}) have got the shape of $\sum q^n/n$ for $q<1$. 
Now one performs a standard elementary computation on this power series: \begin{eqnarray}\label{geoharmo1} \sum_{n=1}^\infty \frac{q^n}{n}&=&\sum_{n=0}^\infty \int_0^q r^n dr= \int_0^q\sum_{n=0}^\infty r^n dr = \int_0^q \frac{1}{1-r} dr=-\ln(1-q) ,\end{eqnarray} which immediately gives \begin{eqnarray*} &&\left(1-e^{-s\left(r+\frac{\mu^2}{2\sigma^2}\right)}\right)^{1/2} \\ &\leq &\exp\left(- \sum_{n=1}^\infty \frac{e^{-rns}}{n}\PP^0\left\{X_{ns}\leq 0\right\}\right)\\ & \leq &\left(1-e^{-s\left(r+\frac{\mu^2}{\sigma^2}\right)}\right)^\frac{1}{2\sqrt{2}} \end{eqnarray*} when applied to equation (\ref{geoharmonic}). Due to de l'Hospital's rule, the differences in the brackets on the left and right hand sides of the last estimate behave like $s$ when $s\downarrow 0$. This is sufficient to prove the estimate in the Theorem for the case of $\mu\geq 0$. \\
{\em Case II: $\mu\leq 0$ and $r>\frac{\mu^2}{ 2\sigma^2 }$.} In that case we employ the estimates $$\forall x\leq 0\quad \forall y\leq 0\quad -\frac{x^2}{2}-\frac{y^2}{2}\leq -\frac{|x-y|^2}{2} \leq -\frac{x^2}{4}+\frac{y^2}{2} $$ and proceed analogously to Case I, to obtain \begin{equation}\label{elestimII}\sqrt{\frac{ {\pi}}{2}} e^\frac{-y^2}{2}\leq \int_{-\infty}^{-y} \exp\left(-\frac{x^2}{2}\right) dx \leq \sqrt{{\pi}}e^{\frac{y^2}{2}} .\end{equation} In the special case of $y:= \frac{\mu}{\sigma}(ns)^{1/2}\leq 0$, this leads to the esimate in the statement of the Theorem via \begin{equation}\label{}\frac{1}{2}\sum_{n=1}^\infty \frac{e^{-ns\left(r+\frac{\mu^2}{2\sigma^2}\right)}}{n} \leq \sum_{n=1}^\infty \frac{e^{-rns}}{n}\PP^0\left\{X_{ns}\leq 0\right\} \leq \frac{1}{\sqrt{2}}\sum_{n=1}^\infty \frac{e^{-ns\left(r-\frac{\mu^2}{2\sigma^2}\right)}}{n} .\end{equation} \\
Therefore in case $\mu=0$ the scaling exponent is exactly $\frac{1}{2}$.
\end{proof}

\begin{cor} \label{phasecor1d} Assume $d=1$ and let, as in the previous Theorem \ref{phase1d}, $(X_t)_{t\geq 0}=\left(X_0+\sigma\cdot B_t+\left(r-\frac{\sigma^2}{2}\right)\right)_{t\geq 0}$, in words: $X$ be the logarithmic price process of the one-dimensional Black-Scholes model with constant volatility $\sigma$ and discount rate $r>0$. Furthermore, suppose $g=K-\exp$ and let $G^s:=(-\infty,\gamma^s]$ denote the optimal exercise region for a (one-dimensional) perpetual Bermudan put option of exercise mesh size $s$ and strike price $K$ on the (one-dimensional) basket $X$. Define $\mu:=r-\frac{\sigma^2}{2}$, $\gamma_0:=\sup_{s>0}\gamma^s$ and $G^0:=(-\infty,\gamma^0]$. Then we have for all $s>0$ $$\lim_{t\downarrow 0}\bar V_{G^0}^t(\gamma^0) - \bar V_{G^s}^s(\gamma^s)=o(\gamma^0-\gamma^s) + K\cdot \left[\exp\left(-\sum_{n=1}^\infty\frac{e^{-rns}}{n} \PP^0\left\{X_{ns}\leq 0\right\}\right) \right]^{t=s}_{t\downarrow 0}.$$ Moreover, there are constants $c_0,C_0,c_1,C_1$, such that if $\mu\geq 0$, $$C_0\cdot s^{\frac{1}{2\sqrt{2}}}\geq\lim_{t\downarrow 0}\bar V_{G^0}^t(\gamma^0) - \bar V_{G^s}^s(\gamma^s)\geq c_0\cdot s^{\frac{1}{2}},$$ and if both $\mu\leq 0$ and $r> \frac{\mu^2}{2\sigma^2}$, $$C_1\cdot s^{\frac{1}{2}}\geq \lim_{t\downarrow 0}\bar V_{G^0}^t(\gamma^0) - \bar V_{G^s}^s(\gamma^s)\geq c_1\cdot s^{ \frac{1}{\sqrt{2}} }$$ for all sufficiently small $s>0$.
\end{cor}
\begin{proof} The first asymptotic identity in the statement of the Corollary follows from equation (\ref{phase1d identity}) in Theorem \ref{phase1d} as soon as we have remarked that $$\lim_{t\downarrow 0}\bar V_{G^0}^t(\gamma^0)=\lim_{t\downarrow 0}V_{G^0}^t(\gamma^0)=g(\gamma^0)$$ and $$\lim_{t\downarrow 0}\bar V_{G^s}^t(\gamma^s)=\lim_{t\downarrow 0} V_{G^s}^t(\gamma^s)=g(\gamma^s),$$ for these equations yield $$\lim_{t\downarrow 0}\bar V_{G^0}^t(\gamma^0)-\bar V_{G^s}^s(\gamma^s)=g(\gamma^0)-g(\gamma^s)+\lim_{t\downarrow 0}\bar V_{G^s}^t(\gamma^s)-\bar V_{G^s}^s(\gamma^s)$$ which by the differentiability of $g$ means $$\lim_{t\downarrow 0}\bar V_{G^0}^t(\gamma^0)-\bar V_{G^s}^s(\gamma^s)=o\left(\gamma^0-\gamma^s\right)+\lim_{t\downarrow 0}\bar V_{G^s}^t(\gamma^s)-\bar V_{G^s}^s(\gamma^s).$$
We can now use results on the exercise boundary for perpetual Bermudan options obtained by Boyarchenko and Levendorskii \cite[equation (5.3)]{BL} who showed $$\gamma^s-\gamma^0\sim s^1$$ for sufficiently small $s$, and the estimates in the Corollary follow directly from the estimates of Theorem \ref{phase1d}.
\end{proof}

\begin{rem} Up to this point, we have derived estimates for the American-Bermundan option price difference at the boundary $\gamma$ of the exercise region $G=(-\infty,\gamma]$ (in case of a put) or $G=[\gamma,+\infty)$ (in case of a call with dividends). We can extend these bounds of the American-Bermudan difference from the exercise boundary to the complement of the exercise region: By continuity, we can even extend the lower bounds or upper bounds, respectively, to a neighbourhood of the exercise boundary: For, if we consider a put for the moment, we get from Lemma \ref{VG pde} that $(x,s)\mapsto \bar V_{G}^s(x)$ is continuous for all $G$, implying that if $x\notin G=(-\infty, \gamma]$ $$\lim_{s\downarrow 0} \frac{\ln\left(\bar V_{G}^0(\gamma)-\bar V_{G}^s(\gamma)\right)}{\ln s}=\lim_{x\uparrow \gamma}\lim_{s\downarrow 0} \frac{\ln\left(\bar V_{G}^0(x)-\bar V_{G}^s(x)\right)}{\ln s},$$ where $\bar V_{G}^0$ is shorthand for $\lim_{t\downarrow 0}\bar V_{G}^t.$ Thus, if $\lim_{s\downarrow 0} \frac{\ln\left(\bar V_{G}^0(\gamma)-\bar V_{G}^s(\gamma)\right)}{\ln s}=:\alpha\in(0,1)$ -- where we have, thanks to Theorem \ref{phase1d} estimates for the limit $\lim_{s\downarrow 0} \frac{\ln\left(\bar V_{G}^0(\gamma)-\bar V_{G}^s(\gamma)\right)}{\ln s}$ -- we will for any $\varepsilon>0$ get a $\delta >0$ such that $$\forall x\in[\gamma,\gamma+\delta)\quad \lim_{s\downarrow 0} \frac{\ln\left(\bar V_{G}^0(x)-\bar V_{G}^s(x)\right)}{\ln s}\in(\alpha-\varepsilon,\alpha+\varepsilon).$$
Analogously, we can proceed to derive bounds for the American-Bermundan call option price difference (for an option on a dividend-paying asset) in a neighbourhood of the exercise boundary.
\end{rem}

\section{One-dimensional continuity corrections outside the Black-Scholes model}

The identity (\ref{phase1d identity}) of Theorem \ref{phase1d} can be used to derive estimates in the spirit of the second half of Theorem \ref{phase1d} in more general situations. We will illustrate this by means of the following example:

\begin{ex}[Merton's jump-diffusion model with positive jumps and ``moderate'' volatility] Suppose the logarithmic price process $X$ is governed by an equation of the form $$ \forall t\geq 0 \quad X_t=X_0+\alpha t+ \beta Z_t+ \sigma B_t $$ where $\alpha\in\RR$, $\beta,\sigma>0$, $Z$ is the Poisson process (thus, in this setting, only positive jumps are allowed for simplicity) and $B$ a normalised one-dimensional Brownian motion, and the stochastic processes $B$ and $Z$ are assumed to be independent. Let $\PP^\cdot$ be an associated family of risk-neutral measures and $r>0$ the discount rate. In order to employ (\ref{phase1d identity}), we shall compute the sum $\sum_{n=0}^\infty\frac{e^{-rns}}{n}\PP^0\left\{X_{ns}\leq 0\right\}$ for all $s>0$. Since $\PP^0\left\{X_{ns}\leq 0\right\}=\PP^0\left\{\frac{X_{ns}}{\beta}\leq 0\right\}$ for arbitrary $n,s$ we may without loss of generality take $\beta =1$. Let us also assume $\alpha\geq 0$; note that since $\exp\left(X_t-rt\right)_{t\geq 0}$ is a martingale -- as $X$ is a logarithmic price process -- , $\sigma$ must be such that $\alpha -r+\frac{\sigma^2}{2}+ \EE^0\left[e^{Z_1}\right]=0$ (if $r>0$ and $\alpha\geq 0$ are given), hence $\alpha\geq 0$ implies $\sigma\leq \sqrt{2\left(r+ \EE^0\left[e^{Z_1}\right]\right)}$. Now, by definition of the Poisson distribution together with the symmetry and scaling invariance of Brownian motion \begin{eqnarray} \label{jumpsum1} \nonumber &&\sum_{n=0}^\infty\frac{e^{-rns}}{n}\PP^0\left\{X_{ns}\leq 0\right\} \\ \nonumber &=& \sum_{n=0}^\infty \sum_{k=0}^\infty e^{-ns}\frac{(ns)^k}{k!}\cdot \frac{e^{-rns}}{n}\PP^0\left\{\sigma B_{ns}\leq -\alpha ns-k\right\} \\ &=& \sum_{n=0}^\infty \sum_{k=0}^\infty e^{-ns}\frac{(ns)^k}{k!}\cdot \frac{e^{-rns}}{n}\underbrace{\PP^0\left\{B_{1}\leq -\frac{\alpha}{\sigma} (ns)^{\frac{1}{2}}-\frac{k}{\sigma}(ns)^{-\frac{1}{2}}\right\}}_{=(2\pi)^{-1/2}\int_{-\infty}^{-\frac{\alpha}{\sigma} (ns)^{\frac{1}{2}}-\frac{k}{\sigma}(ns)^{-\frac{1}{2}}} \exp\left(\frac{-x^2}{2}\right)dx}  \end{eqnarray} (with the convention that $0^0=1$). Now let us first of all try and find estimates for the probability in the last line. By equation (\ref{elestim}) applied to $y:=-\frac{\alpha}{\sigma} (ns)^{\frac{1}{2}}-\frac{k}{\sigma}(ns)^{-\frac{1}{2}}\leq 0$, \begin{eqnarray*}&&\frac{ e^{ -\left(\frac{\alpha}{\sigma} (ns)^{\frac{1}{2}}+\frac{k}{\sigma}(ns)^{-\frac{1}{2}} \right)^2 } }{2\sqrt{2}} \\ &\leq& \PP^0\left\{B_{1}\leq -\frac{\alpha}{\sigma} (ns)^{\frac{1}{2}}-\frac{k}{\sigma}(ns)^{-\frac{1}{2}}\right\} \\ &\leq& \frac{ e^\frac{-\left(\frac{\alpha}{\sigma} (ns)^{\frac{1}{2}}+\frac{k}{\sigma}(ns)^{-\frac{1}{2}} \right)^2}{2} }{{2}} \end{eqnarray*} which yields, using the abbreviation $\alpha':=\frac{\alpha}{\sigma}$, \begin{eqnarray*}&&\frac{ e^{ -{\alpha'}^2 ns -2 \frac{\alpha'}{\sigma}k - \frac{k^2}{\sigma^2ns}  } }{2\sqrt{2}} \\ &\leq& \PP^0\left\{B_{1}\leq -\frac{\alpha}{\sigma} (ns)^{\frac{1}{2}}-\frac{k}{\sigma}(ns)^{-\frac{1}{2}}\right\} \\ &\leq& \frac{ e^{ -\frac{{\alpha'}^2}{2} ns -\frac{\alpha '}{\sigma}k -\frac{k^2}{2\sigma^2ns} } }{2} , \end{eqnarray*} so \begin{eqnarray} \label{lowerjump1} &&\frac{ e^{ -{\alpha'}^2 ns -k\left(2 \frac{\alpha'}{\sigma} + \frac{k}{\sigma^2ns}\right) } }{2\sqrt{2}} \\ \nonumber &\leq& \PP^0\left\{B_{1}\leq -\frac{\alpha}{\sigma} (ns)^{\frac{1}{2}}-\frac{k}{\sigma}(ns)^{-\frac{1}{2}}\right\} \\ \nonumber &\leq& \frac{ e^{ -\frac{{\alpha'}^2}{2} ns -\frac{\alpha '}{\sigma}k } }{2}. \end{eqnarray} Thus, we can perform the following estimates to derive an upper bound of the sum in (\ref{jumpsum1}): \begin{eqnarray*} && \sum_{n=0}^\infty \sum_{k=0}^\infty e^{-ns}\frac{(ns)^k}{k!}\cdot \frac{e^{-rns}}{n}\PP^0\left\{B_{1}\leq -\frac{\alpha}{\sigma} (ns)^{\frac{1}{2}}-\frac{k}{\sigma}(ns)^{-\frac{1}{2}}\right\} \\ &\leq& \frac{1}{2} \sum_{n=0}^\infty \frac{e^{-ns\left( 1+r+\frac{{\alpha'}^2}{2}\right)}}{n}\sum_{k=0}^\infty\frac{1}{k!}\left(e^\frac{-\alpha '}{\sigma}\cdot ns\right)^k \\ &=& \frac{1}{2} \sum_{n=0}^\infty \frac{e^{-ns\left( 1+r+\frac{{\alpha'}^2}{2}\right)} }{n} e^{e^\frac{-\alpha '}{\sigma}\cdot ns} \\ &=& \frac{1}{2} \sum_{n=0}^\infty \frac{e^{-ns\left( 1+r+\frac{{\alpha'}^2}{2}-e^\frac{-\alpha '}{\sigma}\right)} }{n} \\ &=& -\frac{1}{2}\ln\left(1-e^{-s\left( 1+r+\frac{{\alpha'}^2}{2}-e^\frac{-\alpha '}{\sigma}\right) }\right) \end{eqnarray*} where the last line uses that $\frac{\alpha'}{\sigma}=\frac{\alpha}{\sigma^2}\geq \alpha\cdot\left(r+\EE^0\left[e^{Z_1}\right]\right)\geq 0$ and we need to impose the condition that $e^{-\alpha\cdot\left(r+\EE^0\left[e^{Z_1}\right]\right)}\leq 1+r+\frac{\alpha^2}{2\sigma^2}$ (which, given $r>0$ and $\alpha$, will be satisfied if $\sigma>0$ is sufficiently small) to employ the identity \begin{equation}\label{geoharmo}\forall q<1\quad \sum_{n=0}^\infty\frac{q^n}{n}=\ln\frac{1}{1-q}.\end{equation} The lower bound follows simply from $$\forall n\in\NN_0\forall s>0\quad \sum_{k=0}^\infty\frac{1}{k!}\left(ns\cdot e^{ -2 \frac{\alpha'}{\sigma} - \frac{k}{\sigma^2ns}}\right)^k\geq 1$$ (for $n=0$ recall that $0^0=1$ in this paragraph by our earlier convention) as this entails (when exploiting the estimate (\ref{lowerjump1}) and finally (\ref{geoharmo}) ): \begin{eqnarray*} && \sum_{n=0}^\infty \sum_{k=0}^\infty e^{-ns}\frac{(ns)^k}{k!}\cdot \frac{e^{-rns}}{n}\PP^0\left\{B_{1}\leq -\frac{\alpha}{\sigma} (ns)^{\frac{1}{2}}+\frac{k}{\sigma}(ns)^{-\frac{1}{2}}\right\} \\ &\geq& \frac{1}{2\sqrt{2}} \sum_{n=0}^\infty \frac{e^{-ns\left( 1+r+{\alpha'}^2\right)}}{n} \sum_{k=0}^\infty\frac{1}{k!}\left(ns\cdot e^{-2 \frac{\alpha'}{\sigma} - \frac{k}{\sigma^2ns}}\right)^k\\ &\geq& \frac{1}{2\sqrt{2}} \sum_{n=0}^\infty \frac{e^{-ns\left( 1+r+{\alpha'}^2\right)}}{n}  \\ &=& -\frac{1}{2\sqrt{2}}\ln\left(1-e^{-s\left( 1+r+{\alpha'}^2\right)}\right)\end{eqnarray*} As a consequence of these estimates and using the Taylor expansion of $\exp$ around $0$, we now get the existence of two constants $c_3>0$ and $C_3>0$ (which can be computed explicitly) such that for all sufficiently small $s$, $$c_3\cdot s^\frac{1}{2}\leq \exp\left(-\sum_{n=0}^\infty\frac{e^{-rns}}{n}\PP^0\left\{X_{ns}\leq 0\right\} \right)\leq C_3\cdot s^\frac{1}{2\sqrt 2}.$$ Finally, we may apply identity (\ref{phase1d identity}) from Theorem \ref{phase1d} -- as this is an immediate consequence of Feller's identity \cite[p. 606, Lemma 3]{F} and our Lemma \ref{phaselemma1} -- and conclude that if $g=K-f$ and $G=(-\infty,\gamma]$, then $$Kc_3\cdot s^\frac{1}{2}\leq \lim_{t\downarrow 0}\bar V_{G}^t(\gamma) - \bar V_{G}^s(\gamma) \leq KC_3\cdot s^\frac{1}{2\sqrt 2} $$ for all sufficiently small $s>0$.

\end{ex}

\section{Continuity corrections in higher dimensions}

The proof of Theorem \ref{phase1d} relies heavily on the use of Feller's result \cite[p. 606, Lemma 3]{F} which in turn is proven by means of elementary Fourier analysis and a so-called ``basic identity'' \cite[p. 600, equation (1.9)]{F}.

Hence, if one aims at generalising Theorem \ref{phase1d} to higher dimensions, one should first of all find a multi-dimensional analogue of the said basic identity.

Indeed, we shall see that this is feasible. Let us for the following fix a stochastic process $X=(X_t)_{t\geq 0}$ on $\RR^d$ with stationary and independent increments.

\begin{lem} Suppose $H$ is a measurable subset of $\RR^d$, and $s>0$. Define for all $n\in\NN$ $$\forall K\in\cB\left(\RR^d\right) \quad R_n(K):=\PP^0\left[\bigcap_{1\leq i< n}\left\{X_{is}\in \complement H\right\}\cap \left\{X_{ns}\in K\cap H\right\}\right], $$ as well as  $$\forall K\in\cB\left(\RR^d\right)\quad Q_n(K):=\PP^0\left[\bigcap_{1\leq i< n}\left\{X_{is}\in \complement H\right\}\cap \left\{X_{ns}\in \complement H\cap K\right\}\right] $$ (in particular $R_0=\delta_{0}\left[\cdot\cap H\right]=0$ and $Q_0=\delta_{0}\left[\cdot\cap \complement H\right]=\delta_0$). Then for all $n\in\NN_0$, $$Q_{n+1}+R_{n+1}=Q_n\ast {\PP^0}_{X_s}.$$
\end{lem}
\begin{proof} Consider a measurable $K\subseteq \RR^d$. Clearly, \begin{equation} \label{Fellergeneral} \left(Q_{n+1}+R_{n+1}\right)(K)=\PP^0\left[\bigcap_{i=1}^{n}\left\{X_{is}\in \complement H\right\}\cap \left\{X_{(n+1)s}\in K\right\}\right]. \end{equation} On the other hand, since $X$ is a Markov process, we have $$Q_n(K) =\left(P_s\left(\chi_{\complement H}\cdot\right)\right)^{\circ n}\chi_K(0)$$ (where $(P_t)_{t\geq 0}:=\left(\PP_{X_t}^0\ast\cdot\right)_{t\geq 0}$ is the translation-invariant Markov semigroup of transition functions for the process $X$ whose increments are stationary and independent), thus $$\int_{\RR^d}f(y)Q_n(dy)= \left(P_s\left(\chi_{\complement H}\cdot\right)\right)^{\circ n}f(0)$$ for all nonnegative measurable functions $f$. But this implies \begin{eqnarray*}\left(Q_{n}\ast \PP^0_{X_s}\right)(K)&=&\int_{\RR^d}\int_{\RR^d}\chi_K\left(z+y\right)\PP_{X_s}^0(dz)Q_n(dy)\\ &=& \left(P_s\left(\chi_{\complement H}\cdot\right)\right)^{\circ n}\left(\int_{\RR^d}\chi_{K-\cdot}(z)\PP_{X_s}^0(dz)\right)(0) \\ &=& \left(P_s\left(\chi_{\complement H}\cdot\right)\right)^{\circ n}\circ\left(\PP_{X_s}^0\ast\chi_K\right)(0)\\ &=& \left(P_s\left(\chi_{\complement H}\cdot\right)\right)^{\circ n}\circ P_s\chi_K(0),
\end{eqnarray*}
and the right hand side of this equation coincides with the one of identity (\ref{Fellergeneral}).
\end{proof}

Applying Fourier transforms we obtain

\begin{cor}\label{1.9general} Let us adopt the notation of the preceding Lemma and define the {\em Fourier transform} of a countable sequence $\left(\mu_n\right)_n$ of finite measures on $\RR^d$, denoted by $\widehat{\left(\mu_n\right)_n}=\widehat{\mu}:(0,1)\times\RR^d\rightarrow \CC$, by $$\forall q\in(0,1)\quad \forall\zeta \in \RR^d\quad \widehat{\left(\mu_n\right)_n}(q,\zeta)=\sum_{n=0}^\infty q^n\int_{\RR^d} e^{i\cdot{{^t}\zeta}y}\mu_n(dy)=\sum_{n=0}^\infty q^n\widehat{\mu_n}(\zeta).$$ Then for all $q\in (0,1)$, and $\zeta\in\RR^d$  the equation $$1-\widehat{\left(R_n\right)_n}(q,\zeta) =\widehat{\left(Q_n\right)_n}(q,\zeta)\left(1-q\widehat{{\PP^0}_{X_s}}(\zeta)\right)$$ holds.
\end{cor}
\begin{proof} The result of the previous Lemma reads $$\forall n\in\NN_0 \quad \widehat{Q_{n+1}}+\widehat{R_{n+1}}=\widehat{Q_n}\widehat{{\PP^0}_{X_s}}$$ when we apply the Fourier transform. After multiplication with $q^{n+1}$ and summing up over $n\in\NN_0$, one arrives at $$\forall q\in(0,1)\quad\forall \zeta\in\RR^d\quad \widehat Q(q,\zeta)- \underbrace{\chi_{\complement H}(0)}_{=\widehat{Q_0}(\zeta)}+\widehat R(q,\zeta) - \underbrace{\chi_{H}(0)}_{=\widehat{R_0}(\zeta)}= q\widehat Q(q,\zeta)\widehat{{\PP^0}_{X_s}}(\zeta) ,$$ hence $$\forall q\in(0,1)\quad \widehat R(q,\cdot)-1 = q\widehat Q(q,\cdot)\widehat{{\PP^0}_{X_s}}(\cdot) -\widehat Q(q,\cdot). $$ This is our claim.
\end{proof}

\begin{Def} A subset $A\subseteq \RR^d$ is called {\em $+$-closed} if and only if $A$ is measurable and $A+A\subseteq A$, that is sums of elements of $A$ are again elements of $A$.
\end{Def}

\begin{lem}[\`a la Feller, Wiener, Hopf] \label{alaWienerHopf}Suppose $H$ is a $+$-closed set and its complement $\complement H$ is a $+$-closed set as well. Assume furthermore $0\not\in H$ (ensuring $R_0 = 0$), and let $\ln$ the main branch of the logarithm on $\CC$. Then $$-\ln\left({1-\widehat{R}(q,\zeta)}\right)=\sum_{n=1}^\infty\frac{ q^n}{n}\int_He^{i\cdot{{^t}\zeta}x}\left({\PP^0}_{X_s}\right)^{\ast n}(dx)$$ for all $(q,\zeta)\in (0,1)\times\RR^d$ such that the left-hand side is well-defined. In general, for all $q\in (0,1)$, one has at least $${1-\widehat{R}(q,0)}=\exp\left(-\sum_{n=1}^\infty\frac{q^n\cdot\left( \widehat{\PP^0_{X_s}\left[\cdot\cap H\right]}(0)\right)^n }{n}\right).$$
\end{lem}
\begin{proof} Let $q\in (0,1)$. According to the previous Corollary \ref{1.9general}, we have \begin{equation}\forall \zeta\in U\quad \ln\frac{1}{1-q\widehat{{\PP^0}_{X_s}}(q,\zeta)}=\ln\frac{1}{1-\widehat{R}(q,\zeta)} - \ln\widehat Q(q,\zeta)\label{alaWienerHopfLogEq} \end{equation} wherever this is defined. Due to the identities $\sum_{n=1}^\infty \frac{r^n}{n}=\ln\frac{1}{1-r}$ for all $r\in B_1(0)\subset \CC$ (cf equation (\ref{geoharmo1}) in the proof of Theorem \ref{phase1d} above) and $\widehat{{\PP^0}_{X_s}}^n=\widehat{{{\PP^0}_{X_s}}^{\ast n}}$ this can also be written as \begin{eqnarray}  &&\label{alafeller}\sum_{n=1}^\infty\frac{ q^n }{n}\int_{\RR^d} e^{i\cdot{{^t}\zeta}x}\left({\PP^0}_{X_s}\right)^{\ast n}(dx) \\ \nonumber &=& \sum_{n=1}^\infty \frac{1}{n} \left(\widehat{R}(q,\zeta)\right)^n + \sum_{n=1}^\infty \frac{(-1)^n}{n}\left(\widehat{Q}(q,\zeta)-1\right)^n.\end{eqnarray} However, at least for $\zeta=0$ and arbitrary choice of $q$, one may still state identity (\ref{alaWienerHopfLogEq}) as this follows from Corollary \ref{1.9general} more or less directly: First we note that \begin{eqnarray*}&&-\ln\left(\left(1-q\widehat{{\PP^0}_{X_s}}(q,0)\right)\cdot \widehat{\left(Q_n\right)_n}(q,0)\right)\\&=&-\ln \widehat{\left(Q_n\right)_n}(q,0)-\ln\left(1-q\widehat{{\PP^0}_{X_s}}(q,0)\right)\end{eqnarray*} (as in these statements the arguments of $\ln$ are positive, hence surely in the domain of $\ln$) and written in series notation \begin{eqnarray*}&& \sum_{n=1}^\infty \frac{1}{n}\left(1-\widehat{\left(Q_n\right)_n}(q,0)\right)^n+ \sum_{n=1}^\infty\frac{q^n}{n}\widehat{{\PP^0}_{X_s}}(q,0)^n \\ &=& \sum_{n=1}^\infty \frac{(-1)^n}{n}\left(\left(1-q\widehat{{\PP^0}_{X_s}}(q,0)\right)\cdot \widehat{\left(Q_n\right)_n}(q,0)-1\right)^n.\end{eqnarray*} But Corollary \ref{1.9general} implies $$ \forall n\in\NN\quad\frac{(-1)^n}{n}\left(\left(1-q\widehat{{\PP^0}_{X_s}}(q,0)\right)\cdot \widehat{\left(Q_n\right)_n}(q,0)-1\right)^n = \frac{1}{n}\widehat{\left(R_n\right)_n}(q,0).$$ Combining these two equations yields (\ref{alaWienerHopfLogEq}). Next, note that $$\mu_{R,q}:=\sum_{n=0}^\infty q^nR_n$$ is still a finite measure -- concentrated on $H$ -- and thus possesses a Fourier transform. Analogously, the measure $\mu_{Q,q}:=\sum_{n=0}^\infty q^nQ_n$ is concentrated on $\complement H$ and also has a Fourier transform as it is finite. Now, for arbitrary $n\in\NN$, the properties of the Fourier transform imply \begin{eqnarray*}\left(\widehat{R}(q,\cdot)\right)^n&=& \left(\widehat{\mu_{R,q}}\right)^n=\widehat{{\mu_{R,q}}^{\ast n}}, \\ \left(\widehat{Q}(q,\cdot)-1\right)^n &=& \left(\widehat{\mu_{Q,q}-\delta_0}\right)^n=\widehat{\left( \left({\mu_{Q,q}-\delta_0}\right)^{\ast n} \right)}.\end{eqnarray*} But since $H$ and $\complement H$ are $+$-closed sets, i.e. $H+H\subseteq H$ and $\complement H+\complement H\subseteq\complement H$ , the measures on the right hand sides of these two equations, ${{\mu_{R,q}}^{\ast n}}$ and $\left({\mu_{Q,q}-\delta_0}\right)^{\ast n}$, have to be (signed) measures on $H$ and $\complement H$, respectively. Let us now split the sum in (\ref{alafeller}) and insert the terms we have previously identified:
\begin{eqnarray} \nonumber&& \forall \zeta\in U \\ &&\sum_{n=1}^\infty\frac{ q^n }{n}\int_{H} e^{i\cdot{{^t}\zeta}x}\left({\PP^0}_{X_s}\right)^{\ast n}(dx) + \sum_{n=1}^\infty\frac{ q^n }{n}\int_{\complement H} e^{i\cdot{{^t}\zeta}x}\left({\PP^0}_{X_s}\right)^{\ast n}(dx) \nonumber \\  &=& \sum_{n=1}^\infty \frac{1}{n} \widehat{\left({\mu_{R,q}}^{\ast n}\right)}(\zeta) + \sum_{n=1}^\infty \frac{(-1)^n}{n}\widehat{\left({\mu_{Q,q}-\delta_0}\right)^{\ast n}}(\zeta).\end{eqnarray} It is the injectivity of the Fourier transform that yields from this \begin{eqnarray} \nonumber &&\sum_{n=1}^\infty\frac{ q^n }{n}\left(\left({\PP^0}_{X_s}\right)^{\ast n}\left(\cdot\cap H\right)\right) + \sum_{n=1}^\infty\frac{ q^n }{n} \left(\left({\PP^0}_{X_s}\right)^{\ast n}\left(\cdot\cap \complement H\right)\right) \\ \nonumber &=& \sum_{n=1}^\infty \frac{1}{n}\left({\mu_{R,q}}^{\ast n}\right) + \sum_{n=1}^\infty \frac{(-1)^n}{n}\left({\mu_{Q,q}-\delta_0}\right)^{\ast n}.\end{eqnarray} Either side of this equation equals the sum of two (signed measures), and we recall that the first measure on the left hand side and first measure on the right hand side are both concentrated on $H$, whilst the second measure on the left hand side as well as the second measure on the right hand side are both concentrated on $\complement H$. The only way for this to be true is that the two measures that are concentrated on each of $H$ or $\complement H$ are equal: $$\sum_{n=1}^\infty\frac{ q^n }{n}\left(\left({\PP^0}_{X_s}\right)^{\ast n}\left(\cdot\cap H\right)\right) = \sum_{n=1}^\infty \frac{1}{n}\left({\mu_{R,q}}^{\ast n}\right),$$ and also $$ \sum_{n=1}^\infty\frac{ q^n }{n} \left(\left({\PP^0}_{X_s}\right)^{\ast n}\left(\cdot\cap \complement H\right)\right) =  \sum_{n=1}^\infty \frac{(-1)^n}{n}\left({\mu_{Q,q}-\delta_0}\right)^{\ast n} ,$$ the former identity being exactly what the statement of the Lemma expresses in the language of Fourier transforms.
\end{proof}

Based on this result, we may partially generalise Theorem \ref{phase1d} to higher dimensions when we require $G$ (the set that we refer to the exercise region) to be $+$-closed set.

\begin{Th}\label{phase} Let us make the assumptions of Lemma \ref{characteriseputdifference}, viz: Let $X$ be the logarithmic price process of the multidimensional Black-Scholes model for independent assets with constant volatility and interest rate, that is $$\forall t\geq 0 \quad X_t=\left(\left(X_0\right)_i+\sigma_i\cdot (B_t)_i + \left(r-\frac{1}{2}{\sigma_i}^2\right) t\right)_{i=1}^d$$ (where $B$ is the $d$-dimensional Wiener process) for some $r>0$ and $\sigma\in{\RR_{>0}}^d$.
Let $g=K-f$, wherein $K\geq 0$ be a real number and $f\geq 0$ be a continuous function that is monotonely increasing in each component and such that $\left(e^{-rt}f(X_t)\right)_{t\geq 0}$ is a martingale.
Finally, consider a measurable set of the shape $G=\gamma-H'$ for some $\gamma\in\RR^d$ and some convex north-east connected set $H'\subseteq {\RR_+}^d$ (making $G$ convex and south-west connected) such that $g=K-f$ is nonnegative on $G$ and $0\not \in H'$. Suppose furthermore that both $H:=-H'$ and $\complement H=-\complement H' $ are $+$-closed.
Then for all $s>0$, \begin{eqnarray*} && \lim_{t\downarrow 0}\bar V_{G}^t(\gamma) - \bar V_{G}^s(\gamma)\\&=& \left.K\cdot\sum_{n=1}^\infty  e^{-rnt} \PP^{\gamma} \left[ \bigcap_{i=1}^{n-1}\left\{ X_{it}\in\complement G\right\}\cap\left\{X_{nt}\in G\right\}\right]\right|^{t\downarrow 0}_{t=s} \\&=& \left.K\cdot\sum_{n=1}^\infty  e^{-rnt} \PP^{0} \left[ \bigcap_{i=1}^{n-1}\left\{ X_{it}\in\complement H\right\}\cap\left\{X_{nt}\in H\right\}\right]\right|^{t\downarrow 0}_{t=s}\\ &=& \left. K\exp\left(-\sum_{n=1}^\infty \frac{e^{-rns}}{n}\PP^{0} \left\{X_{ns}\in G-\gamma \right\} \right)\right|^{t=s}_{t\downarrow 0} .\end{eqnarray*} 
\end{Th}
\begin{proof} The existence of $\lim_{t\downarrow 0}\bar V_{G}^t(\gamma)$ is a consequence of Lemma \ref{dyadic limits generic}. The subsequent identity follows directly from Lemma \ref{alaWienerHopf} (just as in the proof of the corresponding equation in Theorem \ref{phase1d}, except that {\em en lieu} of Lemma \ref{alaWienerHopf}, the proof of Theorem \ref{phase1d} makes use of Feller's original result \cite[p. 606, Lemma 3]{F}): For, the second equation in Lemma \ref{alaWienerHopf} may be read $${1-\sum_{n=0^\infty}q^nR_n\left[\RR^d\right] } = \exp\left(-\sum_{n=1}^\infty\frac{q^n\cdot\left( \PP^0_{X_s}\left[H\right]\right)^n }{n}\right)$$ that is \begin{eqnarray*} &&  1-\sum_{n=0^\infty}q^n\PP^0\left[\bigcap_{1\leq i< n}\left\{X_{is}\in \complement H\right\}\cap \left\{X_{ns}\in H\right\}\right]\\ & = & \exp\left(-\sum_{n=1}^\infty\frac{q^n\cdot\left( \PP^0 \left\{X_s\in H\right\}\right)^n }{n}\right) \end{eqnarray*} for all $q\in (0,1)$, in particular for $q=e^{-rs}$.

\end{proof}

Again, an analogous result can be accomplished when the function $g$ is replaced by $f-K$ and the set $H$ by $-H$:

\begin{rem}

In special cases, one can find estimates for $\PP^{0} \left\{X_{s}\in H\right\}=\PP^{0} \left\{X_{s}\in G-\gamma \right\}$ that are strong enough to establish multi-dimensional generalisations of the estimates in the second half of Theorem \ref{phase1d}. We shall give a few examples.

In general our results can be used for the extrapolation from (multi-dimensional) Bermudan to American barrier knock-in option prices when the barrier regions and their complements are, up to a constant factor, closed with respect to multiplication, and when, in addition, the barrier region is convex as well as south-west connected (in the case of put options) or north-east connected (in the case of call options), and is contained in the immediate exercise region of the corresponding American option.

\end{rem}

\begin{ex} Let $d$ be any natural number. Consider the convex, north-east connected, $+$-closed set $H:={\RR_{<0}}^d$ (whose complement is also $+$-closed) and set $G=\gamma+H$. Let us impose the same assumptions on $X$, $f$, and $g$ as in the statement of the previous Theorem \ref{phase}. Then one has, due to the independence of the components of $X$, the following bounds for all $s>0$:
$$c s^\alpha\leq\exp\left(-\sum_{n=1}^\infty \frac{e^{-rns}}{n}\PP^{0} \left\{X_{ns}\in H\right\} \right) \leq  C s^\beta,$$ where $c, C, \alpha,\beta$ depend on $\mu$. $\alpha=2^{-d}$ and $\beta=2^{-d/2}$ if $\mu \geq 0$, $\alpha=2^{-3d/2}$ and $\beta = 2^{-d}$ if $\mu\leq 0$. Hence also, \begin{eqnarray*}&& K\cdot \exp\left(-\sum_{n=1}^\infty \frac{e^{-rns}}{n}\PP^{0} \left\{X_{ns}\in H\right\} \right) \\ &=& K\cdot \exp\left(-\sum_{n=1}^\infty \frac{e^{-rns}}{n}\PP^{\gamma} \left\{X_{ns}\in G\right\} \right)\\ &=& K\cdot \left[\exp\left(-\sum_{n=1}^\infty \frac{e^{-rnt}}{n}\PP^{\gamma} \left\{X_{nt}\in G\right\} \right)\right]^{t\downarrow 0}_{t=s} \\ &=& \lim_{t\downarrow 0}\bar V_G^t(\gamma)- \bar V_G^s(\gamma) \geq \lim_{t\downarrow 0}\bar V_G^t(x)- \bar V_G^s(x) \end{eqnarray*} for arbitrary $s>0$ and $x\geq \gamma$ componentwise. 
\end{ex}

\begin{ex}\label{Ex2} Suppose $d=2$ and let again $X$, $f$, $K$ and $g$ be as in Theorem \ref{phase}, though we will later on have to impose the condition of $\mu\geq 0$ (componentwise). Furthermore consider the convex south-west connected $+$-closed set $$H:=\left\{(x,y)\in \RR^2 \ : \ x+cy\leq 0\right\}$$ whose complement is also $+$-closed for $c>0$ and set $G=\gamma+ H$. Note that in this situation $\exp G$ (which one might refer to as the non-logarithmic exercise region) equals $$\exp G= e^\gamma\cdot\left\{(u,v)\in\RR^2 \ : \ u,v>0,\quad v\leq u^{-c}\right\}$$ ($\exp$ denoting componentwise exponentiation as before). Then we get for all $t>0$ and arbitrary $d\in\NN$, \begin{eqnarray} \PP^\gamma\left\{ X_t\in G \right\}&= & \PP^0\left\{ X_t\in H \right\}\\ &=& \PP^0\left\{ B_t\in H -t\left(\frac{\mu_i}{\sigma_i}\right)_{i=1}^d \right\} \nonumber \\ &=&\PP^0\left\{ t^{-1/2}B_t\in t^{-1/2}H -t^{1/2}\left(\frac{\mu_i}{\sigma_i}\right)_{i=1}^2\right\}\nonumber  \\ &=&\PP^0\left\{ B_1\in H -t^{1/2}\left(\frac{\mu_i}{\sigma_i}\right)_{i=1}^d\right\} \nonumber \\ &=&\PP^0\left\{ -B_1\in H -t^{1/2}\left(\frac{\mu_i}{\sigma_i}\right)_{i=1}^d\right\} \nonumber  \\&=& {\nu_{0,1}}^d\left\{(x,y)\in\RR^d  : \left(x-\frac{\mu_1}{\sigma_1}t^{1/2},y-\frac{\mu_2}{\sigma_2}t^{1/2}\right)\in -H\right\} \nonumber \\ &=& {\nu_{0,1}}^d\left\{(x,y)\in\RR^2  : x+cy \geq \frac{\mu_1}{\sigma_1}t^{1/2}+ c \frac{\mu_2}{\sigma_2}t^{1/2}\right\}
\label{X_tinGprobability}\end{eqnarray} (where $\nu_{0,1}$ is the normal Gaussian measure on $\RR$ and ${\nu_{0,1}}^d={\nu_{0,1}}^{\otimes d}$).
 
Next observe that for any real number $\alpha\geq 0$, from rotating the set $\left\{(x,y)\in{\RR_+}^2  : x+cy\geq \alpha\right\}$ by $\frac{\pi}{2}$, $\pi$ and $\frac{3\pi}{2}$, we obtain, via exploiting the translation-invariance of the two-dimensional normal Gaussian measure ${\nu_{0,1}}^2$, the relation $$4{\nu_{0,1}}^2\left\{(x,y)\in{\RR_+}^2  : x+cy\geq \alpha\right\} \geq 2{\nu_{0,1}}^2\left\{(x,y)\in\RR^2  : x+cy\geq \alpha\right\}.$$ Using the trivial estimate ${\nu_{0,1}}^2\left\{(x,y)\in\RR^2  : x+cy\geq \alpha\right\}\geq {\nu_{0,1}}^2\left\{(x,y)\in{\RR_+}^2  : x+cy\geq \alpha\right\}$, we arrive at \begin{eqnarray*} {\nu_{0,1}}^2\left\{(x,y)\in{\RR_+}^2  : x+cy\geq \alpha\right\}& \leq& {\nu_{0,1}}^2\left\{(x,y)\in{\RR}^2  : x+cy\geq \alpha\right\}\\ &\leq& 2\cdot {\nu_{0,1}}^2\left\{(x,y)\in{\RR}^2  : x+cy\geq \alpha\right\}.\end{eqnarray*} But of course, by a change of coordinates, viz $z:=\frac{y}{c}$, one has $${\nu_{0,1}}^2\left\{(x,y)\in{\RR}^2  : x+cy\geq \alpha\right\} = c{\nu_{0,1}}^2\left\{(x,z)\in{\RR}^2  : x+z\geq \alpha\right\}$$ and $${\nu_{0,1}}^2\left\{(x,y)\in{\RR_+}^2  : x+cy\geq \alpha\right\} = c{\nu_{0,1}}^2\left\{(x,z)\in{\RR_+}^2  : x+z\geq \alpha\right\}.$$ 

Applying this to the equation (\ref{X_tinGprobability}) for $\PP^\gamma\left\{X_t\in G\right\}$ and using the assumption $\mu\leq 0$ (componentwise) yields \begin{eqnarray} c{\nu_{0,1}}^2\left\{(x,z)\in{\RR_+}^2  : x+z\geq \alpha\right\}&\leq &\PP^\gamma\left\{X_t\in G\right\} \nonumber \\ &\leq& 2c{\nu_{0,1}}^2\left\{(x,z)\in{\RR_+}^2  : x+z\geq \alpha\right\} \label{XtinG probability 2}\end{eqnarray} for $$\alpha:=\frac{\mu_1}{\sigma_1}t^{1/2} +c \frac{\mu_2}{\sigma_2}t^{1/2}.$$

We can find the following bounds for the measure in the previous estimate:

\begin{lem} For all $\alpha\geq 0$ and $t>0$, 
$$\frac{1}{4\sqrt{2}} e^{-\alpha^2} \leq \nu_{0,1}\left\{(x,z)\in\RR^2 \ : \ x+z\geq \alpha, \quad x,z\geq 0\right\}\leq \frac{1}{\sqrt{2}}e^{-\alpha^2/4}$$
\end{lem}
\begin{proof} The elementary proof has two parts. Firstly, we observe that for all $\alpha,z\geq 0$ and $x\in[0,\alpha]$, $$ |z+\alpha-x|^2\geq z^2+(\alpha-x)^2 $$ which implies \begin{eqnarray} 0&\leq& \nonumber\int_0^\alpha e^{-x^2/2}\int_{\alpha - x}^\infty e^{-z^2/2}dz \ dx\\ &\leq &\int_0 ^\alpha e^{-\left(\frac{x^2}{2}-\alpha x+\frac{\alpha^2}{2}\right)-\frac{x^2}{2}}dx\cdot \int_0^\infty e^{-z^2/2}dz\nonumber \\&=& \int_0 ^\alpha e^{-\left(x-\frac{\alpha}{2}\right)^2 -\frac{\alpha^2}{4}}dx\cdot \int_0^\infty e^{-z^2/2}dz\nonumber \\&=& \int_{-\alpha/2}^{\alpha/2}e^{-x^2}dx\cdot e^{-\alpha^2/4}\sqrt{\frac{\pi}{2}}\nonumber  \\&\leq & \int_{-\infty}^{\infty}e^{-x^2}dx\cdot e^{-\alpha^2/4}\sqrt{\frac{\pi}{2}}=\frac{\pi}{\sqrt{2}}e^{-\alpha^2/4} \label{xbelowalphaintegral}. \end{eqnarray} Secondly, we have for all $\alpha\geq 0$ and $x\geq 0$, $$2x^2+2\alpha^2\geq |x+\alpha|^2\geq x^2+\alpha^2.$$ Thus, \begin{eqnarray*}&&\frac{\sqrt{\pi}}{2}e^{-\alpha^2}=\int_0^\infty e^{-x^2}{dx} \cdot e^{-\alpha^2} \\ &\leq& \int_0^\infty e^{-\left|x+\alpha\right|^2/2} \ dx \\ &\leq& \int_0^\infty e^{-x^2/2}dx\cdot e^{-\alpha^2/2}=\sqrt{\frac{{\pi}}{2}}e^{-\alpha^2/2}\end{eqnarray*} which via $$\int_\alpha^\infty e^{-x^2/2}\int_0^\infty e^{-z^2/2}dz \ dx= \sqrt{\frac{{\pi}}{2}}\int_0^\infty e^{-\left|x+\alpha\right|^2/2} \ dx $$ and (\ref{xbelowalphaintegral}) gives 
\begin{eqnarray*}&& \frac{\pi}{2\sqrt{2}}e^{-\alpha^2}+0\\ &\leq& \int_\alpha^\infty e^{-x^2/2}\int_0^\infty e^{-z^2/2}dz \ dx + \int_0^\alpha e^{-x^2/2}\int_{\alpha - x}^\infty e^{-z^2/2}dz \ dx\\ &\leq & \frac{\pi}{2} e^{-\alpha^2/2} + \frac{\pi}{\sqrt{2}}e^{-\alpha^2/4} \leq \pi\sqrt{2}e^{-\alpha^2/4}.\end{eqnarray*}  But 
\begin{eqnarray*} && \int_\alpha^\infty e^{-x^2/2}\int_0^\infty e^{-z^2/2}dz \ dx + \int_0^\alpha e^{-x^2/2}\int_{\alpha - x}^\infty e^{-z^2/2}dz \\ &=& \int_0^\infty e^{-x^2/2}\int_{\left(\alpha - x\right)\vee 0}^\infty e^{-z^2/2}dz \ dx  \\ &=& {2\pi}\cdot \nu_{0,1}\left\{(x,z)\in\RR^2 \ : \ x+z\geq \alpha, \quad x,z\geq 0\right\}\end{eqnarray*} from which the Lemma follows.

\end{proof}

This Lemma's inequalities admit by means of identity (\ref{XtinG probability 2}) the following conlusion: $$\forall t>0 \quad \frac{c}{4\sqrt{2}} e^{-\left(\frac{\mu_1}{\sigma_1}+c \frac{\mu_2}{\sigma_2}\right)^2t} \leq \PP^\gamma\left\{X_t\in G\right\} \leq  2\cdot \frac{c}{\sqrt{2}}e^{-\left(\frac{\mu_1}{\sigma_1}+c \frac{\mu_2}{\sigma_2}\right)^2\frac{t}{4}}.$$

By Theorem \ref{phase} and the formula $\sum_n q^n/n = \ln\frac{1}{1-q}$ for all $q\in(0,1)$ , we conclude, analogously to the deliberations in the proof of Theorem \ref{phase1d} that \begin{eqnarray*} && K\left(1-e^{-\left(r+\frac{\left(\frac{\mu_1}{\sigma_1}+c \frac{\mu_2}{\sigma_2}\right)^2}{4}\right) t}\right)^{{c}{\sqrt{2}}} \\ &\leq &\lim_{t\downarrow 0}\bar V_{G}^t(\gamma)-\bar V_{G}^s(\gamma) \\ &\leq & K\left(1-e^{-\left(r+\left(\frac{\mu_1}{\sigma_1}+c \frac{\mu_2}{\sigma_2}\right)^2\right) t}\right)^{\frac{c}{4\sqrt{2}}}.\end{eqnarray*} After applying de l'Hospitals rule to the bases of the powers on each side of this estimate, we get constants $C_0,C_1>0$ such that for all sufficiently small $s$, $$C_0 \cdot s^{{c}{\sqrt{2}}} \leq \bar V_{G}^0(\gamma)-\bar V_{G}^s(\gamma)\leq C_1\cdot s^{\frac{c}{4\sqrt{2}}}.$$

\end{ex}

\begin{ex}[a special Extended Black-Scholes Model] In this example we do not assume a multi-dimensional Black-Scholes model, but we presume the discounted price process vector $\tilde S$ to satisfy the stochastic differential equation $$d\tilde S_t=C\cdot \tilde S_t dt+D\cdot\tilde S_td B_t,$$ where $C,D:\Omega\rightarrow \RR^{d\times d}$ are mutually commuting symmetric random matrices and $B$ is a one-dimensional Brownian motion, subject to the initial condition $$\tilde S_0=e^x.$$ Then, due to Albeverio and Steblovskaya \cite[Proposition 4]{AS02}, we have got an explicit solution of that stochastic differential equation, given by \begin{eqnarray} \forall t\in[0,T]\quad \tilde S_t&=&\exp\left(t\cdot\left(C-\frac{1}{2}D^2+ D\cdot B_t\right)\right)\cdot e^x\nonumber \\ &=& \exp\left(t\left(C-\frac{1}{2}D^2\right)\right)\exp\left(tB_t\cdot D\right)\cdot e^x.\label{AS}\end{eqnarray} Then $\ln \tilde S$ and hence the logarithmic non-discounted process $\ln \tilde S +r\cdot$ are L\'evy processes, and thus Theorem \ref{phase} applies. In this setting we can compute the expression in the last line (\ref{AS}) by applying the Spectral Theorem to the symmetric matrices $C-\frac{1}{2}n D^2$ and $tB_t\cdot D$.

\end{ex}

\chapter{From perpetual to non-perpetual Bermudan barrier options}

Recall how the function $U^t(T):=U^t_G(T):\RR^d\rightarrow\RR$, the expected payoff of a non-perpetual Bermudan option on a Feller basket with validity $T$, log-price payoff function $g$ and exercise mesh $h$ as a function of the logarithmic start price vector, given that the option is exercised on the first entry into $G\subset \RR^d$, was defined: $$U^t_G(T):x \mapsto \EE^x\left[e^{-r\left(\tau_G^h\wedge T\right)}g\left(X_{\tau_G^h\wedge T}\right)\right].$$

The purpose of the following Lemma \ref{UVdifference} is to see see that for all $y\not\in G$, the limiting behaviour of the difference $U_G^t(Nt)(y)- U_G^s(Ns)(y)$ as $t$ tends to zero whilst $Nt=Ns$ remains constant must be the same as the one of the difference $V_G^t(y)-V_G^s(y)$. In words: In the (sub-optimal case) of a non-stationary exercise policy for a non-perpetual option, the American-Bermudan barrier option price difference has the same limiting behaviour as the American-Bermudan difference for the corresponding perpetual barrier options.

\begin{lem} \label{UVdifference} Suppose $X$ is a $d$-dimensional Feller basket with $\PP^\cdot $ and $r>0$ being an associated family of risk-neutral probability measures and discount rate, respectively. If we define $$\forall h>0 \forall k\in\NN \quad H_{k,h} :=\left\{\tau_G^h=kh\right\}=\left\{X_{kh}\in G \right\}\cap \bigcap_{\ell =0}^{k-1}\left\{X_\ell\not\in G\right\},$$ then one has for all $h>0$, $N\in\NN$, measurable $G\subset \RR^d$ and $x\in\complement G$, \begin{eqnarray*} &&\bar V_G^h(x) -U_G^h(Nh)(x)\\ &=& {e^{-rNh}}\PP^x\left[\bigcap_{j=0}^N\left\{X_{jh}\not\in G\right\}\right]\cdot\int_{\complement G}\bar V_G^h(y) \PP_{X_{Nh}}^x(dy)\\ &&- e^{-rNh}\PP^x\left[\complement \bigcup_{k=0}^N H_{k,h}\right] \cdot \EE^x\left[ g\left(X_{Nh}\right) \left|\left\{X_{Nh}\in\complement G\right\} \right.\right].\end{eqnarray*}
\end{lem}

\begin{proof} For all $N\in\NN$, $x\in\RR^d$, $h>0$ we can use the Markov property of the Feller process $X$ and the definition of the sequence of events $\left(H_{k,h}\right)_{k\in\NN_0}$ to obtain the following expressions for $U_G$ and $V_G$:
\begin{eqnarray*} U_G^h(Nh)(x) &=& \sum_{k=0}^N e^{-rkh}\PP^x H_{k,h}\cdot\EE^x\left[ \left. g\left(X_{kh}\right) \right| H_{k,h} \right] \\ && + e^{-rNh}\PP^x\left[\complement \bigcup_{k=0}^N H_{k,h}\right] \cdot\EE^x\left[ g\left(X_{Nh}\right) \left| \complement \bigcup_{k=0}^N H_{k,h} \right.\right] \\ &=& \sum_{k=0}^N e^{-rkh}\PP^x H_{k,h}\cdot\EE^x\left[ \left. g\left(X_{kh}\right) \right| \left\{X_{kh} \in G \right\}\right] \\ && + e^{-rNh}\PP^x\left[\complement \bigcup_{k=0}^N H_{k,h}\right] \cdot\EE^x\left[ \left. g\left(X_{Nh}\right) \right| \left\{X_{Nh} \in \complement G \right\}\right] ,\end{eqnarray*} 
as well as \begin{eqnarray*} \bar V_G^h(x) &=& \sum_{k=0}^\infty e^{-rkh}\PP^x \left[H_{k,h}\right]\cdot\EE^x\left[ \left. g\left(X_{kh}\right) \right| H_{k,h} \right] \\ &=& \sum_{k=0}^\infty e^{-rkh}\PP^x \left[H_{k,h}\right]\cdot\EE^x\left[ \left. g\left(X_{kh}\right) \right| \left\{X_{kh} \in G \right\}\right] .\end{eqnarray*} 
Then immediately for all $h$, $N$, $x$, \begin{eqnarray*} \bar V_G^h(x) &=& U_G^h(Nh)(x) - e^{-rNh}\PP^x\left[\complement \bigcup_{k=0}^N H_{k,h}\right] \cdot \EE^x\left[ \left. g\left(X_{Nh}\right) \right| \left\{X_{Nh} \in \complement G \right\}\right] \\ && +\sum_{k=N+1}^\infty e^{-rkh}\PP^x \left[H_{k,h}\right]\cdot\EE^x\left[ \left. g\left(X_{kh}\right) \right| \left\{X_{kh} \in G \right\}\right] \end{eqnarray*} 
Regarding the difference between $\bar V_G$ and $U_G$, observe that again for all $h>0$, $N\in\NN$ and $x\in\RR^d$, \begin{eqnarray*} && \sum_{k=N+1}^\infty e^{-rkh}\PP^x \left[H_{k,h}\right]\cdot\EE^x\left[ \left. g\left(X_{kh}\right) \right| \left\{X_{kh} \in G \right\}\right] \\ &=& \sum_{k=N+1}^\infty e^{-rkh}\cdot\EE^x\left[ \left. g\left(X_{kh}\right) \right., H_{k,h}\right] \\ &=& e^{-rNh}\sum_{\ell=1}^\infty e^{-r\ell h}\cdot\EE^x\left[ \left. g\left(X_{(\ell+N)h}\right) \right., H_{(\ell+N),h}\right] \\ &=& e^{-rNh}\sum_{\ell=1}^\infty e^{-r\ell h}\cdot\EE^x\left[ \left. g\left(X_{\ell h}\circ\theta_{Nh} \right) \cdot\chi_{H_{\ell,h}}\circ\theta_{Nh} \right. , \bigcap_{j=0}^N\left\{X_{jh}\not\in G\right\}\right]  \\ &=& e^{-rNh}\sum_{\ell=1}^\infty e^{-r\ell h}\cdot\EE^{X_{Nh}}\left[ \left. g\left(X_{\ell h}\right) \cdot\chi_{H_{\ell,h}} \right., \bigcap_{j=0}^N\left\{X_{jh}\not\in G\right\}\right] \\ &=& e^{-rNh}\EE^x\left[ \left. \bar V_G^h\left({X_{Nh}}\right)\right., \bigcap_{j=0}^N\left\{X_{jh}\not\in G\right\}\right] \\ &=& e^{-rNh}\PP^x\left[\bigcap_{j=0}^N\left\{X_{jh}\not\in G\right\}\right]\cdot \EE^x\left[ \left. \bar V_G^h\left({X_{Nh}}\right)\right| \left\{X_{Nh}\not\in G\right\}\right]
,\end{eqnarray*}where again we have exploited several times the Markov property of $X$, and -- in order to interchange $\EE$ and $\sum$ -- the assumption that $g\geq 0$ on $G$.

\end{proof}

This proves

\begin{cor} \label{UVdifferenceCor} Under the assumptions of Lemma \ref{UVdifference} as well as $g=K-f$ for nonnegative $f:\RR^d\rightarrow\RR$ and $g\geq 0$ on $G$, we have the identity \begin{eqnarray*} &&\lim_{t\downarrow 0} U_G^{t}(T)(x)- U_G^{T\cdot{2^{-n}}}(T)(x)\\ &=& \lim_{t\downarrow 0} \bar V_G^t(x) - \bar V_G^{T\cdot{2^{-n}}}(x) \\ && - {e^{-rT}} \PP^x\left[\bigcap_{j=0}^{2^n}\left\{X_{j\cdot 2^{-n}T}\not\in G\right\}\right] \int_{\complement G}\left(\lim_{t\downarrow 0} \bar V_G^t(y) - \bar V_G^{T\cdot{2^{-n}}}(y)\right) \PP_{X_{T}}^x(dy)\\ &&+ \left[\PP^x\left[\bigcap_{j=0}^{2^m}\left\{X_{j\cdot 2^{-m}T}\not\in G\right\}\right]\right]_{m=n}^{m\rightarrow\infty} \int_{\complement G}\lim_{t\downarrow 0} \bar V_G^t(y)\PP_{X_{T}}^x(dy) \\ && + e^{-rT}\left[\PP^x\left[\complement \bigcup_{k=0}^{2^m} H_{k,2^{-m}T}\right]\right]_{m=n}^{m\rightarrow\infty} \cdot \underbrace{\EE^x\left[ g\left(X_{T}\right) \left|\left\{X_T\in\complement G\right\} \right.\right]}_{\leq K}.\end{eqnarray*}
\end{cor}

\begin{rem} Put informally, this Corollary \ref{UVdifferenceCor} means that as soon as one has established order estimates (in $n$) on the difference $\lim_{t\downarrow 0} \bar V_G^t(x) - \bar V_G^{T\cdot{2^{-n}}}(x)$ (for instance the ones from Theorem \ref{phase1d}), one only needs to find estimates on the probabilities $\left[\PP^x\left[\bigcap_{j=0}^{2^m}\left\{X_{j\cdot 2^{-m}T}\not\in G\right\}\right]\right]_{m=n}^{m\rightarrow\infty} $ and $\left[\PP^x\left[\complement \bigcup_{k=0}^{2^m} H_{k,2^{-m}T}\right]\right]_{m=n}^{m\rightarrow\infty}$ to obtain order estimates on the difference $\lim_{t\downarrow 0} U_G^{t}(T)(x)- U_G^{T\cdot{2^{-n}}}(T)(x)$.
\end{rem}

\part{Convergence of some approximate pricing algorithms}

\chapter{Bermudan option pricing based on piecewise harmonic interpolation and the r\'eduite}

\section{Introduction}

We intend to approximate the function that assigns the value of a Bermudan option with payoff function $g$ and no dividends to the logarithmic start prices of the underlying assets by piecewise harmonic functions. In the first step, we will compute a piecewise harmonic approximation to the function that assigns the European option price associated with $g$ and the Bermudan's maturity $T>0$ to the logarithmic asset prices at the penultimate time $T-t$ where exercise is possible. Then we iteratively compute the expectation of this function after time $t$, discount, take the maximum with the payoff function $g$, and perform a r\'eduite-based interpolation (in the one-dimensional setting: a piecewise harmonic interpolation).

Now we would like to answer the following questions: Given the stationarity of perpetual Bermudan option prices, can we prove that there exists a minimal fixed point of the iteration step described above (which would then be an approximation to the perpetual Bermudan price)? If so, can we characterise it explicitly? Is the iteration step monotone?

First, we will discuss these questions in the one-dimensional setting -- very little knowledge of potential theory has to be assumed for the proofs in that section. Second, we shall generalise that approach to higher dimensions; this will entail a few technical subtleties.

\section{Piecewise harmonic Bermudan option pricing for options on one asset}

Consider $\{a_0,\dots,a_m\}\subset\RR$, the {\em set of (mutually distinct) support abscissas}, and let $L:C^\infty(\RR,\RR)\rightarrow C^\infty(\RR,\RR)$ be the infinitesimal generator of a Markov semigroup of operators on Lebesgue measurable functions from $\RR$ to $\RR$. We call a function $f\in C^\infty(\RR,\RR)$ {\em $P$-harmonic} (or shorter: {\em harmonic}, if no ambiguity can arise) if and only if $Lf=0$. Let $(P_t)_{t\geq 0}$ denote the semigroup generated by $L$.

For the following, assume $L$ to be a {\em second-order differential operator}, that is, there are constants $\beta_1,\beta_2\in\RR$ such that $$L:f\mapsto \beta_1 f' + \beta_2 f''.$$

A function $g:\RR\rightarrow \RR$ is said to be {\em subharmonic} ({\em superharmonic}) if and only if $g$ is right- and left-differentiable (thus, letting $Lg$ become well-defined as a function from $\RR$ to $\RR\cup\{\pm\infty\}$) and $Lg\geq 0$ ($Lg\leq 0$, respectively). 
 
In particular, the supremum (infimum) of countably many harmonic functions is subharmonic (superharmonic). 

\begin{lem} \label{uniqueinterpol} Given two support abscissas and ordinates, there is a unique harmonic interpolation, provided $L$ is a second-order differential operator with a non-trivial second-order part (i e $\beta_2\neq 0$) or a non-zero first-order part (i e $\beta_2\neq 0$).
\end{lem}
\begin{proof}[Proof sketch] The uniqueness is a consequence of the maximum principle for harmonic functions. The existence follows (in our one-dimensional setting) by distinguishing the cases delineated in the statement of the Lemma. If $L$ is a second-order operator and it has only a non-zero term of second order, then the space of solutions are all affine-linear functions from $\RR$ to $\RR$. This space is two-dimensional. If there are terms of different order, the space of solutions will have basis elements of the form $\exp(\alpha\cdot)$ and we have to solve a linear or quadratic equation to find the $\alpha$ (or $\alpha$'s) satisfying this linear or quadratic equation. Since $L$ is sub-Markovian, there will be at least one real solution to this equation for $\alpha$. 
\end{proof}

The Lemma implies 

\begin{cor}\label{2dimharmonic}There cannot be more than two linearly independent harmonic functions: There is a canonical monomorphism from the space of functions to the -- two-dimensional -- space of pairs of subordinates.
\end{cor}

\begin{lem} \label{threezeroes}A subharmonic function from $\RR$ to $\RR$ is constantly zero if it has three zeros.
\end{lem}
\begin{proof} The left- and right-differentiablility of subharmonic functions entail that for all subharmonic $g$, $Lg$ will be defined as a function from $\RR$ to $\RR\cup\{\pm\infty\}$.
\end{proof}

If there is only a first order non-zero term, the space of harmonic functions will just coincide with the space of constant functions.

\begin{lem}\label{subinterpol} 
\begin{enumerate}
\item Piecewise harmonic interpolation with respect to the support abscissas $\{a_0,\dots,a_m\}$ preserves subharmonicity on $[a_0,a_m]$: The interpolating function dominates the interpolated function on $A:=[a_0,a_m]$, and if the interpolating function $f$ equals the harmonic function $f_i$ on $[a_i,a_{i+1}]$ for all $i<m$, then we have $f=\sup\{f_0,\dots,f_{m-1}\}$. 
\item The interpolating function $f$ is strictly dominated by the interpolated function $\cI(f)$ on the intervals $(-\infty,a_0)$ and $(a_m,+\infty)$.
\end{enumerate}
\end{lem}
\begin{proof}[Proof sketch]
\begin{enumerate} 
\item The domination part follows from the maximum principle for harmonic functions. From the maximum principle, we also get for all $i<m-1$ that if $f_i\neq f_{i+1}$, then $$\left\{f_i=f_{i+1}\right\}=\left\{a_{i+1}\right\}.$$ Now there are two possibilities: either $f_i< f_{i+1}$ on $(-\infty,a_{i+1})$ and $f_i> f_{i+1}$ on $(a_{i+1},+\infty)$ or the other way round $f_i> f_{i+1}$ on $(-\infty,a_{i+1})$ and $f_i< f_{i+1}$ on $(a_{i+1},+\infty)$. However, in the former case, the interpolating function would equal $f_i\wedge f_{i+1}$ on $[a_i,a_{i+2}]$, which is superharmonic, and it would also dominate the subharmonic interpolated function $\eta$ on $[a_i,a_{i+2}]$. Then, $\eta-\left(f_i\wedge f_{i+1}\right)$ would be nonpositive and subharmonic on $[a_i,a_{i+2}]$ and it would have three zeroes, in $a_i$, $a_{i+1}$ and $a_{i+2}$. By Lemma \ref{threezeroes}, this can only be true if $\eta-\left(f_i\wedge f_{i+1}\right)=0$ on $[a_i,a_{i+2}]$. Thus, $\eta=f_i\wedge f_{i+1}$ on $[a_i,a_{i+2}]$. Since $\eta$ is subharmonic on $[a_i,a_{i+2}]$, so must be $f_i\wedge f_{i+1}$ then, and therefore, $f_i\wedge f_{i+1}$ is harmonic on $[a_i,a_{i+2}]$. This means $f_i=f_{i+1}$ (as both $f_i$ and $f_{i+1}$ are harmonic) which contradicts our assumption that $f_i\neq f_{i+1}$. Therefore, $f_i> f_{i+1}$ on $(-\infty,a_{i+1})$ and $f_i< f_{i+1}$ on $(a_{i+1},+\infty)$ for all $i<m$. 

Inductively, this yields $f\geq f_i$ on $[a_0,a_m]$ for all $i<m$, hence $f=\sup\{f_0,\dots,f_{m-1}\}$ on $[a_0,a_m]$.
\item The function $f-\cI(f)$ is subharmonic on $(-\infty,a_1]$ and it has two zeroes in $a_0$ and $a_1$. Moreover, it is nonpositive on $(a_0,a_1)$. Because of Lemma \ref{threezeroes}, then $f-\cI(f)$ has to be positive or negative on $(-\infty,a_0)$. In the former case, we are done. In the latter case, due to the maximum principle, $f-\cI(f)$ must be decreasing and therefore in $a_0$ we would have $L\left(f-\cI(f)\right)(a_0)<0$, which is absurd.
A symmetric argument works for the proof of the domination of $\cI(f)$ by $f$ on the interval $(a_m,+\infty)$.
\end{enumerate}
\end{proof}

\begin{lem}\label{monointerpol} Piecewise harmonic interpolation to a set of support absicssas $\{a_0,\dots,a_m\}$ is monotone on $[a_0,a_m]$ in the sense that if $f\leq g$ on $[a_0,a_m]$, then the piecewise harmonic interpolation of $f$ will be dominated by the piecewise harmonic interpolation of $g$ on $[a_0,a_m]$.
\end{lem}
\begin{proof} Use the maximum principle on each of the intervals $[a_i,a_{i+1}]$ for $i<m$.
\end{proof}

\begin{lem}\label{monointerpolsub} Let $\cI:\RR^{[a_0,a_m]}\rightarrow \RR^\RR$ denote the operator of piecewise harmonic interpolation with respect to the set of support abscissas $\{a_0,\dots,a_m\}$. Let $f:\RR\rightarrow\RR$ be subharmonic on $\RR$. Consider a harmonic function $h:\RR\rightarrow\RR$, assumed to dominate $f$: $f\leq h$ on $\RR$. Then $\cI(f)\leq h$ on $\RR$.
\end{lem}
\begin{proof} From the previous Lemma \ref{monointerpol}, we already know that $\cI(f)(x)\leq\cI(h)(x)$ holds for all $x\in [a_0,a_m]$. However, $\cI(h)=h$, hence $\cI(f)\leq h$ on $[a_0,a_m]$ and from Lemma \ref{subinterpol}, we conclude that $h\geq f\geq\cI(f)$ on the intervals $(-\infty,a_0)$ and $(a_m,+\infty)$.
\end{proof}

\begin{Th}\label{fixexist} Let $\cI:\RR^{[a_0,a_m]}\rightarrow \RR^\RR$ again denote the operator of piecewise harmonic interpolation with respect to the set of support abscissas $\{a_0,\dots,a_m\}$. Let $g$ be a subharmonic function, let $c$ be nonnegative and subharmonic, and let $h$ be harmonic. Let $c$ be, moreover, harmonic on each of the intervals $[a_i,a_{i+1}]$ for $i<m$. Suppose $c\leq g$ on $[a_0,a_m]$ and $c,g\leq h$ on $\RR$, $r>0$ and let $t>0$. Now define $$\cK:f\mapsto\cI\left(e^{-rt}P_t\left(\cI(f)\vee c\right)\vee g\right)\restriction [a_0,a_m]$$ as well as $$Q:=\left\{f\restriction [a_0,a_m]\ : \ \begin{array}{c} f:\RR\rightarrow\RR \text{ subharmonic}, \quad f\geq c \text{ on }[a_0,a_m], \\ \forall i\in\{1,\dots,m-2\} \quad f\text{ harmonic on }[a_i,a_{i+1}], \\ f\text{ harmonic on }(-\infty,a_1),(a_{m-1},+\infty), \quad f\leq h \end{array}\right\}.$$ Then $\cK$ maps the convex and bounded subset $Q$ of $C^0[a_0,a_m]$ continuously to itself. Moreover, due to Lemma \ref{uniqueinterpol}, $Q$ is a subset of a finite-dimensional subspace of $C^0[a_0,a_m]$ (this subspace being the space of all functions from $[a_0,a_m]$ that are harmonic on each of the intervals $[a_i,a_{i+1}]$ for $i<m$. By Brouwer's Fixed Point Theorem, $\cK$ has got a fixed point in $Q$. Finally, $\cK$ is a composition of monotone functions on $[a_0,a_m]$ and therefore monotone as well.
\end{Th}
\begin{proof}[Proof sketch] We can divide the proof for $\cK(Q)\subseteq Q$ into three parts:
\begin{enumerate}
\item The cone of subharmonic functions is closed under $\vee$, under $P_t$, under multiplication by constants and under piecewise harmonic interpolation $\cI$ (cf Lemma \ref{subinterpol}), therefore the image of $Q$ under $\cK$ can only consist of subharmonic functions. 
\item The upper bound on the elements of the image $\cK(Q)$ follows from the monotonicity of $P_t$ and $\cI$ (Lemma \ref{monointerpol}), combined with the equations $P_th=h$ and $\cI(h)=h$ as well as the Lemma \ref{monointerpolsub}: First, we may state $e^{-rt}P_t\left(\cI(f)\vee c\right)\vee g\leq h $ for all $f\leq h$,  which by Lemma \ref{monointerpolsub} allows us to deduce $$\cI\left(e^{-rt}P_t\left(\cI(f)\vee c\right)\vee g\right)\leq h$$ for all $f\in Q$.
\item The lower bound follows again from the monotonicity of $\cI$, but this time only by exploiting $c\leq g$ on $[a_0,a_m]$ and employing the fact that the space of those functions that are harmonic on each of the intervals $[a_i,a_{i+1}]$ for $i<m$ is invariant under the composition of $\cI$ with the restriction to $[a_0,a_m]$ (yielding $\cI(c)=c$ on $[a_0,a_m]$). 
\end{enumerate}
Since $c$ is nonnegative, we get that $Q$ is bounded by $\sup_{[a_0,a_m]}h\geq 0$ as a subset of $C^0[a_0,a_m]$, and because $Q$ is finite-dimensional, we may apply Schauder's Theorem, provided we are given the continuity of $\cK$. However, this last assertion follows from the maximum principle.
\end{proof}

The existence of a minimal fixed point for $\cK$ can be proven constructively as well:

\begin{cor} Let us adopt the notation of the previous Theorem. Then the sequence $\left(\cK^n(g\vee 0)\right)_{n\in\NN_0}$ is monotone on $[a_0,a_m]$, bounded and dominated by $h$. Therefore we have the existence of a limit on $[a_0,a_m]$ given by $$\forall x\in[a_0,a_m]\quad q(x):=\lim_{n\rightarrow \infty}\cK^n(g\vee 0)(x)=\sup_{n\in\NN_0}\cK^n(g\vee 0)(x).$$ This limit is an element of $Q$ and therefore can be canonically extended to the whole of $\RR$. By the continuity of $\cK$, $q$ is a fixed point of $\cK$. On $[a_0,a_m]$, the convergence in the last equation will be uniform.
\end{cor}
\begin{proof} The only part of the Corollary that does not follow directly from the preceding Theorem \ref{fixexist} is the uniformity of the convergence and that $q$ will be harmonic on each of the intervals $[a_i,a_{i+1}]$ for $i<m$. However, monotone convergence on compact sets preserves harmonicity and is always uniform (cf e g Meyer \cite{M} -- or, more directly, Port and Stone \cite[Theorem 3.9]{PS} if $P$ is the Brownian semigroup).
\end{proof}

\begin{lem}\label{Lemma7.6}In the preceding Corollary's notation, $q$ is the minimal nonnegative fixed point of $\cK$.
\end{lem}
\begin{proof} The proof partly copies the one for Lemma \ref{soundmonotone}. Any nonnegative fixed point $p$ of $\cK$ must be greater or equal $g$ on $[a_0,a_m]$. Therefore the monotonicity of $\cK$ on $[a_0,a_m]$, implies $$\forall n\in \NN_0 \quad p=\cK^np\geq \cK^n(g\vee 0) \text{ on }[a_0,a_m],$$ yielding $$p\geq \sup_{n\in\NN_0}\cK^n(g\vee 0)=q\text{ on }[a_0,a_m].$$
\end{proof}

\begin{ex}[Bermudan vanilla call on a dividend-paying asset in a special Black-Scholes model]Assume $$P:=(P_t)_{t\geq 0}:=\left(\nu_{\mu t, \sigma^2t}\ast\cdot\right)_{t\geq 0},$$ where $$\sigma>0,\quad \mu:=r-\frac{\sigma^2}{2},$$ thus $P$ can be perceived as the semigroup associated to the logarithmic price process under the risk-neutral measure in the one-dimensional Black-Scholes model). We will assume that (possibly after a linear change of the time scale) $\sigma=1$ and we assume that $r$ has been cut to discount dividends. Define $$g:=\exp-K$$ (the payoff on exercise of a one-dimensional call option with strike price $K$). The infinitesimal generator of the Markov semigroup $P$ is $$L=\frac{1}{2}\Delta + \mu\nabla=\frac{1}{2}\Delta+\left(r-\frac{1}{2}\right)\nabla.$$ Thus we obtain $$Lg=r\dot \exp \geq 0,$$ hence $g$ is, $P$-subharmonic. We can find the $P$-harmonic functions for $\mu\neq 0$ (otherwise they are simply the affine linear functions) by observing that for all $\alpha\in\RR$, \begin{eqnarray*}&& 0=L\exp(\alpha\cdot)=\frac{1}{2}\left(\alpha^2+2\mu\alpha\right)\exp(\alpha)\\ &\Leftrightarrow& \alpha\in\left\{ 0, -2\mu\right\}. \end{eqnarray*} If $\mu\neq 0$, the functions $1:x\mapsto 1$ and $\exp(-2\mu\cdot):x\mapsto e^{-2\mu x}$ are two linearly independent harmonic functions, thus by Corollary \ref{2dimharmonic}, we have already found a basis for the space of harmonic functions. If $\mu=0$, the harmonic functions are exactly the affine linear functions. In order to obtain the setting of Theorem \ref{fixexist}, we will assume $\mu\leq \frac{1}{2}$ such that the sum $h$ of $\exp(2\mu\cdot)$ and a sufficiently large positive constant is a harmonic function dominating $g=\exp-K$. In order to satisfy the conditions on $c$ we could simply take $c=-K$ for instance. 
\end{ex}

\section{R\'eduite-based approximation of Bermudan option prices}

Suppose $P$ is a Markov semigroup on $\RR^d$ ($d\in\NN$) and $L$ is the infinitesimal generator of $P$. We will call a function $f:\RR^d\rightarrow \RR$ {\em subharmonic} if and only if $$\forall t> 0 \quad P_tf \geq f$$ holds pointwise. A function $f:\RR^d\rightarrow \RR$ will be called {\em superharmonic} if and only if $-f$ is subharmonic, and $f:\RR^d\rightarrow \RR$ will be called {\em harmomic} if it is both super- and subharmonic.

Let $\cU$ denote the operator of upper-semicontinuous regularisation, that is, for all functions $f:\RR^d\rightarrow \RR$, $$\cU f=\inf\left\{\ell \geq f \ : \ \ell:\RR^d\rightarrow\RR \text{ subharmonic}\right\}$$ (of course, this is a priori only defined as a function taking values in $\RR\cup\{--\infty\}$). Consider a harmonic function $h:\RR^d\rightarrow \RR$ and a closed (and therefore $F_{\sigma}$) set $B$ and define the {\em r\'eduite} operator $\cR=\cR_{h,B}$ on the set of all subharmonic functions $f:\RR^d\rightarrow \RR$ dominated by $h$ via $$\cR f:=\cU\left(\sup\left\{\ell \leq h \ : \ \ell:\RR^d\rightarrow\RR \text{ subharmonic}, \quad \ell \leq f \text{ on } B\right\}\right).$$ It is a well-known  result from potential theory (cf e g the work of Paul-Andr\'e Meyer \cite[Th\'eor\`eme T22]{M}) that there will be a greatest subharmonic function dominated by $f$ on $B$ and that this function will be equal to $\cR f$. Moreover, we have that $f=\cR f$ on $B$ except on a set of potential zero, in probabilistic/potential-theoretic jargon $$f=\cR f \text{ q.e. on }B,$$ where ``q.e.'' is, as usual, short-hand for ``quasi-everywhere''. Now define $$Q:=\left\{f\leq h \ : \ f:\RR^d\rightarrow\RR \text{ subharmonic}\right\}.$$ Then our definition of the r\'eduite operator $\cR$ implies $\cR f\leq h$ (as $h$ is dominating the function whose upper-semicntinuous regularisation is, according to our definition, the r\'eduite $\cR f$ of $f$) and our potential-theoretic characterisation of the r\'eduite -- as the greatest subharmonic function dominated by $f$ on $B$ -- ensures the subharmonicity of $\cR f$. Therefore, $$\cR:Q\rightarrow Q.$$ We also have that $\cU$ is monotone (in the sense that for all $f_0\leq f_1$, $\cU f_0\leq \cU f_1$) so that $\cR$ must be monotone as well (from the $\subseteq$-monotonicity of $\sup$ and the definition of $\cR$).

Hence

\begin{lem}\label{propR} Adopting the notation of the preceding paragaph, $\cR:Q\rightarrow Q$ and whenever $f_0\leq f_1$, $\cR f_0\leq \cR f_1$. 
\end{lem}

Let $g:\RR^d\rightarrow \RR$ be a subharmonic function such that $g\leq h$ and let $r>0$. The next step is going to be the consideration of the following family of operators: $$\phi_t:f\mapsto e^{-rt}P_tf\vee g$$ for $t\geq 0$. If $f\leq h$, $P_tf\leq P_t h=h$ for all $t\geq 0$, since the operators $P_t$ are positive and linear, and $h$ was assumed to be harmonic. Thus, since $g\leq h$ and $r>0$, one must have $\phi_t f\leq h$ for all $f\leq h$ and $t\geq 0$. Moreover, the operators $P_t$ preserve subharmonicity and the maximum of two subharmonic functions is subharmonic again, therefore $\phi_t f$ must be subharmonic for all subharmonic $f$. Finally, since $P_t$ is monotone, $\phi_t$ has to be monotone for all $t\geq 0$ Summarising this, we obtain

\begin{lem}\label{propphi} Using the notation introduced previously, $\phi_t:Q\rightarrow Q$ and whenever $f_0\leq f_1$, $\phi_t f_0\leq \phi_t f_1$ for all $t\geq 0$. 
\end{lem}

As a consequence, we derive from the two Lemmas \ref{propR} and \ref{propphi} the following:

\begin{cor} If we define $\cK_t:=\cR\circ\phi_t$ (adopting the notation of the previous paragraph), we have $\cK_t:Q\rightarrow Q$ and whenever $f_0\leq f_1$, $\cK_t f_0\leq \cK_t f_1$. 
\end{cor}

\begin{cor} The map $f\mapsto\cK_tf\vee 0$ is a sound iterative Bermudan option pricing algorithm for the payoff function $g\vee 0$ (in the sense of Definition \ref{algorithmclasses}).
\end{cor}

This already suffices to prove the following

\begin{Th} \label{Theorem7.2}Let $t\geq0$. Then for all $n\in\NN_0$, \begin{equation}\label{mono}{\cK_t}^{n+1}(g\vee 0)\geq {\cK_t}^{n}(g\vee 0).\end{equation} Furthermore, $$q:=\sup_{n\in \NN_0} {\cK_t}^{n}(g\vee 0)$$ (which a priori is only defined as a function with range in $\RR\cup\{+\infty\}$) is an element of $Q$ and indeed is the least nonnegative fixed point of $\cK_t$.
\end{Th}
\begin{proof}[Proof] \begin{enumerate}
\item Relation (\ref{mono}) follows from the fact that $\cK_t$ is a sound algorithm and Remark \ref{soundmonotone}.
\item Since $\cK_t$ maps $Q$ to itself, the whole sequence $\left({\cK_t}^{n}(g\vee 0)\right)_{n\in\NN_0}$ is bounded by $h$. This entails $q\leq h$ as well. Applying Beppo Levi's Theorem on swapping $\sup$ and $\int\cdot d\mu$ -- for bounded monotonely increasing sequences of measurable nonnegative functions and an arbitrary measure $\mu$ -- to the measures $P_t(\cdot, x)$, $x\in\RR^d$ and the sequence $\left({\cK_t}^{n}(g\vee 0)\right)_{n\in\NN_0}$, we can exploit the subharmonicity of the functions ${\cK_t}^{n}(g\vee 0)$, ${n\in\NN_0}$, to deduce \begin{eqnarray*}\forall x\in\RR^d\quad P_tq(x)&=&\sup_{n\in\NN_0} P_t\left({\cK_t}^{n}(g\vee 0)\right)(x)\\ &\geq& \sup_{n\in\NN_0}{\cK_t}^{n}(g\vee 0)(x) =q(x),\end{eqnarray*} which is the subharmonocity of $q$. As we have already seen, $q\leq h$, so $q\in Q$.
\item If we employ Beppo Levi's Theorem again, we can show that $\cK_t$ and $\sup_{n\in\NN_0}$ commute for bounded monotonely increasing sequences of functions. Thereby $$\cK_tq=\sup_{n\in\NN_0}\cK_t {\cK_t}^{n}(g\vee 0)=\sup_{n\in\NN}{\cK_t}^{n}(g\vee 0)=q.$$
\item That $q$ is the least nonnegative fixed point is seen as in the proof of Lemma \ref{inffixedpoint}. Any nonnegative fixed point $p$ of $\cK_t$ must be greater or equal $g\vee 0$. Therefore by the monotonicity of $\sup$ and $\cK_t$, $$\sup_{n\in\NN_0}{\cK_t}^{n}p\geq \sup_{n\in\NN_0}{\cK_t}^{n}(g\vee 0)=q.$$
\end{enumerate}
\end{proof}

\begin{ex}[Bermudan call option with equidistant exercise times in $t\cdot \NN_0$ on the weighted arithmetic average of a basket in a special Black-Scholes model] Let $\beta_1,\dots,\beta_d\in[0,1]$ be a convex combination and for simplicity, assume that the assets in the basket are independent and each follow the Black-Scholes model with one and the same volatility $\sigma_1=\dots=\sigma_d=:\sigma$, and let $r>0$ be the interest rate of the bond. We may assume that, possibly after a linear change of the time-scale, $\sigma=1$. Then $\left(P_t\right)_{t\geq 0}=\left(\nu_{{^t}\left(r-\frac{1}{2}\right)_{i=1}^dt,t}\ast\cdot\right)_{t\geq 0}$ is the semigroup of this Markov (even L\'evy) basket. Then one has $$L=\frac{1}{2}\Delta + \left(r-\frac{1}{2}\right)_{i=1}^d\cdot\nabla$$ (cf e g Revuz and Yor's exposition \cite{RY}), and for $$g:x\mapsto \sum_{i=1}^d \beta_i\exp\left(x_i\right)-K$$ we obtain $$Lg=\sum_{i=1}^d\left(\frac{{\beta_i}^2}{2} +\left(r-\frac{\beta_i}{2}\right)\right)\exp \left((\cdot)_i\right)$$ which is pointwise nonnegative if and only if  $$r\geq \frac{\max_{i\in\{1,\dots,d\}}{\beta_i}}{2}=\frac{1}{2d}+\frac{\sum_{i=1}^d {\beta_i}^2}{2d}. $$ Hence, if $r$ is sufficiently large, $g$ is subharmonic and we can apply the theory developed earlier in this Chapter, in particular Theorem \ref{Theorem7.2}.
\end{ex}

\chapter{Soundness and convergence rate of perpetual Bermudan option pricing via cubature}

When Nicolas Victoir studied ``asymmetric cubature formulae with few points'' \cite{V} for symmetric measures such as the Gaussian measure, the idea of (non-perpetual) Bermudan option pricing via cubature in the log-price space was born. In the following, we will discuss the soundness and convergence rate of this approach when used to price perpetual Bermudan options.

Consider a convex combination $(\alpha_1,\dots,\alpha_m)\in[0,1]^d$ (that is, $\sum_{k=1}^m \alpha_k = 1$) and $x_1,\dots,x_m\in\RR^d$. Then there is a canonical {\em weighted arithmetic average operator} $A$ associated with $\vec \alpha, \vec x$ given by $$\forall f\in\RR^\RR\quad Af=\sum_{k=1}^m\alpha_k f(\cdot-x_k).$$ Now suppose $c\in(0,1)$, $g,h:\RR^d\rightarrow \RR$, $Ag\geq g$, $Ah=h$ and $0\vee g\leq h$. Define an operator $\cD$ on the cone of nonnegative measurable functions by $$\cD:f\mapsto \left(c\cdot Af\right)\vee g.$$ Since $A$ is positive and linear, thus monotone (in the sense that for all $f_0\leq f_1$, $A f_0\leq A f_1$), it follows that $\cD$ must be monotone as well. Furthermore, whenever $Af\geq f$, we have that $A\cD f\geq \cD f$, as the linearity and positivity of $A$ combined with our assumption on $g$ imply $$A\cD f \geq cA^2f\vee g\geq \left(c\cdot Af\right)\vee g =\cD f.$$ Finally, due to our assumptions on $h$ and $g$, we have for all nonnegative $f\leq h$, $$\cD f\leq cAh\leq h.$$ Summarising this, we are entitled to state

\begin{lem} \label{discreteLemma}Adopting the previous paragraph's notation and setting $$Q:=\left\{f\leq h \ : \ Af\geq h\right\},$$ we have that $$\cD:Q\rightarrow Q,$$ $\cD$ is monotone (i e order-preserving), and $A\cD -\cD$ is nonnegative.
\end{lem}

\begin{cor} The map $f\mapsto\cD f\vee 0$ is a sound iterative Bermudan option pricing algorithm for the payoff function $g\vee 0$ (in the sense of Definition \ref{algorithmclasses}).
\end{cor}

This is sufficient to prove 

\begin{Th} \label{monoconv}For all $n\in\NN_0$, \begin{equation}\label{discretemono}\cD^{n+1}(g\vee 0)\geq \cD^n(g\vee 0)=:q_n.\end{equation} Furthermore, $$q:=\lim_{n\rightarrow \infty}\cD^n(g\vee 0)=\sup_{n\in\NN_0}\cD^n(g\vee 0)\in Q$$ and $q$ is the smallest nonnegative fixed point of $\cD$.
\end{Th}
\begin{proof}
\begin{enumerate}
\item Relation (\ref{mono}) follows from the soundness of $\cD$ is a sound algorithm and Remark \ref{soundmonotone}.
\item Since $\cD$ maps $Q$ itself, the whole sequence $\left({\cD}^{n}(g\vee 0)\right)_{n\in\NN_0}$ is bounded by $h$. This entails $q\leq h$ as well. Using the linearity of $\sup$ and our previous observation that $A\cD-\cD\geq 0$ (Lemma \ref{discreteLemma}), we can show \begin{eqnarray*}\forall y\in\RR^d\quad Aq(y)&= &\sup_{n\in\NN_0} A\left({\cD}^{n}(g\vee 0)\right)(y)\\ &\geq& \sup_{n\in\NN_0}{\cD}^{n}(g\vee 0)(y) =q(y),\end{eqnarray*} which means $Aq\geq g$. As we have already seen, $q\leq h$, so $q\in Q$.
\item Again, due to the linearity of $\sup$ and the special shape of $\cD$ that is based on a weighted arithmetic average operator, $\cD$ and $\sup_{n\in\NN_0}$ commute for bounded monotonely increasing sequences of functions. Thereby $$\cD q=\sup_{n\in\NN_0}\cD{\cD}^{n}(g\vee 0)=\sup_{n\in\NN} {\cD}^{n}(g\vee 0)=q.$$
\item Just as in the proof of Lemma \ref{inffixedpoint}, we see that $q$ is the minimal nonnegative fixed point. For, any nonnegative fixed point $p$ of $\cD$ must be greater or equal $g\vee 0$. Thus, by the monotonicity of $\sup$ and $\cD$, $$\sup_{n\in\NN_0}{\cD}^{n}p\geq \sup_{n\in\NN_0}{\cD}^{n}(g\vee 0)=q.$$
\end{enumerate}
\end{proof}

\begin{lem} \label{gLemma} Using the previous Theorem's notation, we have for all $x\in\RR^d$ and $n\in\NN_0$, if $q_{n+1}(x)=g(x)$, then $q_n(x)=g(x)$.
\end{lem}
\begin{proof} By the monotonicity of the sequence $(q_n)_{n\in\NN_0}$ (Theorem \ref{monoconv}), we have $$g(x)\leq q_0(x)\leq q_n(x)\leq q_{n+1}(x).$$
\end{proof}

\begin{Th} For all $n\in\NN$, $$\left\|q_{n+1}-q_n\right\|_{C^0\left(\RR^d,\RR\right)}\leq c\cdot \left\|q_{n}-q_{n-1}\right\|_{C^0\left(\RR^d,\RR\right)}.$$
\end{Th}

\begin{proof} The preceding Lemma \ref{gLemma} yields \begin{eqnarray*} \left\|q_{n+1}-q_n\right\|_{C^0\left(\RR^d,\RR\right)}&=& \left\|q_{n+1}-q_n\right\|_{C^0\left(\left\{c\cdot Aq_{n}>g\right\},\RR\right)}\\ &=&\left\|c\cdot Aq_{n}-\left(\left(c\cdot A q_{n-1}\right)\vee g\right)\right\|_{C^0\left(\left\{c\cdot Aq_{n}>g\right\},\RR\right)}\end{eqnarray*} via the definition of $q_{i+1}$ as $\left(cAq_i\right)\vee g$ for $i=n$ and $i=n+1$. But the last equality implies \begin{eqnarray*}\left\|q_{n+1}-q_n\right\|_{C^0\left(\RR^d,\RR\right)}&\leq&\left\|c\cdot Aq_{n}-c\cdot A q_{n-1}\right\|_{C^0\left(\left\{c\cdot Aq_{n}>g\right\},\RR\right)}\\ &\leq&\left\|c\cdot Aq_{n}-c\cdot A q_{n-1}\right\|_{C^0\left(\RR^d,\RR\right)}.\end{eqnarray*} Since $A$ is linear as well as an $L^\infty$-contraction (and therefore a $C^0$-contraction, too), we finally obtain $$ \left\|q_{n+1}-q_n\right\|_{C^0\left(\RR^d,\RR\right)}\leq c\left\|A\left(q_{n}-q_{n-1}\right)\right\|_{C^0\left(\RR^d,\RR\right)}\leq c\left\|q_{n}-q_{n-1}\right\|_{C^0\left(\RR^d,\RR\right)}.$$
\end{proof}

\begin{ex}[Bermudan put option with equidistant exercise times in $t\cdot \NN_0$ on the weighted arithmetic average of a basket in a discrete Markov model with a discount factor $c=e^{-rt}$ for $r>0$] Let $\beta_1,\dots,\beta_d\in[0,1]$ be a convex combination and assume that $A$ is such that \begin{equation}\label{discretecondition}\forall i\in \{1,\dots,d\}\quad \sum_{k=1}^m\alpha_ke^{-(x_k)_i}=1,\end{equation} then the functions $$g:x\mapsto K-\sum_{i=1}^d\beta_i \exp\left(x_i\right)$$ and $h:=K$ (where $K\geq 0$) satisfy the equations $Ah=h$ and $Ag=g$, respectively. Moreover, by definition $g\leq h$. Then we know that the (perpetual) Bermudan option pricing algorithm that iteratively applies $\cD$ to the payoff function $g\vee 0$ on the $\log$-price space, will increase monotonely and will have a limit which is the smallest nonnegative fixed point of $\cD$. Moreover, the convergence is linear and the contraction rate can be bounded by $c$.

The condition (\ref{discretecondition}) can be achieved by a change of the time scale (which ultimately leads to different cubature points for the distribution of the asset price)
\end{ex}

One might also be interested in determining the convergence rate for the approximation of non-perpetual American option prices based on non-perpetual Bermudan option pricing via cubature. After proving a series of Lemmas we will end up with a Theorem that asserts linear convergence and also provides bounds for the convergence factor. 

From now on, $c$ and $A$ will no longer be fixed but their r\^ole will be played by $e^{-rt}$ and $P_t$ (for $t\in s\NN_0$ where $s>0$ shall be fixed) respectively, where $r>0$ and $\left(P_{s\cdot m}\right)_{m\in\NN_0}$ describes a Markov chain on $\RR^d$ (By the Chapman-Komogorov equation this is tantamount to $\forall s,t\geq 0\quad P_sP_t=P_{s+t}$).

\chapter{Some convergence estimates for non-perpetual American option pricing based on cubature}
\label{cubature_dyadic}

For this Chapter, let us consider an arbitrary but fixed translation-invariant finite-state Markov chain $P:=\left(P_t\right)_{t\in I}$ with state space $\RR^d$ (for $d\in\NN$) where $I= h\NN_0$ for some real number $h>0$, as well as a real number $T>h>0$ (the time horizon, or maturity), a real number $r>0$, a continuous function $\bar f:\RR^d\rightarrow [0,+\infty)$ that is monotone in each coordinate, a nonnegative real number $K\geq 0$ and let us set $$g:=K- \bar f$$ as well as defining a family of maps $B_t:L^0\left(\RR^d,[0,+\infty)\right)\rightarrow L^0\left(\RR^d,[0,+\infty)\right)$, $t\geq 0$, by $$\forall t\geq 0 \quad B_t:f\mapsto \max\left\{e^{-rt}P_tf,g\right\}= \left(e^{-rt}P_tf\right)\vee g.$$ (Note that $B_tf$ will always be nonnegative for $f\geq 0$ -- hence, for all $f\geq 0$, $B_tf\geq g\vee 0$.) Furthermore, we shall denote by $\left\{x_1^{(t)},\dots, x_{m^{\frac{t}{h}}}^{(t)}\right\}$ the set of (distinct) states at time $t$ after starting the process at time $0$ in $0\in\RR^d$ and by $\left\{\alpha^{(t)}_1,\dots,\alpha^{(t)}_{m^\frac{t}{h}}\right\}\subset(0,1]$ the weights for each of these states, thereby imposing on the sets $\left\{\alpha^{(t)}_1,\dots,\alpha^{(t)}_{m^{\frac{t}{h}}}\right\}$ for $t\in I$, in addition to it being a subset of $(0,1]$, the condition that they be a convex combination, viz. $$\forall t\in I\quad \sum_{i=1}^{m^{\frac{t}{h}}}\alpha_{i}^{(t)}=1.$$ Summarising this, we write $$\forall t\in I=h\NN_0\quad P_t:f\mapsto \sum_{i=1}^{m^{\frac{t}{h}}}\alpha_{i}^{(t)}f\left(\cdot-x_i^{(t)}\right).$$ For the whole of this section, the Lebesgue measure on $\RR$ shall be denoted by $\lambda$, and $\lambda^d$ will be shorthand for the measure-theoretic power $\lambda^{\otimes d}$.

The operators $\max$ and $\min$ when applied to subsets of $\RR^d$ will be understood to be taken componentwise. Analogously, we will interpret the relations $\leq$ and $\geq $ componentwise on $\RR^d\times\RR^d$.

For convenience, we allow all $L^p$-norms (including the $C^0$ norm) of measurable functions to take values in the interval $[0,+\infty]$, thereby extending the domain for each of the $L^p$-norm to $L^0$, the vector lattice of measurable functions. Furthermore, any functions occurring in this Chapter will be assumed to be measurable. Thus, eg the relation $f_0\geq f_1$ should be read as shorthand for $f_0\in L^0\left(\lambda^d\right)\cap\left\{\cdot \geq f_1\right\}$ for all functions $f_0,f_1$; analogously for the relation $f_0\leq f_1$.

Finally, we will use the operation $\vee$ in such a way that it is applied prior to $+$, but only after $P_s$ and multiplication with other functions or constants have taken place: $$C\cdot P_sf_0\vee f_3\cdot f_1 +f_2 =\max\left\{C\cdot P_sf_0,f_3\cdot f_1\right\}+f_2$$

In this Chapter we are aiming to understand the convergence behaviour of the sequence $\left(\left(B_{T\cdot 2^{-n}}\right)^{\circ 2^n}(g\vee 0)\right)_{n}$. We will start by noting that this sequence is monotonely increasing:

\begin{lem} The sequence $\left(\left(B_{T\cdot 2^{-n}}\right)^{\circ 2^n}f\right)_{n\in\NN_0\cap\left\{T\cdot 2^{-\cdot}\in I\right\} }$ is monotonely increasing for all functions $f:\RR^d\rightarrow \RR$. Furthermore, if there exists a function $\tilde g\geq g\vee 0$ such that $\tilde g$ is {\em $e^{-r\cdot}P_\cdot$-harmonic} (ie $e^{-rh}P_h\tilde g=\tilde g$) and $f\leq \tilde g$, then for all $n\in\NN_0$, $\left(B_{T\cdot 2^{-n}}\right)^{\circ 2^n}f\leq \tilde g$.
\end{lem}
\begin{proof} Consider $f:\RR^d\rightarrow \RR$ and $n\in\NN_0$ such that $T\cdot 2^{-(n+1)}\in I=h\NN_0$. Then $$\left(B_{T\cdot 2^{-(n+1)}}\right)^{\circ 2^{n+1}}f=\left(\left(B_{T\cdot 2^{-n+1}}\right)^{\circ 2}\right)^{\circ 2^n}f$$ and by the monotonicity of the operators $P_s$ for $s\in I$, \begin{eqnarray*}\left(B_{T\cdot 2^{-n+1}}\right)^{\circ 2}&=&e^{-rT \cdot 2^{-(n+1)}}P_{T \cdot 2^{-(n+1)}}\left(e^{-rT \cdot 2^{-(n+1)}}P_{T \cdot 2^{-(n+1)}}(\cdot)\vee g\right)  \vee g\\ &\geq &e^{-rT \cdot 2^{-(n+1)}}P_{T \cdot 2^{-(n+1)}}\left(e^{-rT \cdot 2^{-(n+1)}}P_{T \cdot 2^{-(n+1)}}(\cdot)\right)\vee g\\ &=& e^{-rT \cdot 2^{-n}}P_{T \cdot 2^{-n}}(\cdot)\vee g=B_{T\cdot 2^{-n}},\end{eqnarray*} where the last line is a consequence of the Chapman-Kolmogorov equation. This completes the proof for the monotonicity of the sequence $\left(\left(B_{T\cdot 2^{-n}}\right)^{\circ 2^n}f\right)_{n\in\NN_0\cap\left\{T\cdot 2^{-\cdot}\in I\right\}}$. 

Now suppose there exists such a function $\tilde g$ as in the statement of the Lemma. Then $e^{-rs}P_s\tilde g=\tilde g$ for all $s\in I$ and therefore $B_s\tilde g=\tilde g$ for all $s\in I$. Also, the map $B_s$ is monotone in the sense that $g_0\leq g_1$ always implies $B_sg_0\leq B_sg_1$ (because it is the composition of two monotone maps: $e^{-rs}P_s$ and $\cdot\vee g$) for all $s\in I$. Thus we see that for all $f\leq\tilde g,$ $$\left(B_{T\cdot 2^{-n}}\right)^{\circ 2^n}f\leq \left(B_{T\cdot 2^{-n}}\right)^{\circ 2^n}\tilde g =\tilde g.$$
\end{proof}

\begin{lem} \label{B_tcontractsifgeq g} For all measurable functions $f_1\geq f_0\geq g\vee 0$, as well as for all $t\in I$ and $p\in\{1,\infty\}$ one has \begin{eqnarray*} \left\|B_tf_1-B_tf_0\right\|_{L^p\left(\lambda^d\right)} &\leq& e^{-rt} \left\|f_1-f_0\right\|_{L^p\left(\lambda^d\left[\left\{e^{-rt}P_tf_1>g\right\}\cap\cdot\right]\right)} \\ &\leq& e^{-rt} \left\|f_1-f_0\right\|_{L^p\left(\lambda^d\right)}  \end{eqnarray*} (with the usual convention that $x\leq +\infty$ for all $x\in\RR\cup\{\pm\infty\}$).
\end{lem}
\begin{proof} The map $B_t$ is monotone. Thus we have \begin{eqnarray*} \left\{B_tf_1=g\right\}&=& \left\{B_tf_0\leq B_tf_1=g\right\} \\ &=&  \left\{g\vee 0 \leq B_tf_0\leq B_tf_1=g\right\} = \left\{B_tf_0=g\right\}\cap  \left\{B_tf_1=g\right\}\\ &\subseteq&  \left\{B_tf_1-B_tf_0=0\right\} \end{eqnarray*} for $f_1\geq f_0\geq g\vee 0$.
Since $B_tf_1\geq g$, this implies \begin{eqnarray*}0\leq B_tf_1-B_tf_0&=& \chi_{\left\{B_tf_1>g\right\}}\left(e^{-rt}P_tf_1\vee g - e^{-rt}P_tf_0\vee g\right) \\ &=& \chi_{\left\{B_tf_1>g\right\}}\left(e^{-rt}P_tf_1-e^{-rt}P_tf_0\vee g\right)\\ &\leq & \chi_{\left\{e^{-rt}P_tf_1>g\right\}}\left(e^{-rt}P_tf_1-e^{-rt}P_tf_0\right) \\ &=& e^{-rt}\chi_{\left\{e^{-rt}P_tf_1>g\right\}}P_t\left(f_1-f_0\right)\end{eqnarray*} which yields the assertion as $P_t$ is an $L^p(\lambda^d)$-contraction (for $p=\infty$ this is immediate and for $p=1$ it follows from the translation-invariance of both $P_t$ and the Lebesgue measure). 
\end{proof}

\begin{lem} \label{ggeq0onCEt}Suppose $$\max_{k\in\{1,\dots,m\}}x_k^{(h)}\leq 0$$ componentwise, implying $x_i^{(s)}\leq 0$ componentwise for all $s\in I=h\NN_0$ and $i\in\left\{1,\dots,m^\frac{s}{h}\right\}$. Then for all $s\in I$, $g$ is nonnegative on $\left\{P_s(g\vee 0)=P_sg\right\}$.
\end{lem}
\begin{proof} Recalling our notational convention that $\leq$ as relation on $\RR^d\times \RR^d$ and $\max$ when applied to subsets of $\RR^d$ are to be interpreted componentwise, we may write $$\forall\ell\in{\left\{1,\dots,m^\frac{s}{h}\right\}}\quad \RR^d\ni 0\geq\max_{i\in\left\{1,\dots,m^\frac{s}{h}\right\}} x_i^{(t)}\geq x_\ell^{(t)},$$ due to the componentwise monotonicity of $g$, yields for all $s\in I$ the inclusion
\begin{eqnarray*}\left\{P_s(g\vee 0)=P_s g\right\}&=& \left\{\forall i\in\left\{1,\dots,m^\frac{s}{h}\right\}\quad g\left(\cdot-x_i^{(s)}\right)\geq 0\right\}\\ &\subseteq& \left\{g\left(\cdot-\max_{i\in\left\{1,\dots,m^\frac{s}{h}\right\}} \right)\geq 0\right\}\\ &\subseteq&\left\{g\geq 0\right\}.\end{eqnarray*}
\end{proof}

\begin{lem} \label{offEformula} Suppose there is a $\gamma_0\geq 1$ (without loss of generality, $\gamma_0\in[1,e^r)$) such that $$P_t\bar f\geq {\gamma_0}^t\bar f$$ for all $t\in(0,T]\cap I$ (where  $I=h\NN_0$ with $h>0$ whence it is sufficient that this estimate holds for $t=h$). In addition, assume that $g\geq 0$ on the subset $\left\{P_t(g\vee 0)=P_tg\right\}$ of $\RR^d$ for all $t\in(0,T]\cap I$ (this assumption being, due to Lemma \ref{ggeq0onCEt}, satisfied in particular if $\max_{i\in\{1,\dots,m\}}x_i^{(h)}\leq 0$). Then for all $t\in(0,T]\cap I$, \begin{eqnarray*} \left\{g\geq e^{-rt}P_t(g\vee 0)\right\} &\supseteq&  \left\{P_t(g\vee 0) = P_tg\right\} \end{eqnarray*}
\end{lem}
\begin{proof} Let $t\in (0,T]\cap I$. Due to our assumption of $g\geq 0$ on $\{P_t(g\vee 0)=P_tg\}$, one has \begin{eqnarray} && \left\{g\geq e^{-rt}P_t(g\vee 0)\right\}\cap \left\{P_t(g\vee 0)=P_tg\right\} \nonumber \\ &=& \left\{g\geq e^{-rt}P_tg\right\}\cap \left\{P_t(g\vee 0)=P_tg\right\}  \nonumber \\ &=&\{g\geq 0\}\cap \left\{\left(\id-e^{-rt}P_t\right) g \geq 0\right\}\cap \left\{P_t(g\vee 0)=P_tg\right\} \nonumber  \\ &=& \{g\geq 0\}\cap \left\{\left(1-e^{-rt}\right)K\geq \left(\id-e^{-rt}P_t\right)\bar f\right\} \nonumber \\ &&\cap \left\{P_t(g\vee 0)=P_tg\right\}. \label{lemma1.2eq}\end{eqnarray}
On the other hand \begin{eqnarray*}&&\left\{\left(1-e^{-rt}\right)K\geq \left(\id-e^{-rt}P_t\right)\bar f\right\} = \left\{\left(1-e^{-rt}\right)K\geq \bar f-e^{-rt}\underbrace{P_t\bar f}_{\geq {\gamma_0}^t\bar f}\right\} \\ &\supseteq&   \left\{\left(1-e^{-rt}\right)K\geq \left(1-e^{-rt}{\gamma_0}^t\right)\bar f\right\}, \end{eqnarray*} that is \begin{eqnarray*}&&\left\{\left(1-e^{-rt}\right)K\geq \left(\id-e^{-rt}P_t\right)\bar f\right\} \\ &\supseteq& \left\{\frac{1-e^{-rt}}{1-e^{-rt}{\gamma_0}^t}K \geq  \bar f\right\}, \end{eqnarray*} where we have exploited $\gamma_0<e^r$. Now $\gamma_0\in[1,e^r)$ gives $$\frac{1-e^{-rt}}{1-e^{-rt}{\gamma_0}^t}K\geq K$$ since $K\geq 0$. Combining this estimate with the previous inclusion, one obtains $$\left\{\left(1-e^{-rt}\right)K \geq \left(\id-e^{-rt}P_t\right)\bar f\right\} \supseteq \{K\geq \bar f\}$$ and hence $$\left\{\left(1-e^{-rt}\right)K\geq \left(\id-e^{-rt}P_t\right)\bar f\right\} \cap \{g\geq 0\} =\{g\geq 0\}.$$ This result, combined with the first equation (\ref{lemma1.2eq}) in this Proof, yields \begin{eqnarray*} && \left\{g\geq e^{-rt}P_t(g\vee 0)\right\}\cap \left\{P_t(g\vee 0)=P_tg\right\} \\ &=& \{g\geq 0\}\cap \left\{P_t(g\vee 0)=P_tg\right\}. \end{eqnarray*} However, one of our assumptions reads $$\left\{P_t(g\vee 0)=P_tg\right\}\subseteq \{g \geq 0\}$$ whence we conclude \begin{eqnarray*} && \left\{g\geq e^{-rt}P_t(g\vee 0)\right\}\cap \left\{P_t(g\vee 0)=P_tg\right\} \\ &=& \left\{P_t(g\vee 0)=P_tg\right\}. \end{eqnarray*} 
\end{proof}

\begin{rem} The assumption of the existence of a $\gamma_0\leq e^r$ such that $P_t\bar f\geq {\gamma_0}^t\bar f$ for all $t\in I$ is natural: If $(X_t)_{t\in I}$ was a Markov process evolving according to $(P_t)_{t\in I}$, the (stronger) condition $$\forall t\in I \quad P_t\bar f=e^{+rt }\bar f$$ simply means that the process $\bar f(X_\cdot)$ is, after discounting, a martingale. Now, if $X_\cdot$ was a Markov model for a vector of logarithmic asset prices (a {\em Markov basket} in our terminology) and $\bar f$ would assign to each vector the arithmetic average of the exponentials of its components, this is by definition true if $P_\cdot$ governs the process $X_\cdot$ under a risk-neutral measure. Furthermore, the said assumption $$\forall t\in I \quad P_t\bar f=e^{+rt }\bar f$$ trivially implies $$\exists \gamma_1>0\quad \forall t\in I\quad P_t\bar f \leq {\gamma_1}^t\bar f$$ and therefore provides us with some vindication for assuming the last assertion in some of the subsequent Lemmas of this Chapter.
\end{rem}

\begin{lem} \label{Bs/2-Bsestimate} Suppose there is a $\gamma_1>0$ such that $$P_t\bar f\leq {\gamma_1}^t\bar f$$ for all $t\in(0,T]\cap I$ (for which in our case of $I=h\NN_0$ with $h>0$ it is sufficient that this estimate holds for $t=h$), and let us assume without loss of generality that this $\gamma_1$ be $\geq e^r$. Then, setting $$R:=K\cdot \sup_{t\in(0,T]\cap I}\frac{{\gamma_1}^t -1}{t},$$ we have found an $R<+\infty$ such that for all $s\in(0,T)\cap (2\cdot I)\subset I$ and measurable $f\geq g\vee 0$, $$\left\|\left(B_{\frac{s}{2}}\right)^{\circ 2}f - B_sf\right\|_{L^\infty(\RR^d)}\leq R\cdot\frac{s}{2}.$$
\end{lem}
\begin{proof} Let $s\in(0,T]\cap (2\cdot I)\subset I$ and consider a measurable $f\geq g\vee 0$. Then by our assumption of $P_t\bar f\leq {\gamma_1}^t\bar f$ for $\gamma_1>0$, we firstly have (inserting $\frac{s}{2}$ for $t$)
\begin{eqnarray*}&& g-e^{-r\frac{s}{2}}P_{\frac{s}{2}}g=K-\bar f-e^{-r\frac{s}{2}}\left(K-P_{\frac{s}{2}}\bar f\right)\\ &=& K\left(1-e^{-r\frac{s}{2}}\right)+e^{-r\frac{s}{2}}P_{\frac{s}{2}}\bar f-\bar f\\ &\leq& K\left(1-e^{-r\frac{s}{2}}\right)+ \left(e^{-r\frac{s}{2}}{\gamma_1}^\frac{s}{2}-1\right)\cdot\bar f \end{eqnarray*}
and therefore (using $f\geq g \vee 0$ and $\gamma_1\geq e^r$ as well as the monotonicity of $P_{\frac{s}{2}}$), 
\begin{eqnarray} 0&\leq& \chi_{\left\{g\geq e^{-r\frac{s}{2}}P_{\frac{s}{2}}f\right\}}\cdot\left(g-e^{-r\frac{s}{2}}P_{\frac{s}{2}}f\right) \nonumber\\ &\leq& \chi_{\left\{g\geq e^{-r\frac{s}{2}}P_{\frac{s}{2}}f\right\}}\cdot\left(g-e^{-r\frac{s}{2}}P_{\frac{s}{2}}g\right) \nonumber\\ &\leq& \chi_{\left\{g\geq 0\right\}}\cdot\left(K\left(1-e^{-r\frac{s}{2}}\right)+ \left(e^{-r\frac{s}{2}}{\gamma_1}^\frac{s}{2}-1\right)\cdot\bar f \right) \nonumber\\ &\leq& \chi_{\left\{\bar f\leq K\right\}}\cdot\left(K\left(1-e^{-r\frac{s}{2}}\right)+ \left(e^{-r\frac{s}{2}}{\gamma_1}^\frac{s}{2}-1\right)\cdot K \right)\nonumber \\ &\leq& K\cdot\left(\left(1-e^{-r\frac{s}{2}}\right)+ \left(e^{-r\frac{s}{2}}{\gamma_1}^\frac{s}{2}-1\right) \right)\nonumber \\ & = & K\cdot e^{-r\frac{s}{2}}\cdot \left({\gamma_1}^\frac{s}{2}-1\right) \label{Rbound} \end{eqnarray} 
Now, $$\sup_{t\in(0,T]\cap I}\frac{{\gamma_1}^t -1}{t}<+\infty$$ since $t\mapsto {\gamma_1}^t-1$ is right-differentiable in zero with derivative $\ln \gamma_1$. Therefore $$R=K\cdot \sup_{t\in(0,T]\cap I}\frac{{\gamma_1}^t -1}{t}<+\infty.$$ Via estimate (\ref{Rbound}), we arrive at \begin{equation}0\leq \chi_{\left\{g\geq e^{-r\frac{s}{2}}P_{\frac{s}{2}}f\right\}}\left(g-e^{-r\frac{s}{2}}P_{\frac{s}{2}}f\right)\leq R\cdot \frac{s}{2}.\label{Chapter8,1.3a}\end{equation} But $$\chi_{\left\{g\geq e^{-r\frac{s}{2}}P_{\frac{s}{2}}f\right\}}\left(g-e^{-r\frac{s}{2}}P_{\frac{s}{2}}f\right)=\left(e^{-r\frac{s}{2}}P_{\frac{s}{2}}f\right)\vee g - e^{-r\frac{s}{2}}P_{\frac{s}{2}}f,$$ and -- in combination with the linearity of $P_{\frac{s}{2}}$ and the Chapman-Kolmogorov equation -- this implies \begin{eqnarray*}&& e^{-r\frac{s}{2}}P_{\frac{s}{2}}\left(\chi_{\left\{g\geq e^{-r\frac{s}{2}}P_{\frac{s}{2}}f\right\}}\left(g-e^{-r\frac{s}{2}}P_{\frac{s}{2}}f\right)\right)\\ &=&e^{-r\frac{s}{2}}P_{\frac{s}{2}}\left(\left(e^{-r\frac{s}{2}}P_{\frac{s}{2}}f\right)\vee g\right) - e^{-rs}P_{{s}}f\end{eqnarray*} hence by equation (\ref{Chapter8,1.3a})\begin{eqnarray}&& R\cdot\frac{s}{2}\nonumber\\ &\geq& e^{-r\frac{s}{2}}P_{\frac{s}{2}}\left(\chi_{\left\{g\geq e^{-r\frac{s}{2}}P_{\frac{s}{2}}f\right\}}\left(g-e^{-r\frac{s}{2}}P_{\frac{s}{2}}f\right)\right)\\ &\geq&e^{-r\frac{s}{2}}P_{\frac{s}{2}}\left(\left(e^{-r\frac{s}{2}}P_{\frac{s}{2}}f\right)\vee g\right) - e^{-rs}P_{{s}}f\vee g \label{Rs/2eqn}\end{eqnarray} (where in (\ref{Rs/2eqn}) we have exploited the fact that $P_\frac{s}{2}$ is an $L^\infty$-contraction).

Now, again by the Chapman-Kolmogorov equation and the monotonicity of $P_t$ for any $t\in I$, we have $$\left(B_{\frac{s}{2}}\right)^{\circ 2}f\geq B_s f$$ and therefore (due to the estimate $B_tf_0\geq g\vee 0$ which holds for arbitrary $t\in I$ and $f_0\geq 0$) \begin{eqnarray*} \left\{\left(B_{\frac{s}{2}}\right)^{\circ 2}f=g\right\}&=& \left\{B_sf\leq \left(B_{\frac{s}{2}}\right)^{\circ 2}f=g\right\} \\ &=&  \left\{g\vee 0 \leq B_sf\leq \left(B_{\frac{s}{2}}\right)^{\circ 2}f=g\right\} \\ &=& \left\{B_sf=g\right\}\cap  \left\{\left(B_{\frac{s}{2}}\right)^{\circ 2}f=g\right\}\\ &\subseteq&  \left\{\left(B_{\frac{s}{2}}\right)^{\circ 2}f-B_sf=0\right\} \end{eqnarray*} 
But $$\left(B_{\frac{s}{2}}\right)^{\circ 2}f\geq g\vee 0\geq g,$$ thus the last inclusion yields \begin{eqnarray*}0&\leq& \left(B_{\frac{s}{2}}\right)^{\circ 2}f-B_sf\\ &=&\left(\left(B_{\frac{s}{2}}\right)^{\circ 2}f-B_sf\right)\chi_{\left\{ \left(B_{\frac{s}{2}}\right)^{\circ 2}f\geq g\right\}}\\ &=& e^{-r\frac{s}{2}}\left(P_{\frac{s}{2}}\left(e^{-r\frac{s}{2}}P_{\frac{s}{2}}f\vee g\right) - \left(e^{-rs}P_{{s}}f\vee g\right)\right)\chi_{\left\{ \left(B_{\frac{s}{2}}\right)^{\circ 2}f\geq g\right\}}\\ &\leq& R\cdot\frac{s}{2}\end{eqnarray*} where the last line has used the estimate (\ref{Rs/2eqn}) derived previously.
\end{proof}

Later on, in Lemma \ref{lowerboundoffE}, we will see that it is impossible to obtain estimates for $\left\|\left(B_{\frac{s}{2}}\right)^{\circ 2}f - B_sf\right\|_{L^\infty(\RR^d)}$ that are both uniform in $f\geq g\vee 0$ and of higher than linear order in $s$.

We can draw from the proof of Lemma \ref{Bs/2-Bsestimate} the following Corollary:

\begin{cor}\label{Bs/2-Bsequation} Suppose there is a $\gamma_1>0$ such that $$P_t\bar f\leq {\gamma_1}^t\bar f$$ for all $t\in(0,T]\cap I$ (for which in case $I=h\NN_0$ with $h>0$ it is sufficient that this estimate holds for $t=h$). Then for all measurable $f\geq g\vee 0$ \begin{eqnarray*}0 &\leq& \left(B_{\frac{s}{2}}\right)^{\circ 2}f-B_sf\\ &=& e^{-r\frac{s}{2}}P_{\frac{s}{2}}\left(\left(e^{-r\frac{s}{2}}P_{\frac{s}{2}}f\right)\vee g - e^{-r\frac{s}{2}}P_{\frac{s}{2}}f\right)\chi_{\left\{ \left(B_{\frac{s}{2}}\right)^{\circ 2}f\geq g\right\}}\\&\leq& e^{-r\frac{s}{2}}P_{\frac{s}{2}}\left(\chi_{\left\{g\geq e^{-r\frac{s}{2}}P_{\frac{s}{2}}(g\vee 0)\right\}}\cdot\left(g-e^{-r\frac{s}{2}}P_{\frac{s}{2}}(g\vee 0)\right) \right)\\ &=& \left(B_{\frac{s}{2}}\right)^{\circ 2}(g\vee 0)-B_s(g\vee 0).\end{eqnarray*}
\end{cor}

We will continue to assume $I=h\NN_0$, $h<T$ and define 

\begin{Def}$$\forall t\in I\quad E^t = \left\{P_t(g\vee 0)>P_tg\right\}.$$\end{Def}

\begin{rem} \label{Etproperties}Equivalent expressions for $E^\cdot$ are:
\begin{eqnarray*} \forall t\in I \quad E^t &=& \left\{P_t(g\vee 0)>P_tg\right\} \\ &=& \complement \left\{P_t(g\vee 0)=P_tg\right\} \\ &=& \left\{P_t(g\wedge 0)<0\right\}  \\ &=& \left\{\exists i\in\{1,\dots,m^{\frac{t}{h}}\}\quad g\left(\cdot-x_i^{(t)}\right)<0\right\}\\ &=& \bigcup_{i\in\{1,\dots,m^{\frac{t}{h}}\}}\left\{K<\bar f\left(\cdot-x_i^{(t)}\right)\right\} .\end{eqnarray*} These formulae for $E^\cdot$ imply, by the monotonicity of $\bar f$, that $E^t$ is north-east connected (that is $E^t\subset E^t+a$ for all $a\leq 0$) for all $t\in I$. Furthermore, if one had for all $t\in I$ and $i\in\left\{1,\dots, m^\frac{t}{h}\right\}$ an index $k=k(i)\in\left\{1,\dots, m^\frac{t}{h}\right\}$ such that $x_i^{(t)}+x_{k(i)}^{(t)}\leq 0$ componentwise (for instance if the set $\left\{x_1^{(h)}, \dots,x_m^{(h)} \right\}$ could be written as the sum of a reflection symmetric subset of $\RR^d$ and a componentwise nonpositive vector), then the north-east connectedness of the $E^t$'s entails for all $t\in I$ and $i$, \begin{eqnarray*}E^t-x_i^{(t)}&\subseteq& E^t+x_i^{(t)}+x_{k(i)}^{(t)}-x_i^{(t)}=E^t+x_{k(i)}^{(t)}\\ &=&\bigcup_{j\in\left\{1,\dots,m^\frac{t}{h}\right\}}\left\{K<\bar f<\left(\cdot-x_j^{(t)}-x_{k(i)}^{(t)}\right)\right\} \\ &\subseteq&\bigcup_{j_0,j_1\in\left\{1,\dots,m^\frac{t}{h}\right\}}\left\{K<\bar f<\left(\cdot-x_{j_0}^{(t)}-x_{j_1}^{(t)}\right)\right\}\\ &=&\bigcup_{\ell\in\left\{1,\dots,m^\frac{2t}{h}\right\}}\left\{K<\bar f<\left(\cdot-x_\ell^{(2t)}\right)\right\} \\ &=& E^{2t}\end{eqnarray*} where for the penultimate line we have used the Chapman-Kolmogorov equation, of course. Therefore $$\chi_{E^t}\left(\cdot+x_i^{(t)}\right)=\chi_{E^t-x_i^{(t)}}\leq\chi_{E^{2t}}.$$ Also, if there exists an $i_0\in\{1,\dots,m\}$ such that $$\forall j\in\{1,\dots, d\}\quad \left(x_{i_0}^{(h)}\right)_j\leq 0 $$ one has -- due to the monotonicity of $\bar f$ in each coordinate -- first of all $\bar f\left(\cdot -x_{i_0}^{(h)}\right)\geq \bar f$ and thence for all $n\in\NN$ the inclusion \begin{eqnarray*} E^{nh} &=& \bigcup_{i\in\{1,\dots,m^n\}}\left\{K<\bar f\left(\cdot-x_i^{(nh)}\right)\right\}  \\ &=& \bigcup_{i_1,\dots,i_n\in\{1,\dots,m\}}\left\{K<\bar f\left(\cdot-x_{i_1}^{(h)}-\dots- x_{i_n}^{(h)}\right)\right\}\\ &\supseteq& \bigcup_{i_1,\dots,i_{n-1}\in\{1,\dots,m\}}\left\{K<\bar f\left(\cdot-x_{i_0}^{(h)}-x_{i_1}^{(h)}-\cdots -x_{i_{n-1}}^{(h)}\right)\right\} \\ &\supseteq& \bigcup_{i_1,\dots,i_{n-1}\in\{1,\dots,m\}}\left\{K<\bar f\left(\cdot-x_{i_1}^{(h)}-\cdots- x_{i_{n-1}}^{(h)}\right)\right\}  \\ &=& E^{(n-1)h} \end{eqnarray*} This means $$E^s\uparrow \text{ as }s\uparrow\infty\text{ in }I$$ and for all $T\in[h,+\infty]$, $$\bigcap_{s\in (0,T]\cap I}E^s = E^h.$$
\end{rem}

The reason for $E^h$ not being the whole space is that the measure $B\mapsto P_t\chi_B$ on the Borel $\sigma$-algebra of $\RR$ has compact support. 

If one interprets $g$ as a logarithmic payoff function (eg $g=K-\exp$, $d=1$ in case of a vanilla one-dimensional put) and $P$ as a Markov chain that models the stochastic evolution of the logarithmic prices of assets in a given portfolio, then the set $E^t$, for $t\in I$ consists of all those vectors of logarithmic start prices where the probability of exercising the option at time $t$ is strictly positive.

\begin{lem} \label{lowerboundoffE} Suppose there is a $\gamma_0>1$ (without loss of generality, $\gamma_0\in\left(1,e^r\right]$) such that $$P_t\bar f\geq {\gamma_0}^t\bar f$$ for all $t\in(0,T]\cap I$ (where  $I=h\NN_0$ with $h>0$ whence it is sufficient that this estimate holds for $t=h$). 
Assume furthermore that $\max_{i\in\{1,\dots, m\}} x_i^{(h)}\leq 0$, implying that $g>0$ on the set $\left\{P_t(g\vee 0)=P_tg\right\}$. Then for all $\varepsilon_1>0$ there is an $\varepsilon_0<T$ independent of $h<T$ such that for all $t\in 2\cdot\left( \left(0,\varepsilon_0\right)\cap I \right)$ and $A\supset \complement E^\frac{t}{2}$ (with positive Lebesgue measure), \begin{eqnarray*}&&\sup_{f\geq g\vee 0}\left\|\left(B_{\frac{s}{2}}\right)^{\circ 2}f - B_sf\right\|_{L^\infty\left(A \right)}\geq \left\|\left(B_{\frac{s}{2}}\right)^{\circ 2}(g\vee 0) - B_s(g\vee 0)\right\|_{L^\infty\left(A \right)}  \\ &\geq& \left(\min_{i\in\left\{1,\dots,m\right\}} \alpha^{\left(h\right)}\right)^\frac{T}{h} e^{-r\frac{t}{2}}K\left(\ln\gamma_0-\varepsilon_1 \right)\cdot \frac{t}{2},\end{eqnarray*} as well as \begin{eqnarray*}&&\sup_{f\geq g\vee 0}\left\|\left(B_{\frac{s}{2}}\right)^{\circ 2}f - B_sf\right\|_{L^1\left(A \right)}\geq \left\|\left(B_{\frac{s}{2}}\right)^{\circ 2}(g\vee 0) - B_s(g\vee 0)\right\|_{L^1\left(A \right)} \\ &\geq& \left(\min_{i\in\left\{1,\dots,m\right\}} \alpha^{\left(h\right)}\right)^\frac{T}{h} \lambda^d\left[\left\{P_{\frac{t}{2}}(g\vee 0)=P_{\frac{t}{2}} g\right\}\right]\cdot e^{-r\frac{t}{2}}K\left(\ln\gamma_0-\varepsilon_1 \right)\cdot \frac{t}{2} \end{eqnarray*} (the left hand side, following the usual convention, being $+\infty$ if $\lambda^d\left[\left\{P_{\frac{t}{2}}(g\vee 0)=P_{\frac{t}{2}} g\right\}\right]=+\infty$, $m>1$ and $\varepsilon_1<\ln\gamma_0$).
\end{lem}
\begin{proof} Let us first remark that, due to Corollary \ref{Bs/2-Bsequation}, we have \begin{eqnarray*} &&\sup_{f\geq g\vee 0} \left\|\left(B_{\frac{s}{2}}\right)^2f-B_{s}f\right\|_{L^\infty\left(A\right)} \\ & = & \left\|\left(B_{\frac{s}{2}}\right)^2(g\vee 0)-B_{s}(g\vee 0)\right\|_{L^\infty\left(A\right)}\\ &=&  \left\|\chi_A e^{-r\frac{t}{2}}P_{\frac{t}{2}}\left(\chi_{\left\{g\geq e^{-r\frac{t}{2}}P_{\frac{t}{2}}(g\vee 0)\right\}}\cdot\left(g-e^{-r\frac{t}{2}}P_{\frac{t}{2}}(g\vee 0)\right)\right) \right\|_{L^\infty\left(A\right)}\end{eqnarray*} 
as well as 
\begin{eqnarray} &&\sup_{f\geq g\vee 0} \left\|\left(B_{\frac{s}{2}}\right)^2f-B_{s}f\right\|_{L^1\left(A\right)} \nonumber\\ & = & \left\|\left(B_{\frac{s}{2}}\right)^2(g\vee 0)-B_{s}(g\vee 0)\right\|_{L^\infty\left(A\right)}\nonumber\\ &=&  \left\|e^{-r\frac{t}{2}}P_{\frac{t}{2}}\left(\chi_{\left\{g\geq e^{-r\frac{t}{2}}P_{\frac{t}{2}}(g\vee 0)\right\}}\cdot\left(g-e^{-r\frac{t}{2}}P_{\frac{t}{2}}(g\vee 0)\right)\right) \right\|_{L^1\left(A\right)}\nonumber\\ &=&  \left\|\left(\chi_A\cdot e^{-r\frac{t}{2}}P_{\frac{t}{2}}\left(\chi_{\left\{g\geq e^{-r\frac{t}{2}}P_{\frac{t}{2}}(g\vee 0)\right\}}\cdot\left(g-e^{-r\frac{t}{2}}P_{\frac{t}{2}}(g\vee 0)\right)\right) \right)\left(\cdot+\frac{t}{2h}x_{i_0}^{(h)}\right)\right\|_{L^1\left(\RR^d\right)}\label{chapter8a,4*}\end{eqnarray} 
for all $i_0\in \left\{1,\dots,m^\frac{t}{2h}\right\}$ (using the translation invariance of $\int_{\RR^d}\cdot\lambda^d$).

Next let us note that by our assumption of $x_{i_0}^{(h)}\leq 0$ componentwise for all $i_0\in \left\{1,\dots,m^\frac{s}{h}\right\}$, combined with the north-east connectedness of $E^s$ (which entails south-west connectedness of $\complement E^s$), we have $$\forall s\in I\quad -\frac{s}{h}x_{i_0}^{(h)}+\complement E^s\supseteq \complement E^s.$$ Therefore we may conclude that for all $t\in 2\cdot I$ and $i_0\in \left\{1,\dots,m^\frac{t}{2h}\right\}$, 
\begin{eqnarray}&&\chi_{A} \cdot e^{-r\frac{t}{2}}P_{\frac{t}{2}}\left(\chi_{\left\{g\geq e^{-r\frac{t}{2}}P_{\frac{t}{2}}(g\vee 0)\right\}}\cdot\left(g-e^{-r\frac{t}{2}}P_{\frac{t}{2}}(g\vee 0)\right)\right) \nonumber \\ &\geq&  \chi_{\complement E^\frac{t}{2}} \cdot e^{-r\frac{t}{2}}P_{\frac{t}{2}}\left(\chi_{\left\{g\geq e^{-r\frac{t}{2}}P_{\frac{t}{2}}(g\vee 0)\right\}}\cdot\left(g-e^{-r\frac{t}{2}}P_{\frac{t}{2}}(g\vee 0)\right)\right) \nonumber \\ &\geq& \left(\chi_{\left\{g\geq e^{-r\frac{t}{2}}P_{\frac{t}{2}}(g\vee 0)\right\}}\cdot\left(g-e^{-r\frac{t}{2}}P_{\frac{t}{2}}(g\vee 0)\right)\right)\left(\cdot-\frac{t}{2h}x_{i_0}^{(h)}\right)\nonumber \\&& \cdot  \min_{i\in\left\{1,\dots,m^\frac{t}{2h}\right\}}{\alpha^{\left(\frac{t}{2}\right)}} \chi_{\complement E^\frac{t}{2}} \nonumber \\ &\geq& \left(\chi_{\left\{g\geq e^{-r\frac{t}{2}}P_{\frac{t}{2}}(g\vee 0)\right\}}\cdot\left(g-e^{-r\frac{t}{2}}P_{\frac{t}{2}}(g\vee 0)\right)\chi_{\complement E^\frac{t}{2}-\frac{t}{2h}x_{i_0}^{(h)}} \right)\left(\cdot-\frac{t}{2h}x_{i_0}^{(h)}\right)\nonumber \\&& \cdot  \left(\min_{i\in\left\{1,\dots,m\right\}} \alpha^{\left(h\right)}\right)^\frac{t}{2h}\nonumber  \\ &\geq& \left(\chi_{\left\{g\geq e^{-r\frac{t}{2}}P_{\frac{t}{2}}(g\vee 0)\right\}}\cdot\left(g-e^{-r\frac{t}{2}}P_{\frac{t}{2}}(g\vee 0)\right)\chi_{\complement E^\frac{t}{2}}\right)\left(\cdot-\frac{t}{2h}x_{i_0}^{(h)}\right)\nonumber \\&& \cdot  \left(\min_{i\in\left\{1,\dots,m\right\}} \alpha^{\left(h\right)}\right)^\frac{t}{2h}\label{chapter8a,1.6} \\ (& \geq& 0)\nonumber \end{eqnarray}

Now, off $E^s$ one has due to Lemma \ref{offEformula} (which may be applied thanks to our assumption $\max_ix_i^{(h)}\leq 0$) the following situation: \begin{eqnarray}&&\chi_{\left\{g\geq e^{-rs}P_{s}(g\vee 0)\right\}}\cdot\left(g-e^{-r{s}}P_{{s}}(g\vee 0)\right) \nonumber \\ &=& \chi_{\left\{P_s(g\vee 0)=P_sg\right\}\cap\{g\geq 0\}}\cdot\left(g-e^{-r{s}}P_{{s}}(g\vee 0)\right) \quad\text{ on }\complement E^s\nonumber  \\ &=& \chi_{\left\{P_s(g\vee 0)=P_sg\right\}\cap\{g\geq 0\}}\cdot\left(g-e^{-r{s}}P_{{s}}g\right) \quad\text{ on }\complement E^s \nonumber  \\ &=& \chi_{\left\{P_s(g\vee 0)=P_sg\right\}\cap\{g\geq 0\}}\cdot\left(K-\bar f-e^{-r{s}}K+e^{-rs}P_s\bar f\right) \quad\text{ on }\complement E^s \nonumber \\ &\geq & \chi_{\left\{P_s(g\vee 0)=P_sg\right\}\cap\{g\geq 0\}}\cdot\left(K-\bar f-e^{-r{s}}K+e^{-rs}{\gamma_0}^s\bar f\right) \quad\text{ on }\complement E^s \label{lowerboundoffEest} \end{eqnarray}
However, one can also perform the calculation 
\begin{eqnarray} && \chi_{\{g\geq 0\}}\left(K-\bar f-e^{-r{s}}K+{\gamma_0}^s e^{-rs}\bar f\right) \nonumber \\ &=& \chi_{\{K-\bar f\geq 0\}}\left(\left(K-\bar f\right)\left(1-{\gamma_0}^se^{-r{s}}\right) +Ke^{-rs}\left({\gamma_0}^s-1\right) \right)\label{lowerboundoffEpenult} \\ &\geq & Ke^{-rs}\left({\gamma_0}^s-1\right) \label{lowerboundoffEfin}\end{eqnarray} (where we have used the assumption $\gamma_0\leq e^r$ to get from (\ref{lowerboundoffEpenult}) to (\ref{lowerboundoffEfin})). Combining estimates (\ref{lowerboundoffEfin}) and (\ref{lowerboundoffEest}), we arrive at \begin{eqnarray*} &&\chi_{\left\{g\geq e^{-rs}P_{s}(g\vee 0)\right\}}\left(g-e^{-r{s}}P_{{s}}(g\vee 0)\right) \\ &\geq& Ke^{-rs}\left({\gamma_0}^s-1\right)  \quad\text{ on }\complement E^s \\ &\geq& Ke^{-rs}\left(\ln{\gamma_0}-\varepsilon_1\right)\cdot s \quad\text{ on }\complement E^s \end{eqnarray*} for every $s<\varepsilon_0$ for some $\varepsilon_0>0$ dependent on $\varepsilon_1>0$ and finally (using estimate (\ref{chapter8a,1.6}), $\min_{i\in\left\{1,\dots,m\right\}} \alpha^{\left(h\right)}\leq 1$ and $T\geq t$) \begin{eqnarray*} && \chi_{A} \cdot e^{-r\frac{t}{2}}P_{\frac{t}{2}}\left(\chi_{\left\{g\geq e^{-r\frac{t}{2}}P_{\frac{t}{2}}(g\vee 0)\right\}}\cdot\left(g-e^{-r\frac{t}{2}}P_{\frac{t}{2}}(g\vee 0)\right)\right) \\ &\geq &   \left(\min_{i\in\left\{1,\dots,m\right\}} \alpha^{\left(h\right)}\right)^\frac{t}{2h}\cdot \chi_{\complement E^\frac{t}{2}}\left(\cdot-\frac{t}{2h}x_{i_0}^{(h)}\right) \cdot Ke^{-r{\frac{t}{2}}}\left(\ln{\gamma_0}-\varepsilon_1\right)\cdot {\frac{t}{2}} \\ & \geq & \left(\min_{i\in\left\{1,\dots,m\right\}} \alpha^{\left(h\right)}\right)^\frac{T}{2h}\cdot \chi_{\complement E^\frac{t}{2}+\frac{t}{2h}x_{i_0}^{(h)}}\cdot Ke^{-r{\frac{t}{2}}}\left(\ln{\gamma_0}-\varepsilon_1\right)\cdot {\frac{t}{2}}\end{eqnarray*} for all $t\in 2\cdot \left(I\cap(0,\varepsilon_0)\right)$ and $i_0\in \left\{1,\dots,m^\frac{t}{2h}\right\}$.

This yields -- due to the translation-invariance of the Lebesgue measure (which gave us estimate (\ref{chapter8a,4*})) -- the first line of the Lemma's $L^1$ norm estimate. It also implies the $L^\infty$ norm estimate of the Lemma since for all $s\in I$ (in particular for $s=\frac{t}{2}$), \begin{eqnarray*}\lambda^d\left[\complement E^{s} +\frac{s}{h}x_{i_0}^{(h)}\right]&=&\lambda^d\left[\complement E^\frac{s}{h}\right] \\ &= &\lambda^d\left\{P_s(g\vee 0)=P_s g\right\}\\ &\geq &\lambda^d\left\{ g\left(\cdot - \frac{s}{h}\min_{i\in\{1,\dots, m\}} x_i^{(h)}\right)\geq 0\right\}>0\end{eqnarray*} (a consequence of the monotonicity of $g$), and therefore \begin{eqnarray*}&&\left\|\left(\min_{i\in\left\{1,\dots,m\right\}} \alpha^{\left(h\right)}\right)^\frac{T}{2h}\cdot \chi_{\complement E^\frac{t}{2}+\frac{t}{2h}x_{i_0}^{(h)}}\cdot Ke^{-r{\frac{t}{2}}}\left(\ln{\gamma_0}-\varepsilon_1\right)\cdot {\frac{t}{2}}\right\|_{L^\infty(\RR^d)}\\&=&\left(\min_{i\in\left\{1,\dots,m\right\}} \alpha^{\left(h\right)}\right)^\frac{T}{2h}\cdot  Ke^{-r{\frac{t}{2}}}\left(\ln{\gamma_0}-\varepsilon_1\right)\cdot {\frac{t}{2}}.\end{eqnarray*}
\end{proof}

\begin{rem} Assume $\bar f$ is not strictly less than $K$, say $\bar f(z)\geq K$ for some $z\in\RR^d$. We can use the property of $\bar f$ being monotonely increasing in each component to see, via Remark \ref{Etproperties} that 
\begin{eqnarray*}\complement E^h&=&\left\{P_h(g\vee 0)=P_hg\right\} \\ &=&\complement\left\{\exists i\in\{1,\dots,m\}\quad g\left(\cdot-x_i^{(h)}\right)<0\right\} \\ &=&\bigcap_{i=1}^m\left\{g\left(\cdot-x_i^{(h)}\right)\geq 0\right\} \\ &=&\bigcap_{i=1}^m\left\{\bar f\left(\cdot-x_i^{(h)}\right) \leq K \right\} \\ &=&\bigcap_{i=1}^m\left\{\bar f\left(\cdot-x_i^{(h)}\right) \leq K \leq \bar f(z)\right\} \\ &\supset &\bigcap_{i=1}^m\left\{\forall j\in\{1,\dots,d\}\quad \left(\cdot-x_i^{(h)}\right)_j \leq z_j \right\} \\ &= &\bigcap_{i=1}^m\left\{\forall j\in\{1,\dots,d\}\quad \left(\cdot\right)_j \leq z_j +\left(x_i^{(h)}\right)_j \right\} \\ &= & \left\{\forall j\in\{1,\dots,d\}\quad \left(\cdot\right)_j \leq z_j +\min_{i\in\{1,\dots,m\}}\left(x_i^{(h)}\right)_j \right\} \\ &= & \bigotimes_{j=1}^d\left(-\infty,z_j+ \min_{i\in\{1,\dots,m\}}\left(x_i^{(h)}\right)_j \right] ,\end{eqnarray*} where the set in the last line has infinite Lebesgue measure. 

Thus, $\lambda^d\left[\complement E^h\right]=+ \infty$ whenever $\bar f< K$ fails to hold.

\end{rem}

Keeping Corollary \ref{Bs/2-Bsequation} in mind, our next step shall consist in proving

\begin{lem} \label{estimateonE} Let $T\in I$. Suppose there is a $\gamma_1>0$ such that $$P_t\bar f\leq {\gamma_1}^t\bar f$$ for all $t\in(0,T]\cap I$ (where  $I=h\NN_0$ with $h>0$ and therefore it is sufficient that this estimate holds for $t=h$). Let us define $$\tilde D:=\chi_{(0,e^r)}\left(\gamma_1\right)\inf_{\bigcup_{t\in(0,T]\cap I}E^t}\bar f +\chi_{[e^r,+\infty)}\left(\gamma_1\right)\sup_{\bigcup_{t\in(0,T]\cap I}E^t}\bar f\geq 0.$$
Then there is a constant $C_0\in \RR$ given by $$C_0:=K\left(\sup_{s\in(0,T]\cap I}\frac{1-e^{-rs}}{s}- r\right)+\left(\sup_{s\in(0,T]\cap I} \frac{{\gamma_1}^se^{-rs}-1}{s}-\ln\gamma_1+ r\right)$$ such that for all $s\in(0,T]\cap I$ and measurable $A$, \begin{eqnarray*}&& \left\|\chi_{\left\{g> e^{-rs}P_{s}(g\vee 0)\right\}}\left(g-e^{-r{s}}P_{{s}}(g\vee 0)\right) \right\|_{L^1\left(\bigcap_{s\in (0,T]\cap I}E^s\cap A\right)} \\ &\leq &\left\|\chi_{\left\{g> e^{-rs}P_{s}(g\vee 0)\right\}}\left(g-e^{-r{s}}P_{{s}}(g\vee 0)\right) \right\|_{L^1(E^s\cap A)} \\ &\leq& \lambda^d\left[\left\{e^{rs}g>P_s(g\vee 0)>P_s g\right\}\cap A\right] \cdot \left(\left(\ln\gamma_1-r\right)\tilde D+rK+ C_0 \right)\cdot s .\end{eqnarray*}

\end{lem}
\begin{proof} For all $s\in(0,T]\cap I$, the following estimates hold on $E^s$: \begin{eqnarray*}0&\leq& \chi_{\left\{g> e^{-rs}P_{s}(g\vee 0)\right\}}\left(g-e^{-r{s}}P_{{s}}(g\vee 0)\right) \\ &=& \chi_{\left\{e^{rs}g>P_s(g\vee 0)>P_sg\right\}}\left(g-e^{-r{s}}P_{{s}}(g\vee 0)\right) \quad\text{ on }E^s  \\ &\leq& \chi_{\left\{e^{rs}g>P_s(g\vee 0)>P_sg\right\}}\left(g-e^{-r{s}}P_{{s}}g\right) \\ &\leq& \chi_{\left\{e^{rs}g>P_s(g\vee 0)>P_sg\right\}}\left(g-e^{-r{s}}P_{{s}}g\right) \\ &=& \chi_{\left\{e^{rs}g>P_s(g\vee 0)>P_sg\right\}}\left(K-\bar f-e^{-rs}K+e^{-rs} P_s\bar f\right) \\ &\leq& \chi_{\left\{e^{rs}g>P_s(g\vee 0)>P_sg\right\}}\left(K\left(1-e^{-rs}\right)+\left({\gamma_1}^se^{-rs}-1\right)\bar f\right) \\ &\leq& \chi_{\left\{e^{rs}g>P_s(g\vee 0)>P_sg\right\}}\left(K\left(1-e^{-rs}\right)+\left({\gamma_1}^se^{-rs}-1\right) \tilde D\right) \\ &\leq& \chi_{\left\{e^{rs}g>P_s(g\vee 0)>P_sg\right\}}\left(rKs+\left(\ln{\gamma_1}-r\right)\tilde D\cdot s +C\cdot s \right)\end{eqnarray*} for some real constant $C>0$ that can be bounded by $$C\leq K\cdot\left(\sup_{s\in(0,T]\cap I}\frac{1-e^{-rs}}{s}-r\right)+ \tilde D\cdot\left(\sup_{s\in(0,T]\cap I} \frac{{\gamma_1}^se^{-rs}-1}{s}-\left(\ln\gamma_1-r\right)\right)=C_0.$$ This gives a uniform pointwise estimate for the nonnegative function $\chi_{\left\{e^{rs}g>P_s(g\vee 0)\right\}}\left(g-e^{-r{s}}P_{{s}}g\right)$ on $E^s$ from which the Lemma's estimate can be derived immediately.
\end{proof}

\begin{cor}\label{corestimateonE}Let $T\in I$. Assume there exists an $i_0\in\{1,\dots,m\}$ such that $$\forall j\in\{1,\dots,d\}\quad \left(x_{i_0}^{(h)}\right)_j\leq 0.$$ Then one has $$\tilde D=\chi_{(0,e^r)}\left(\gamma_1\right)\inf_{E^T}\bar f +\chi_{[e^r,+\infty)}\left(\gamma_1\right)\sup_{E^T}\bar f$$ and for all $s$ and measurable $A$,
\begin{eqnarray*}&& \left\|\chi_{\left\{g> e^{-rs}P_{s}(g\vee 0)\right\}}\left(g-e^{-r{s}}P_{{s}}(g\vee 0)\right) \right\|_{L^1\left(E^h\cap A\right)} \\ &\leq& \lambda^d\left[\left\{e^{rs}g>P_s(g\vee 0)>P_s g\right\}\cap A\right]\cdot \left(\left(\ln\gamma_1-r\right)\tilde D +rK+ C_0  \right)\cdot s \end{eqnarray*}
\end{cor}
\begin{proof} The assumption about $i_0$ implies that $E^t\uparrow$ as $t\uparrow\infty$ by Remark \ref{Etproperties}), hence $\bigcup_{t\in(0,T]\cap I}E^t=E^T$ which suffices to prove the Corollary. 
\end{proof}

\begin{lem}If there is an $i_1\in\{1,\dots,m\}$ such that $x_{i_1}^{(h)}\geq 0$ componentwise (which entails $\{g>0\}\subseteq \bigcap_{s\in (0,T]\cap I}\left\{P_s(g\vee 0)>0\right\}$ by the monotonicity of $g$ in each component), one will have the following upper bound for the measure of the set occuring in the preceding Lemma \ref{estimateonE}:
\begin{eqnarray*}&&   \lambda^d\left[\left\{e^{rs}g>P_s(g\vee 0)>P_s g\right\}\cap A\right]\\ &\leq& \lambda^d\left[\left\{P_s(g\vee 0)>P_s g>0\right\}\cap A\right] + \lambda^d\left[\left\{P_s(g\vee 0)>0\right\}\cap\left\{P_s g\leq 0\right\}\cap A\right]\end{eqnarray*} for all measurable $A\subseteq\RR^d$.
\end{lem}

\begin{proof} We shall establish an upper bound for the set $\left\{e^{rs}g>P_s(g\vee 0)>P_s\right\}$. Since by our assumption $$\left\{g>0\right\}\subset \left\{P_s(g\vee 0)>0\right\}$$ for all $s\in(0,T]\cap I$, we may, once again exploiting $P_s(g\vee 0)\geq 0$, derive \begin{eqnarray*} && \left\{e^{rs}g>P_s(g\vee 0)>P_sg\right\}\cap \left\{P_sg\leq 0\right\} \\ &\subseteq & \{g>0\}\cap \left\{P_s(g\vee 0)>P_sg\right\}\cap \left\{P_sg\leq 0\right\}\\ &\subseteq& \left\{P_s(g\vee 0)>0\right\}\cap \left\{P_sg\leq 0\right\}\cap \left\{P_s(g\vee 0)>P_sg\right\} \\ &=& \left\{P_s(g\vee 0)>0\right\}\cap \left\{P_sg\leq 0\right\}\end{eqnarray*} This implies \begin{eqnarray*} && \left\{e^{rs}g>P_s(g\vee 0)>P_sg\right\}\\ &\subseteq& \left\{P_s(g\vee 0)>P_sg>0\right\}\cup \left(\left\{P_s(g\vee 0)>0\right\}\cap \left\{P_sg\leq 0\right\}\right). \end{eqnarray*}
\end{proof}

This and Lemma \ref{estimateonE} readily yield, via Corollary \ref{Bs/2-Bsequation}, the following

\begin{lem} \label{chapter8a,Lemma1.9}Suppose $T\in I$ and $\max_{i} x_i^{\left(h\right)}\leq 0$ componentwise. Assume furthermore that there exists a real number $\gamma_1>0$ such that $$P_h\bar f\leq{\gamma_1}^h\bar f.$$ Then for all $s\in(0,T]\cap (2\cdot I)$, \begin{eqnarray*} &&\left\|\left(B_{\frac{s}{2}}\right)^{\circ 2}f-B_sf\right\|_{L^1\left(E^h \cap \bigcap_{k=1}^{m^\frac{s}{2h}} \left(A+ x_k^{\left(\frac{s}{2}\right)}\right)\right)} \\ & \leq & \left\|\left(B_{\frac{s}{2}}\right)^{\circ 2}f-B_sf\right\|_{L^1\left(E^\frac{s}{2}\cap \bigcap_{k=1}^{m^\frac{s}{2h}} \left(A+ x_k^{\left(\frac{s}{2}\right)}\right)\right)}\\ &\leq& \lambda^d\left[\left\{e^{rs/2}g>P_{s/2}(g\vee 0)>P_{s/2} g\right\}\cap A\right]\cdot e^{-r\frac{s}{2}}\\ &&\cdot \left(\left(\ln\gamma_1-r\right)\tilde D +rK+ C_0  \right)\cdot \frac{s}{2}, \end{eqnarray*} wherein $$\tilde D=\chi_{(0,e^r)}\left(\gamma_1\right)\inf_{E^T}\bar f +\chi_{[e^r,+\infty)}\left(\gamma_1\right)\sup_{E^T}\bar f\geq 0 $$ and $$C_0=K\left(\sup_{s\in(0,T]\cap I}\frac{1-e^{-rs}}{s}- r\right)+\left(\sup_{s\in(0,T]\cap I} \frac{{\gamma_1}^se^{-rs}-1}{s}-\ln\gamma_1+ r\right)$$ 
\end{lem}
\begin{proof} Consider $t\in I$. Via our assumption of $\max_{i\in\{1,\dots,m\}} x_i^{\left(h\right)}\leq 0$ componentwise, one has $$\max_k x_k^{\left(t\right)}=\frac{t}{h}\max_{i} x_i^{\left(h\right)}\leq 0$$ componentwise. Since the set $E^t$ is north-east connected, this yields ${E^t+x_k^{\left(t\right)}}\supseteq {E^t}$ for all $k\in\left\{1,\dots,m^{\frac{t}{h}}\right\}$, which in turn -- via $\chi_{E^t}\left(\cdot{-x_k^{\left(t\right)}}\right)=\chi_{E^t+x_k^{\left(t\right)}}\geq \chi_{E^t}$ for all $k\in\left\{1,\dots, m^\frac{t}{h}\right\}$ -- gives 
\begin{eqnarray*} P_{t}\left(\chi_{E^t}f\right)&=& \sum_{k=1}^{m^\frac{t}{h}}\alpha^{(t)}_k\chi_{E^t} \left(\cdot-x_k^{(t)}\right)f\left(\cdot-x_k^{(t)}\right)\\ &\geq& \sum_{k=1}^{m^\frac{t}{h}}\alpha^{(t)}_k\chi_{E^ t} \left(\cdot\right)f\left(\cdot-x_k^{(t)}\right)\\ &=& \chi_{E^t}{P_{t}}f\end{eqnarray*} for all $f\geq 0$. This yields, replacing $f$ by $f\chi_A$,
\begin{eqnarray*}P_{t}\left(\chi_{E^t\cap A}f\right) = P_{t}\left(\chi_{E^t}\cdot\chi_A f\right)&\geq& \chi_{E^t}\cdot{P_t}\left(\chi_Af\right) \\&\geq& \chi_{E^t}\cdot\left(\min_k\chi_{A}\left(\cdot-x_k^{\left(t\right)}\right)\right){P_t}f\\&\geq & \chi_{E^t}\chi_{\bigcap_k\left(A+x_k^{\left(t\right)}\right)} {P_t}f\end{eqnarray*} for all $f\geq 0$. Therefore -- using in addition the translation-invariance of $P_{t}$ and $\lambda^d$ (which makes $P_t$ a map that preserves the $L^1\left(\lambda^d\right)$-norm of nonnegative measurable functions) -- we deduce that for all measurable $f\geq 0$, 
\begin{eqnarray*} \left\|{P_t}f\right\|_{L^1\left({E^t}\cap \bigcap_{k=1}^{m^\frac{t}{h}} \left(A+ x_k^{\left(t\right)}\right)\right)}&\leq& \left\| P_t\left(f\chi_{E^t\cap A}\right) \right\|_{L^1\left(\RR^d\right)}\\ &=&\left\| f\chi_{E^t\cap A} \right\|_{L^1\left(\RR^d\right)}\leq \left\| f \right\|_{L^1\left(E^t\cap A\right)}.\end{eqnarray*} 
From this, using Corollary \ref{Bs/2-Bsequation}, we derive
\begin{eqnarray*} &&\left\|\left(B_{\frac{s}{2}}\right)^{\circ 2}f-B_sf\right\|_{L^1\left(E^\frac{s}{2}\cap \bigcap_k\left(A+x_k^{\left(\frac{s}{2}\right)}\right)\right)}\\ &\leq& \left\|e^{-r\frac{s}{2}}P_{\frac{s}{2}}\left(\chi_{\left\{g>e^{-r\frac{s}{2}}P_{\frac{s}{2}}(g\vee 0)\right\}}\cdot\left(g-e^{-r\frac{s}{2}}P_{\frac{s}{2}}(g\vee 0)\right) \right)\right\|_{L^1\left(E^\frac{s}{2}\cap \bigcap_k\left(A+x_k^{\left(\frac{s}{2}\right)}\right) \right)} \\ &\leq& e^{-r\frac{s}{2}} \left\| \chi_{\left\{g>e^{-r\frac{s}{2}}P_{\frac{s}{2}}(g\vee 0)\right\}}\cdot\left(g-e^{-r\frac{s}{2}}P_{\frac{s}{2}}(g\vee 0)\right) \right\|_{L^1\left(E^\frac{s}{2}\cap A\right)} \end{eqnarray*}
This is enough to prove the Lemma once one takes advantage of Lemma \ref{estimateonE} and Corollary \ref{corestimateonE}.
\end{proof}

\begin{rem} Let the translation-invariant Markov semigroup $P$ be derived from a cubature formula for the Gaussian measure with points $\{y_1,\dots,y_m\}$ in such a way that a geometric Brownian motion with logarithmic drift $\mu=\left(r-\frac{{\sigma_k}^2}{2}\right)_{k\in\{1,\dots,d\}}$ ($r>0$ and $\sigma\in{\RR_+}^d$ being the interest rate of the price process and the volatility vector, respectively) shall be approximated, that is to say $$\forall i\in\{1,\dots,m\}\forall j\in\{1,\dots,d\}\quad \left(x_i^{(h)}\right)_k=\mu_k h+\sigma_k h^{\frac{1}{2}}\left(y_i\right)_k.$$ Then the assumption that all the $x_i^{(h)}$ be componentwise nonpositive for $i\in\{1,\dots,m\}$ reads $$\max_{i} x_i^{\left(h\right)}=\mu h+\left(\max_i\left(y_i\right)_k\cdot \sigma_k h^\frac{1}{2}\right)_{k\in\{1,\dots,d\}}\leq 0$$ and therefore simply means that $\mu h$ is componentwise at least as small or even smaller than $-\left(\max_i\left(y_i\right)_k\cdot \sigma_k h^\frac{1}{2}\right)_{k\in\{1,\dots,d\}}$ which, needless to say, equals \linebreak $\left(\min_i\left(y_i\right)_k\cdot \sigma_k h^\frac{1}{2}\right)_{k\in\{1,\dots,d\}}$ in case of an axis-symmetric cubature formula for the Gaussian measure. This assumption is tantamount to $$\forall k\in\{1,\dots,d\}\quad {\sigma_k}^2-2h^{-\frac{1}{2}}\max_i\left(y_i\right)_k\cdot\sigma_k -2\frac{r}{h}\geq 0, $$ that is \begin{eqnarray*}{\sigma_k}&\in& \RR_+\setminus\left( \begin{array}{c} h^{-\frac{1}{2}}\max_i\left(y_i\right)_k-\sqrt{h^{-1}\cdot\left(\max_i\left(y_i\right)_k\right)^2 + 2\frac{r}{h}}, \\ h^{-\frac{1}{2}}\max_i\left(y_i\right)_k +\sqrt{h^{-1}\cdot\left(\max_i\left(y_i\right)_k\right)^2 + 2\frac{r}{h}}\end{array}\right)\\ &=& \RR_+\setminus\left( \begin{array}{c} h^{-\frac{1}{2}}\left(\max_i\left(y_i\right)_k-\sqrt{\cdot\left(\max_i\left(y_i\right)_k\right)^2 + 2{r}}\right), \\ h^{-\frac{1}{2}}\left(\max_i\left(y_i\right)_k +\sqrt{\cdot\left(\max_i\left(y_i\right)_k\right)^2 + 2r}\right)\end{array}\right)\end{eqnarray*} for all $k\in\{1,\dots,d\}$, entailing that $P$ models a basket of logarithmic asset prices whose volatilities are bounded below by the positive number \linebreak $h^{-\frac{1}{2}}\left(\max_i\left(y_i\right)_k +\sqrt{\cdot\left(\max_i\left(y_i\right)_k\right)^2 + 2r}\right)$.
\end{rem}

Now, emphasising again that our investigations are only concerned with discrete translation-invariant Markov chains $\left(P_t\right)_{t\in I}$ (Markov chains which are derived from cubature formulae, for instance), we can use rather elementary inequalities to find upper bounds on the subsets of $\RR^d$ occurring in the estimates of Lemma \ref{estimateonE}.

We will start with the simple, nevertheless practically important, example of a one-dimensional American vanilla put:

\begin{lem} \label{chapter8a,1dvanilla}Suppose $d=1$ and $\bar f=\exp$. Under these assumptions there exists a $\gamma_1>0$ such that $P_t\bar f={\gamma_1}^t\bar f$, and furthermore, one has for all $s\in I$, $$\left\{e^{rs}g>P_s(g\vee 0)>P_s g\right\}\subseteq\ln K+ \left( \frac{1}{h}\cdot\min_{i\in\left\{1,\dots,m \right\}}x_i^{(h)}, 0\right)\cdot s.$$
\end{lem}
\begin{proof} The real number $\gamma_1$ is given by the relation $${\gamma_1}^h=\sum_{i=1}^m \alpha_i^{(h)}e^{-x_i^{(h)}},$$ that is $$\gamma_1=e^\frac{\ln \left(\sum_{i=1}^m \alpha_i^{(h)}e^{-x_i^{(h)}}\right)}{h}.$$
Next we observe that on the one hand by Remark \ref{Etproperties} \begin{eqnarray*}\left\{P_s(g\vee 0)>P_s g\right\}&=& \left\{\exists k\in\left\{1,\dots,m^\frac{s}{h}\right\}\quad g\left(\cdot-x_k^{(s)}\right)<0\right\} \\  &= &\left\{\min_{k\in\left\{1,\dots,m^\frac{s}{h}\right\}} g\left(\cdot-x_k^{(s)}\right)<0\right\} \\ &=& \left\{K-\max_{k\in\left\{1,\dots,m^\frac{s}{h}\right\}} \exp\left(\cdot-x_k^{(s)}\right)<0\right\}\\ &=&\left\{K- \exp\left(\cdot-\min_{k\in\left\{1,\dots,m^\frac{s}{h}\right\}}x_k^{(s)}\right)<0\right\}\\ &=& \left\{\ln K< \cdot-\min_{k\in\left\{1,\dots,m^\frac{s}{h}\right\}}x_k^{(s)}\right\} \\&=&\left(\ln K + \min_{k\in\left\{1,\dots,m^\frac{s}{h}\right\}}x_k^{(s)}, +\infty\right) \\&=&\left(\ln K + \frac{s}{h}\cdot\min_{i\in\left\{1,\dots,m \right\}}x_i^{(h)}, +\infty\right) \end{eqnarray*} and secondly 
$$\{g>0\}=\left\{K>\exp\right\}=\left(-\infty, \ln K\right),$$ thus \begin{eqnarray*}\left\{e^{rs}g>P_s(g\vee 0)>P_s g\right\}&\subseteq &\left(\ln K + \frac{s}{h}\cdot\min_{i\in\left\{1,\dots,m \right\}}x_i^{(h)}, \ln K\right)\\ &=& \ln K+ \left( \frac{1}{h}\cdot\min_{i\in\left\{1,\dots,m \right\}}x_i^{(h)}, 0\right)\cdot s \end{eqnarray*}
\end{proof}

Applying the preceding Lemmas and using Corollary \ref{Bs/2-Bsequation}, we conclude by stating

\begin{Th}\label{estimateBs/2-BsonE} Suppose $d=1$ and $\bar f=\exp$. Under these assumptions there is a $\gamma_1>0$ such that $P_t\bar f={\gamma_1}^t\bar f$ for all $t\in I$. Assume, moreover, that $$\forall {i\in\{1,\dots,m\}}\forall j\in\{1,\dots, d\}\quad \left(x_i^{\left(h\right)}\right)_j\leq 0.$$ Then there is a real constant $D$ such that for all $s\in(0,T]\cap (2\cdot I)$ and for all $f\geq g\vee 0$, $$\left\|\left(B_{\frac{s}{2}}\right)^{\circ 2}f - B_sf\right\|_{L^1\left(E^h\right)}\leq \left\|\left(B_{\frac{s}{2}}\right)^{\circ 2}f - B_sf\right\|_{L^1\left(E^\frac{s}{2}\right)}\leq\frac{D}{2}\cdot{s}^{2}.$$ We can compute $D$ explicitly as $$D=\left(\left(\ln\gamma_1-r\right)\tilde{D} +rK+ C_0  \right)\cdot \frac{\min_i x_i^{(h)}}{h} $$
\end{Th}
\begin{proof} One only has to apply Lemma \ref{chapter8a,Lemma1.9} for $A=\RR^d$, which one is entitled to by Lemma \ref{chapter8a,1dvanilla}.
\end{proof}

Now we shall proceed to establish convergence estimates for the sequence $\left(B_{T\cdot 2^{-n}}f\right)_{n\in\NN}$ in the $L^1(E^h\cap A)$-norm, for all measurable $f\geq g\vee 0$ and measurable $A\subseteq\RR^d$.

\begin{lem}\label{B_t-B_s_estimate} Suppose $d=1$ and $\bar f=\exp$. Under these assumptions there is a $\gamma_1$ such that $P_t\bar f={\gamma_1}^t\bar f$. Assume, moreover, that $$\forall {i\in\{1,\dots,m\}}\forall j\in\{1,\dots, d\}\quad \left(x_i^{\left(h\right)}\right)_j\leq 0.$$ Under these assumptions there exists a real number $D>0$ (the same as in Theorem \ref{estimateBs/2-BsonE}) such that for all $k\in\NN_0$, $s\in(0,T]\cap \left(2^{k+1}\cdot I\right)$ and measurable $f\geq g\vee 0$, one has $$\left\|\left(B_{s\cdot 2^{-(k+1)}}\right)^{\circ\left(2^{k+1}\right)}f - \left(B_{s\cdot 2^{-k}}\right)^{\circ\left(2^k\right)}f\right\|_{L^1\left(\lambda^1\left[E^h\cap\cdot\right]\right)}\leq D\cdot {s}^2\cdot{2}^{-(k+1)}.$$ 
\end{lem}

The proof is contrived inductively, the base step being Theorem \ref{estimateBs/2-BsonE}, and the induction step being the first part of Lemma \ref{from_k=1_to_kinN}. However, the second and more general part of Lemma \ref{from_k=1_to_kinN} -- which we will need later on in this Chapter when we study options on multiple assets -- requires the following auxiliary result.

\begin{lem}\label{B_tcontractsL1Eh} Let $t\in I$, $A\subseteq\RR^d$ measurable, and assume $$\max_{i\in\{1,\dots,m\}}x^{(h)}\leq 0$$ (which due to $I=h\NN_0$ is equivalent to $\max_{i\in\left\{1,\dots,m^\frac{s}{h}\right\}}x^{(s)}\leq 0$ for all $s\in I$). Then $\bigcap_{s\in (0,T]\cap I}E^s=E^h$ by Remark \ref{Etproperties}, and for all $f_1\geq f_0\geq g\vee 0$ and $p\in\{1,+\infty\}$, \begin{eqnarray*} && \left\|B_tf_1-B_tf_0\right\|_{L^p\left(\lambda^d\left[E^h\cap \bigcap_{k\in\left\{1,\dots, m^\frac{t}{h}\right\}}\left(A+x_k^{(t)}\right)\cap \cdot\right]\right)}\\ &\leq &e^{-rt}\left\|f_1-f_0\right\|_{L^p\left(\lambda^d\left[E^h\cap A\cap\cdot\right]\right)}\end{eqnarray*}
\end{lem}
\begin{proof} Consider a measurable set $A\subset \RR^d$ and measurable functions $f_0,f_1\geq g\vee 0$. Similarly to the proof of Lemma \ref{B_tcontractsifgeq g}, we observe that due to the monotonicity of $B_t$ and the fact that $B_tf\geq g\vee 0\geq g$ for all $f\geq 0$, \begin{eqnarray*}\left\{B_tf_1=g\right\}&=&\left\{B_tf_0\leq B_tf_1=g\right\} \\ &=& \left\{g\leq B_tf_0\leq B_tf_1=g\right\}\\ &=& \left\{B_tf_1=g\right\}\cap \left\{B_tf_0=g\right\}\\ &\subseteq& \left\{B_tf_1-B_tf_0=0\right\},\end{eqnarray*} that is $$\left\{B_tf_1-B_tf_0\neq 0\right\} \subseteq \left\{B_tf_1\neq g\right\} = \left\{B_tf_1>g\right\}$$ Combining this with the monotonicity of $P_t$ as well as the fact that $E^t-x_i^{(t)}\subseteq E^t$ for all $i$ (which in turn is a consequence of the north-east connectedness of $E^t$ -- cf Remark \ref{Etproperties} -- and the assumption that $x_i^{(t)}\leq 0$ for all $i$), we obtain
\begin{eqnarray}&& 0\leq\left(B_tf_1-B_tf_0\right)\chi_{E^h \cap \bigcap_k\left(A+x_k^{(t)}\right)}\nonumber  \\ &= & \left(B_tf_1-B_tf_0\right)\chi_{E^h \cap \bigcap_k\left(A+x_k^{(t)}\right)\cap \left\{B_tf_1>g\right\}}\nonumber  \\ &= & \chi_{E^h \cap  \bigcap_k\left(A+x_k^{(t)}\right)\cap\left\{B_tf_1>g\right\}}\left(e^{-rt}P_tf_1\vee g-e^{-rt}P_tf_0\vee g\right) \nonumber \\ &=& \chi_{E^h \cap  \bigcap_k\left(A+x_k^{(t)}\right)\cap\left\{e^{-rt}P_tf_1>g\right\}}\left(e^{-rt}P_tf_1-e^{-rt}P_tf_0\vee g\right) \nonumber \\ &\leq&\chi_{E^h \cap \bigcap_k\left(A+x_k^{(t)}\right)\cap \left\{e^{-rt}P_tf_1>g\right\}}\left(e^{-rt}P_tf_1-e^{-rt}P_tf_0\right) \nonumber \\ &\leq& e^{-rt}\chi_{E^h \cap \bigcap_k\left(A+x_k^{(t)}\right)}\left(P_tf_1-P_tf_0\right) \nonumber \\ &=&  e^{-rt}\sum_{i=1}^{m^{\frac{t}{h}}}\alpha_i^{(t)}\chi_{E^h \cap \bigcap_k\left(A+x_k^{(t)}\right)}\left(f_1\left(\cdot-x_i^{(t)}\right)-f_0\left(\cdot-x_i^{(t)}\right)\right) \nonumber \\ &=&  e^{-rt}\sum_{i=1}^{m^{\frac{t}{h}}}\alpha_i^{(t)}\left(\chi_{\left(E^h \cap \bigcap_k\left(A+x_k^{(t)}\right)\right)-x_i^{(t)}}\left(f_1-f_0\right)\right)\left(\cdot-x_i^{(t)}\right) \nonumber \\ &=&  e^{-rt}\sum_{i=1}^{m^{\frac{t}{h}}}\alpha_i^{(t)}\left(\chi_{\left(E^h-x_i^{(t)}\right) \cap \bigcap_k\left(A+x_k^{(t)}-x_i^{(t)}\right)}\left(f_1-f_0\right)\right)\left(\cdot-x_i^{(t)}\right) \nonumber\end{eqnarray}

Now, since $$\bigcap_{k\in\left\{1,\dots,m^\frac{t}{h}\right\}}\left(A+x_k^{(t)}-x_i^{(t)}\right) \subseteq A$$ for all $i\in\left\{1,\dots,m^\frac{t}{h}\right\}$ and $f_1-f_0\geq 0$, this means
\begin{eqnarray} && \left(B_tf_1-B_tf_0\right)_{E^h \cap \bigcap_k\left(A+x_k^{(t)}\right)}\nonumber  \\ &\leq& e^{-rt}\sum_{i=1}^{m^{\frac{t}{h}}}\alpha_i^{(t)}\left(\chi_{\left(E^h-x_i^{(t)}\right) \cap A}\left(f_1-f_0\right)\right)\left(\cdot-x_i^{(t)}\right) \nonumber 
\end{eqnarray} 

Combining this pointwise estimate with the translation-invariance of the Lebesgue measure yields
\begin{eqnarray}&& \left\|B_tf_1-B_tf_0\right\|_{L^1\left(\lambda^d\left[E^h \cap \bigcap_k\left(A+x_k^{(t)}\right)\cap\cdot\right]\right)}\nonumber  \\ &= &\int_{\RR^d} \left(B_tf_1-B_tf_0\right)\chi_{E^h \cap \bigcap_k\left(A+x_k^{(t)}\right)}d\lambda^d \nonumber \\ &\leq& e^{-rt}\sum_{i=1}^{m^{\frac{t}{h}}}\alpha_i^{(t)}\int_{\RR^d}\left(\chi_{\left(E^h-x_i^{(t)}\right) \cap A}\left(f_1-f_0\right)\right)\left(\cdot-x_i^{(t)}\right) d\lambda^d \nonumber \\ &=& e^{-rt}\sum_{i=1}^{m^{\frac{t}{h}}}\alpha_i^{(t)}\int_{\RR^d}\left(\chi_{\left(E^h-x_i^{(t)}\right) \cap A}\left(f_1-f_0\right)\right) d\lambda^d \nonumber \\ &\leq& e^{-rt}\sum_{i=1}^{m^{\frac{t}{h}}}\alpha_i^{(t)}\int_{E^h\cap A}\left(f_1-f_0\right) d\lambda^d\nonumber \\ &\leq& \int_{\RR^d} e^{-rt}\chi_{E^h \cap A}\left(f_1-f_0\right)  d\lambda^d\nonumber \\ &=& e^{-rt}\left\| f_1-f_0\right\|_{L^1\left(E^h \cap A\right)}\nonumber, \end{eqnarray}
where we have used the inclusion $E^s-x_k^{(t)}\subseteq E^s$ which -- owing to the north-east connectedness of the sets $E^s$ and our assumption $\max_{i\in\{1,\dots,m\}}x^{(h)}\leq 0$ -- holds for arbitrary $k\in\left\{1,\dots,m^\frac{t}{h}\right\}$ and $s,t\in I$ as well as the assumption $f_1-f_0\geq0$.
Similarly, the translation-invariance and the sub-linearity of the $\esssup_{\RR^d}$-norm imply  
\begin{eqnarray}&& \left\|B_tf_1-B_tf_0\right\|_{L^\infty\left(\lambda^d\left[E^h \cap \bigcap_k\left(A+x_k^{(t)}\right)\cap\cdot\right]\right)}\nonumber  \\ &= &\esssup_{\RR^d} \left[ \left(B_tf_1-B_tf_0\right)\chi_{E^h \cap \bigcap_k\left(A+x_k^{(t)}\right)}\right] \nonumber  \\ &\leq& \esssup_{\RR^d}\left[ e^{-rt}\sum_{i=1}^{m^{\frac{t}{h}}}\alpha_i^{(t)}\left(\chi_{\left(E^h-x_i^{(t)}\right) \cap A}\left(f_1-f_0\right)\right)\left(\cdot-x_i^{(t)}\right) \right]\nonumber \\ &\leq& e^{-rt}\sum_{i=1}^{m^{\frac{t}{h}}}\alpha_i^{(t)}\esssup_{\RR^d}\left[\left(\chi_{\left(E^h-x_i^{(t)}\right) \cap A} \left(f_1-f_0\right)\right)\left(\cdot-x_i^{(t)}\right) \right]\nonumber \\ &\leq& e^{-rt}\sum_{i=1}^{m^{\frac{t}{h}}}\alpha_i^{(t)}\esssup_{\RR^d}\left[\chi_{\left(E^h-x_i^{(t)}\right) \cap A} \left(f_1-f_0\right)\right] \nonumber \\ &\leq& e^{-rt}\sum_{i=1}^{m^{\frac{t}{h}}}\alpha_i^{(t)}\esssup_{\RR^d}\left[\chi_{E^h \cap A} \left(f_1-f_0\right)\right]\nonumber \\ &=& e^{-rt}\sum_{i=1}^{m^{\frac{t}{h}}}\alpha_i^{(t)}\esssup_{E^h \cap A} \left(f_1-f_0\right) \nonumber \\ &=& e^{-rt}\esssup_{E^h \cap A} \left(f_1-f_0\right) \nonumber  \end{eqnarray}
where again one has exploited the inclusion $E^s-x_k^{(t)}\subseteq E^s$ that holds for any $k\in\left\{1,\dots,m^\frac{t}{h}\right\}$ and $s,t\in I$.

\end{proof}

\begin{lem}\label{from_k=1_to_kinN} Let $T\in I$ and $p\in\{1,+\infty\}$. Consider a real number $D'>0$ and a measurable set $C\subseteq\RR^d$. Suppose one has an estimate of the kind \begin{eqnarray*}&&\forall f\geq g\vee 0\forall s\in (2\cdot I)\cap(0,T)\\ && \left\|\left(B_{\frac{s}{2}}\right)^{\circ 2}f - B_sf\right\|_{L^p\left(\lambda^d\left[C \cap\cdot\right]\right)}\leq \frac{D'}{2} \cdot{s}^{2}.\end{eqnarray*} Assume, moreover, $\max_{i\in\{1,\dots,m\}}x^{(h)}\leq 0$ (which by the Chapman-Kolmogorov equation is firstly equivalent to $x_k^{(t)}\leq 0$ for all $k\in\left\{1,\dots,m^\frac{t}{h}\right\}$ and $t\in I$ and secondly also entails $\bigcap_{t\in (0,T]\cap I}E^t=E^h$). Then we get for all measurable $f\geq g\vee 0$ and for all $k\in\NN_0$, $s>0$  such that $s\in(0,T)\cap \left(2^{k+1}\cdot I\right)$, the estimate \begin{eqnarray*}&&\left\|\left(B_{s\cdot 2^{-(k+1)}}\right)^{\circ\left(2^{k+1}\right)}f - \left(B_{s\cdot 2^{-k}}\right)^{\circ\left(2^k\right)}f\right\|_{L^p\left(\lambda^d\left[C \cap\cdot\right]\right)}\\&\leq& D'\cdot {s}^2\cdot{2}^{-(k+1)}.\end{eqnarray*} Furthermore, if one assumes in addition $$0\in\left\{ x_i^{(h)}\ : \ i\in\{1,\dots,m\}\right\},$$ then one has a related implication for $L^p\left(E^h\cap \bigcap_{i\in\left\{1,\dots,m^{\frac{s}{2}}\right\}}\left(A+x_i^{\left(\frac{s}{2}\right)}\right)\right)$ instead of $L^p\left(C \right)$ for all measurable $A\subset\RR^d$: If under these assumptions the assertion \begin{eqnarray*}&&\forall f\geq g\vee 0\forall s\in (2\cdot I)\cap(0,T)\\ &&\left\|\left(B_{\frac{s}{2}}\right)^{\circ 2}f - B_sf\right\|_{L^p\left(E^h\cap \bigcap_{i}\left(A+x_i^{\left(s\right)}\right)\right)} \leq \frac{D'}{2} \cdot{s}^{2}\end{eqnarray*} holds, then the estimate \begin{eqnarray*}&& \forall f\geq g\vee 0\\ && \left\|\left(B_{s\cdot 2^{-(k+1)}}\right)^{\circ\left(2^{k+1}\right)}f - \left(B_{s\cdot 2^{-k}}\right)^{\circ\left(2^k\right)}f\right\|_{L^p\left(E^h \cap \bigcap_{i}\left(A+x_i^{\left(s\right)}\right)\right)}\\ &\leq &D'\cdot {s}^2\cdot{2}^{-(k+1)}\end{eqnarray*} holds for all $k\in\NN_0$ and $s>0$ such that $s\in(0,T)\cap \left(2^{k+1}\cdot I\right)$.

\end{lem}

\begin{proof} For both parts of the Lemma, we will conduct an induction in $k\in\NN_0$, the initial (or base) step being tautological each time. We have for all $s\in(0,T)\cap \left(2^{k+1}\cdot I\right)$ and $f\geq g\vee 0$ the estimate \begin{eqnarray*} &&\left(B_{s\cdot 2^{-(k+1)}}\right)^{\circ\left(2^{k+1}\right)}f - \left(B_{s\cdot 2^{-k}}\right)^{\circ2^k}f\\ &=& \left(B_{s\cdot 2^{-(k+1)}}\right)^{\circ\left(2^{k}\right)}\circ\left(B_{s\cdot 2^{-(k+1)}}\right)^{\circ\left(2^{k}\right)}f \\ && -\left(B_{s\cdot 2^{-(k+1)}}\right)^{\circ\left(2^{k}\right)}\circ\left(B_{s\cdot 2^{-k}}\right)^{\circ\left(2^{k-1}\right)}f \\ &&  + \left(B_{s\cdot 2^{-(k+1)}}\right)^{\circ\left(2^{k}\right)}\circ\left(B_{s\cdot 2^{-k}}\right)^{\circ\left(2^{k-1}\right)}f \\ && - \left(B_{s\cdot 2^{-k}}\right)^{\circ\left(2^{k-1}\right)}\circ \left(B_{s\cdot 2^{-k}}\right)^{\circ\left(2^{k-1}\right)}f \\ &=& \left(B_{s\cdot 2^{-(k+1)}}\right)^{\circ\left(2^{k}\right)}\circ\left(\left(B_{\frac{s}{2}\cdot 2^{-k}}\right)^{\circ\left(2^{k}\right)}-\left(B_{\frac{s}{2}\cdot 2^{-(k-1)}}\right)^{\circ\left(2^{k-1}\right)}\right)f \\ &&  + \left(\left(B_{\frac{s}{2}\cdot 2^{-k}}\right)^{\circ\left(2^{k}\right)} - \left(B_{ \frac{s}{2} \cdot 2^{-(k-1)}}\right)^{\circ\left(2^{k-1}\right)}\right)\circ \left(B_{s\cdot 2^{-k}}\right)^{\circ\left(2^{k-1}\right)}f \end{eqnarray*} which plays a crucial part in both the first and the second part of the Lemma.
For, we can first of all note that the induction hypothesis in the situation of the first part of the Lemma reads \begin{eqnarray}&&\forall f\geq g\vee 0\forall t\in (2^k\cdot I)\cap(0,T)\nonumber \\ &&  \left\| \left(B_{t\cdot 2^{-k}}\right)^{\circ\left(2^{k}\right)}f - \left(B_{ t\cdot 2^{-(k-1)}}\right)^{\circ\left(2^{k-1}\right)}f \right\|_{L^p\left(C\right)}\nonumber\\ &\leq& D'\cdot t^2\cdot 2^{-k}.\label{chapter8a,inductionhypothesis1}
\end{eqnarray} 
And if one now applies this induction hypothesis (\ref{chapter8a,inductionhypothesis1}) for $t=\frac{s}{2}$ (recalling that by assumption $s\in 2^{k+1}\cdot I$, thus $\frac{s}{2}\in 2^k\cdot I$) to the previous two equations and uses Lemma \ref{B_tcontractsL1Eh}, then one gets by the triangle inequality for the ${L^p\left(C\right)}$-norm,\begin{eqnarray*} && \left\|\left(B_{s\cdot 2^{-(k+1)}}\right)^{\circ\left(2^{k+1}\right)}f - \left(B_{s\cdot 2^{-k}}\right)^{\circ2^k}f\right\|_{L^p\left(C\right)} \\ &\leq & \left\|\left(\left(B_{\frac{s}{2}\cdot 2^{-k}}\right)^{\circ\left(2^{k}\right)} - \left(B_{ \frac{s}{2} \cdot 2^{-(k-1)}}\right)^{\circ\left(2^{k-1}\right)}\right)\circ \left(B_{s\cdot 2^{-k}}\right)^{\circ\left(2^{k-1}\right)}f \right\|_{L^p\left(C\right)}\\ && + \left\| \left(B_{s\cdot 2^{-(k+1)}}\right)^{\circ\left(2^{k}\right)}\circ\left(\left(B_{\frac{s}{2}\cdot 2^{-k}}\right)^{\circ\left(2^{k}\right)}-\left(B_{\frac{s}{2}\cdot 2^{-(k-1)}}\right)^{\circ\left(2^{k-1}\right)}\right)f \right\|_{L^p\left(C\right)}\\ &\leq & \left\|\left(\left(B_{\frac{s}{2}\cdot 2^{-k}}\right)^{\circ\left(2^{k}\right)} - \left(B_{ \frac{s}{2} \cdot 2^{-(k-1)}}\right)^{\circ\left(2^{k-1}\right)}\right)\circ \left(B_{s\cdot 2^{-k}}\right)^{\circ\left(2^{k-1}\right)}f \right\|_{L^p\left(C\right)} \\ && + \left\| \left(\left(B_{\frac{s}{2}\cdot 2^{-k}}\right)^{\circ\left(2^{k}\right)}-\left(B_{\frac{s}{2}\cdot 2^{-(k-1)}}\right)^{\circ\left(2^{k-1}\right)}\right)f \right\|_{L^p\left(C\right)} \\ &\leq & D'\cdot\frac{s^2}{4}\cdot 2^{-k}+D'\cdot\frac{s^2}{4}\cdot 2^{-k} = D'\cdot{s}^2\cdot 2^{-(k+1)}.\end{eqnarray*} In order to be entitled to apply Lemma \ref{B_tcontractsL1Eh} in this situation we have successively used the fact that $$\forall t\in I \forall \ell\geq g\vee 0 \quad B_t\ell\geq g\vee 0.$$ This completes the induction step for the first part of the Lemma. 

Turning to the proof of the second assertion in the Lemma (where $\max_{i\in\{1,\dots,d\}}x_i^{(h)}\leq 0$ is assumed), we remark that \begin{eqnarray}&&\forall t\in I\forall k\leq\frac{\ln t-\ln h}{\ln 2} \nonumber \\ A(t)&:=& \bigcap_{i\in\left\{1,\dots,m^{\frac{t}{h}}\right\}} \left(A+x_i^{\left({t}\right)}\right)\nonumber \\&=& \bigcap_{\ell,k\in\left\{1,\dots,m^{\frac{t}{h}}\right\}} \left(A+x_k^{\left(\frac{t}{2}\right)}+x_\ell^{\left(\frac{t}{2}\right)}\right)\nonumber \\ &=& \left(A\left(\frac{t}{2}\right)\right)\left(\frac{t}{2}\right)\\&=& \label{A(t)properties} \bigcap_{i_1,\dots,i_{2^k}\in\left\{1,\dots,m^{\frac{t}{h}\cdot2^{-k-1}}\right\}} \left(A+x_{i_1}^{\left( t \cdot 2^{-(k+1)}\right)}+\cdots+x_{i_{2^k}}^{\left( {t} \cdot 2^{-(k+1)}\right)}\right).\end{eqnarray} In particular, if $0\in\left\{ x_i^{(h)}\ : \ i\in\{1,\dots,m\}\right\}$, $A(s)$ is decreasing in $s$: $$\forall s,t\in I\left(s\leq t\Rightarrow A(s)\supseteq A(t)\right).$$ Similarly to proof of the first part of the present Lemma, we deduce \begin{eqnarray*} && \left\|\left(B_{s\cdot 2^{-(k+1)}}\right)^{\circ\left(2^{k+1}\right)}f - \left(B_{s\cdot 2^{-k}}\right)^{\circ2^k}f\right\|_{L^p\left(E^h \cap A\left({s}\right) \right)} \\ &\leq & \left\|\left(B_{s\cdot 2^{-(k+1)}}\right)^{\circ\left(2^{k+1}\right)}f - \left(B_{s\cdot 2^{-k}}\right)^{\circ2^k}f\right\|_{L^p\left(E^h \cap A\left(s\right) \right)} \\ &\leq & \left\|\left(\left(B_{\frac{s}{2}\cdot 2^{-k}}\right)^{\circ\left(2^{k}\right)} - \left(B_{ \frac{s}{2} \cdot 2^{-(k-1)}}\right)^{\circ\left(2^{k-1}\right)}\right)\circ \left(B_{s\cdot 2^{-k}}\right)^{\circ\left(2^{k-1}\right)}f \right\|_{L^p\left(E^h \cap A\left(s\right) \right)} \\ && + \left\| \left(B_{s\cdot 2^{-(k+1)}}\right)^{\circ\left(2^{k}\right)}\circ\left(\left(B_{\frac{s}{2}\cdot 2^{-k}}\right)^{\circ\left(2^{k}\right)}-\left(B_{\frac{s}{2}\cdot 2^{-(k-1)}}\right)^{\circ\left(2^{k-1}\right)}\right)f \right\|_{L^p\left(E^h \cap A\left(2^{k+1}\frac{s}{2}\right) \right)} \\ &\leq & \left\|\left(\left(B_{\frac{s}{2}\cdot 2^{-k}}\right)^{\circ\left(2^{k}\right)} - \left(B_{ \frac{s}{2} \cdot 2^{-(k-1)}}\right)^{\circ\left(2^{k-1}\right)}\right)\circ \left(B_{s\cdot 2^{-k}}\right)^{\circ\left(2^{k-1}\right)}f \right\|_{L^p\left(E^h \cap A\left(s\right) \right)} \\ && + \left\| \left(B_{s\cdot 2^{-(k+1)}}\right)^{\circ\left(2^{k}\right)}\circ\left(\left(B_{\frac{s}{2}\cdot 2^{-k}}\right)^{\circ\left(2^{k}\right)}-\left(B_{\frac{s}{2}\cdot 2^{-(k-1)}}\right)^{\circ\left(2^{k-1}\right)}\right)f \right\|_{L^p\left(E^h \cap A\left(s\right)\right)} \end{eqnarray*} from the triangle inequality.
But by a successive application of Lemma \ref{B_tcontractsL1Eh}, combined with the properties (\ref{A(t)properties}) of $A(\cdot)$, we have for all $f_1\geq f_0\geq g\vee 0$,
\begin{eqnarray*} && \left\|\left(B_{s\cdot 2^{-(k+1)}}\right)^{\circ\left(2^{k}\right)}\circ\left(f_1-f_0\right)\right\|_{L^p\left(E^h\cap A(s)\right)}\\ &\leq& \left\|\left(B_{s\cdot 2^{-(k+1)}}\right)^{\circ\left(2^{k}\right)}\circ\left(f_1-f_0\right)\right\|_{L^p\left(E^h\cap \bigcap_{\ell\in\left\{1,\dots,m^\frac{s}{2h}\right\}} \left( A\left(\frac{s}{2}\right) +x_\ell^{\left(\frac{s}{2}\right)}\right)\right)} \\ &=&  \left\|\left(B_{s\cdot 2^{-(k+1)}}\right)^{\circ\left(2^{k}\right)}\circ\left(f_1-f_0\right)\right\|_{L^p\left(E^h\cap \bigcap_{i_1,\dots,i_{2^k}\in\left\{1,\dots,m^\frac{s}{2^{k+1}h}\right\}} \left(A\left(\frac{s}{2}\right)+x_{i_{1}}^{\left(\frac{s}{2^{k+1}}\right)}+\dots+x_{i_{2^k}}^{\left(\frac{s}{2^{k+1}}\right)}\right)\right)}\\ &\leq& \left\|\left(B_{s\cdot 2^{-(k+1)}}\right)^{\circ\left(2^{k}-1\right)}\circ\left(f_1-f_0\right)\right\|_{L^p\left(E^h\cap \bigcap_{i_1,\dots,i_{2^k}\in\left\{1,\dots,m^\frac{s}{2^{k+1}h}\right\}} \left(A\left(\frac{s}{2}\right)+x_{i_{1}}^{\left(\frac{s}{2^{k+1}}\right)}+\dots+x_{i_{2^k-1}}^{\left(\frac{s}{2^{k+1}}\right)}\right)\right)}\\ &\leq& \\ &\vdots&\\ &\leq& \left\|f_1-f_0\right\|_{L^p\left(E^h\cap A\left(\frac{s}{2}\right)\right)}.\end{eqnarray*} In light of the inclusion $A(s)\supseteq A\left(\frac{s}{2}\right)$, we also have
\begin{eqnarray*} && \left\|\left(B_{s\cdot 2^{-(k+1)}}\right)^{\circ\left(2^{k+1}\right)}f - \left(B_{s\cdot 2^{-k}}\right)^{\circ2^k}f\right\|_{L^p\left(E^h \cap A\left({s}\right) \right)}\\ &\leq & \left\|\left(\begin{array}{c}\left(B_{\frac{s}{2}\cdot 2^{-k}}\right)^{\circ 2^{k}} -\\ \left(B_{ \frac{s}{2} \cdot 2^{-(k-1)}}\right)^{\circ 2^{k-1}}\end{array}\right)\circ \left(B_{s\cdot 2^{-k}}\right)^{\circ 2^{k-1}}f \right\|_{L^p\left(E^h \cap A\left(s\right) \right)} \\ && + \left\| \left(\left(B_{\frac{s}{2}\cdot 2^{-k}}\right)^{\circ\left(2^{k}\right)}-\left(B_{\frac{s}{2}\cdot 2^{-(k-1)}}\right)^{\circ\left(2^{k-1}\right)}\right)f \right\|_{L^p\left(E^h \cap A\left(\frac{s}{2}\right) \right)}\end{eqnarray*} Combining the previous two sets of estimates leads to 
\begin{eqnarray*} && \left\|\left(B_{s\cdot 2^{-(k+1)}}\right)^{\circ\left(2^{k}\right)}\circ\left(f_1-f_0\right)\right\|_{L^p\left(E^h\cap A(s)\right)}\\ &\leq& \left\|\left(\begin{array}{c}\left(B_{\frac{s}{2}\cdot 2^{-k}}\right)^{\circ 2^{k}} -\\ \left(B_{ \frac{s}{2} \cdot 2^{-(k-1)}}\right)^{\circ 2^{k-1}}\end{array}\right)\circ \left(B_{s\cdot 2^{-k}}\right)^{\circ 2^{k-1}}f \right\|_{L^p\left(E^h \cap A\left(\frac{s}{2}\right) \right)} \\ && + \left\| \left(\left(B_{\frac{s}{2}\cdot 2^{-k}}\right)^{\circ\left(2^{k}\right)}-\left(B_{\frac{s}{2}\cdot 2^{-(k-1)}}\right)^{\circ\left(2^{k-1}\right)}\right)f \right\|_{L^p\left(E^h \cap A\left(\frac{s}{2}\right) \right)} \\ &\leq & D'\cdot\frac{s^2}{4}\cdot 2^{-k}+D'\cdot\frac{s^2}{4}\cdot 2^{-k} = D'\cdot{s}^2\cdot 2^{-(k+1)},\end{eqnarray*} where in the last line we have taken advantage of the induction hypothesis \begin{eqnarray*}&&\forall k\in\NN_0\forall f\geq g\vee 0\forall t>0 \\ &&\left(\begin{array}{c}t\in (2^k\cdot I)\cap(0,T)\Rightarrow\\\left\| \left(B_{t\cdot 2^{-k}}\right)^{\circ\left(2^{k}\right)}f - \left(B_{ t\cdot 2^{-(k-1)}}\right)^{\circ\left(2^{k-1}\right)}f \right\|_{L^p\left(E^h \cap \bigcap_{i}\left(A+x_i^{\left(2^{k} {t} \right)}\right)\right)}\\ \leq D'\cdot t^2\cdot 2^{-k} \end{array}\right)\end{eqnarray*} for the special case $t=\frac{s}{2}$
\end{proof}

The assumption of $0\in\left\{ x_i^{(h)}\ : \ i\in\{1,\dots,m\}\right\}\subseteq\RR^d$ while $\max_{i\in\{1,\cdots,m\}}x_i^{(h)}\leq 0$ componentwise corresponds to the volatility attaining a certain critical value:

\begin{rem}\label{scaling remark} Consider a cubature formula for the one-dimensional Gaussian measure with cubature points $\{z_1,\dots,z_m\}$ which will then give rise to a new Markov chain via $$\forall i\in\{1,\dots,m\} \quad x_i^{(h)} = \left(r-\frac{{\sigma_0}^2}{2}\right)h+ z_i\sigma_0 h^\frac{1}{2}$$ (if simply $\{z_1,z_2\}=\{\pm 1\}$, then this was a discrete model for a logarithmic asset price evolution that converges weakly to the Black-Scholes model with volatility $\sigma_0$ and discount rate $r>0$ when $h\downarrow 0$). In this setting, the set of pairs $(r,\sigma_0)\in\RR^2$ such that $$ \max_i x_i^{(h)}=\left(r-\frac{{\sigma_0}^2}{2}\right)h + \sigma_0h^\frac{1}{2}\cdot\max_iz_i=0$$ has at most two elements, ie it is a Lebesgue null set. However, in practice, we will not have the exact values of the volatility $\sigma_0$ (and if the maturity is sufficiently large, one will not even have an exact value for the interest rate $r$), but we will only know that $\sigma_0\in \left(\tilde \sigma-\varepsilon,\tilde \sigma+\varepsilon\right)$ for some $\varepsilon>0$. So, given $r>0$ the set of volatility parameters that both fit the model and allow for the previous Lemma to be applied will equal $$\left\{ \sigma_0\in \left(\tilde \sigma-\varepsilon,\tilde \sigma+\varepsilon\right)\ :\ \left(r-\frac{{\sigma_0}^2}{2}\right)h + \sigma_0h^\frac{1}{2}\cdot\max_iz_i=0\right\}. $$ If $\varepsilon>0$ and the equation characterising this set has a solution $\sigma'\in \left(\tilde \sigma-\varepsilon,\tilde \sigma+\varepsilon\right)$, this set will at least have positive Lebesgue measure, so that there is some hope that our condition of $ \max_i x_i^{(h)}=0$ (which we had to impose in the second part of the previous Lemma \ref{from_k=1_to_kinN}) can be satisfied in practice at least occasionally.
\end{rem}

With the first half of Lemma \ref{from_k=1_to_kinN}, we have completed the proof of Lemma \ref{B_t-B_s_estimate}. We shall now apply this result to finally get to a convergence bound for $\left(B_{T\cdot 2^{-n}}(g\vee 0)\right)_n$ -- which can be conceived of as a sequence of non-perpetual Bermudan option prices when successively halving the exercise mesh size.

\begin{lem}\label{from_k,k+1_to_M,N} Let $p\in[1,+\infty]$. Consider a real constant $D>0$ as well as a measurable set $C$ and a set $E$ of nonnegative measurable functions, and suppose one has an estimate of the kind \begin{eqnarray*}&&\forall k\in\NN_0\forall f\in E\forall s\in (2^{k+1}\cdot I)\cap(0,T)\\ &&  \left\|\left(B_{s\cdot 2^{-(k+1)}}\right)^{\circ\left(2^{k+1}\right)}f - \left(B_{s\cdot 2^{-k}}\right)^{\circ\left(2^k\right)}f\right\|_{L^p\left(\lambda^d\left[C\cap \cdot\right]\right)}\leq D \cdot {s}^2\cdot{2}^{-(k+1)} .\end{eqnarray*}
Then for all $N>M\in\NN$, $s\in(0,T)\cap \left(2^N\cdot I\right)$ and $f\in E$, the estimate \begin{eqnarray*}\left\|\left(B_{s\cdot 2^{-N}}\right)^{\circ\left(2^{N}\right)}f - \left(B_{s\cdot 2^{-M}}\right)^{\circ\left(2^M\right)}f\right\|_{L^p(C)}&\leq& D\cdot {s^2}\cdot{2}^{-M}\left(1- 2^{-(N-M-1)}\right) \\ &\leq& D\cdot {s^2}\cdot{2}^{-M}\end{eqnarray*}
holds.
\end{lem}
\begin{proof} With $M,N$, $s$, $f$ as in the statement of the Theorem, we obtain by the triangle inequality \begin{eqnarray*} && \left\|\left(B_{s\cdot 2^{-N}}\right)^{\circ\left(2^{N}\right)}f - \left(B_{s\cdot 2^{-M}}\right)^{\circ\left(2^M\right)}f\right\|_{L^p(C)}\\ &=&\left\|\sum_{k=M}^{N-1} \left(\left(B_{s\cdot 2^{-(k+1)}}\right)^{\circ\left(2^{k+1}\right)}f - \left(B_{s\cdot 2^{-k}}\right)^{\circ\left(2^k\right)}f\right)\right\|_{L^p(C)} \\ &\leq&\sum_{k=M}^{N-1}\left\| \left(\left(B_{s\cdot 2^{-(k+1)}}\right)^{\circ\left(2^{k+1}\right)}f - \left(B_{s\cdot 2^{-k}}\right)^{\circ\left(2^k\right)}f\right)\right\|_{L^p(C)} \\ &\leq & \sum_{k=M}^{N-1}D\cdot s^2\cdot 2^{-k-1}=D\cdot\frac{s^2}{2}\cdot\sum_{k=0}^{N-1-M}2^{-k}2^{-M}\\ &=& D\cdot\frac{s^2}{2}\cdot 2^{-M}\cdot \frac{1-2^{-(N-M-1)}}{1-2^{-1}}\leq D\cdot\frac{s^2}{2}\cdot 2^{-M}\cdot 2.\end{eqnarray*}
This suffices to prove the Theorem.
\end{proof}

Thus, if we combine this last Lemma \ref{from_k,k+1_to_M,N} with Lemma \ref{B_t-B_s_estimate} we arrive at 

\begin{Th} \label{cubatureBermudanconv} Suppose, as before, $d=1$ and $\bar f=\exp$. Under these assumptions there is a $\gamma_1$ such that $P_t\bar f={\gamma_1}^t\bar f$, and let us suppose this $\gamma_1\in(0,e^r]$. Assume, moreover, that $$\forall {i\in\{1,\dots,m\}}\forall j\in\{1,\dots, d\}\quad \left(x_i^{\left(h\right)}\right)_j\leq 0.$$  Under these assumptions there exists a real number $D>0$ such that for all $N>M\in\NN$, $s\in(0,T]\cap \left(2^N\cdot I\right)$ and monotonely decreasing $f\geq g\vee 0$, one has \begin{eqnarray*}&& \left\|\left(B_{s\cdot 2^{-N}}\right)^{\circ\left(2^{N}\right)}f - \left(B_{s\cdot 2^{-M}}\right)^{\circ\left(2^M\right)}f\right\|_{L^1\left(E^h\right)}\\ &\leq& D\cdot {s}^2\cdot{2}^{-M}\left(1- 2^{-(N-M-1)}\right) \\ &\leq& D\cdot {s}^2\cdot{2}^{-M}.\end{eqnarray*} 
\end{Th}

Analogously, we may proceed to prove convergence of higher order in $s$ for $\bar f=\sum_{j=1}^d w_j\exp\left((\cdot)_j\right)$, where $w_1,\dots,w_d$ is a convex combination (the weights for a weighted average of the components/assets in a $d$-dimensional basket), as well as for the choices $\bar f=\min_{j\in\{1,\dots,d\}}\exp\left((\cdot)_j\right)$ and $\bar f=\max_{j\in\{1,\dots,d\}}\exp\left((\cdot)_j\right)$. However, this time, we shall employ different norms: $L^1\left(\lambda^d\left[E^h\cap A\cap\cdot\right]\right)$ for a compact subset $A\subset \RR^d$ such that $\lambda^d\left[\bigcap_{s\in (0,T]\cap I}E^s \cap A\right]\in (0,+\infty)$.

The first part of this endeavour will be to prove certain generalisations of Lemmas \ref{estimateonE} and \ref{B_tcontractsL1Eh}.

\begin{lem} If $\bar f=\sum_{j=1}^d w_j\exp\left((\cdot)_j\right)$, then $$P_s\bar f\leq {\gamma_1}^s\bar f $$ where $$\gamma_1:=\left(\underbrace{ \max_{j\in\{1,\dots,d\}} \sum_{i=1}^m \alpha_i^{(h)} e^{\left(-x_i^{(h)}\right)_j} }_{>0} \right) ^{\frac{1}{h}}.$$
\end{lem}

\begin{proof} We have for all $s\in I$ the estimate
\begin{eqnarray*}P_s\bar f&=&\sum_{i=1}^{m^\frac{s}{h}}\alpha_i^{(s)} \bar f\left(\cdot-x_i^{(s)}\right)\\ &=& \sum_{j=1}^d \sum_{i=1}^{m^\frac{s}{h}} w_j\alpha_i^{(s)} \exp\left(\left(\cdot-x_i^{(s)}\right)_j\right)\\ &=& \sum_{j=1}^d \sum_{i=1}^{m^\frac{s}{h}} \alpha_i^{(s)} e^{\left(-x_i^{(s)}\right)_j} w_j\exp\left(\left(\cdot\right)_j\right)\\ &\leq& \underbrace{ \left(\max_{j\in\{1,\dots,d\}} \sum_{i=1}^{m^\frac{s}{h}} \alpha_i^{(s)} e^{\left(-x_i^{(s)}\right)_j} \right) }_{>0}\sum_{j=1}^d w_j\exp\left(\left(\cdot\right)_j\right), \end{eqnarray*} in particular this estimate holds for $s=h$. But this is to say $$P_h\bar f\leq {\gamma_1}^h\bar f ,$$ hence we have proven the estimate in the Lemma for $s=h$. This readily suffices to prove the Lemma's assertion in its full generality, as $\left(P_s\right)_{s\in I}$ is a Markov semigroup and by applying the Chapman-Komogorov equation (and the monotonicity of $P_h$) inductively, $$\forall n\in\NN\quad P_{nh}\bar f=\underbrace{P_{h}\cdots P_{h}}_{n}\bar f\leq \underbrace{{\gamma_1}^h\cdots \gamma_1^h}_{n}\bar f={\gamma_1}^{hn}\bar f.$$
\end{proof}

\begin{lem}\label{barf=min_gamma1} If $\bar f=\min_{j\in\{1,\dots,d\}}\exp\left((\cdot)_j\right)$, then $$P_s\bar f\leq {\gamma_1}^s\bar f $$ where $$\gamma_1:=\left(\underbrace{ \sum_{i=1}^{m^\frac{s}{h}} \alpha_i^{(s)} \left(\max_{\ell\in\{1,\dots,d\}} e^{\left(-x_i^{(s)}\right)_\ell}\right) }_{>0} \right) ^{\frac{1}{h}}.$$
\end{lem}

\begin{proof} We have for all $s\in I$ the estimate
\begin{eqnarray*}P_s\bar f&=&\sum_{i=1}^{m^\frac{s}{h}} \alpha_i^{(s)} \bar f\left(\cdot-x_i^{(s)}\right)\\ &=& \sum_{i=1}^{m^\frac{s}{h}} \alpha_i^{(s)} \min_{j\in\{1,\dots,d\}}\exp\left(\left(\cdot-x_i^{(s)}\right)_j\right) \\ &\leq& \sum_{i=1}^{m^\frac{s}{h}} \alpha_i^{(s)}\ \chi_{\left\{k\in\{1,\dots,d\} \ : \ {\tiny \begin{array}{c}\min_{j\in\{1,\dots,d\}} \exp\left(\left(\cdot\right)_j\right) = \exp\left(\left(\cdot\right)_k\right),\\ \forall n> k\quad\min_{j\in\{1,\dots,d\}} \exp\left(\left(\cdot\right)_j\right) < \exp\left(\left(\cdot\right)_n\right)\end{array}} \right\}}(\ell) \\&&\cdot e^{\left(-x_i^{(s)}\right)_\ell} \exp\left(\left(\cdot\right)_\ell\right) \\ &\leq& \underbrace{ \sum_{i=1}^{m^\frac{s}{h}} \alpha_i^{(s)} \left(\max_{\ell\in\{1,\dots,d\}} e^{\left(-x_i^{(s)}\right)_\ell}\right) }_{>0} \min_{j\in\{1,\dots,d\}} \exp\left(\left(\cdot\right)_j\right)\end{eqnarray*} in particular the estimate holds for $s=h$ again. But this means $$P_h\bar f\leq {\gamma_1}^h\bar f ,$$ hence we arrive at the estimate of the Lemma for $s=h$. This readily suffices to prove the Lemma's assertion in its full strength, as $\left(P_s\right)_{s\in I}$ is a Markov semigroup and by applying the Chapman-Komogorov equation inductively, $$\forall n\in\NN\quad P_{nh}\bar f =\underbrace{P_{h}\cdots P_{h}}_{n}\bar f\leq \underbrace{{\gamma_1}^h\cdots {\gamma_1}^h}_{n}\bar f={\gamma_1}^{hn}\bar f.$$
\end{proof}

\begin{lem} \label{barf=max_gamma1}If this time $\bar f=\max_{j\in\{1,\dots,d\}}\exp\left((\cdot)_j\right)$, then $$P_s\bar f\leq {\gamma_1}^s\bar f $$ where $$\gamma_1:=\left(\underbrace{ \left( \sum_{i=1}^{m^\frac{s}{h}} \alpha_i^{(s)} \max_{j\in\{1,\dots,d\}} e^{\left(-x_i^{(s)}\right)_j} \right) }_{>0} \right) ^{\frac{1}{h}}.$$
\end{lem}

\begin{proof} We have for all $s\in I$ the estimate
\begin{eqnarray*}P_s\bar f&=&\sum_{i=1}^{m^\frac{s}{h}} \alpha_i^{(s)} \bar f\left(\cdot-x_i^{(s)}\right)\\ &=& \sum_{i=1}^{m^\frac{s}{h}} \alpha_i^{(s)} \max_{j\in\{1,\dots,d\}}\exp\left(\left(\cdot-x_i^{(s)}\right)_j\right)\\ &\leq& \sum_{i=1}^{m^\frac{s}{h}} \alpha_i^{(s)} \max_{\ell\in\{1,\dots,d\}} e^{\left(-x_i^{(s)}\right)_\ell} \max_{j\in\{1,\dots,d\}} \exp\left(\left(\cdot\right)_j\right)\\ &=& \underbrace{ \sum_{i=1}^{m^\frac{s}{h}} \alpha_i^{(s)} \max_{\ell\in\{1,\dots,d\}} e^{\left(-x_i^{(s)}\right)_\ell} }_{>0}\max_j\exp\left(\left(\cdot\right)_j\right), \end{eqnarray*} in particular the estimate holds for $s=h$. But this is -- as it was in the proofs of the two preceding Lemmas -- to say $$P_h\bar f\leq {\gamma_1}^h\bar f ,$$ hence we have proven the estimate in the Lemma for $s=h$. This readily suffices to prove the Lemma's assertion, as $\left(P_s\right)_{s\in I}$ is a Markov semigroup and by applying the Chapman-Komogorov equation inductively, $$\forall n\in\NN\quad P_{nh} \bar f = \underbrace{P_{h}\cdots P_{h}}_{n}\bar f\leq \underbrace{{\gamma_1}^h\cdots {\gamma_1}^h}_{n}\bar f={\gamma_1}^{hn}\bar f.$$
\end{proof}

Thus at least for certain choices of $\bar f$ -- viz. weighted arithmetic average of the exponential components, minimum of the exponential components and maximum of the exponential components -- we can apply Lemma \ref{estimateonE}.

Therefore we shall next turn our attention to deriving upper bounds for the measures of the sets in the estimates of Lemma \ref{estimateonE} for the said examples of $\bar f=\sum_{j=1}^d w_j\exp\left((\cdot)_j\right)$, $\bar f=\min_{j\in\{1,\dots,d\}}\exp\left((\cdot)_j\right)$ and $\bar f=\max_{j\in\{1,\dots,d\}}\exp\left((\cdot)_j\right)$. We continue to use the notation $I=h\NN_0$ and $$P_s: f\mapsto\sum_{i=1}^{m^\frac{s}{h}}\alpha_i^{(s)} f\left(\cdot-x_i^{(s)}\right),$$ where $\left(P_s\right)_{s\in I}=\left(P_{nh}\right)_{n\in \NN_0}=\left(\underbrace{P_{h}\cdots P_{h}}_{n}\right)_{n\in \NN_0}$ is the Markov chain generated by $P_h$.

\begin{lem} If $\bar f=\min_{j\in\{1,\dots,d\}}\exp\left((\cdot)_j\right)$, then for all $s\in I$, \begin{eqnarray*}&&\left\{\forall i\in\left\{1,\dots,m^\frac{s}{h}\right\}\quad g\left(\cdot-x_i^{(s)}\right) \leq 0\right\}\\ &=& \bigotimes_{j=1}^d\left[\ln K+ \frac{s}{h} \max_{i\in\{1,\dots,m\}} \left(x_i^{(h)}\right)_j,  +\infty\right)\end{eqnarray*} as well as \begin{eqnarray*}&& \left\{\forall i\in\left\{1,\dots,m^\frac{s}{h}\right\}\quad g\left(\cdot-x_i^{(s)}\right)\geq 0\right\}\\ &\subset& \bigcup_{j=1}^d \left( \underbrace{ \RR\times\cdots\times\RR }_{j-1} \times \left(-\infty,\ln K+ \frac{s}{h} \max_{i\in\{1,\dots,m\}} \left(x_i^{(h)}\right)_j \right]\times \underbrace{ \RR\times\cdots\times\RR }_{d-j} \right)\end{eqnarray*}
\end{lem}
\begin{proof} Let $s\in I$. Then 
\begin{eqnarray*}&& \left\{\forall i\in\left\{1,\dots,m^\frac{s}{h}\right\}\quad g\left(\cdot-x_i^{(s)}\right)\leq 0\right\} \\ &= &\left\{\max_{i\in\left\{1,\dots,m^\frac{s}{h}\right\}} g\left(\cdot-x_i^{(s)}\right)\leq 0\right\} \\ &=& \left\{K-\min_{i\in\left\{1,\dots,m^\frac{s}{h}\right\}} \min_{j\in\{1,\dots,d\}} \exp\left(\left(\cdot-x_i^{(s)}\right)_j\right)\leq 0\right\}\\  &=& \left\{K-\min_{j\in\{1,\dots,d\}}\min_{i\in\left\{1,\dots,m^\frac{s}{h}\right\} } \exp\left(\left(\cdot-x_i^{(s)}\right)_j\right)\leq 0\right\}\\  &=& \left\{K- \min_{j\in\{1,\dots,d\}}\exp\left(\left(\cdot\right)_j -\max_{i\in\left\{1,\dots,m^\frac{s}{h}\right\}} \left(x_i^{(s)}\right)_j\right) \leq 0\right\} \\ &=& \left\{\ln K\leq \min_{j\in\{1,\dots,d\}} \left(\left(\cdot\right)_j - \max_{i\in\left\{1,\dots,m^\frac{s}{h}\right\}} \left(x_i^{(s)}\right)_j \right)\right\} \\ &=& \bigcap_{j=1}^d\left\{\ln K+\max_{i\in\left\{1,\dots,m^\frac{s}{h}\right\}} \left(x_i^{(s)}\right)_j\leq (\cdot)_j \right\}\\ &=& \bigcap_{j=1}^d\left\{\ln K+\frac{s}{h}\max_{i\in\left\{1,\dots,m\right\}} \left(x_i^{(h)}\right)_j\leq (\cdot)_j \right\}, \end{eqnarray*} 
and also 
\begin{eqnarray*}&& \left\{\forall i\in\{1,\dots,m\}\quad g\left(\cdot-x_i^{(s)}\right)\geq 0\right\} \\ &= &\left\{\min_{i\in\left\{1,\dots,m^\frac{s}{h}\right\}} g\left(\cdot-x_i^{(s)}\right)\geq 0\right\} \\ &=& \left\{K-\max_{i\in\left\{1,\dots,m^\frac{s}{h}\right\}} \min_{j\in\{1,\dots,d\}} \exp\left(\left(\cdot-x_i^{(s)}\right)_j\right)\geq  0\right\}\\  &\subset& \left\{K-\min_{i\in\left\{1,\dots,m^\frac{s}{h}\right\}} \min_{j\in\{1,\dots,d\}} \exp\left(\left(\cdot-x_i^{(s)}\right)_j\right)\geq  0\right\}\\  &=& \left\{K- \min_{j\in\{1,\dots,d\}}\min_{i\in\left\{1,\dots,m^\frac{s}{h}\right\}} \exp\left(\left(\cdot-x_i^{(s)}\right)_j\right)\geq  0\right\}\\  &=& \left\{K- \min_{j\in\{1,\dots,d\}}\exp\left(\left(\cdot\right)_j -\max_{i\in\left\{1,\dots,m^\frac{s}{h}\right\}} \left(x_i^{(s)}\right)_j\right) \geq 0\right\} \\ &=& \left\{\ln K \geq  \min_{j\in\{1,\dots,d\}} \left(\cdot\right)_j - \max_{i\in\left\{1,\dots,m^\frac{s}{h}\right\}} \left(x_i^{(s)}\right)_j \right\} \\ &=& \bigcup_{j=1}^d\left\{\ln K+\max_{i\in\left\{1,\dots,m^\frac{s}{h}\right\}} \left(x_i^{(s)}\right)_j\geq  (\cdot)_j \right\}\\ &=& \bigcup_{j=1}^d\left\{\ln K+\frac{s}{h}\max_{i\in\left\{1,\dots,m\right\}} \left(x_i^{(h)}\right)_j\geq (\cdot)_j \right\}.\end{eqnarray*} 
\end{proof}

\begin{cor}\label{barf=min_estimate} If $\bar f=\min_{j\in\{1,\dots,d\}}\exp\left((\cdot)_j\right)$, then for all $s\in I$,
\begin{eqnarray*}&&\left\{P_s(g\vee 0)>0\right\}\\ &=& \bigcup_{j=1}^d \left( \underbrace{ \RR\times\cdots\times\RR }_{j-1} \times \left(-\infty,\ln K+ \frac{s}{h} \max_{i\in\{1,\dots,m\}} \left(x_i^{(h)}\right)_j \right]\times \underbrace{ \RR\times\cdots\times\RR }_{d-j} \right)\end{eqnarray*}
\end{cor}
\begin{proof} Let $s\in I$. We simply remark that
\begin{eqnarray*}&&\left\{P_s(g\vee 0)>0\right\}\\ &=&\left\{\exists i\in\left\{1,\dots,m^\frac{s}{h}\right\}\quad g\left(\cdot-x_i^{(s)}\right)>0\right\} \\ &=&\complement \left\{\forall i\in\left\{1,\dots,m^\frac{s}{h}\right\}\quad g\left(\cdot-x_i^{(s)}\right) \leq 0\right\} \end{eqnarray*} 
\end{proof}

\begin{lem} If $\bar f=\max_{j\in\{1,\dots,d\}}\exp\left((\cdot)_j\right)$, then for all $s\in I$, \begin{eqnarray*}&&\left\{\forall i\in\left\{1,\dots,m^\frac{s}{h}\right\}\quad g\left(\cdot-x_i^{(s)}\right)\leq 0\right\}\\ &\subset&\bigcup_{j=1}^d \left( \underbrace{ \RR\times\cdots\times\RR }_{j-1} \times \left[\ln K+ \frac{s}{h} \min_{i\in\{1,\dots,m\}} \left(x_i^{(h)}\right)_j ,+\infty\right) \times \underbrace{ \RR\times\cdots\times\RR }_{d-j} \right)\end{eqnarray*} as well as \begin{eqnarray*}&& \left\{\forall i\in\left\{1,\dots,m^\frac{s}{h}\right\}\quad g\left(\cdot-x_i^{(s)}\right)\geq 0\right\}\\ &=& \bigotimes_{j=1}^d\left(-\infty,\ln K+ \frac{s}{h} \min_{i\in\{1,\dots,m\}} \left(x_i^{(h)}\right)_j\right]\end{eqnarray*}
\end{lem}
\begin{proof} Let $s\in I$. Then 
\begin{eqnarray*}&& \left\{\forall i\in\{1,\dots,m\}\quad g\left(\cdot-x_i^{(s)}\right)\leq 0\right\} \\ &= &\left\{\max_{i\in\left\{1,\dots,m^\frac{s}{h}\right\}} g\left(\cdot-x_i^{(s)}\right)\leq 0\right\} \\ &=& \left\{K-\min_{i\in\left\{1,\dots,m^\frac{s}{h}\right\} } \max_{j\in\{1,\dots,d\}} \exp\left(\left(\cdot-x_i^{(s)}\right)_j\right)\leq 0\right\}\\  &\subset& \left\{K-\max_{i\in\left\{1,\dots,m^\frac{s}{h}\right\}} \max_{j\in\{1,\dots,d\}} \exp\left(\left(\cdot-x_i^{(s)}\right)_j\right)\leq 0\right\}\\  &=& \left\{K-\max_{j\in\{1,\dots,d\}}\max_{i\in\left\{1,\dots,m^\frac{s}{h}\right\} } \exp\left(\left(\cdot-x_i^{(s)}\right)_j\right)\leq 0\right\}\\  &=& \left\{K- \max_{j\in\{1,\dots,d\}}\exp\left(\left(\cdot\right)_j -\min_{i\in\left\{1,\dots,m^\frac{s}{h}\right\}} \left(x_i^{(s)}\right)_j\right) \leq 0\right\} \\ &=& \left\{\ln K\leq \max_{j\in\{1,\dots,d\}} \left(\cdot\right)_j - \min_{i\in\left\{1,\dots,m^\frac{s}{h}\right\}} \left(x_i^{(s)}\right)_j \right\} \\ &=& \bigcup_{j=1}^d\left\{\ln K+\min_{i\in\left\{1,\dots,m^\frac{s}{h}\right\}} \left(x_i^{(s)}\right)_j\leq (\cdot)_j \right\}\\ &=& \bigcup_{j=1}^d\left\{\ln K+\frac{s}{h}\min_{i\in\left\{1,\dots,m\right\}} \left(x_i^{(h)}\right)_j\leq (\cdot)_j \right\}, \end{eqnarray*} 
and also 
\begin{eqnarray*}&& \left\{\forall i\in\{1,\dots,m\}\quad g\left(\cdot-x_i^{(s)}\right)\geq 0\right\} \\ &= &\left\{\min_{i\in\left\{1,\dots,m^\frac{s}{h}\right\}} g\left(\cdot-x_i^{(s)}\right)\geq 0\right\} \\ &=& \left\{K-\max_{i\in\left\{1,\dots,m^\frac{s}{h}\right\}} \max_{j\in\{1,\dots,d\}} \exp\left(\left(\cdot-x_i^{(s)}\right)_j\right)\geq 0\right\} \\  &=& \left\{K- \max_{j\in\{1,\dots,d\}}\max_{i\in\left\{1,\dots,m^\frac{s}{h}\right\}} \exp\left(\left(\cdot-x_i^{(s)}\right)_j\right)\geq  0\right\}\\  &=& \left\{K- \max_{j\in\{1,\dots,d\}}\exp\left(\left(\cdot\right)_j -\min_{i\in\left\{1,\dots,m^\frac{s}{h}\right\}} \left(x_i^{(s)}\right)_j\right)\geq 0\right\} \\  &=& \left\{K\geq \exp\left(\max_{j\in\{1,\dots,d\}}\left(\left(\cdot\right)_j -\min_{i\in\left\{1,\dots,m^\frac{s}{h}\right\}} \left(x_i^{(s)}\right)_j\right)\right)\right\} \\ &=& \left\{\ln K \geq  \max_{j\in\{1,\dots,d\}} \left(\left(\cdot\right)_j - \min_{i\in\left\{1,\dots,m^\frac{s}{h}\right\}} \left(x_i^{(s)}\right)_j \right)\right\} \\ &=& \bigcap_{j=1}^d\left\{\ln K+\min_{i\in\left\{1,\dots,m^\frac{s}{h}\right\}} \left(x_i^{(s)}\right)_j\geq  (\cdot)_j \right\}\\ &=& \bigcap_{j=1}^d\left\{\ln K+\frac{s}{h}\min_{i\in\left\{1,\dots,m\right\}} \left(x_i^{(h)}\right)_j\geq (\cdot)_j \right\}. \end{eqnarray*} 
\end{proof}

\begin{cor}\label{Psgv0>Psg.max} If $\bar f=\max_{j\in\{1,\dots,d\}}\exp\left((\cdot)_j\right)$, then for all $s\in I$,
\begin{eqnarray*}&&\left\{P_s(g\vee 0)>P_sg\right\} \\ &=& \bigcup_{j=1}^d \left( \underbrace{ \RR\times\cdots\times\RR }_{j-1} \times \left(\ln K+ \frac{s}{h} \min_{i\in\{1,\dots,m\}} \left(x_i^{(h)}\right)_j ,+\infty\right) \times \underbrace{ \RR\times\cdots\times\RR }_{d-j} \right)\end{eqnarray*}
\end{cor}
\begin{proof} Let $s\in I$. We simply remark that
\begin{eqnarray*}&&\left\{P_s(g\vee 0)>P_sg\right\}\\ &=&\left\{\exists i\in\left\{1,\dots,m^\frac{s}{h}\right\}\quad g\left(\cdot-x_i^{(s)}\right)<0\right\} \\ &=&\complement \left\{\forall i\in\left\{1,\dots,m^\frac{s}{h}\right\}\quad g\left(\cdot-x_i^{(s)}\right)\geq 0\right\} \end{eqnarray*} 
\end{proof}

These estimates lead to the following Corollary that will enable us -- under the assumption of $$\max_{i\in\{1,\dots,m\}}x_{i}^{(h)}=0$$ (in order to be entitled to apply eg Lemma \ref{from_k=1_to_kinN}) -- to prove an $L^1$-convergence estimate (on a particular subset of $\RR^d$) for $\left(B_{T\cdot 2^{-n}}f\right)_{n\in\NN}$ for any measurable $f\geq g\vee 0$.

\begin{cor} Suppose $\bar f=\max_{j\in\{1,\dots,d\}} \exp\left((\cdot)_j\right)$ and consider any compact set $B\subset \RR^d$. Then \begin{eqnarray*}&&\lambda^d\left[\left\{e^{rs}g>P_s(g\vee 0)>P_sg\right\}\cap B\right] \\&\leq& s\cdot R^{d-1}\frac{1}{h} \sum_{j=1}^d \left(-\min_{i\in\{1,\dots,m\}} \left(x_i^{(h)}\right)_j\vee 0\right)\end{eqnarray*} for all $s\in I$.
\end{cor}
\begin{proof} Let $s\in I$. Since $$\left\{e^{rs}g>P_s(g\vee 0)>P_sg\right\}\subseteq \left\{P_s(g\vee 0)>P_sg\right\}\cap \{g>0\}$$ by the monotonicity of $P_s$, we only have to observe that \begin{eqnarray*} \{g>0\} &=& \left\{K>\max_{j\in\{1,\dots,d\}} \exp\left((\cdot)_j\right)\right\} \\ &=&\bigcap_{j=1}^d\left\{K>\exp\left((\cdot)_j\right)\right\}\\ &=&\bigotimes_{j=1}^d \left(-\infty,\ln K\right)\end{eqnarray*} to arrive -- after taking advantage of the preceding Corollary \ref{Psgv0>Psg.max} -- at \begin{eqnarray*}&& \left\{e^{rs}g>P_s(g\vee 0)>P_sg\right\}\\ &\subseteq & \bigcup_{j=1}^d \left( \begin{array}{c}\underbrace{ \left(-\infty,\ln K\right)\times\cdots\times\left(-\infty,\ln K\right) }_{j-1} \\\times \left(\ln K+ \frac{s}{h} \min_{i\in\{1,\dots,m\}} \left(x_i^{(h)}\right)_j ,\ln K\right)\\ \times \underbrace{ \left(-\infty,\ln K\right)\times\cdots\times\left(-\infty,\ln K\right) }_{d-j}\end{array} \right) \\ &=& \bigcup_{j=1}^d \left( \underbrace{ \RR_{<0}\times\cdots\times\RR_{<0} }_{j-1} \times \left(\frac{1}{h} \min_{i\in\{1,\dots,m\}} \left(x_i^{(h)}\right)_j,0 \right)\cdot s\times \underbrace{ \RR_{<0}\times\cdots\times\RR_{<0} }_{d-j} \right) \\ && +(\ln K)_{j=1}^d.\end{eqnarray*}

However, by our assumption that $B$ be compact, there is some $R>0$ such that $B-(\ln K)_{j=1}^d\subset\left[-{R},R\right]^d$. Thus \begin{eqnarray*}&&\left\{e^{rs}g>P_s(g\vee 0)>P_sg\right\}\cap B \\ &\subseteq& \bigcup_{j=1}^d \left( [-R,0)^{j-1} \times \left(\frac{1}{h} \min_{i\in\{1,\dots,m\}} \left(x_i^{(h)}\right)_j,0 \right)\cdot s\times \RR^{d-j} \right) \\ && +(\ln K)_{j=1}^d,\end{eqnarray*} and from this inclusion we may deduce the estimate given in the Lemma.
\end{proof}

The inequality we have just derived implies that the \linebreak $\lambda^d\left[E^h \cap \bigcap_{\ell}\left(B+x_\ell^{\left(s\right)}\right)\cap\cdot\right]$-volume of the set occurring in Lemma \ref{estimateonE} is of order $s$ for any compact $B$ and for $\bar f=\max_{j\in\{1,\dots,d\}} \exp\left((\cdot)_j\right)$. Hence again by Lemma \ref{estimateonE} (which is applicable because of Lemma \ref{barf=max_gamma1}) we obtain that the difference $B_{s/2}f-B_sf$ is of order $s^2$ (this time, however in the $L^1\left(E^h \cap \bigcap_{\ell}\left(B+x_\ell^{(s)}\right)\right)$-norm). This estimate on the norm of $B_{s/2}f-B_sf$ leads, via Lemmas \ref{from_k,k+1_to_M,N} and \ref{from_k=1_to_kinN} to the result that the analogon of the difference in Theorem \ref{cubatureBermudanconv} is of order $s^2\cdot 2^{-M}$, too:

\begin{Th} \label{cubatureBermudanconv_barf=min} Suppose $\bar f=\max_{j\in\{1,\dots,d\}}\exp\left((\cdot)_j\right)$ and consider a compact set $B$. Assume that $$\max_{i\in\{1,\dots,m\}}x_{i}^{(h)}=0.$$ Under these assumptions there exists a real number $D>0$ such that for all $N>M\in\NN$, $s\in(0,T)\cap \left(2^N\cdot I\right)$ and $f\geq g\vee 0$, one has 
\begin{eqnarray*}&&\left\|\left(B_{s\cdot 2^{-N}}\right)^{\circ\left(2^{N}\right)}f - \left(B_{s\cdot 2^{-M}}\right)^{\circ\left(2^M\right)}f\right\|_{L^1\left(E^h \cap \bigcap_{\ell}\left(B+x_\ell^{\left(s\right)}\right)\right)} \\ &\leq& D\cdot {s}^2\cdot{2}^{-M}\left(1- 2^{-(N-M-1)}\right) \\ &\leq& D\cdot {s}^2\cdot{2}^{-M}. \end{eqnarray*}
\end{Th}

\begin{rem} This $L^1$-convergence result has some (however, because of our assumption $\max_{i\in\{1,\dots,m\}}x_{i}^{(h)}=0$, fairly limited) practical interest, as in practice quite frequently the exact start price of the (multiple) asset on which an option is issued, is unknown. Instead, one will have the logarithmic start price vector $x\in\RR^d$ a short time $\delta>0$ before the actual option contract becomes valid. Now, asuming that $\PP^x_{X_\delta}$ has a continuous density $\frac{\PP^x_{X_\delta}}{\lambda^d}$, this function $\frac{\PP^x_{X_\delta}}{\lambda^d}$ will be bounded on $E^h \cap \bigcap_{\ell}\left(B+x_\ell^{\left(s\right)}\right)$ by some constant $$C:=\sup_{E^h \cap \bigcap_{\ell}\left(B+x_\ell^{\left(s\right)}\right)} \frac{\PP^x_{X_\delta}}{\lambda^d}<+\infty.$$ One will therefore have for all $f\geq g\vee 0$, $s\in(0,T]\cap I$ and $N>M\in\NN$, \begin{eqnarray*}&&\EE^x\left[\begin{array}{c}\left|\left(B_{s\cdot 2^{-N}}\right)^{\circ\left(2^{N}\right)}f - \left(B_{s\cdot 2^{-M}}\right)^{\circ\left(2^M\right)}f\right|\left(X_\delta\right),\\ X_\delta\in E^h \cap \bigcap_{\ell}\left(B+x_\ell^{\left(s\right)}\right)\end{array}\right]\\&\leq&
\left\|\left(B_{s\cdot 2^{-N}}\right)^{\circ\left(2^{N}\right)}f - \left(B_{s\cdot 2^{-M}}\right)^{\circ\left(2^M\right)}f\right\|_{L^1\left(\PP^x_{X_\delta}\left[E^h \cap \bigcap_{\ell}\left(B+x_\ell^{\left(s\right)}\right)\cap\cdot\right]\right)} \\ &\leq& C\cdot\left\|\left(B_{s\cdot 2^{-N}}\right)^{\circ\left(2^{N}\right)}f - \left(B_{s\cdot 2^{-M}}\right)^{\circ\left(2^M\right)}f\right\|_{L^1\left(\lambda^d\left[E^h \cap \bigcap_{\ell}\left(B+x_\ell^{\left(s\right)}\right)\cap\cdot\right]\right)} \\ &\leq& C\cdot D\cdot {s}^2\cdot{2}^{-M}. \end{eqnarray*}

\end{rem}

\part{Numerical analysis of cubature-based American pricing}

\chapter{Motivating Bermudan pricing based on cubature}

\label{treerecursion}

\section{The general setting}

This Chapter has been designed to elaborate the idea of Bermudan option pricing via cubature and to put it into a context of other Bermudan option prcing algorithms. Given its informal character, it can also be perceived as another introductory chapter.

Consider a basket of $d$ assets. A {\em $d$-dimensional Bermudan option} is an option that can be exercised at a discrete set of exercise times, yielding payoff $g(x)$ if $x\in\RR^d$ is the vector of logarithmic prices (of stocks in the $d$-dimensional basket) at that time, for some $g:\RR^d\rightarrow\RR$ which will be called the {\em payoff function}. In the case of a one-dimensional put option with strike price $K$ for example, one would have $g=(K-\exp)^+=(K-\exp)\vee 0$. In case of a call on a stock index, $g:\RR^d\rightarrow \RR$ would be the positive part of the difference between a weighted sum of exponential functions of the coordinate entries and the strike price. Unless specified otherwise, we will from now on assume the exercise times to be equidistant with an exercise mesh size $h>0$. Adopting the notation of Chapter \ref{formalintro}, this is to say $J=h\NN_0$.

We regard such a Bermudan option as a binary tree of European options. This means that at each exercise time one has to decide whether it is more rewarding to keep the option or to exercise it -- in other words, whether the payoff at that exercise time is less than the value of the (European) option to exercise at the next exercise time. A recursive algorithm is thus obtained. Note however that this binary tree of European options has continuum many nodes at each level, one for each price vector at the subsequent exercise time. 

Let us now describe this recursion in detail. Suppose the option is non-perpetual, i.e. it has a maturity time $T\in(0,\infty)$, and assume furthermore that $T=hN$ for some $N\in\NN$. Then the Bermudan option price for a start price vector $x_0\in\RR^d$ will be $V_0(x_0)$ where the $V_k, k\in\{0,\dots,N\}$ are computed according to the following backward recursion (where for simplicity we assume the logarithmic discount rate to be a constant $r>0$ and $(Y_t)_{t}$ denotes the ${\RR_+}^d$-valued process of vectors comprised of the asset prices in the basket): $$V_N=g,$$ \begin{equation}\label{recurs}\forall k\in\{1,\dots,N\}\quad \forall x\in\RR^d\quad V_{k-1}(x)=\max\left\{g(x), \EE^x e^{-rh}V_k\left(Y_{h\cdot k}\right)\right\}.\end{equation} 

Recursion formulae of this kind for the pricing of Bermudan options are fairly standard and can be found for example in textbooks such as Hull's \cite{Hull} or Wilmott, Howison and Dewynne's \cite{WHD}. To use these recursion formulae practically, one needs a way to summarise or approximate the state $V_k$ in a way that permits the equivalent summarisation or approximation for $V_{k-1}$ (this for all positive integer $k$). One method of achieving this in a one-dimensional setting is the application of Fourier-Hermite expansions to the functions $V_k$, as studied in the paper by Chiarella, el-Hassan and Kucera \cite{CHK}. Our goal is to develop their approach; in particular we will employ cubature formulae for symmetric measures. These methods of approximating integrals by weighted (finite) averages can be computationally efficient, and with increasing dimension may be superior to other approaches. Victoir \cite{V} introduced a vital improvement by constructing sequences that scale well. This route to high-dimensional Bermudan and American option pricing was proposed for the first time by my supervisor \cite{L}.

\section{Application to the Black-Scholes model}

Let us in this section work within the situation of the multi-dimensional Black-Scholes model, that is to say that the logarithmic price processes of the $d$ assets in the basket are independent Brownian motions with drift. Let us assume the volatilities of the assets to be constants $\sigma_i$ for $i\in\{1,\dots,d\}$.

We set $$\mu_i:=r-\frac{\sigma_i^2}{2}$$ for $i\in \{1,\dots,d\}$ and denote by $\nu_{\alpha, \sigma^2}$ for $\alpha\in\RR^d$, $\sigma>0$ the $d$-dimensional Gaussian probability measure of variance $\sigma^2$ centered at $\alpha$. We will assume that -- possibly after an appropriate change of the time scale and the discount rate $r$ (by a linear transformation from the left) -- we have  $\sigma_i=1/\sqrt{d}$ for all $i\in\{1,\dots,d\}$. $\mu$ will change according to its definition.

Then It\^o's Lemma implies that the logarithmic price process in the $i$-th coordinate is -- with repect to the risk-neutral measure -- just a Brownian motion with drift $\mu_i$ and volatility $1/\sqrt{d}$, thus the process of logarithmic prices of the assets in the basket is a Brownian motion with drift $\mu$ and volatility $1$. Therefore our recursion fomula (\ref{recurs}) becomes \begin{eqnarray*}&&  \forall k\in\NN\cap[0,T]\quad \forall\xi\in\RR^d \\ V_{k-1}(\xi)&=&\max\left\{e^{-rh}\int_{\RR^d}  V_{k}d\nu_{\xi +\mu h, h}, g(\xi)\right\}\\ &=&\max\left\{e^{-rh}\int_{\RR^d}V_{k}(\xi+\mu h+h^{1/2}\cdot)d\nu_{0,1}, g(\xi)\right\} \end{eqnarray*}  Now, if the points $\xi_j\in\RR^d$ with respective weights $\alpha_j>0$, $j\in\{1,\dots,m\}$, determine cubature formulae for the standard Gaussian measure $\nu_{0,1}$, we can approximate the previous recursion by the following formula: \begin{eqnarray}\label{recursion} \nonumber && \forall k\in\NN\cap[0,T]\forall\eta\in\RR^d \\ \tilde V_{k-1}(\eta)&=&\max\left\{e^{-rh}\sum_{j=1}^m\alpha_j \tilde V_{k}(\eta+\mu h+h^{1/2}\xi_j), g(\eta)\right\}. \end{eqnarray}

\section{Exploiting combinatorial aspects of Gaussian cubature}

Thanks to the work of Nicolas Victoir (which has later been extended by Christian Litterer), there are ``cubature formulae with few points'' \cite{V} for the integration of polynomials with respect to the standard Gaussian measure up to a certain degree. Although ``asymmetric'' \cite{V}, their shape is quite regular and uniform. Since the recursion following the previous recursion formula amounts to the evaluation of payoff functions at (modified) sums of these cubature points and we therefore desire recombination of these sums, this will turn out to be a computationally palpable advantage.

The commutativity of $(\RR^d,+)$ and the equidistance of the exercise times already enable us to perform a geometric argument based on the regular and uniform shape of the cubature points, which results in

\begin{Th} \label{treeispolynomial}Let $d=3k-2$ for some $k\in\NN$. The recursion according to (\ref{recursion}), using the cubature formula for the integration of degree 5 polynomials with respect to a standard Gaussian measure from Victoir's example \cite[5.1.1]{V}, is polynomial in $\frac{1}{h}$.
\end{Th}
\begin{proof} The cubature points of the cubature formulae referred to in the Theorem form a finite subset of $\sqrt[4]{3}\{0,\pm 1\}^d$. Sums of length $\frac{1}{h}$ (provided this fraction is an integer) of the cubature points are therefore always elements of $\sqrt[4]{3}\left(\ZZ^d\cap\left\{\left|\cdot\right|\leq \frac{1}{h}\right\}\right)$ (and this set has only $\left(2\frac{1}{h}\right)^d$ elements), and the points used in the recusion formula stated above are comprised of a subset of $h^{1/2}\cdot\sqrt[4]{3}\left(\ZZ^d\cap\left\{\left|\cdot\right|\leq \frac{1}{h}\right\}\right)+ h \cdot\mu \left\{0,\dots,\frac{1}{h}\right\}+\xi_0$.
\end{proof}

However, this is not the only recombination that can be accomplished in the case where $d=3k-2$ :

\begin{rem} Let us look at the tree obtained from starting at some point $\xi_0\in\RR^d$ and then at each node letting exactly $\left|\{x_0\}\cup\sqrt[4]{3}\cdot\cG_3X_1\right|$ branches leave (where $\{x_0\}\cup\sqrt[4]{3}\cdot\cG_3X_1$ in Victoir's notation is the set of cubature points he uses in the example \cite[5.1.1]{V} we are referring to), exactly one branch for each element of the set $\xi_0+\left(\{x_0\}\cup\sqrt[4]{3}\cdot\cG_3X_1\right)$.

If we intend to find and eliminate the branches of the tree that are computed ``wastefully'', it is reasonable to divide the sums (of length $\frac{1}{h}$) of the cubature points by $\sqrt[4]{3}$ and consider them coordinate-wise modulo $2$. Then one is looking at elements of the vector space $\left(\ZZ/2\ZZ\right)^d$. For the sake of simplicity, let $d=7$, that is $k=3$ in the notation of the previous Theorem \ref{treeispolynomial}. The coordinate-wise projection of the $\frac{1}{\sqrt[4]{3}}$-multiple of our set of cubature points $\{x_0\}\cup \sqrt[4]{3}\cdot\cG_3X_1$, where $$X_1:=\left\{x_{1,1},\dots, x_{1,7}\right\}$$ and $$x_0:={ 0}, \quad \cG_3=\left(\{\pm 1 \}^3,\ast\right)$$ (in the notation of \cite[5.1.1]{V}), into the vector space $\left(\ZZ/2\ZZ\right)^7$ now contains only eight points (instead of $57$ as before). 

Thus, using basic linear algebra in a $7$-dimensional $\ZZ /2\ZZ$-vector space, we are easily able to classify the non-trivial zero representations from elements of the projected cubature points. 

Perceiving $X_1$ as a $7$-element subset of $\left(\ZZ /2\ZZ\right)^7$, we see that $(x)_{x\in X_1}$ is an invertible $\left(\ZZ /2\ZZ\right)^{7\times 7}$-matrix. Therefore we cannot expect any recombination from representations of zero by nontrivial linear combinations of elements of $X_1\subset \left(\ZZ /2\ZZ\right)^7$. Moreover, the fact that $A:=(x)_{x\in X_1}$ is invertible, shows that $x_0 ={ 0}$ can only be written trivially as a sum of elements of $X_1$. Hence we have shown that we exploit symmetries optimally if we use: (i) the commutativity of $(\RR^d, +)$; (ii) the obvious symmetries due to the construction of the cubature formulae by means of the action of a reflection group on certain points; (iii) the fact that addition of $x_0$ does nothing at all.
\end{rem}

\chapter{Numerical results}

In this Chapter we shall present some numerical results. We have decided to choose a $7$-dimensional example, since (1) most previous research has stopped short of numerically tackling American options on baskets with more than $5$ assets, (2) it is the smallest dimension $d\geq 5$ in which Victoir's cubature formulae for the normal Gaussian measure of dimensions $d=3k-2$ (where $k\in\NN_0$) \cite[Example 5.1.1]{V} hold. 

We shall assume that the basket $X$ as a logarithmic price process follows the Black-Scholes model for independent assets with discount rate $r>0$ and volatilities $\sigma_i>0$, $i\in\{1,\dots,d\}$, that is $$\forall t\geq 0 \quad X_t=\left(\left(X_0\right)_i+\sigma_i\cdot (B_t)_i + \left(r-\frac{1}{2}{\sigma_i}^2\right) t\right)_{i=1}^d$$ (where $B$ is the $d$-dimenional Wiener process). Given a payoff function, a strike price, a maturity $T>0$, and logarithmic start price vector $x\in\RR^d$ we shall vary the exercise mesh size $h>0$ (say $h\in\{h_0,\dots,h_n\}$) and compute approximate Bermudan prices $U^{h\NN_0}(T)(x)$. Then we will extrapolate the function $h\mapsto U^{h\NN_0}(T)(x)$ to $h=0$ by assuming assume that $h\mapsto U^{h\NN_0}(T)(x)$ is a polynomial of degree $n$ in $h^\alpha$ for a given $\alpha>0$ that finally shall be varied as well. 

Unfortunately, it is difficult to find data on American option prices for dimension $d=7$. However, one can of course use our algorithm sub-optimally for $d=5$ through letting the payoff function only depend on the first five coordinates. Then a comparison with the numerical value computed by the 50S algorithm (as stated in Rogers \cite{R}) sadly yields a 3.64 \% difference after 9.87 seconds of computations on a 1.4 GHz Personal Computer (whereas 50S needed 14 seconds on a 600 MHz PC). 

More extensive numerical experiments (on computers of better performance) may find, however, that a cubature-based algorithm is superior to a Monte-Carlo routine when higher dimensions than $d=7$ are considered. On a different note, recall that in practice for the vast majoriy of derivative options, pricing algorithms are only used as part of hedging programs -- and with hedging, the accuracy of the prices computed is of lesser importance than the processor time the algorithms actually requires.

We conclude this Chapter by stating some Bermudan and American option prices computed through our cubature-based algorithm.

For a min-put on a basket of seven independent assets with discount rate $r=0.06$, maturity at time $0.5$ and strike price $K=100$, one will get the following numerical results. (Here, extrapolation I is the extrapolation of $U^{h\NN_0}(T)(x)$ from $h\in\left\{\frac{T}{1},\frac{T}{2},\frac{T}{3}\right\}$ to $h=0$ with scaling exponent $\alpha=1.0$, and extrapolation II is the corresponding extrapolation with scaling exponent $\alpha=\frac{1}{2}$. The amount of time elapsed during each computation is given in seconds.) 

If all volatilities $\sigma_1,\dots,\sigma_d$ are equal to $0.4$, then 

\begin{tabular}{*{5}{c}} Start prices & \begin{tabular}{c}Bermudan \\ ($\frac{T}{h}=3$)\end{tabular} & Extrapolation I & Extrapolation II & Time\\&&&&\\

$80,\dots,80$ & $42.7981$ & $44.4427$ & $45.7474$ & $12.35$ \\
$90,\dots,90$ & $36.0074$ & $37.8231$ & $39.2082$ & $11.52$\\
$100,\dots,100$ & $29.2172$ & $31.2058$ & $32.6757$ & $11.81$ \\
$110,\dots,110$ & $23.2196$ & $24.2340$ & $25.3323$ & $12.07$ \\
$120,\dots,120$ & $17.3103$ & $17.2211$ & $17.8944$ & $11.78$ \\
\begin{tabular}{c}$80,90,90,100$,\\$110,110,120$\end{tabular} & $34.1646$ & $34.1625$ & $34.2569$ & $12.89$ 
\end{tabular}

In case the volatilities $\sigma_1,\dots,\sigma_d$ are all equal to $0.6$, then

\begin{tabular}{*{5}{c}} Start prices & \begin{tabular}{c}Bermudan \\ ($\frac{T}{h}=3$)\end{tabular} & Extrapolation I & Extrapolation II & Time\\&&&&\\

$80,\dots,80$ & $53.0963$ & $55.1235$ & $56.2195$ & $11.94$\\
$90,\dots,90$ & $47.5948$ & $49.8475$ & $51.0137$ & $12.49$\\
$100,\dots,100$ & $42.0951$ & $44.5798$ & $45.8316$ & $12.18$ \\
$110,\dots,110$ & $36.9143$ & $39.1641$ & $40.3082$ & $11.90$ \\
$120,\dots,120$ & $32.2508$ & $33.5128$ & $34.2440$ & $11.60$\\
\begin{tabular}{c}$80,90,90,100$,\\$110,110,120$\end{tabular} & $45.9639$ & $45.8861$ & $ 44.4914$ & $11.56$
\end{tabular}

Finally, if $\sigma_1=0.3$, $\sigma_2=0.4$, $\sigma_3=0.5$, $\sigma_4=0.6$, $\sigma_5=0.7$, $\sigma_6=0.8$, $\sigma_7=0.9$, one has the following figures:

\begin{tabular}{*{5}{c}} Start prices & \begin{tabular}{c}Bermudan \\ ($\frac{T}{h}=3$)\end{tabular} & Extrapolation I & Extrapolation II & Time\\&&&&\\

$80,\dots,80$ & $56.6740$ & $57.1258$ & $55.2092$ & $11.50$\\
$90,\dots,90$ & $51.6227$ & $52.1136$ & $49.9160$ & $13.04$ \\
$100,\dots,100$ & $46.5720$ & $47.1037$ & $44.6290$ & $11.51$ \\
$110,\dots,110$ & $41.5507$ & $42.1357$ & $39.3828$ & $11.81$ \\
$120,\dots,120$ & $36.7872$ & $37.2093$ & $33.9543$ & $11.52$ \\
\begin{tabular}{c}$80,90,90,100$,\\$110,110,120$\end{tabular} & $45.1416$ & $44.8038$ & $42.5739$ & $11.51$

\end{tabular}

These data suggest that the optimal scaling exponent $\alpha$ for the extrapolation from Bermudan to American min-put prices will have to depend on both the volatility vector $\left(\sigma_1,\dots,\sigma_d\right)$ and the vector of (logarithmic) start prices $\left(x_1,\dots,x_d\right)$.

Our second example concerns itself with the pricing of Bermudan and American put options on the arithmetic average of a basket of independent assets.

If all volatilities $\sigma_1,\dots,\sigma_d$ are equal to $0.4$, then 

\begin{tabular}{*{5}{c}} Start prices & \begin{tabular}{c}Bermudan \\ ($\frac{T}{h}=3$)\end{tabular} & Extrapolation I & Extrapolation II & Time\\&&&&\\

$90,\dots,90$ & $10.0000$ & $10.0000$ & $10.0000$ & $10.52$\\
$100,\dots,100$ & $3.15446$ & $3.13173$ & $2.44775$ & $10.68$ \\
$110,\dots,110$ & $0.715959$ & $0.681342$ & $0.793472$ & $10.24$ \\
$120,\dots,120$ & $0.113121$ & $0.149574$ & $0.214042$ & $10.86$ \\
\begin{tabular}{c}$80,90,90,100$,\\$110,110,120$\end{tabular} & $3.17743$ & $3.22211$ & $2.70252$ & $10.44$ 
\end{tabular}

In case the volatilities $\sigma_1,\dots,\sigma_d$ are all equal to $0.6$, then

\begin{tabular}{*{5}{c}} Start prices & \begin{tabular}{c}Bermudan \\ ($\frac{T}{h}=3$)\end{tabular} & Extrapolation I & Extrapolation II & Time\\&&&&\\

$80,\dots,80$ & $20.0000$ & $20.0000$ & $20.0000$ & $10.18$\\
$90,\dots,90$ & $10.8834$ & $10.8347$ & $10.4006$ & $10.62$\\
$100,\dots,100$ & $5.34809$ & $5.16546$ & $4.54056$ & $10.30$ \\
$110,\dots,110$ & $2.22539$ & $2.35572$ & $2.18241$ & $11.09$ \\
$120,\dots,120$ & $0.881495$ & $0.758719$ & $0.791531$ & $10.85$\\
\begin{tabular}{c}$80,90,90,100$,\\$110,110,120$\end{tabular} & $5.38701$ & $5.21894$ & $4.58038$ & $10.67$
\end{tabular}

And if $\sigma_1=0.3$, $\sigma_2=0.4$, $\sigma_3=0.5$, $\sigma_4=0.6$, $\sigma_5=0.7$, $\sigma_6=0.8$, $\sigma_7=0.9$, we obtain the following results:

\begin{tabular}{*{5}{c}} Start prices & \begin{tabular}{c}Bermudan \\ ($\frac{T}{h}=3$)\end{tabular} & Extrapolation I & Extrapolation II & Time\\&&&&\\

$80,\dots,80$ & $20.0000$ & $20.0000$ & $20.0000$ & $11.32$\\
$90,\dots,90$ & $11.3241$ & $11.4565$ & $11.4368$ & $10.55$ \\
$100,\dots,100$ & $5.69733$ & $5.67916$ & $5.47655$ & $10.98$ \\
$110,\dots,110$ & $2.44459$ & $2.52824$ & $2.35541$ & $10.52$ \\
$120,\dots,120$ & $0.937346$ & $0.915058$ & $1.10870$ & $11.04$ \\
\begin{tabular}{c}$80,90,90,100$,\\$110,110,120$\end{tabular} & $6.24348$ & $6.26401$ & $6.26811$ & $10.33$

\end{tabular}

The first line of each of these sets of figures of course simply means that immediate exercise is optimal if the start price of each asset is at $80$ or below (and, in case $\sigma_i=0.4$ for all $i\in\{1,\dots,7\}$, even if each asset start price is at $90$).

Again it is apparent from these numerical data that the scaling exponent needed for the extrapolation from Bermudan to American prices put-on-the-average option prices has to be varied with the vector of asset start prices and possibly the volatilities of the underlying assets (otherwise the American price computed by extrapolation would be at times very significantly below the approximate price of a Bermudan option on the same basket and with the same payoff function).

\part{High-dimensional approximate $\Delta$-hedging}

\chapter{Hedging options on multiple assets -- a suggestion for further research}

\section{Theoretical suggestions}

Up to this point, the subject of our investigation has been the {\em pricing} of high-dmensional American and Bermudan options. In practice, there is at least as much (if not even significantly more) interest in the heding of such options as in finding out their price -- the latter task often being simply left to the markets. To this extent, any pricing algorithm gains much of its practical interest merely from being employable as a subroutine of a hedging algorithm.

The canonical way of hedging -- that is replicating a portfolio, ideally without risk -- that does not need to introduce utility functions for portfolios which  sometimes may not be that easy to justify themselves is $\Delta$-hedging. Unfortunately, however, there is no straightforward multi-dimensional generalisation of $\Delta$-hedging in the discrete binomial model (in the sense of eg Hull \cite{Hull} or Wilmott, Howison, Dewynne \cite{WHD}). For, if the price processes of all of the assets in a portfolio of $d$ different types of shares each follow the binomial model, then at each time step $n$ where the vector of current asset prices equals $x_n\in\RR^d$ there are $2^d$ possible states of the market that may be encountered at the next time step (given by a set of the form $$\left\{\left(\left(x_n\right)_i\cdot y_i\right)_{i\in\{1,\dots,d\}}\in\RR^d \ : \ \forall i\in\{1,\dots,d\}  \quad y_i\in\left\{\alpha_i,\beta_i\right\}\right\},$$ as each asset $i\in\{1,\dots,d\}$ is assumed to move either by a factor $\alpha_i>0$ or by a factor $\beta_i>0$ where without loss of generality one may assume $\alpha_i\neq\beta_i$ for all $i$) compared with only $d+1$ elements in the portfolio (including the bond). Assuming translation-invariance of the Markov chain, we introduce the notation $$p\left(z_1,\dots,z_d\right):=\PP\left[\left.\left\{X_{n+1}=\left(z_i\cdot x_i\right)_{i=1}^d\right\}\right| \left\{X_n=x\right\} \right]$$

Then the volatility $\sigma_i$ of the $i$-th asset is defined to be the square root of the variance of the one-dimensional random walk with steps $\alpha_i, \beta_i$ and transition probabilities $$p_i:=\sum_{z_1,\dots,z_{i-1},z_{i+1},\dots z_d}p\left(z_1,\dots,z_{i-1},\alpha_i,z_{i+1},\dots z_d\right)$$ and $1-p_i$ respectively, on the set  $\alpha_i\NN_0 +\beta_i\NN_0$. This is to say, $$\sigma_i:=\sqrt{{\alpha_i}^2\cdot p_i+{\beta_i}^2\cdot \left(1- p_i\right)-\left(\alpha_i\cdot p_i+\beta_i\cdot \left(1- p_i\right)\right)^2}.$$

We shall define for each such vector $y$ in the set $$Y:=\left\{z\in\RR^d\ : \forall i\in\{1,\dots,d\}  \quad z_i\in\left\{\alpha_i,\beta_i\right\}\right\}$$ the {\em overall absolute correlation} by $$\rho(y):=\prod_{i=1}^d \frac{\left|y_i- \left(\alpha_i\cdot p_i+\beta_i\cdot \left(1- p_i\right)\right)\right|}{\prod_{i=1}^d{\sigma_i}}.$$

One can now think of various {\em approximate} $\Delta$-hedging algorithms -- previsible transaction policies that whilst being unable to eliminate the $\Delta$ altogether, reduce it significantly. At least three classes of such algorithms come to one's mind:

\begin{enumerate}

\item \label{subsetchange} At each time step $\Delta$-hedging of proper subsets of the portfolio, possibly changing the subset with time.
\item \label{correlationhedge} Removing $2^d-(d+1)$ of the elements of $Y$ (thereby making the market model a ``$d+1$-nomial'' one) via a correlation analysis (cf Section \ref{correl}).
\item \label{cubaturehedge} The use of cubature formulae to achieve this elimination.

\end{enumerate}

A natural method of comparing these hedging algorithms will be to look at the $\ell^p$-norms of the resulting sequences of $\Delta$'s. The algortihms of \ref{subsetchange}. and \ref{cubaturehedge}. are straightforward modifications of standard $\Delta$-hedging algorithms for $d$-component portfolios in market models that only allow for $d+1$ possible states of the market at the respective subsequent time step. 

We will therefore dedicate the rest of this short Chapter to the algorithm suggested in \ref{correlationhedge}.

\section{A $\Delta$-hedging algorithm based on a correlation analysis}
\label{correl}

Following suggestion \ref{correlationhedge}. of the preceding paragraph, we will now propose an algorithm that constructs a new set $Y'$ with cardinality $d+1$ from $Y$ (in the notation of the previous Section). It will then be possible to apply $\Delta$-hedging to the Markov chain market model that is given by \begin{eqnarray*}\forall z\in Y'\forall n\in\NN_0 \forall x\in\RR^d &&\\ P\left[\left.\left\{X_{n+1}=\left(z_i\cdot x_i\right)_{i=1}^d\right\}\right| \left\{X_n=x\right\} \right] &=& \frac{p\left(z_1,\dots,z_d\right)}{\sum_{\bar z\in Y'} p\left(\bar z_1,\dots, \bar z_d\right)}.\end{eqnarray*} 

In order to describe the algorithm that produces $Y'$ from $Y$, two cases according to the size of $d$ have to be distinguished. 

{\em Case I}: $d\leq 3$.

In this situation, $$ \left| Y\setminus Y'\right| =2^d-(d+1)\leq d+1= \left| Y'\right|$$ and one will determine the set of those $2^d-(d+1)$ vectors in the set $Y\subseteq \RR^d$ which will be removed, that is the elements of $Y\setminus Y'$ (as opposed to finding the elements of $Y'$ themselves). 

The set $Y\setminus Y'$ will comprise exactly the $d+1$ elements $y$ of $Y$ with the smallest overall absolute correlation $\rho(y)$. 

Note that {\em a priori} there can be $y\neq z\in Y$ such that $\rho(y)=\rho(z)$. Therefore, in order to get the procedure of constructing $Y\setminus Y'$ from $Y$ well-defined, it is necessary to first define a well-ordering $\prec$ on $Y$ and to define a linear order $\leq_\rho$ on $Y$ by $$\forall y,z\in Y\quad y\leq_\rho z:\Leftrightarrow \rho(y)\leq \rho(z).$$ Then, the product $\left(\leq_\rho\times\prec\right)$ will be a well-ordering and we will define $Y\setminus Y'$ to be the set of the $\left(\leq_\rho\times\prec\right)$-smallest $2^d-(d+1)$ elements of $Y$.

{\em Case II}: $d>3$.

In that case $$ \left| Y\setminus Y'\right|=2^{d}-(d+1)>d+1= \left| Y'\right|.$$ For this reason it is faster to single out the elements of $Y'$ directly, rather than determining the elements of its complement $Y\setminus Y'$ first.

Using the well-ordering defined above, one will thus choose the $\left(\leq_\rho\times\prec\right)$-greatest $d+1$ elements of the finite set $Y$ (the elements of $Y$ with the largest overall absolute correlation $\rho$, that is).

\part{Appendix}

\appendix

\chapter{Re-formulation of the perpetual Bermudan pricing problem in $L^1$ and $L^2$}

\section{Non-applicability of the $L^2(\RR^d)$ Spectral Theorem}

Consider a $d$-dimensional L\'evy  basket $X$ with associated family of risk-neutral probability measures $\PP^\cdot$ and discount rate $r>0$.

Fixing $h>0$ and defining $$P:=\pi_{L^2\left(\RR^d\setminus G\right)},\quad A:=A^h:= \II - 
e^{-rh}\PP_{X_0-X_h}\ast\cdot=\left(\delta_0 - e^{-rh}\PP_{X_0-X_h}\right)\ast\cdot,$$ we can rewrite the result of Lemma \ref{convolveproject} as follows:
\begin{equation}\label{1.1} PA\left(V_G^h-g_1\right)=-PAg_1\end{equation} where we assume that $g$ has a square-integrable extension from $G$ to the whole of $\RR^d$; given this assumption, the $g_1\in L^2(\RR^d)$ of the previous identity can be any such extension.

We will suppress the superscript of $A$ for the rest of this paragraph.

Also, without loss of generality, we will assume in this Chapter that the components of the basket $X$ when following the Black-Scholes model all have volatility $1$.

\begin{lem} Let $X$ be a L\'evy basket with associated family of risk-neutral probability measures $\PP^\cdot$ and discount rate $r>0$. Then $A$ and $PA\restriction L^2(\RR^d\setminus G)$ are invertible. Furthermore, the $L^2$ norm of $A$ is bounded by $\left(1+e^{-rh}\right)^\frac{1}{2}$, if $X=\mu\cdot+B$ (thus $\mu=\left(r-\frac{1}{2}\right)_{i=1}^d$) where $B$ is a standard Brownian motion. Moreover, $A$ is a contraction if $\mu=0$. 
\end{lem}
\begin{proof} Suppose $0\neq u\in L^2(\RR^d\setminus G)$ and $u$ is bounded. Then $\esssup |u|\neq 0$ and we may choose a set $H\subset G$ of positive Lebesgue measure such that $e^{-rh}\esssup|u| <|u(x)|$ for all $x\in H$ (this is possible because $r,h>0$ and therefore $e^{-rh}<1$), we deduce \begin{eqnarray*}\forall x\in H\quad \left|e^{-rh}\PP_{X_0-X_h}\ast u(x)\right| &\leq& e^{-rh} \esssup|u|\\&<& u(x), \end{eqnarray*} which means that $$\forall x\in H \quad PAu(x)=u(x)-e^{-rh}\PP_{X_0-X_h}\ast u(x)\neq 0,$$ hence $PAu\neq 0$ (for $H$ has positive Lebesgue measure). So $$\ker PA\restriction L^2(\RR^d\setminus G)=\ker PA\cap L^2(\RR^d\setminus G)=\{0\}$$ and we are done for the invertibility of $PA\restriction L^2(\RR^d\setminus G)$. Similarly, one can prove the invertibility of $A$. Finally, $A$ is seen to be a contraction by application of the Fourier transform: The Fourier transform is an $L^2$ isometry (by Plancherel's Theorem), thus \begin{eqnarray*}&&\left\|\left(\delta_0 - e^{-rh}\PP_{X_0-X_h}\right)\ast f\right\|_{L^2\left(\RR^d\right)}\\ &=&\left\|\left(\left(\delta_0 - e^{-rh}\PP_{X_0-X_h}\right)\ast f\right)^{\widehat{}}\right\|_{L^2\left(\RR^d\right)} \\ &=& \left\|\left(1 - e^{-rh}\widehat{\PP_{X_0-X_h}}\right)\cdot \widehat f\right\|_{L^2\left(\RR^d\right)} \\ &=& \left\|\left(1 - e^{-rh}e^{ih{^t}\mu\cdot}e^{-|\cdot|^2h/2}\right) \widehat f\right\|_{L^2\left(\RR^d\right)}\\ &=& \left(\int_{\RR^d}\left|1-e^{-rh+ih{^t}\mu x- \frac{|x|^2h}{2}}\right|^2\left|\widehat f(x)\right|^2 \ dx\right)^\frac{1}{2}. \end{eqnarray*}Now, the factor in front of $\left|\widehat f(x)\right|^2$ in the last line can be bounded by $\left(1+e^{-rh}\right)^{2}$, and it is strictly less than one for $\mu=0$. Using Plancherel's Theorem again, this yields the result. 
\end{proof}

Now, this is sufficient to apply a Wiener-Hopf factorisation (for a general treatment of this kind of factorisations, one may consult eg Speck \cite{Sp85}, our application uses in particular \cite[1.1, Theorem 1]{Sp85}) and state

\begin{Th} \label{l2whf} Let $G\subseteq\RR^d$ and let $X$ be a L\'evy basket with associated family of risk-neutral probability measures $\PP^\cdot$ and discount rate $r>0$. Then $V_G^h$, the expected payoff of a perpetual Bermudan option for $G$ with exercise mesh size $h>0$ and payoff function $g$, is -- using the above notation -- given by $$V_G^h=g_1-\left(PA\restriction L^2(\RR^d\setminus G)\right)^{-1}PAg_1=g_1-A_+^{-1}PA_-^{-1}PAg_1$$ where $A=A_-A_+$ is a Wiener-Hopf factorisation of $A$.
\end{Th}

We observe

\begin{lem} The Hilbert space operator $A:L^2(\RR^d,\CC)\rightarrow L^2(\RR^d,\CC)$ is normal.
\end{lem}
\begin{proof} We define $p:=e^{-rh} \frac{d\PP_{X_0-X_h}}{d\lambda^d}$ (where $\lambda^d$ is the $d$-dimensional Lebesgue measure) and via the Fubini Theorem one has for every $f,g\in L^2(\RR^d)$
\begin{eqnarray*} \langle Af,g\rangle&=& \int_{\RR^d}\int_{\RR^d}\left(\delta_0-p\right)  (x-y)f(y)dy \bar g(x) dx \\ &=& \int_{\RR^d}\int_{\RR^d}f(y)\left(\delta_0-p\right)  (x-y) \bar g(x) dx dy\\ &=& \langle f, \overline{\left(\delta - \bar p\circ(-\II)\ast g\right) }\rangle,\end{eqnarray*} that is $$A^*=\left(\delta_0-\bar p\circ(-\II)\right)\ast\cdot $$ But since the convolution is associative and commutative, this implies  \begin{eqnarray*}A^*A&=& \left(\delta_0- \bar p\circ(-\II)\right)\ast\left(\delta_0-p\right)\ast\cdot \\ &=& \left(\delta_0-p\right)\ast \left(\delta_0-\bar p\circ(-\II)\right)\ast\cdot \\&=& AA^*.\end{eqnarray*}
\end{proof}

However, it will not be possible to find a basic system of eigenvectors and eigenvalues for this operator, since 

\begin{lem} The operator $A$ fails to be compact.
\end{lem} 
\begin{proof} Any normalised basis provides a counterexample for the compactness assertion.
\end{proof}

Therefore, the equation (\ref{1.1}) cannot easily be applied to compute the expected option payoff by means of a spectral analysis. Thus, our examination of the Hilbert space approach in the second part of this Chapter has led to a negative outcome.

However, one can also conceive of the operators $A^h$ as operators on the Banach space $L^1\left(\RR^d\right)$:

\section{The $L^1$ operator equation: analyticity in the exercise mesh size}

From now on, $h$ will no longer be fixed and we will therefore write $A^h$ instead of $A$.

If we now assume $g_1$ to be an integrable extension of $G$ to the complement of $\complement G$ as an element of Quite similarly to \ref{l2whf}, we can prove

\begin{Th} Let $G\subseteq\RR^d$ and let $X$ be a L\'evy basket with associated family of risk-neutral probability measures $\PP^\cdot$ and discount rate $r>0$. Then $V_G^h$, the expected payoff of a perpetual Bermudan option for $G$ with exercise mesh size $h>0$ and payoff function $g$, is -- using the above notation -- given by $$V_G^h=g_1-\left(PA\restriction L^1(\RR^d\setminus G)\right)^{-1}PAg_1=g_1-A_+^{-1}PA_-^{-1}PAg_1,$$ where $A=A_-A_+$ is a Wiener-Hopf factorisation of $A$.
\end{Th}

It suffices to observe that $A^h$ is -- due to the $L^1$ norm estimate for the convolution of two integrable functions (as the product of the norms of the convolved functions) -- also a bounded operator on $L^1\left(\RR^d\right)$.

We shall now identify $g$ and $g_1$.

\begin{Th} With the notation previously introduced, we define $E$ to be the semigroup $$(E_t)_{t\geq 0}=\left(e^{-rt}\nu_{\mu t,t}\ast\cdot\right)_{t\geq 0}=\left( e^{-rt}g_{\mu t,t}\ast\cdot\right)_{t\geq 0}\in L\left(L^1(\RR^d),L^1(\RR^d)\right)^{[0,+\infty)},$$ where $$g_{\mu t,t}=(2\pi t)^{-d/2}e^{-|\mu t -\cdot|^2/(2t)}$$ is the distribution of the logarithmic price vector at time $t$. Suppose $x\not\in G$ and, with the notation from the previous chapters, $X$ is a (normalised) Brownian motion with (possibly zero) drift ({\em Black-Scholes model}). Then $h\mapsto V_G^h(x)$ is real analytic in $h$ on $(0,+\infty)$ as  function with range in the Banach space $L^1(\RR^d)$. 
\end{Th}

\begin{proof} It is obvious that $E$ is a semigroup. According to \cite[Theorem 1.48]{D}, the set $$\cE:=\left\{f\in L^1(\RR^d) \ : \ t\mapsto E_tf \in L^1(\RR^d)\text{ entire} \right\}$$ is dense in $L^1(\RR^d)$. Hence it is possible to approximate every $g$ by a sequence $\{g_k\}_k\subset\cE$ in $L^1(\RR^d)$. Since $$\forall t\geq 0 \quad \left\|E_t\right\|_{L^1(\RR^d)}\leq e^{-rt}\leq 1,$$ we obtain $E_tg_k\rightarrow E_tg$ for $k\rightarrow\infty$ uniformly in $t$ on $\RR_+$, where $E_tg_k$ is entire for every $t>0$ and $k\in\NN$. Thus, $t\mapsto E_tg$, and thereby $t\mapsto Pg-PE_tg$, is an analytic function on $(0,+\infty)$ taking values in the Banach space $L^1(\RR^d)$. 
Now observe that for arbitrary open $U\subset\subset\RR_+$ (the symbol ``$\subset\subset$'' indicating that $U$ is contained in a compact subset of $\RR_{>0}$) the following equations hold: \begin{eqnarray}\nonumber  \forall t\in U \quad V_G^t&=&g-\left(\II-PE_t\restriction PL^1(\RR^d)\right)^{-1}\left(Pg-PE_tg\right) \\ \nonumber &=& g-\left(\sum_{k=0}^\infty(PE_t)^k\right)\left(Pg-PE_tg\right) \\ \nonumber &=& g -\sum_{k=0}^\infty (PE_t)^{k}Pg +\sum_{k=0}^\infty (PE_t)^{k+1}g  \\  \label{harmonic} &=& \sum_{k=0}^\infty (PE_t)^{k}\left(g-Pg\right),\end{eqnarray} since the sums converge uniformly in $t$ on $U\subset\subset \RR_+$, yielding the analyticity of $t\mapsto V_G^t$ as a function whose range lies in the Banach space $L^1(\RR^d)$.
\end{proof}

\begin{lem} \label{deriv} Let $u>0$, $n\in\NN$. Then the equation $$\frac{d^n}{du^n}E_u = \left(-r+\frac{1}{2}\Delta-{{^t}\mu} \nabla\right)^n g_{t\mu ,t}\ast \cdot $$ holds (where ${^t}y$ denotes the transpose of a vector $y$). In particular, if $\frac{1}{2}\Delta f-\mu\nabla f= \lambda f$ for some $\lambda, f$, $$\frac{d^n}{du^n}E_u f= \left(-r+\lambda\right)^n g_{t\mu ,t}\ast f.$$ 
\end{lem}
\begin{proof} According to Davies \cite[Proof of Theorem 2.39]{D}, we have \begin{equation} \label{infgen} \frac{d^n}{du^n}E_u=\left(ZE_{u/n}\right)^n,\end{equation} where $Z$ denotes the infinitesimal generator of the semigroup $E$. Now, define $C$ to be the convolution operator semigroup $\left(g_{t\mu,t}\ast\cdot\right)_{t\geq 0}$ of (normalised) Brownian motion with drift $\mu$ (as before denoting by $g_{z,\sigma^2}$ the Lebesgue density of the Gaussian distribution centered around $z$ of variance $\sigma^2$ for all $z\in\RR^d$ and $\sigma>0$). It is well-known (cf e g \cite[p. 352]{RY}) that the infinitesimal generator of this semigroup $C$ is $$L:=\frac{1}{2}\Delta +{{^t}\mu} \nabla.$$ By our requirements on $f$, $Lf=0$ on $U$. Furthermore, $L$ and $C$ commute: $$\forall t\geq 0 \quad C_tL=LC_t.$$ 
Thus, \begin{eqnarray*} \forall t\geq 0 \quad ZE_t &=& \frac{d}{d}E_t =\frac{d}{dt}\left(e^{-rt}\cdot C_t\right)\\ &=& -re^{-rt}C_t+e^{-rt}\frac{d}{dt}C_t \\&=& e^{-rt}C_t(-r+L) ,\end{eqnarray*} which due to equation (\ref{infgen}) already suffices for the proof of the Lemma in the general case. And if $f$ is an eigenfunction of $L$ for the eigenvalue $\lambda$, one has $(-r+L)^nf=(-r+\lambda)^nf$.

\end{proof}

\begin{Th}\label{taylor} The Taylor series for the expected payoff of a perpetual Bermudan option as a function of the exercise mesh with respect to a fixed exercise region $G$ is for all $s>0$: \begin{eqnarray*} \forall t>0 \quad V_G^t &=& \sum_{k=0}^\infty (t-s)^k \sum_{m=1}^\infty e^{-rms}\sum_{\begin{array}{c} l_1+\dots+l_m=n \\ (l_1,\dots,l_m)\in{\NN_0}^m \end{array} } \left(\prod_{i=1}^m\frac{1}{l_i !}\right) \\ && \left( \chi_{\RR^d\setminus G}\cdot \left(g_{s\mu,s}\ast \cdot \right)\left(-r+\frac{1}{2}\Delta+{{^t}\mu}\nabla\right)^{\circ l_i}\right)\left(\chi_Gg\right),\end{eqnarray*} where, in order to avoid confusion with pointwise exponentiation, $A^{\circ k}$ denotes $A^k$ for any operator $A$.
\end{Th}
\begin{proof} We know about the real analyticity of $t\mapsto E_t$ on $\RR^{>0}$ and even, thanks to the previous Lemma, the explicit Taylor series. Thereby we also have the Taylor series for $t\mapsto PE_t$. So we can use equation (\ref{harmonic}) and see by means of a binomial expansion
\begin{eqnarray*}\sum_{k=0}^\infty \left(PE_t\right)^k &=&\sum_{k=0}^\infty \left(P\sum_{\ell= 0}^\infty \frac{(t-s)^\ell}{\ell !}\left(e^{-rs}\left(-r+L\right)\right)^\ell C_s\right)^k\\ &=& \sum_{n=0}^\infty \sum_{m=1}^\infty \sum_{\begin{array}{c} l_1+\dots+l_m=n \\ (l_1,\dots,l_m)\in{\NN_0}^m \end{array} } \prod_{i=1}^m \frac{(t-s)^{l_i}}{l_i!}e^{-rs}P(-r+L)^{l_i}C_s \\&=& \sum_{n=0}^\infty (t-s)^n \sum_{m=1}^\infty e^{-rms} \\&& \sum_{\begin{array}{c} l_1+\dots+l_m=n \\ (l_1,\dots,l_m)\in{\NN_0}^m \end{array} } \prod_{i=1}^m \frac{1}{l_i!} P(-r+L)^{l_i}C_s.\end{eqnarray*} 
\end{proof}

This Taylor series fails to provide any straightforward possibility for the computation of $V_G$. Instead we state the following immediate Corollary of equation (\ref{harmonic}):

\begin{cor}\label{firstderiv} With the notation as in the previous Theorem, \begin{eqnarray*} \forall s>0\quad \frac{d}{ds}V_G^s &=& \frac{d}{ds}\sum_{m=1}^\infty e^{-rms}\left( \chi_{\RR^d\setminus G}\cdot \left(g_{s\mu,s}\ast \cdot \right)\right)^{\circ m} \left(\chi_G\cdot g\right).\end{eqnarray*} 
\end{cor}

\chapter{An algebraic perpetual Bermudan pricing method and its natural scaling}

In this Chapter, we will, in a more algebraically flavoured way, present an approach that approximates the perpetual Bermudan option price as the fixed point of some map on the space of polynomials that is defined by means of not only the max operator, but also interpolation with respect to a given, fixed, set of interpolation points, as well as convolution with one and the same Gaussian (not necessarily normalised) measure. This set of interpolation points could, for example, be a set of cubature points for the distribution of the time $h$ increment of the logarithmic price process (if this process is assumed to be Gaussian with stationary increments).

Let us, for this purpose, adopt Victoir's notation \cite{V} and denote the space of all polynomials of degree $m$ and degree at most $m$ by $\RR_{m}[X]$ and $\RR_{\leq m}[X]$, respectively, for all $m\in\NN$. We will write polynomials in the form $p[X]$ and denote by $p(\cdot)$ the associated polynomial function from $\RR$ to $\RR$.

Now, let $m$ be any positive integer. We introduce the {\em interpolation map} $$\cI:\RR^{\RR}\times \RR^{m+1}\rightarrow \RR_{\leq m }[X]$$ that assigns to each pair $(f,\vec x)$ of a real-valued function $f$ and a vector $\vec x$  of interpolation points the well-defined (see e g Stoer and Bulirsch \cite{StB}) polynomial $p(X)\in\RR_m[X]$ satisfying $$\forall j\in\{0,\cdots,m\} \quad p\left(x_j\right)=f\left(x_j\right).$$

Now let $r>0$ be a positive real number (interpreted to be the discout rate), $\mu:=\frac{1}{2}-r\in\RR$ (which is the drift of the logarithmic price process in a Black-Scholes model with normalised time scale -- ensuring the volatility $\sigma$ to be equal to one), and $g:\RR\rightarrow\RR$ a continuous function (interpreted to be the payoff function of an option defined on the space $\RR$ of logarithmic underlying asset prices). Consider a vector-valued function $\vec y:\RR_{>0}\rightarrow \RR^{m+1}$ such that the range of the function $\vec y $ consists exclusively of vectors with mutually distinct entries. The elements of the range of $\vec y$ can in this case serve as sets of interpolation points (these points also called {\em support abscissas}). Hence using the notation of previous paragraphs, we may define another map $$\cH_h:\polm \rightarrow \polm$$ for all $h>0$ by $$\cH_h:p[X] \mapsto \cI\left(\left(p(\cdot)\ast \nu_{-\mu h, h}\cdot e^{-rh}\right)\vee g,\vec y(h)\right).$$

Note that $$\nu_{-\mu h,h}\ast\cdot:\RR_m[X]\rightarrow\RR_m[X].$$ This can be shown using the linearity of the convolution and the fact that for all $k\in\NN$, the convolution $\nu_{-\mu h,h}\ast(\cdot)^k$ of the measure $\nu_{-\mu h,h}$ with the function $x\mapsto x^k$ is again a polyonmial function of degree $k$. Writing down an explicit formula for the map $\nu_{-\mu h,h}\ast\cdot:\RR_m[X]\rightarrow\RR_m[X]$, we see that this function is continuous with respect to the Euclidean topology on the $(m+1)$-dimensional real vector space $\RR_{\leq m }[X]$.

This ushers in the proof of the following Lemma which is one of the first observations leading to the fixed point equation mentioned at the beginning of this Chapter.

\begin{lem} Consider any $h>0$ and arbitrary $p[X]\in\polm$. $\cH_h$ is continuous with respect to the Euclidean topology on the $(m+1)$-dimensional real vector space $\RR_{\leq m }[X]$. Also, if and only if $\cH_hp[X]=p[X]$, that is, $p[X]$ is a fixed point of $\cH_h$, there will exist a polynomial $q[X]\in \polm$ such that $p[X]=\lim_{n\rightarrow\infty}\left(\cH_h\right)^nq[X]$ (in the Euclidean toplogy of the $(m+1)$-dimensional vector space $\polm$). 
\end{lem}
\begin{proof} The continuity of $\cH_h$ is a consequence of the continuity of the map $\nu_{-\mu h,h}\ast\cdot:\RR_m[X]\rightarrow\RR_m[X]$. For the second part of the Lemma observe that provided the existence of such a $q[X]$ as in the statement of the Lemma, we can deduce the fixed point equation of the Lemma's statement from the continuity of the map $\cH_h$. For the converse implication, simply take $q[X]=p[X]$.
\end{proof}

\begin{Not} For any vector $\vec q={^t}\left(q_0,\cdots,q_m\right)\in\RR^{m+1}$, $q[X]$ will be understood to be the polynomial $q[X]=\sum_{k=0}^m q_k\cdot X^k\in\RR_{\leq m}[X]$, and for all $q[X]\in\RR_{\leq m}[X]$, $q(\cdot)$ shall be understood to denote the polynomial function $$ q(\cdot):x\mapsto \sum_{k=0}^m q_k\cdot x^k,$$ and .
\end{Not}

\begin{lem}\label{convpol} For all $h>0$ and each $p[X]\in\polm$, the function $(x,h)\mapsto p(\cdot)\ast \nu_{-\mu h,h}(x)$ is a polynomial function in both $x$ and $h$. Its degree in $x$ is $m$, its leading coefficient in $x$ being the leading coefficient of $p[X]$. Furthermore, the function $h\mapsto p(\cdot)\ast \nu_{-\mu h,h}(x)$ is $o(h)+ x^m$ for all $x\in\RR$. Moreover, if $x= o\left(h^\frac{1}{2}\right)$, then the first two leading terms of $(\cdot)^k\ast \nu_{-\mu h, h}(x)$ as a function of $h$ are $x^k= o\left(h^\frac{k}{2}\right)$ and a term of order $o\left(h^\frac{k+1}{2}\right)$, respectively, if $k$ is odd -- and if $k$ is even, the first two leading terms of $h\mapsto (\cdot)^k\ast \nu_{-\mu h, h}(x)$ are $x^k= o\left(h^\frac{k}{2}\right)$ and some term of order $o\left(h^\frac{k}{2}\right)$, respectively.
\end{lem}
\begin{proof} First, let us once again remark that $\nu_{-\mu h, h}\ast \cdot $ is, as a convolution operator, linear and that we therefore may restrict our attention to the functions $\nu_{-\mu h, h}\ast(x\mapsto x^k)$ for $k\leq m$. Consider any $x_0\in\RR$. Using the transformation $y\mapsto \frac{x+\mu h-y}{\sqrt{2h}}$ we find that \begin{eqnarray*} \nu_{-\mu h, h}\ast(x\mapsto x^k)(x_0) &=& \frac{1}{\sqrt{\pi}}\int_\RR\left(x_0+\mu h-z\sqrt{2h}\right)^k e^{-z^2}\ dz \\ &=& x^k+ \text{polynomial terms in ${x_0}^0,\cdots, {x_0}^{k-1},h$}, \end{eqnarray*} where the last line follows from expanding the binomial and using the identity $$\forall \ell\in\NN_0\quad \int_\RR z^{2\ell+1}e^{-z^2 } \ dz = 0$$ (an immediate consequence of the ``oddness'' of the integrand), which entails that all odd terms in $h^{1/2}$ will be cancelled out (as they have to be odd terms in $z$ as well).
\end{proof}

\begin{rem} \label{cubaturepoints o(h^1/2)} The assumption of $x= o\left(h^\frac{1}{2}\right)$ will hold in particular for any non-zero element of a set of cubature points for the measure $\nu_{-\mu h,h}$ derived from a cubature formula for $\nu_{0,1}$.
\end{rem}

\begin{Def} For all $h>0$, $\sigma\in\cS_{m+1}$ (the symmetric group of $\{0,\dots,m\}$), and $n\leq m$, we define the $(m+1)\times(m+1)$-matrix 
\begin{eqnarray*} && A_n^\sigma(h):= \\ && \left(\begin{array}{*{3}{c}}y_{\sigma (0)}(h)^0 & \cdots & y_{\sigma (0)}(h)^m \\ \vdots & \cdots & \vdots \\ y_{\sigma (n)}(h)^0 & \cdots  & y_{\sigma (n)}(h)^m \\ \begin{array}{r}-e^{-rh}(\cdot)^0\ast\nu_{-\mu h,h}\left(y_{\sigma (n+1)}(h)\right)\\ + y_{\sigma (n+1)}(h)^0 \end{array}&\cdots & \begin{array}{r} -e^{-rh}(\cdot)^{m}\ast\nu_{-\mu h,h}\left(y_{\sigma (n+1)}(h)\right)\\ + y_{\sigma (n+1)}(h)^m \end{array} \\ \vdots & \cdots & \vdots \\ \begin{array}{r} -e^{-rh} (\cdot)^0\ast\nu_{-\mu h,h}\left(y_{\sigma (m)}(h)\right)\\ + y_{\sigma (m)}(h)^0 \end{array} &\cdots & \begin{array}{r} -e^{-rh} (\cdot)^{m}\ast\nu_{-\mu h,h}\left(y_{\sigma (m)}(h)\right)\\ + y_{\sigma (m)}(h)^m \end{array}\end{array}\right)\end{eqnarray*} and the following vector: $$ \left( \begin{array}{c} g\left(y_{\sigma(0)}(h)\right) \\ \vdots \\ g\left(y_{\sigma(n)}(h)\right)\\ 0 \\ \vdots \\ 0\end{array}\right)=:\vec v\left(\vec y, n,h, \sigma\right).$$ We also define $\vec p\left(\vec y, n,h, \sigma\right)$ to be the solution $\vec q$ of $$A_n^\sigma(h)\cdot \vec q=\vec v\left(\vec y, n,h, \sigma\right),$$ provided there exists a unique solution to this equation.
\end{Def}

From the matrix formulation of the interpolation problem (again cf Stoer and Bulirsch \cite{StB}), the following Lemma is immediate:

\begin{lem} \label{fixalg} The polynomial $q[X]=\sum_{k=0}q_kX^k\in\polm$ is a fixed point of $\cH_h$ if and only if there is an $n\leq m$ and a $\sigma\in\cS_{m+1}$ (the symmetric group of $\{0,\dots,m\}$) such that $$A_n^\sigma(h)\cdot {^t}(q_0,\dots,q_m)= \vec v\left(\vec y, n,h, \sigma\right)$$ and, in addition, \begin{eqnarray*} &&\forall j\in\{0,\ldots,n\} \\ && g\left(y_{\sigma(j)}(h)\right)\vee e^{-rh}\left( q(h)(\cdot)\ast\nu_{-\mu h,h}\right)\left(y_{\sigma(j)}(h)\right) \\ &=& g\left(y_{\sigma(j)}(h)\right)\end{eqnarray*} as well as \begin{eqnarray*}&& \forall j\in\{n+1,\ldots,m\} \\  && g\left(y_{\sigma(j)}(h)\right)\vee e^{-rh}\left( q(h)(\cdot)\ast\nu_{-\mu h,h}\right)\left(y_{\sigma(j)}(h)\right)\\ & =& e^{-rh}\left( q(h)(\cdot)\ast\nu_{-\mu h,h}\right)\left(y_{\sigma(j)}(h)\right) .\end{eqnarray*}
\end{lem}

\begin{lem} \label{detneq0} Suppose there is an $\alpha\in(0,1)$ such that for all $j\in\{0,\dots,m\}$, $h\mapsto y_j(h)$ is a non-constant polynomial in $ h^\alpha$, except for possibly one $j_0\in\{0,\dots,m\}$ where $y_{j_0}(h)=0$ for all $h>0$. Then for all sufficiently small $h$, $\det A_n^\sigma(h)\neq 0$. The upper bound in $\RR_{>0}\cup\{+\infty\}$ on all those $h$ that satisfy the previous inequality $\det A_n^\sigma(h)\neq 0$ for all $n\leq m$ and $\sigma\in\cS_{m+1}$ shall be denoted by $h_0\left(\vec y\right)$.
\end{lem}

\begin{rem} According to Remark \ref{cubaturepoints o(h^1/2)}, the assumption of $h\mapsto y_j(h)$ being non-constannt polynomial in $ h^\alpha$ for all $j\in\{0,\dots,m\}$ (apart from possibly one zero coordinate) for some $\alpha$ holds in particular for any set of cubature points for the measures $\nu_{-\mu h,h}$ that is derived from a cubature formula for the normalised Gaussian measure (in that case $\alpha=\frac{1}{2}$). Similar assertions hold if one replaces $\nu_{-\mu h,h}$ by $m_h$ where $(m_t)_{t\geq 0}$ is the convolution semigroup associated to some other symmetric stable process.
\end{rem}

\begin{proof}[Proof sketch for Lemma \ref{detneq0}] The function $h\mapsto\det A_n^\sigma(h)$ is, by our assumptions on the functions $y_j$ on the one hand polynomial in $h^\alpha$, as one can see exactly as in the proof of Lemma \ref{convpol}. On the other hand, one can show, using the polyonmiality in $h$ and the assumpion that all the entries of $\vec y(h)$ are mutually distinct for all $h>0$, that the function $h\mapsto\det A_n^\sigma(h)$ is non-constant. Hence, $h\mapsto\det A_n^\sigma(h)$ is a non-constant analytic function in $h^{\alpha}$, and note that $h\mapsto h^\alpha$ is a bijection on the unit interval $(0,l)$. Therefore, $h\mapsto\det A_n^\sigma(h)$ cannot be constantly zero on the open unit interval $(0,1)$, but it can also only have finitely many critical points on that interval. Hence there must be an $\bar h>0$ such that either $\det A_n^\sigma(h)< 0$ for all $h\in\left(0,\bar h\right)$ or $\det A_n^\sigma(h)< 0$ for all $h\in\left(0,\bar h\right)$.
\end{proof}

On the other hand, we have the following

\begin{lem}\label{constsig,n} Suppose the real-valued function $h\mapsto y_j(h)$ is a non-constant polynomial in $h\mapsto h^\alpha$ for all $j\in\{0,\dots,m\}$ for some $\alpha\in(0,1)$, except for possibly one $j_0\in\{0,\dots,m\}$ where $y_{j_0}(h)=0$ for all $h>0$. Furthermore, take $g$ to be analytic. Then for all $n\leq m$ and all $\sigma\in\cS_{m+1}$ there exists a vector $\cR(n,\sigma)\in\{\leq,>\}^{m+1}$ of relations such that for all sufficiently small $h>0$, \begin{eqnarray*}&&\forall i\in\{0,\dots,m\} \\ && \left( g-e^{-rh}p\left(\vec y, n,h, \sigma\right)(\cdot)\ast\nu_{-\mu h, h}\right)\left(y_{i}(h)\right) \ \cR(\sigma,n)_{i} \ 0.\end{eqnarray*} The upper bound in $\RR_{>0}\cup\{+\infty\}$ on all those $h_1$ such that for all $h<h_0\left(\vec y\right)\wedge h_1$ the relations in the previous line hold for all $n\leq m$ and $\sigma\in\cS_{m+1}$ shall be denoted by $h_1\left(\vec y\right)$. (Here $h_0\left(\vec y\right)$ is the strictly positive constant of Lemma \ref{detneq0}.) Thus, $h_1\left(\vec y\right)\leq h_0\left(\vec y\right)$.
\end{lem}

\begin{proof} We have already defined $\vec p\left(\vec y, n,h, \sigma\right)$ to be the solution $\vec q$ of $$A_n^\sigma(h) \cdot {^t}\left(q_{\sigma(0)},\dots,q_{\sigma(m)}\right)=\vec v\left(\vec y, n,h, \sigma\right).$$ From Cramer's rule, Lemma \ref{convpol}, and our assumptions on $\vec y$ (coordinatewise polynomial in $(\cdot)^\alpha$) as well as $g$ (analyticity), we derive that $h\mapsto \left(g-e^{-rh}p\left(\vec y, n,h, \sigma\right)(\cdot)\ast\nu_{-\mu h, h}\right)\left(y_j(h)\right)$ is analytic in $(\cdot)^\alpha$ and therefore has only finitely many critical points on $(0,1)$ for all $j\in\{0,\dots, m\}$ and arbitrary choice of $\vec y,n,\sigma$. Thus, when approaching zero, these functions must eventually stay on either side of nought. Put more formally, there must be for all $\sigma\in\cS_{m+1}$ and $n\leq m$ a vector $\cR(n,\sigma)\in\{\leq,>\}^{m+1}$ such that for all sufficiently small $h$, and for all $i\in\{0,\dots, m\}$, $$\left( g-e^{-rh}p\left(\vec y, n,h, \sigma\right)(\cdot)\ast\nu_{-\mu h, h}\right)\left(y_{i}(h)\right) \ \cR(\sigma,n)_{i} \ 0.$$ 
\end{proof}

\begin{cor} \label{existfix} Let the assumptions of the previous Lemma hold. If there is a fixed point of $\cH_{h_1}$ for an $h_1<h_1\left(\vec y\right)$ (the strictly positive constant of Lemma \ref{constsig,n}), all $\cH_h$ with positive $h\leq h_1$ must have fixed points as well. There exist a permutation $\sigma_0\in\cS_{m+1}$ as well as a natural number $n_0\leq m$ such that the coefficient vectors $\vec q={^t}(q_0,\dots,q_m)$ to all these fixed points $ q[X]=\sum_{k=0}^m q_kX^k$ are solutions $\vec q$ to the linear equation $$A_{n_0}^{\sigma_0}(h)\cdot \vec q= \vec v\left(\vec y, n_0,h, \sigma_0\right)$$ and satisfy \begin{eqnarray*} && \forall j\in\{0,\ldots,n_0\} \\ && g\left(y_{\sigma_0(j)}(h)\right)> e^{-rh}\left( q(h)(\cdot)\ast\nu_{-\mu h,h}\right)\left(y_{\sigma_0(j)}(h)\right) \end{eqnarray*} as well as \begin{eqnarray*}&& \forall j\in\{n_0+1,\ldots,m\} \\  && g\left(y_{\sigma_0(j)}(h)\right)\leq e^{-rh}\left( q(h)(\cdot)\ast\nu_{-\mu h,h}\right)\left(y_{\sigma_0(j)}(h)\right) .\end{eqnarray*}
\end{cor}
\begin{proof} By virtue of Lemma \ref{fixalg}, a fixed point is a polynomial $q[X]=\sum_{k=0}^m q_k X^k$ whose coordinate vector $\vec q={^t}(q_0,\dots,q_m)$ solves the linear equation $$A_n^\sigma(h)\cdot \vec q= \vec v\left(\vec y, n,h, \sigma\right)$$ and satisfies, moreover, inequalities of the form $$\left( g-e^{-rh}q(\cdot)\ast\nu_{-\mu h, h}\right)\left(y_{i}(h)\right) \ \cR_{i} \ 0$$ for some vector of relations $\cR\in\{\leq, >\}^m$. Now apply the previous Lemma \ref{constsig,n}.

\end{proof}

\begin{Th}\label{algscale} Suppose there is a vector $\vec\xi\in\RR^{m+1}$ such that $y_j(h)=-\mu h + \xi_j\cdot h^{1/2}$ for all $j\in\{0,\dots,m\}$ (where the mutual distinctness of the entries of $\vec y(h)$ for all $h>0$ entails that $\xi_0,\dots,\xi_m$ are mutually distinct as well), and assume furthermore that $g$ is analytic and satisfies $g(0)\neq 0$. Suppose, moreover, that there exists a fixed point of $\cH_{h_2}$ for some $h_2\in\left(0,h_1\left(\vec y\right)\right)$ (where $h_1\left(\vec y\right)$ is the strictly positive constant from Lemma \ref{constsig,n}). Then by Corollary \ref{existfix} the maps $\cH_h$ do have a fixed point for every $h<h_2$, and let $n,\sigma$ be the natural number $\leq m$ and the permutation whose existence is stated in Corollary \ref{existfix}, respectively. Also write $p(h)[X]:=p\left(\vec y, n,h, \sigma\right)[X]$ for all $h\in(0,h_2)$. Then there is a polynomial $p^0[X]$ such that $$p(h)[X]=p^0[X]+ o\left(h^\frac{1}{2}\right)$$ componentwise.
\end{Th}

\begin{rem} Such a $\xi\in\RR^{m+1}$ exists in particular whenever $\left\{y_0,\dots,y_m\right\} $ is a set of cubature points for the measure $\nu_{-\mu h,h}$ derived from a cubature formula for $\nu_{0,1}$.
\end{rem}

\begin{proof}[Proof of Theorem \ref{algscale}] Observe that by Lemma \ref{convpol}, \begin{eqnarray*} && A_n^\sigma(h) \\ &=& \left(\begin{array}{*{4}{c}} 1& o\left(h^\frac{1}{2}\right) & \cdots& o\left(h^\frac{m}{2}\right)\\ \vdots &\vdots &\cdots & \vdots \\ 1& o\left(h^\frac{1}{2}\right) & \cdots& o\left(h^\frac{m}{2}\right) \\ \left(1-e^{-rh}\right) &\left(1-e^{-rh}\right)\cdot  o\left(h^\frac{1}{2}\right) & \cdots& \left(1-e^{-rh}\right) \cdot o\left(h^\frac{m}{2}\right) \\ \vdots & \vdots &\cdots & \vdots  \\ \left(1-e^{-rh}\right) &\left(1-e^{-rh}\right)\cdot  o\left(h^\frac{1}{2}\right) & \cdots& \left(1-e^{-rh}\right) \cdot o\left(h^\frac{m}{2}\right) \end{array}\right) \\ &=& \left(\begin{array}{*{4}{c}} 1& o\left(h^\frac{1}{2}\right) & \cdots& o\left(h^\frac{m}{2}\right)\\ \vdots &\vdots &\cdots & \vdots \\ 1& o\left(h^\frac{1}{2}\right) & \cdots& o\left(h^\frac{m}{2}\right) \\ o\left(h^\frac{2}{2}\right) & o\left(h^\frac{3}{2}\right) & \cdots& o\left(h^\frac{m+2}{2}\right) \\ \vdots & \vdots &\cdots & \vdots  \\ o\left(h^\frac{2}{2}\right) & o\left(h^\frac{3}{2}\right) & \cdots& o\left(h^\frac{m+2}{2}\right) \end{array}\right) \end{eqnarray*}

Next, we will use Cramer's rule to determine if there is a limit for the solution of $$ A_{n}^{\sigma}(h)\cdot\vec p(h)= \vec v\left(\vec y, n,h, \sigma\right)$$ as $h$ tends to zero and if so, what the convergence rate will be. For this purpose, we have to consider the determinant of the matrix $A_{n,i}^\sigma(h)$ which is defined to be the matrix coinciding with $A_n^\sigma(h)$ in the columns $0,\ldots,i-1,i+1,\ldots m$ and having the vector $$\left( \begin{array}{c} g\left(y_0(h)\right) \\ \vdots \\ g\left(y_n(h)\right)\\ 0 \\ \vdots \\ 0\end{array}\right)$$ as its $i$-th column. Then, since $g$ is right-differentiable in $0$ and by assumption $g(0)\neq 0$, \begin{eqnarray*} && A_{n,i}^\sigma(h) = \\ && \left(\begin{array}{*{7}{c}} o\left(h^\frac{0}{2}\right) & \cdots& o\left(h^\frac{i-1}{2}\right) &  g(0)+o\left(h^\frac{1}{2}\right)& o\left(h^\frac{i+1}{2}\right)& \cdots & o\left(h^\frac{m}{2}\right)\\ \vdots &\cdots & \vdots & \vdots & \vdots  &\cdots & \vdots \\ o\left(h^\frac{0}{2}\right) & \cdots& o\left(h^\frac{i-1}{2}\right) &  g(0)+o\left(h^\frac{1}{2}\right)& o\left(h^\frac{i+1}{2}\right)& \cdots & o\left(h^\frac{m}{2}\right) \\ o\left(h^\frac{2}{2}\right)&\cdots & o\left(h^\frac{i+1}{2}\right) & 0 & o\left(h^\frac{i+3}{2}\right) & \cdots& o\left(h^\frac{m+2}{2}\right) \\ \vdots & \cdots & \vdots &\vdots & \vdots &\cdots &\vdots   \\ o\left(h^\frac{2}{2}\right)&\cdots & o\left(h^\frac{i+1}{2}\right) & 0 & o\left(h^\frac{i+3}{2}\right) & \cdots& o\left(h^\frac{m+2}{2}\right) \end{array}\right) \\ && =\left(\begin{array}{*{7}{c}} o\left(h^\frac{0}{2}\right) & \cdots& o\left(h^\frac{i-1}{2}\right) &  o(1) & o\left(h^\frac{i+1}{2}\right)& \cdots & o\left(h^\frac{m}{2}\right)\\ \vdots &\cdots & \vdots & \vdots & \vdots  &\cdots & \vdots \\ o\left(h^\frac{0}{2}\right) & \cdots& o\left(h^\frac{i-1}{2}\right) &  o(1) & o\left(h^\frac{i+1}{2}\right)& \cdots & o\left(h^\frac{m}{2}\right) \\ o\left(h^\frac{2}{2}\right)&\cdots & o\left(h^\frac{i+1}{2}\right) & 0 & o\left(h^\frac{i+3}{2}\right) & \cdots& o\left(h^\frac{m+2}{2}\right) \\ \vdots & \cdots & \vdots &\vdots & \vdots &\cdots &\vdots   \\ o\left(h^\frac{2}{2}\right)&\cdots & o\left(h^\frac{i+1}{2}\right) & 0 & o\left(h^\frac{i+3}{2}\right) & \cdots& o\left(h^\frac{m+2}{2}\right) \end{array}\right).\end{eqnarray*} From this, we can conclude that for all $i\in\{0,\ldots,m\}$ there is a constant $p^0_i\in\RR$ such that $$p_i(h)=\frac{\det A_{n,i}^\sigma(h) }{\det A_{n}^\sigma(h)}=p_i^0+ o\left(h^\frac{1}{2}\right).$$
\end{proof}

\begin{ex}[$m=2$, the quadratic case, when $\mu=0$] We have \begin{eqnarray*} \forall\mu\in\RR\quad  \forall x\in\RR\quad && (\cdot)^0\ast\nu_{-\mu h,h}(x)=1\\ && (\cdot)^1\ast\nu_{-\mu h,h}(x)=x + \mu h\\ && (\cdot)^2\ast\nu_{-\mu h,h}(x)=x^2+ 2\mu h x + h +h^2,\end{eqnarray*} implying \begin{eqnarray*} \forall x\in\RR\quad && (\cdot)^0\ast\nu_{0,h}(x)=1\\ && (\cdot)^1\ast\nu_{0,h}(x)=x \\ && (\cdot)^2\ast\nu_{0,h}(x)=x^2 + h +h^2 \end{eqnarray*} and $$\forall h>0 \quad \vec y(h) = \left(\begin{array}{c}\xi_0\cdot h^\frac{1}{2}\\ \xi_1\cdot h^\frac{1}{2}\\ \xi_2\cdot h^\frac{1}{2}\end{array}\right).$$ Suppose the parameters from Corollary \ref{existfix} are in our example $n_0=1$ and $\sigma_0=\id$. Furthermore, in our case, for all $h>0$ $$A_1^\id(h)=\left(\begin{array}{*{3}{c}} 1 & \xi_0h^\frac{1}{2} & {\xi_0}^2h\\ 1 & \xi_1h^\frac{1}{2} & {\xi_1}^2h\\ 1-e^{-rh} & \left(1-e^{-rh}\right)\xi_2h^\frac{1}{2} & \left(1-e^{-rh}\right){\xi_2}^2h -e^{-rh}(h+h^2) \end{array}\right).$$ Therefore for all $h>0$, \begin{eqnarray*}\det A_1^\id(h)&=& \xi_1 h^\frac{3}{2}\left(\left(1-e^{-rh}\right){\xi_2}^2 -e^{-rh}(1+h)\right) \\&& + \xi_0{\xi_1}^2h^\frac{3}{2} \left(1-e^{-rh}\right) + \left(1-e^{-rh}\right){\xi_0}^2\xi_2h^\frac{3}{2}\\ &&-\xi_0h^\frac{3}{2}\left(\left(1-e^{-rh}\right){\xi_2}^2-e^{-rh}(1+h)\right)\\&&-{\xi_1}^2\left(1-e^{-rh}\right)\xi_2h^\frac{3}{2}-{\xi_0}^2\xi_1 h^\frac{3}{2} \left(1-e^{-rh}\right) \\ &=& \left(1-e^{-rh}\right)h^\frac{3}{2}\left(\xi_1{\xi_2}^2 + \xi_0{\xi_1}^2 + {\xi_0}^2\xi_2 - \xi_0{\xi_2}^2-{\xi_1}^2\xi_2 -{\xi_0}^2\xi_1\right) \\ && + h^\frac{5}{2}\cdot e^{-rh}\left(\xi_0-\xi_1\right)+h^\frac{3}{2}\cdot e^{-rh}\left(\xi_0-\xi_1\right).\end{eqnarray*} Note that $$\xi_1{\xi_2}^2 + \xi_0{\xi_1}^2 + {\xi_0}^2\xi_2 - \xi_0{\xi_2}^2-{\xi_1}^2\xi_2 -{\xi_0}^2\xi_1=\det B\left(\vec \xi\right),$$ where $B\left(\vec \xi\right)$ is the matrix of the interpolation problem with support abscissas $\xi_0,\xi_1,\xi_2$. Due to the unique solvability of the interpolation problem (see again e g Stoer and Bulirsch \cite{StB}), this determinant $\det B\left(\vec \xi\right)$ never vanishes unless the support abscissas $\xi_0,\xi_1,\xi_2$ fail to be mutually distinct.
\end{ex}

\begin{rem} Similarly one can prove that the function $h\mapsto p^h[X]$ is differentiable in $0$ if one assumes the support abscissas to be polynomial in $h$ rather than $h^\frac{1}{2}$.
\end{rem}

\end{document}